\documentclass[reqno]{amsart}
\numberwithin{equation}{section}
\newtheorem{theorem}{\sc Theorem}[section]
\newtheorem{lemma}{\sc Lemma}[section]
\newtheorem{proposition}{\sc Proposition}[section]
\newtheorem{remark}{\sc Remark}
\newtheorem{definition}{\sc Definition}[section]

\def\bbR{{\mathbb R}}
\def\id{{\text{Id}}}
\def\cf{{\text{Cof}}}
\def\sk{\sqrt{\kappa}}
\def\ks{{\kappa}}
\def\nk{{n^\kappa}}
\def\pk{{\Pi^\kappa}}
\def\vk{{v^\kappa}}
\def\sE{\sup_{t\in[0,T]} E_\kappa(t)}
\def\skE{\sup_{t\in[0,T]} E_\kappa(t)}
\def\psE{P(\sup_{t\in[0,T]} E_\kappa(t))}
\def\pskE{P(\sup_{t\in[0,T]} E_\kappa(t))}

\def\ek{{\eta^{\kappa}}}
\def\vk{{v^\kappa}}
\def\uk{{u^\kappa}}
\def\pk{{p^\kappa}}

\def\ty{{\tilde y}}
\def\tz{{\tilde z}}

\def\qk{{q^\kappa}}
\def\nk{{n^\kappa}}
\def\NK{{N^\kappa}}

\def\o{\Omega}

\def\s{\star}
\def\R{\mathbb R}

\def\n{\nonumber}
\def\d{\displaystyle}
\def\w{\tilde w}

\def\a{\bar a}
\def\b{\tilde b}
\def\v{\tilde v}
\def\q{\tilde q}
\def\e{\tilde\eta}

\def\u{\tilde u}
\def\ur{\tilde u^{r}}
\def\p{\tilde p}

\def\E{\tilde{E}}

\def\ci{\check{I}}

\def\o{\operatorname}
\def\hf{\hfill\break}
\def\un{{\frac{1}{2}}}

\def\ci{{\frac{5}{2}}}
\def\se{{\frac{7}{2}}}
\def\ud{{\frac{1}{2}}}
\def\o{\Omega}

\def\s{\star}
\def\R{\mathbb R}

\def\n{\nonumber}
\def\d{\displaystyle}
\def\w{\tilde w}

\def\en{\overline{\eta}}
\def\uk{{\bar{u}_{\kappa}}}
\def\nk{{\bar{n}_{\kappa}}}
\def\ek{{{\bar\eta_{\kappa}}}}
\def\an{{\bar{a}^\kappa}}
\def\ak{{\bar{a}^k}}
\def\un{{\bar{u}^n}}
\def\a{\tilde a}
\def\b{\tilde b}
\def\p{\tilde p}
\def\q{\tilde q}
\def\e{\tilde\eta}

\def\w{\tilde w}
\def\we{{w_{\epsilon}}}
\def\qe{{q_{\epsilon}}}
\def\ue{{u_{\epsilon}}}
\def\pe{{p_{\epsilon}}}

\def\E{\tilde{E}}

\def\ci{\check{I}}

\def\ud{{\frac{1}{2}}}
\def\td{{\frac{3}{2}}}
\def\cd{{\frac{5}{2}}}
\def\sd{{\frac{7}{2}}}
\def\div{\operatorname{div}}
\def\curl{\operatorname{curl}}
\def\bak{{{\bar a}_\kappa}}
\def\bbk{{{\bar b^l}_\kappa}}

\title[Well-posedness of the free-surface Euler equations]
{Well-posedness of the free-surface incompressible Euler  equations
       with or without surface tension}

\author[D. Coutand and S. Shkoller]{}


\email{coutand@math.ucdavis.edu}
\email{shkoller@math.ucdavis.edu}

\begin{document}
\maketitle

\centerline{\scshape   Daniel Coutand and Steve Shkoller}
 \medskip

  {\footnotesize \centerline{ Department of Mathematics }
 \centerline{ University of California at Davis } \centerline{Davis, CA
   95616 } }

\begin{abstract}
We provide a new method for treating free boundary problems in perfect
fluids, and prove local-in-time well-posedness in Sobolev spaces for the 
free-surface incompressible 3D Euler equations with or without surface tension 
for arbitrary initial data, and without any irrotationality assumption on
the fluid.   This is a free boundary problem for the motion of an 
incompressible perfect liquid in vacuum, wherein the motion of the fluid
interacts with the motion of the free-surface at highest-order.  
\end{abstract}

\tableofcontents

\section{Introduction}
\label{sec_introduction}

\subsection{The problem statement and  background}
For $\sigma\ge 0$ and for arbitrary initial data,
we prove local existence and uniqueness of solutions in Sobolev spaces to the 
free boundary incompressible Euler equations in vacuum:
\begin{subequations}
  \label{euler}
\begin{alignat}{2}
\partial_t u+ \nabla_uu + \nabla p&=0  &&\text{in} \ \ Q \,, \label{euler.a}\\
   \div u &= 0     &&\text{in} \ \ Q \,, \label{euler.b}\\
p &= \sigma \, H \ \ &&\text{on} \ \ \partial Q \,, \label{euler.c}\\
(\partial_t+\nabla_u)|_{\partial Q}&\in T(\partial Q)\,, &&\ \ \label{euler.d}\\
   u &= u_0  &&\text{on} \ \ Q_{t=0} \,, \label{euler.e}\\
   Q_{t=0} &= \Omega\,,  && \label{euler.f}
\end{alignat}
\end{subequations}
where $Q= \cup_{0\le t\le T}\{t\}\times \Omega(t)$, $\Omega(t) \subset 
{\mathbb R}^n$, $n=2$ or $3$,
$\partial Q= \cup_{0\le t\le T}\{t\}\times \partial\Omega(t)$, 
$\nabla_uu = u^j\partial u^i/\partial x^j$, 
and where Einstein's summation convention is employed.  The vector field $u$
is the Eulerian or spatial velocity field defined on the time-dependent
domain $\Omega(t)$, $p$ denotes the pressure function, $H$ is twice the
mean curvature of the boundary of the fluid $\partial \Omega(t)$, and
$\sigma$ is the surface tension.  Equation (\ref{euler.a}) is the conservation
of momentum, (\ref{euler.b}) is the conservation of mass, (\ref{euler.c}) is
the well-known Laplace-Young boundary condition for the pressure function,
(\ref{euler.d}) states that the free boundary moves with the velocity of the 
fluid, (\ref{euler.e}) specifies the initial velocity, and (\ref{euler.f}) fixes
the initial domain $\Omega$.  

Almost all prior well-posedness results were focused on {\it irrotational}
fluids (potential flow), wherein the additional constraint 
$\operatorname{curl}u=0$ is imposed; 
with the irrotationality constraint, the Euler equations
(\ref{euler}) reduce to the well-known water-waves equations, wherein the 
motion of the interface is decoupled from the rest of the fluid and is governed
by singular boundary integrals that arise from the use of complex variables
and the equivalence of incompressibility and irrotationality with the
Cauchy-Riemann equations.
For 2D fluids (and hence 1D interfaces), the earliest local existence results
were obtained by Nalimov \cite{Na1974}, Yosihara \cite{Yo1982}, and 
Craig \cite{Cr1985} for initial data near equilibrium.  Beale, Hou, \&
Lowengrub \cite{BeHoLo1993} proved that the linearization  of the 2D water wave 
problem is well-posed if a Taylor sign condition is added to the problem 
formulation, thus preventing Rayleigh-Taylor instabilities.  Using the
Taylor sign condition, Wu \cite{Wu1997} proved local existence for the 2D 
water wave problem for arbitrary (sufficiently smooth) initial data.  Later
Ambrose \cite{Am2003} and Ambrose \& Masmoudi \cite{AmMa2005}, proved local
well-posedness of the 2D water wave problem with  surface tension
on the boundary replacing the Taylor sign condition.

In 3D, Wu \cite{Wu1999} used Clifford analysis to prove local existence of the
full water wave problem with {\it infinite depth}, showing that the Taylor sign 
condition is always satisfied in the irrotational case by virtue of the maximum
principle holding for the potential flow.  Lannes \cite{La2005} provided a proof
for the {\it finite depth case with varying bottom} by implementing a 
Nash-Moser iteration.  The first well-posedness result for
the full Euler equations with zero surface tension, $\sigma=0$, is due to 
Lindblad \cite{Li2004} with the additional ``physical condition'' that
\begin{equation}\label{lindblad}
\nabla p \cdot n < 0 \text{  on  } \partial Q,
\end{equation}
where $n$ denotes the exterior unit normal to $\partial \Omega(t)$.   The
condition (\ref{lindblad}) is equivalent to the Taylor sign condition, and
provided Christodoulou \& Lindblad \cite{ChLi2000} with enough
boundary regularity to establish {\it a priori} estimates for smooth solutions
to (\ref{euler}) together with (\ref{lindblad}) and $\sigma=0$.  (Ebin
\cite{Eb1987} provided a counterexample to well-posedness when (\ref{lindblad})
is not satisfied.)  
Nevertheless, local existence did not
follow in \cite{ChLi2000}, as finding approximations of the Euler equations for
which existence and uniqueness is known and which retain the transport-type
structure of the Euler equations is highly non-trivial, and this geometric
transport-type 
structure is crucial for the a priori estimates.  In \cite{Li2003}, Lindblad
proved well-posedness of the linearized Euler equations, but the estimates
were not sufficient for well-posedness of the nonlinear problem.  The
estimates were improved
in \cite{Li2004}, wherein Lindblad implemented a Nash-Moser iteration to
deal with the manifest loss of regularity in his linearized model and thus
established the well-posedness result in the case that (\ref{lindblad}) holds
and $\sigma=0$.

Local existence for the case of positive surface tension, $\sigma >0$, 
remained open, and although the Laplace-Young condition (\ref{euler.c}) 
provides improved regularity for the boundary, the required nonlinear estimates
are more difficult to close due to the complexity of the mean curvature 
operator, and the need to study time-differentiated problems, which do not
arise in the $\sigma=0$ case.
It appears that the use of the time-differentiated problem in Lindblad's 
paper \cite{Li2004} is due to the use of certain tangential projection 
operators, but this is not necessary.  We
note that our energy function is different from that in \cite{Li2004}, and
provides better control of the Lagrangian coordinate.

After completing this work, we were informed of the paper of Schweizer
\cite{Sch2006} who studies the Euler equations for $\sigma>0$  in the case that
the free-surface is a graph over the two-torus.  In that paper, 
he obtains a priori estimates
under a smallness assumption for the initial surface; well-posedness follows
under the additional assumption that there is no vorticity on the boundary.
We also learned of the paper by
Shatah and Zeng \cite{ShZe2006} who establish a priori estimates for
both the $\sigma=0$ and $\sigma>0$ cases without any restrictions on the 
initial data.

\subsection{Main results}
We prove two main theorems concerning the well-posedness of (\ref{euler}).  
The first theorem, for the case of positive surface tension $\sigma>0$, is new; for our second
theorem, corresponding to the zero surface tension case, we present a new proof
that does not require a Nash-Moser procedure, and has optimal regularity.

\begin{theorem}[Well-posedness with surface tension]\label{theorem1}
Suppose that $\sigma>0$, $\Gamma$  is 
of class $H^{5.5}$, and $u_0 \in H^{4.5}(\Omega)$.  Then,
there exists $T>0$, and a solution ($u(t)$,$p(t)$,$\Omega(t)$) of
(\ref{euler}) with $u \in L^\infty(0,T; H^{4.5}(\Omega(t)))$,
$p \in L^\infty(0,T; H^{4}(\Omega(t)))$, and
$\Gamma(t) \in H^{5.5}$.  The solution is unique if $u_0 \in H^{5.5}$
and $\Gamma \in H^{6.5}$.
\end{theorem}

\begin{theorem}[Well-posedness with Taylor sign condition]\label{theorem2}
Suppose that $\sigma=0$, $\partial \Omega$  is 
of class $H^{3}$, and $u_0 \in H^{3}(\Omega)$ and condition
(\ref{lindblad}) holds at $t=0$.   Then,
there exists $T>0$, and a unique solution ($u(t)$,$p(t)$,$\Omega(t)$) of
(\ref{euler}) with $u \in L^\infty(0,T; H^{3}(\Omega(t)))$,
$p \in L^\infty(0,T; H^{3.5}(\Omega(t)))$, and
$\partial\Omega(t) \in H^{3}$.
\end{theorem}

\subsection{Lagrangian representation of the Euler equations.}
The Eulerian problem (\ref{euler}), set on the moving domain $\Omega(t)$, is 
converted to a PDE on the fixed domain $\Omega$, by the use of Lagrangian
variables.
Let $\eta(\cdot, t):\Omega \rightarrow \Omega(t)$ be the solution of
$$
\partial_t\eta(x,t) = u(\eta(x,t),t), \ \ \ \eta(x,0)=\text{Id}\,.
$$
and set
$$v(x,t):= u(\eta(x,t),t), \ \ q(x,t):=p(\eta(x,t),t), \ \ \text{and} \ \
a(x,t) := [\nabla \eta(x,t)]^{-1}\,.$$
The variables $v$, $q$ and $a$ are functions of the fixed domain $\Omega$
and denote the material velocity, pressure, and pull-back, respectively. 
Thus, on the fixed domain, (\ref{euler}) transforms to
\begin{subequations}
  \label{leuler}
\begin{alignat}{2}
\eta&=\text{Id} + \int_0^t v\ \ \   &&\text{in} \ \Omega \times (0,T]\,, \label{leuler.a}\\
\partial_t v+ a\,\nabla q&=0  &&\text{in} \ \Omega \times (0,T]\,, \label{leuler.b}\\
   \text{Tr}(a \, \nabla v)  &= 0     &&\text{in} \   \Omega \times (0,T]
\,, \label{leuler.c}\\
q\, a^T N/|a^TN| &=-\sigma \, \Delta_g(\eta) \ \ &&\text{on} \  \Gamma \times (0,T] \,, \label{leuler.d}\\
   (\eta,v) &= (\text{Id}, u_0)  &&\text{on} \ \Omega\times\{t=0\} 
                                                 \,, \label{leuler.e} 
\end{alignat}
\end{subequations}
where $N$ denotes the unit normal to $\Gamma$, and
$\Delta_g$ is the surface Laplacian with respect to the induced metric 
$g$ on
$\Gamma$, written in local coordinates as 
\begin{equation}\label{gstuff}
\Delta_g = \sqrt{g}^{-1} \partial_\alpha [\sqrt{g} g^{\alpha\beta} 
\partial_\beta ]\,, \  g^{\alpha \beta} = [g_{\alpha\beta}]^{-1}\,,  \
g_{\alpha\beta} = \eta_{,\alpha}\cdot \eta_{,\beta}\,, \text{ and }
\sqrt{g} = \sqrt{\det{g}} \,.
\end{equation}

\begin{theorem}[$\sigma>0$]\label{ltheorem1}
Suppose that $\sigma>0$, $\partial \Omega$  is 
of class $H^{5.5}$, and $u_0 \in H^{4.5}(\Omega)$ with $\div u_0=0$.  Then,
there exists $T>0$, and a solution ($v$,$q$) of
(\ref{leuler}) with $v \in L^\infty(0,T; H^{4.5}(\Omega))$,
$q \in L^\infty(0,T; H^{4}(\Omega))$,  
and $\Gamma(t) \in H^{5.5}$. The solution satisfies
$$
\sup_{t\in[0,T]}
\left( |\partial\Omega(t)|_{5.5}^2 +
 \sum_{k=0}^3 \|\partial_t^k v(t)\|_{4.5-1.5k}^2
+\sum_{k=0}^2 \|\partial_t^k q(t)\|_{4-1.5k}^2
\right) \le \tilde M_0
$$
where $\tilde M_0$ denotes a polynomial function of $\|\Gamma\|_{5.5}$ and
$\|u_0\|_{4.5}$.  The solution is unique of $u_0 \in H^{5.5}(\Omega)$ and
$\Gamma \in H^{6.5}$.
\end{theorem}


\begin{remark}
Our theorem is stated for a fluid in vacuum, but the analogous theorem holds for
a vortex sheet, i.e.,
for the motion of the interface separating two inviscid immiscible incompressible 
fluids;  the boundary condition (\ref{euler.c}) is replaced by 
$[p]_\pm = \sigma H$, where $[p]_\pm$ denotes the jump in pressure across the
interface.
\end{remark}

For the zero-surface-tension case, we have 

\begin{theorem}[$\sigma=0$ and condition (\ref{lindblad})]\label{ltheorem2}
Suppose that $\sigma=0$, $\Gamma $  is 
of class $H^{3}$, $u_0 \in H^{3}(\Omega)$, and condition
(\ref{lindblad}) holds at $t=0$.   Then,
there exists $T>0$, and a unique solution ($v$,$q$) of
(\ref{leuler}) with $v \in L^\infty(0,T; H^{3}(\Omega))$,
$q \in L^\infty(0,T; H^{3}(\Omega))$, and
$\Gamma(t) \in H^{3}$.
\end{theorem}
Because of the regularity of the solutions, 
Theorems \ref{ltheorem1} and \ref{ltheorem2} imply
Theorems \ref{theorem1} and \ref{theorem2}, respectively.

\begin{remark}  Note that in 3D, we require
less regularity on the initial data than \cite{Li2004}.
\end{remark}  

\begin{remark}  Since the vorticity satisfies the equation 
$\partial_t \operatorname{curl}u + \pounds _u \operatorname{curl} u=0$,
where $\pounds_u$ denotes the Lie derivative in the direction $u$, it follows
that if $\operatorname{curl}u_0=0$, then $\curl u(t)=0$.  Thus our result
also covers the simplified case of irrotational flow.
In particular, Theorem \ref{ltheorem1} shows that the 3D irrotational 
water-wave problem with 
surface  tension is well-posed.  In the zero surface tension case, our
result improves the regularity of the data required by Wu \cite{Wu1999}.
\end{remark}

\subsection{General methodology and outline of the paper}

\subsubsection{Artificial viscosity and the smoothed $\kappa$-problem}
Our methodology begins with the introduction of a smoothed or {\it approximate}
problem (\ref{smooth}), wherein two basic ideas are implemented: first, we 
smooth the transport velocity using a new tool which we call 
{\it horizontal convolution by layers}; second, we introduce an 
{\it artificial viscosity} term in the Laplace-Young boundary condition 
($\sigma >0 $) which simultaneously preserves the 
transport-type structure of the Euler equations, provides a PDE for which 
we can prove the existence of unique smooth solutions, and for which there exist
a priori estimates which are independent of the artificial viscosity 
parameter $\kappa$.  With the addition of the artificial viscosity term, the 
dispersive boundary condition is converted into a parabolic-type boundary
condition, and thus finding solutions of the smoothed problem becomes
an easier matter.  On the other hand, the a priori estimates for the $\kappa$
problem are more difficult than the formal estimates for the Euler equations.

The horizontal convolution is defined in Section \ref{sec_convolution}. The
domain $\Omega$ is partitioned into coordinate charts, each the image of
the unit cube in ${\mathbb R}^3$.  A double convolution is performed in
the horizontal direction  only (this is equivalent to the tangential direction
in coordinate patches near the boundary).  While there is no smoothing in
the vertical direction, our horizontal convolution commutes with the trace
operator, and avoids the need to introduce an extension operator, the latter
destroying the natural transport structure.   The development of the
horizontal convolution by layers is absolutely crucial in proving the
regularity of the weak solutions that we discuss below.  Furthermore,
it is precisely this tool which enables us to prove Theorem \ref{theorem2}
without the use of Nash-Moser iteration.  To reiterate, this horizontal
smoothing operator preserves the essential transport-type structure of the 
Euler equations.

\subsubsection{Weak solutions in a variational framework and a fixed-point,
$\sigma>0$}
The solution to the smoothed $\kappa$-problem (\ref{smooth}) is obtained via a topological
fixed-point procedure, founded upon the analysis of the linear problem 
(\ref{smoothlinear}).  To solve the linear problem, we introduce  a few
new ideas.  First, we penalize the pressure function; in particular,
with $\epsilon>0$ the penalization parameter, we introduce the penalized
pressure function $q_\epsilon ={\frac{1}{\epsilon}} \text{Tr}( a \, \nabla w)$.
Second, we find a new class of $[H^{\frac{3}{2}}(\Omega)]'$-weak solutions
of the penalized and linearized smoothed $\kappa$-problem in a variational
formulation.   The penalization allows
us to perform difference quotient analysis in order to prove regularity
of our weak solutions; without penalization, difference quotients of
weak solutions do not satisfy the ``divergence-free'' constraint and as such
cannot be used as test functions.  Furthermore, the penalization of the 
pressure function avoids the need to analyze the highest-order time-derivative
of the pressure, which would otherwise be highly problematic.  In the
setting of the penalized problem, we crucially rely on the horizontal
convolution by layers to establish regularity of our weak  penalized solution.
Third, we introduce the Lagrange multiplier 
lemmas, which associate a pressure function to the weak solution of a
variational problem  for which the test-functions satisfy the incompressibility 
constraint.  These lemmas allow us to pass to the limit as the penalization
parameter tends to zero, and thus, together with the Tychonoff fixed-point 
theorem, establish solutions to the smoothed problem (\ref{smooth}).
At this stage, however, the time interval of existence and the bounds for
the solution depend on the parameter $\kappa$.

\subsubsection{Solutions of the $\kappa$-problem for $\sigma=0$ via transport}
For the $\sigma=0$ problem, we use horizontal convolution to smooth the 
transport velocity as well as  the moving domain.   Existence and uniqueness of
this smoothed $\kappa$ problem (\ref{smoothl}) is found using simple 
transport-type arguments that rely on the pressure gaining regularity just as in
the fixed-domain case.  Once again, the time interval of existence and
the bounds for the solution a priori depend on $\kappa$.

\subsubsection{A priori estimates and $\kappa$-asymptotics}
We develop a priori estimates which show that the energy function $E_\kappa(t)$
 in Definition
\ref{kenergy} associated to our smoothed problem (\ref{smooth}) is
bounded by a constant depending only on the initial data and not on
$\kappa$.  The
estimates rely on the Hodge decomposition elliptic estimate (\ref{divcurl0}).

In Section \ref{sec_divcurl}, we obtain estimates for the divergence and
curl of $\eta$, $v$ and there space and time derivatives.  The main novelty
lies in the curl estimate for $\eta$.
The remaining portion of the energy is obtained by studying boundary
regularity via energy estimates.

These nonlinear boundary estimates for the surface tension case $\sigma>0$ 
are more
complicated than the ones for the $\sigma=0$ case with the Taylor sign condition 
(\ref{lindblad}) since it is necessary to analyze the time-differentiated Euler
equations, which is not essential in the $\sigma=0$ case (unless optimal
regularity is sought).


We note that the use of the smoothing operator
in Definition \ref{def_convolution}, where a double convolution is
employed, is necessary in order to find exact (or perfect) derivatives 
for the highest-order error terms.  The idea is that one of the
convolution operators is moved onto a function which is a priori not smoothed,
and commutation-type lemmas are developed for this purpose.

We obtain the a priori estimate
$$
\sup_{t\in[0,T]} E_\kappa(t) \le M_0 + T P(\sup_{t\in[0,T]}E_\kappa(t))\,, 
$$ 
where $M_0$  depends only on the data, and $P$ is a polynomial.
The addition of the artificial viscosity term allows us to prove that
$E_\kappa(t)$ is continuous; thus,
following the development in \cite{CoSh2005b}, there exists a sufficiently
small time $T$,
which is independent of $\kappa$, 
such that
$\sup_{t\in[0,T]}E_\kappa(t) < \tilde M_0$ for $\tilde M_0 > M_0$.  

We then find $\kappa$-independent nonlinear estimates for the $\sigma=0$
case for the energy function (\ref{Elin}).

\vspace{.1 in}

\noindent
{\bf Outline.} Sections \ref{sec_convolution}--\ref{uniqueness} are
devoted to the case of positive surface tension $\sigma >0$.  Sections
\ref{L1}--\ref{L13} concern the problem with zero surface tension $\sigma=0$
together with the Taylor sign condition (\ref{lindblad}) imposed.

\subsection{Notation} \label{1}

Throughout the paper, we shall use the Einstein convention with respect to repeated indices or exponents. 
We specify here our notation for certain vector and matrix operations.
\begin{itemize}
\item[] We write the Euclidean inner-product between two vectors $x$ and $y$ 
as $x\cdot y$, so that $x\cdot y=x^i\ y^i$.
\item[] The transpose of a matrix $A$ will be denoted by $A^T$, {\it i.e.}, 
$(A^T)^i_j=A^j_i$.
\item[] We write the product of a matrix $A$ and a vector $b$ as $A\ b$, 
{\it i.e}, $(A\ b)^i=A^i_j b^j$.
\item[] The product of two matrices $A$ and $S$ will be denoted by 
$A\cdot S$, {\it i.e.}, $(A\cdot S)^i_j=A^i_k\ S^k_j$.
\item[] The trace of the product of two matrices $A$ and $S$ will be denoted by 
$A : S$, {\it i.e.}, $(A: S)^i_j=A^i_k\ S^k_i$.
\end{itemize}

For $\Omega$, a domain of class $H^s$ ($s\ge 2$),  there exists a well-defined
extension operator that we shall make use of later.
\begin{lemma}
There exists $E(\Omega)$, a linear and continuous operator from 
$H^r(\Omega)$ into
$H^r({\mathbb R}^3)$ ($0\le r\le s$), such that for any 
$v\in H^r(\Omega)$ ($0\le r\le s$), $E(\Omega)(v)=v$ in $\Omega$.
\end{lemma} 

We will use the notation $H^s(\Omega)$ to denote either
$H^s(\Omega;\mathbb{R})$ (for a pressure function, for instance) 
or $H^s(\Omega;\mathbb{R}^3)$ (for a velocity vector field) and we 
denote the standard norm of $H^s(\Omega)$ ($s\ge 0$) by 
$\|\cdot\|_s$. 
The $H^s(\Omega)$ inner-product will be denoted $(\cdot,\cdot)_s$.

We shall use the following notation for derivatives:
$\partial_t$ or $(\cdot)_t$ denotes the partial time derivative, 
$\partial$ denotes
the tangential derivative on $\Gamma$ (or in a small enough neighborhood
of $\Gamma$), and $\nabla$ 
denotes the three-dimensional gradient.  

Letting $(x^1,x^2)$ denote a local coordinate system on $\Gamma$, for 
$\alpha=1,2$, we let either $\partial_\alpha$ or $(\cdot),_\alpha$
denote $\frac{\partial}{\partial x^\alpha}$.  We define:
$$
\partial^\alpha:= g_0^{\alpha\beta} \partial_\beta\,, \ \ 
|\partial^k \phi|^2 = 
\partial^{\alpha_1} \partial^{\alpha_2}\cdot\cdot\cdot\partial^{\alpha_k}
\partial_{\alpha_1} \partial_{\alpha_2}\cdot\cdot\cdot\partial_{\alpha_k}
$$
for integers $k\ge 0$, 
where $g_0 = g_{t=0}$ is the (induced) metric on $\Gamma$.
In particular, $|\partial^0 \phi| = |\phi|$, $|\partial^1 \phi|^2 = 
|\partial \phi|^2 = \partial^\alpha \phi \partial_\alpha\phi$.  $\partial^k\phi$
will mean any $k$th tangential derivative of $\phi$.

The area element on $\Gamma$ in local coordinates is 
$dS_0 = \sqrt{g_0} dx^1\wedge dx^2$ and the pull-back of the area element 
$dS$ on $\Gamma(t)=\eta(\Gamma)$ is given by $\eta^*(dS)=\sqrt{g} dS_0$.  
Let $\{U_i\}_{i=1}^K$ denote an open covering of $\Gamma$,
and let $\{\xi_i\}_{i=1}^K$ denote the partition of unity subordinate to this
cover.
The $L^2(\Gamma)$ norm is
$$
|\phi|_0:= \|\phi\|_{L^2(\Gamma)} =\left( \int_{\Gamma} \phi^2
dS_0 \right)^{\frac{1}{2}}\,,
$$
and the $H^k(\Gamma)$ norm  for integers $k\ge 1$ is
$$
| \phi|_k:= \|\phi\|_{H^k(\Gamma)} = 
\left(\sum_{i=1}^k \sum_{l=1}^K|\xi_l \partial^i\phi|^2_0 \right)^2\,.
$$
Similarly, for the Hilbert space inner-products, we use
$$
[\phi,\psi]_0:=[\phi,\psi]_{L^2(\Gamma)} = \int_\Gamma \phi\, \psi \, dS_0, \ \
\ \
[\phi,\psi]_k:=[\phi,\psi]_{H^k(\Gamma)} = 
[\phi,\psi]_0 + \sum_{i=1}^k \sum_{l=1}^K [\xi_l\partial^i \phi, \xi_l\partial_i \psi]_0 \,.
$$
Fractional-order spaces are defined via interpolation using the trace
spaces of Lions (see, for example, \cite{Adams1978}).

The dual of a Banach space $X$ is denoted by $X'$, and the corresponding norm
in $X'$ will be denoted $\|\cdot\|_{X'}$. For
$L\in H^s(\Omega)'$ and $v\in H^s(\Omega)$, the duality pairing between $L$ 
and $v$ is  denoted  by $\langle L,v\rangle_s$.

Throughout the paper, we shall use
$C$ to denote a generic constant, which may possibly depend on the coefficient
$\sigma$,  or on the initial geometry given by $\Omega$ (such as a Sobolev 
constant or an elliptic constant),
and we use $P(\cdot)$ to denote a generic polynomial function of $(\cdot)$.
For the sake of notational convenience, we will often write
$u(t)$ for $u(t,\cdot)$.

\section{Convolution by horizontal layers and the smoothed transport velocity}
\label{sec_convolution}

Let $\Omega \subset {\mathbb R}^n$ denote an open subset of class $H^6$, and
let $\{U_i\}_{i=1}^K$ denote an open covering of $\Gamma:= \partial \Omega$,
such that for each $i \in \{1,2,...,K\}$, 
\begin{gather}
\theta_i :(0,1)^2\times (-1,1) \rightarrow U_i
\ \text{ is an $H^6$ diffeomorphism}\,, \nonumber\\
U_i \cap \Omega = \theta_i ( (0,1)^3 )
\ \text{ and } \ U_i \cap \Gamma = \theta_i ( (0,1)^2\times \{0\} ) \,, 
\nonumber\\
\theta_i(x_1,x_2,x_3)=(x_1,x_2,\psi_i(x_1,x_2)+x_3)  \text{ and } 
\det\nabla \theta_i =1 \text{ in } (0,1)^3\,.
\nonumber
\end{gather}
Next, for $L > K$, let $\{U_i\}_{i=K+1}^{L}$ denote a family of open sets 
contained in $\Omega$ such that 
$\{U_i\}_{i=1}^{L}$ is an open cover of $\Omega$.  Let $\{\alpha_i\}_{i=1}^{L}$
denote the partition of unity subordinate to this covering.

Thus, each coordinate patch is locally represented by the unit cube
$(0,1)^3$ and for the first $K$ patches (near the boundary), the tangential
(or horizontal) direction is represented by $(0,1)^2 \times \{0\}$.

\begin{definition}[Horizontal convolution]\label{def_convolution}
Let $0\le\rho\in\mathcal D((0,1)^2)$ denote an even Friederich's mollifier,
normalized so that $\d\int_{(0,1)^2}\rho=1$,  with corresponding
dilated function 
$$\d\rho_{\frac{1}{\delta}}(x)=\frac{1}{\delta^2}\rho\left(\frac{x}{\delta}\right)\,,
\  \ \ \delta >0 .$$ 
For $w\in H^1((0,1)^3)$ such that
$\d\text{supp}(w)\subset [\delta,1-\delta]^2\times(0,1)$, set
\begin{align*}
\rho_{\frac{1}{\delta}}\star_h w(x_H,x_3)=\int_{\R^2} 
\rho_{\frac{1}{\delta}}(x_H-y_H) w(y_h,x_3) dy_H \,, \ \ 
y_H=(y_1,y_2)\,.
\end{align*}
\end{definition}

We then have the following tangential integration by parts formula
\begin{align*}
\rho_{\frac{1}{\delta}}\star_h w,_\alpha(x_H,x_3)=\int_{\R^2} 
\rho_{\frac{1}{\delta}},_\alpha(x_H-y_H) w(y_h,x_3) dy_H\,, \ \ \alpha=1,2\,,
\end{align*}
while
\begin{align*}
\rho_{\frac{1}{\delta}}\star_h w,_3(x_H,x_3)=\int_{\R^2} \rho_{\frac{1}{\delta}}(x_H-y_H) w,_3(y_h,x_3) dy_H \,.
\end{align*}
It should be clear that $\star_h$ smooths $w$ in the horizontal directions, 
but not in the vertical direction.
Fubini's theorem ensures that 
\begin{equation}\label{fubini}
\|\rho_{\frac{1}{\delta}}\star_h w\|_{s,(0,1)^3}\le C_s \|w\|_{s,(0,1)^3} 
\text{ for any } s\ge 0 \,, 
\end{equation}
and we shall often make implicit use of this inequality.

\begin{remark}
The horizontal convolution $\star_h w$ does not smooth $w$ in the 
vertical
direction, however, it does commute with the trace operator, so that
$$ \left.\left(\rho_{\frac{1}{\delta}}\star_h w\right)\right|_{(0,1)^2\times\{0\}} 
=\rho_{\frac{1}{\delta}}\star_h  \left.w\right|_{(0,1)^2\times\{0\}}  \,,
$$
which  is essential for our
methodology.  Also, note that $\star_h$
smooths without the introduction of an extension operator, required by 
standard convolution operators on bounded domains; the extension to the full
space would indeed be problematic for the transport structure of the 
divergence and curl of solutions to the Euler-type PDEs that we introduce.
\end{remark}

\begin{definition}[Smoothing the velocity field]
\label{smoothv}
For $v \in L^2(\Omega)$ and any $\ks \in (0,\frac{\ks_0}{2})$ with
$$
\ks_0 = \min_{i=1}^K \, \text{dist}\left(
\d\text{supp}(\alpha_i \circ \theta_i) \,, \, [(0,1)^2\times\{0\}]^c\cap\partial [0,1]^3
\right) \,,
$$
set
$$
v_\kappa = \sum_{i=1}^K \sqrt{\alpha_i} \left[ \rho_{ \frac{1}{\ks} }\star_h
[\rho_{ \frac{1}{\ks} } \star_h (( \sqrt{\alpha_i} v)\circ \theta_i)] \right]
\circ \theta_i^{-1} + \sum_{i=K+1}^L \alpha_i v \,.
$$
\end{definition}
It follows from (\ref{fubini}) that there exists a constant $C>0$ which
is independent of $\kappa$ such that for any $v \in H^s(\Omega)$ for $s\ge 0$,
\begin{equation}\label{neednow}
\|v_\kappa\|_s \le C \|v\|_s \ \text{ and } \ |v_\kappa|_{s-1/2} \le C |v|_{s-1/2} \,.
\end{equation}

The smoothed particle displacement field is given by
\begin{equation}\label{etak}
\eta_\kappa = \text{Id} + \int_0^t v_\kappa \,.
\end{equation}
For each $x \in U_i$, let $\tilde x = \theta^{-1}_i(x)$.  
The difference of the velocity field and its smoothed counterpart along
the boundary $\Gamma$ then takes the form
\begin{equation}\label{diffv}
v_\kappa(x) - v(x) =  \sum_{i=1}^K \int\int_{B(0,{\ks})^2} 
\zeta_i(x)
\rho_{\frac{1}{\ks}}(\ty)\rho_{\frac{1}{\ks}}(\tz)
\left[
(\zeta_i v)(\theta_i( \tilde x -(\tilde y+\tilde z) )) - (\zeta_i v)(\theta_i(\tilde x))
\right] d\tz \, d\ty \,,
\end{equation}
where
$\zeta_i(x)=\sqrt{\alpha_i(\theta_i(\tilde x))}$.
Combining (\ref{euler.a}), (\ref{etak}), and (\ref{diffv}), 
\begin{equation}\label{diffe}
\eta_\kappa(x) - \eta(x) =  \sum_{i=1}^K \int\int_{B(0,{\ks})^2} 
\zeta_i(x)
\rho_{\frac{1}{\ks}}(\ty)\rho_{\frac{1}{\ks}}(\tz)
\left[
(\zeta_i \eta)(\theta_i( \tilde x -(\tilde y+\tilde z) )) - (\zeta_i \eta)(\theta_i(\tilde x))
\right] d\tz \, d\ty \,,
\end{equation}
For any $u \in H^{1.5}(\Gamma)$, and for $y \in B(x, \ks)$, where
$B(x,\ks)$ denotes the disk of radius $\ks$ centered at $x$,
the mean value theorem shows that
$$
|u(y)-u(x)| \le C |r^{-1}|_{L^q(B(x,\ks))}|\partial u|_{L^p(B(x,\ks))},
\qquad r=\text{radial coordinate},
$$
so that in particular, with $p=4$ and $q=\frac{4}{3}$,
$$
|u(y)-u(x)| \le C \sk |\partial u|_{L^4} \le  C\kappa | u|_{1.5} \,,
$$
the last inequality following from the Sobolev embedding theorem.
Hence, for $U \in H^{1.5}(\Gamma)$,
\begin{equation}\label{Linf_est}
|U_\kappa(x) - U(x)|_{L^\infty} 
\le C \sk |U|_{1.5} \,.
\end{equation}
Note that the constant $C$ depends on $\max_{i\in\{1,...,K\}}|\theta_i|_{5.5}$.

Letting $\zeta_i = \sqrt{\alpha_i}$,  and $R=(0,1)^2$,
we also have that for any $\phi \in L^2(\Gamma)$,
\begin{align}
\int_\Gamma v_\kappa \, \phi
&
= \sum_{i=1}^K \int_R \rho_{\frac{1}{\ks}} \star_h \rho_{\frac{1}{\ks}} \star_h \zeta_i v(x)\,  \zeta_i \phi(x)
= \sum_{i=1}^K \int_R \rho_{\frac{1}{\ks}} \star_h \zeta_i v(x)\, \rho_{\frac{1}{\ks}}\star_h \zeta_i \phi(x) \nonumber \\
&\qquad\qquad
= \int_\Gamma \sum_{i=1}^K 
[\rho_{\frac{1}{\ks}} \star_h (\zeta_i v \circ \theta_i)] \circ \theta_i^{-1}\, 
[\rho_{\frac{1}{\ks}}\star_h (\zeta_i \phi \circ \theta_i)] \circ \theta_i^{-1} \,.
\label{selfadjoint}
\end{align}
Finally, we need the following
\begin{lemma}[Commutation-type lemma] \label{commutator}
Suppose that $g\in L^2(\Gamma)$ satisfies $\text{dist} (\text{supp}(g), 
\partial R) < \kappa_0)$ and that $f \in H^s(\Gamma)$ for $s> 1$.  Then
independently of $\ks \in (0,\kappa_0)$, there exists a constant $C>0$
such that
$$
\left|
\rho_{\frac{1}{\ks}} \star_h [fg]  - f \rho_{\frac{1}{\ks}} \star_h g 
\right|_{0,R} \le C \, \ks  |f|_{s+1,R} \, |g|_{0,R} \,.
$$
We also have
$$
\left\|
\rho_{\frac{1}{\ks}} \star_h [fg]  - f \rho_{\frac{1}{\ks}} \star_h g 
\right\|_{0,[0,1]^3} \le C \, \ks  \|f\|_{s+\frac{3}{2},[0,1]^3} \, 
\|g\|_{0,[0,1]^3} 
$$
whenever $g\in L^2(\Omega)$, $f\in H^s(\Omega)$ and
$$\frac{\kappa}{2}<\min(\hbox{dist}({supp}\ fg, \{1\}\times [0,1]^2), \text{dist}(\text{supp}\ fg, \{0\}\times [0,1]^2)).$$
\end{lemma}
\begin{proof}
Let $\triangle=
\rho_{\frac{1}{\ks}} \star_h [fg]  - f \rho_{\frac{1}{\ks}} \star_h g $. Then
\begin{align*}
|\triangle (x)| &= \left|\int_{B(x,\ks)} \rho_{\frac{1}{\ks}}(x -y) [f(y)-f(x)]g(y) dy
\right| 
\le C\, \ks |f|_{s+1,R} \int_{B(x,\ks)} \rho_{\frac{1}{\ks}}(x-y) |g(y)| dy\,,
\end{align*}
so that
\begin{align*}
|\triangle |_{0,R} & \le C\, \ks |f|_{s+1,R} \, 
\left| \rho_{\frac{1}{\ks}}\star_h|g|\, \right|_{0,R} 
 \le C\, \ks |f|_{s+1,R} \, | g |_{0,R} \,.
\end{align*}

The inequality on $[0,1]^3$ follows the identical argument with an additional
integration over the vertical coordinate.  The hypothesis on the support of
$fg$ makes the integral well-defined.
\end{proof}
\begin{remark}
Higher-order commutation-type lemmas will be developed for the case of
zero surface tension in Section \ref{L7}.
\end{remark}

\section{Closed convex set used for the fixed-point for $\sigma >0$} \label{2}

In order to construct solutions for  our approximate model (\ref{smooth}), we
use a topological fixed-point argument which necessitates the use of high-regularity
Sobolev spaces.  In particular, we shall assume that the initial velocity
$u_0$ is in $H^{13.5}(\Omega)$ and that $\Omega$ is of class $C^\infty$; after establishing
our result for the smoothed initial domain and velocity, we will show that both $\Omega$ and
$u_0$ can be taken with the optimal regularity stated in Theorem \ref{ltheorem1}.

For $T>0$, we define the following closed convex set of the Hilbert space $L^2(0,T;H^{13.5}(\Omega))$:  
\begin{align*}
C_T=\{ v \in L^2(0,T;H^{13.5}(\Omega))| \  \sup_{[0,T]} \|v\|_{13.5}\le 2 \|u_0\|_{13.5}+1 \},
\end{align*}
 It is clear that $C_T$ is  non-empty, since it contains the constant (in time)
 function $u_0$, and is a convex, bounded and closed subset of the separable 
Hilbert space $L^2(0,T;H^{13.5}(\Omega))$.  

Let $v\in C_T$ be given, and define $\eta$ by (\ref{leuler.a}),
the Bochner integral being taken in the separable Hilbert space 
$H^{13.5}(\Omega)$. 
 
Henceforth, we assume that $T>0$ is given such that
independently of the choice of $v\in C_T$, we have the injectivity of 
$\eta(t)$ on $\overline\Omega$, the existence of a normal vector to $\eta(\Omega,t)$
at any point of $\eta(\Gamma,t)$, and the invertibility of $\nabla\eta(t)$ for 
any point of $\overline\Omega$  and for any $t\in [0,T]$. Such a condition can be 
achieved by selecting $T$ small enough so that
\begin{align}
\|\nabla\eta-\text{Id}\|_{L^\infty(0,T;H^{13.5}(\Omega))}&\le \epsilon_0\,,
\label{eta} 
\end{align}
for $\epsilon_0>0$ taken sufficiently small. Condition (\ref{eta}) holds
if $T\|\nabla u_0\|_{H^{2}}\le\epsilon_0$.
Thus,
\begin{equation}\label{aequation}
a=[\nabla\eta]^{-1}
\end{equation}
is well-defined.

Then choosing $T>0$ even smaller, if necessary, there  exists $\kappa_0>0$ such
that for any $\kappa\in (0,\frac{\kappa_0}{2})$, we have the injectivity of 
$\eta_\kappa(t)$ on $\Omega$ for any $t\in [0,T]$; furthermore, 
$\nabla\eta_\kappa$ satisfies the condition (\ref{eta}) with $\eta_\kappa$
replacing $\eta$. 
We let
$n_\kappa(\eta_\kappa(x))$ denote the exterior unit normal to $\eta_\kappa(\Omega)$ 
at $\eta_\kappa(x)$ with $x\in\Gamma$.  

Our notational convention will be as follows: if we choose $\bar v \in C_T$,
then $\bar \eta$ is the flow map coming from (\ref{leuler.a}), and $\bar a$
is the associated pull-back, $\bar a = [\nabla\bar \eta]^{-1}$.  
Thus, a bar over the velocity field, will
imply a bar over the Lagrangian variable and the associated pull-back.

For a given $v_\kappa$, our notation is as follows:
\begin{gather*}
\eta_\kappa(t)=
\operatorname{Id}+ \int_0^t v_\kappa\ \text{ and } \ \eta_\kappa(0)=\operatorname{Id} \,, \\
a_\kappa = \operatorname{Cof}  \nabla \eta_\kappa\,, \ \ J_\kappa = 
\det \nabla \eta_\kappa \,, \ \  {g_\kappa}_{\alpha\beta} =
\partial_\alpha\eta_\kappa \cdot \partial_\beta \eta_\kappa\,.
\end{gather*}

We take $T$ (which a priori depends on $\kappa$) 
even smaller if necessary to ensure that for $t\in[0,T]$,
\begin{subequations}
\label{deteta}
\begin{align}
\sqrt{ g(t)}^{-1} &\le 2 \sqrt{ g_0}^{-1} \,, \label{deteta.a} \\
\sqrt{ g_\kappa(t)}^{-1} &\le 2 \sqrt {g_0}^{-1} \,, \label{deteta.b}\\
{\frac{1}{2}}\le J_\kappa(t) &\le {\frac{3}{2}} 
\,.\label{deteta.c}
\end{align}
\end{subequations}
 
\begin{lemma} \label{smoothvbisk}
For $v\in C_T$, 
and for any $s\ge 0$, we have independently of the choice of $v\in C_T$ that
\begin{align*}
\sup_{[0,T]}|v_\kappa|_{s} \le C_{\kappa,s} \, P(\|u_0\|_{13.5})\,.
\end{align*}
\end{lemma}
\begin{proof}
By the standard properties of the convolution a.e in $[0,T]$:
\begin{equation}
\label{bis1}
|v_\kappa|_s\le [\frac{C}{\kappa^{s-13}}+1] |v|_{13}\le 
[\frac{C}{\kappa^{s-13}}+1] [2 \|u_0\|_{13.5}+1],
\end{equation}
where we have used the definition of $C_T$ for the second inequality.
\end{proof}

Recall that $\{\theta_i\}_{i=1}^K$ is our open cover of $\Gamma$.
Given $\bar v \in C_T$, we define the
matrix $\bbk=[\nabla(\ek\circ\theta_l)]^{-1}$, and 
assume that $T>0$ is sufficiently small so that independently of 
$\bar v\in C_T$, 
we have the following determinant-type condition for $\bbk$:
\begin{equation}
\label{matrix}
\ud\le (\bbk)_3^3\sum_{i=1}^3 [(\bbk)_i^3]^2,\ \text{in}\ (0,1)^3.
\end{equation}
Such a condition is indeed possible since at time $t=0$ we have 
$(\bbk)_3^3\sum_{i=1}^3 [(\bbk)_i^3]^2=1+\psi_l,_1^2+\psi_l,_2^2.$

\section{The smoothed $\kappa$-problem and its linear fixed-point formulation}
\label{3}

Unlike the case of  zero surface tension, for $\sigma>0$ there does not 
appear to be a simple sequence of approximate problems for the Euler equations 
(\ref{euler}) which can be solved only with simple transport-type
arguments. For the surface tension case, the problem is crucially variational 
in nature, and the addition of an artificial viscosity term on the boundary 
$\Gamma$ seems unavoidable in order to be able to construct a sequence
of approximate or smoothed solutions.

As we shall make precise below, our construction of the approximating sequence
of problems is based on smoothing the transport velocity by use
of the horizontal convolution by layers (see Definition \ref{smoothv}), 
and hence
smoothing the Lagrangian flow map and associated pull-back.  Simultaneously,
we introduce a new type of 
parabolic-type artificial viscosity boundary operator  on $\Gamma$ 
(of the same order in space as the surface tension operator).
Note that unlike the case of interface motion in the fluid-structure interaction
problem
that we studied in \cite{CoSh2005b}, there is not a unique choice of
the artificial viscosity term; in particular, 
other choices of artificial viscosity are possible for the asymptotic limit
as the artificial viscosity is taken to zero. 

We can now define our sequence of smoothed $\kappa$-problems. For our artificial
viscosity parameter $\kappa \in (0, \frac{\kappa_0}{2})$, let $(v,q)$ be the solution
of
\begin{subequations}
  \label{smooth}
\begin{alignat}{2}
\eta&=\text{Id} + \int_0^t v\ \ \   &&\text{in} \ \Omega \times (0,T]\,, \label{smooth.a}\\
\partial_t v+ J_\kappa^{-1}a_\kappa\,\nabla q&=0  &&\text{in} \ \Omega \times (0,T]\,, \label{smooth.b}\\
   \text{Tr}(a_\kappa \, \nabla v)  &= 0     &&\text{in} \   \Omega \times (0,T]
\,, \label{smooth.c}\\
-\sigma \,\frac{\sqrt{g}}{\sqrt{g_\kappa}} \Delta_g(\eta)\cdot n_k(\eta_\kappa) \, n_k(\eta_\kappa)
- \kappa \Delta_0[ v \cdot n_\kappa(\eta_\kappa)]n_\kappa(\eta_\kappa)
&= 
q\, n_\kappa(\eta_\kappa)
\ \ &&\text{on} \  \Gamma \times (0,T] \,, \label{smooth.d}\\
   (\eta,v) &= (\text{Id}, u_0)  &&\text{on} \ \Omega\times\{t=0\} 
                                                 \,, \label{smooth.e} 
\end{alignat}
\end{subequations}
where $n_\kappa(\eta_\kappa) = \frac{(a_\kappa)^T N}{|(a_\kappa)^T N|}$, and
$\Delta_0= \sqrt{g_\kappa}^{-1} \partial_\alpha [\sqrt{g_0} g_0^{\alpha\beta} 
\partial_\beta ]$. Note that on $\Gamma$, $ \sqrt{g_\kappa} = 
|a_\kappa^TN|$, and that $(g_\kappa)_{\alpha\beta} = \eta_\kappa,_\alpha 
\cdot \eta_\kappa,_\beta$.

In order to obtain solutions to the sequence of approximate $\kappa$-problems 
(\ref{smooth}), we study a linear problem whose fixed-point will provide 
the desired solutions. If we denote by $\bar v$ an arbitrary element of $C_T$, 
and $\ek$, $\bak$, and $\bar J_\kappa$ are  the associated smoothed Lagrangian 
variables given by Definition \ref{smoothv}, then we 
define $w$ to be the solution of
\begin{subequations}
  \label{smoothlinear}
\begin{alignat}{2}
\partial_t w+ \bar J_\kappa^{-1}\bar a_\kappa\,\nabla q&=0  &&\text{in} \ \Omega \times (0,T]\,, \label{smoothlinear.b}\\
   \text{Tr}(\bar a_\kappa \, \nabla w)  &= 0     &&\text{in} \   \Omega \times (0,T]
\,, \label{smoothlinear.c}\\
-\sigma \frac{\sqrt{\bar g}}{\sqrt{\bar g_\kappa}}
[\Delta_{\bar g}(\bar\eta)\cdot \nk(\ek)] \, \nk(\ek)
- \kappa \Delta_{\bar 0}[ w \cdot \nk(\ek)]\, \nk(\ek)
&= 
q\, \nk(\ek)
\ \ &&\text{on} \  \Gamma \times (0,T] \,, \label{smoothlinear.d}\\
   (\eta,w) &= (\text{Id}, u_0)  &&\text{on} \ \Omega\times\{t=0\} 
                                                 \,, \label{smoothlinear.e} 
\end{alignat}
\end{subequations}
where $\bar g_{\alpha\beta}= \bar\eta,_\alpha \cdot \bar\eta,_\beta$, and
$\Delta_{\bar 0}= \sqrt{\bar g_\kappa}^{-1} 
\partial_\alpha [\sqrt{g_0} g_0^{\alpha\beta} \partial_\beta ]$.

For a solution $w$ to (\ref{smoothlinear}),
a fixed point of the map $\bar v \mapsto w$ provides a solution of our
smoothed problem
(\ref{smooth}). 

In the following sections, we assume that $\bar v\in C_T$ is given, 
and $\kappa$ is in $(0,\frac{\kappa_0}{2})$. 
Until Section \ref{sec_divcurl},
wherein we study the asymptotic behavior of the problem (\ref{smooth}) as 
$\kappa\rightarrow 0$, the parameter $\kappa$ is fixed.

\section{Hodge decomposition elliptic estimates}
Our estimates are based on the following standard elliptic estimate:
\begin{proposition}\label{prop1}
For an $H^r$ domain $\Omega$, $r \ge 3$, 
if $v \in L^2(\Omega)$ with $\operatorname{curl}v \in H^{s-1}(\Omega)$,
${\operatorname{div}}v\in H^{s-1}(\Omega)$, and $v \cdot N|_{\Gamma} \in
H^{s -{\frac{1}{2}}}(\Gamma)$ for $1 \le s \le r$, then there exists a
constant $C>0$ depending only on $\Omega$ such that
\begin{equation}
\begin{array}{c}
\|v\|_s \le C\left( \|v\|_0 + \|\operatorname{curl} v\|_{s-1}
+ \|\operatorname{div} v\|_{s-1} + |v \cdot N|_{s-{\frac{1}{2}}}\right)\,, \\
\|v\|_s \le C\left( \|v\|_0 + \|\operatorname{curl} v\|_{s-1}
+ \|\operatorname{div} v\|_{s-1} + |v \cdot T_\alpha|_{s-{\frac{1}{2}}}\right)\,,
\end{array}
\label{divcurl0}
\end{equation}
where $T_\alpha$, $\alpha=1,2$ are the tangent vectors to $\Gamma$.
\end{proposition}
The first estimate with $V \cdot N$ is standard (see, for example, 
\cite{Taylor1996}), while the second with
$V\cdot T_\alpha$ follows from the fact that $T_\alpha \cdot N=0$.

\section{Weak solutions for the penalized problem and their regularity} 
\label{4}

The aim of this section is 
to establish the existence of the solution $w_{\epsilon}$ to the penalized 
version (of the divergence-free condition) of linearized and smoothed
$\kappa$-problem (\ref{smoothlinear}).  In particular, we study the weak
form of this problem with the pressure function $q$, approximated by the
penalized pressure 
$$q^\epsilon =-{\frac{1}{\epsilon}} \text{Tr}( \bar a_\kappa \nabla w)
\ \ \text{ for }  0 < \epsilon <<1\,.$$

In this section, as well as in Sections \ref{6} and \ref{7}, we let
\begin{equation} \label{Nu0}
N(u_0,x,y) = P(\|u_0\|_{13.5}, x,y)
\end{equation}
denote a generic polynomial function of $\|u_0\|_{13.5}$, $x$, and $y$,
where $x$ and $y$ will typically denote norms of various quantities.

\subsection{Step 1. Galerkin sequence.}

By introducing a basis $(e_l)_{l=1}^{\infty}$ of $H^1 (\Omega)$ 
and $L^2(\Omega)$, and taking the approximation at rank 
$l\ge 2$ under the form $\displaystyle w_l (t,x)
=\sum_{k=1}^l y_k (t)\ e_k (x)\ ,$ satisfying on $[0,T]$, the system of
ordinary differential equations
\begin{align*}
&\text{(i)}\  (\bar J_\kappa\ {w_l}_{t}, \phi)_{0}
+\kappa [w_l\cdot\nk(\ek),\phi\cdot\nk(\ek)]_1 
-\sigma [L_{\bar g}\bar\eta\cdot\nk(\ek),\phi\cdot\nk(\ek)]_0 \n\\
&\qquad\qquad - ((\bar a_\kappa)_i^j q_l,\ \phi^i,_j)_{0}=0, \ \forall \phi\in 
span(e_1,...,e_l)\,,\\
&\text{(ii)}\  w_l (0)=(u_0)_l,\ \text{in}\ \Omega\ ,
\end{align*}
where
$L_{\bar g} = \frac{\sqrt{\bar g}}{\sqrt{g_0}}  \Delta_{\bar g}$,
$\displaystyle q_l= -\frac{1}{\epsilon} (\bak)_i^j w_l^i,_j$,
and 
$(u_0)_l$ denotes the $L^2(\Omega)$ 
projection of $u_0$ on $span(e_1,...,e_l)$, we see that the 
Cauchy-Lipschitz theorem gives us the local well-posedness for $w_l$ on some $[0,T_{max}]$. The use of
the test function $w_l$ in this system of ODEs (which is allowed as
it belongs to $span(e_1,...,e_l)$) gives us in turn the energy law 
for any $t\in(0,T_{\max})$,
\begin{align*}
  \frac{1}{2}  \|{\bar J_\kappa}^{\frac{1}{2}} {w_l}(t)\|^2_{0} + 
\kappa \int_0^t &[w_l\cdot\nk(\ek),w_l\cdot\nk(\ek)]_1
+ \epsilon \int_0^t\|{q_l}\|_0^2 \n\\
-\ud\int_0^t ((\bar J_\kappa)_t w_l, w_l)_0& = \frac{1}{2}  \|{(u_0)_l}\|^2_{0} +\sigma\int_0^t 
[L_{\bar g}\bar\eta\cdot\nk(\ek),w_l\cdot\nk(\ek)]_0,
\end{align*}
which, with the control of $\bar\eta^\kappa$ provided by the definition of 
$C_T$, gives the bound 
\begin{align}
  \frac{1}{4}  \|{w_l}(t)\|^2_{0} &+ 
C\kappa \int_0^t |w_l\cdot\nk(\ek)|_1^2+ \epsilon \int_0^t\|{q_l}\|_0^2 \le C N(u_0).
\label{g1}
\end{align}

\subsection{Step 2. Weak solution $w_{\epsilon}$ of the penalized problem.}

We then infer from (\ref{g1}) that $w_l$ is defined on $[0,T]$, and that 
there is a subsequence (still denoted with the subscript $l$) satisfying
\begin{subequations}
\label{g2}
\begin{align}
w_{l}&\rightharpoonup   w_{\epsilon}\ \ \text{in}\ L^2(0,T; L^2(\Omega)),
\\
{q_l}&\rightharpoonup {q_{\epsilon}}\ \ \text{in}\ L^2(0,T; L^2(\Omega))\ ,
\end{align}
\end{subequations}
where 
\begin{equation}
q_{\epsilon}=-\frac{1}{\epsilon} (\bak)_i^j {w_{\epsilon}},_j^i\ .
\end{equation}
We can also rewrite (\ref{g2}) as
\begin{subequations}
\label{g3}
\begin{align}
w_{l}&\rightharpoonup   w_{\epsilon}\ \ \text{in}\ L^2(0,T; L^2(\Omega)),
\\
\operatorname{div}( w_l\circ\ek^{-1})(\ek)&\rightharpoonup \operatorname{div} (w_\epsilon\circ\ek^{-1})(\ek)\ \ \text{in}\ L^2(0,T; L^2(\Omega))\ \label{g3.b},
\end{align}
\end{subequations}
which with the bound (\ref{g1}) and the definition of the normal $\nk$ provides
\begin{equation}
\label{g4}
w_l\cdot \nk(\ek)\rightharpoonup w_{\epsilon}\cdot \nk(\ek)\ \ \text{in}\ L^2(0,T; H^1(\Gamma)).
\end{equation}

It follows from standard arguments and the ODE defining $w_l$, that 
$\we_t\in L^2(0,T;H^\td(\Omega)')$, $\we\in \mathcal{C}^0([0,T];H^\td(\Omega)')$
with $\we(0)=u_0$, and that for $\phi\in L^2(0,T;H^\td(\Omega))$,
\begin{align}
  \int_0^T \langle\bar J_\kappa\, \we_{t},\phi\rangle_{\td}  
&+\kappa \int_0^T [\partial ({\we}\cdot\nk(\ek)),\partial (\phi\cdot\nk(\ek))]_0 
\n\\
&-\int_0^T \langle \qe,(\bak)_i^j \phi^i,_j\rangle_\ud=  \sigma \int_0^T [L_{\bar g}
\bar\eta\cdot\nk(\ek),\phi\cdot\nk(\ek)]_0 \ .
\label{wepsilonvar}
\end{align}
Since by definition $\bak=\text{Cof}\nabla\bar\eta_\kappa$, this implies that 
in $\Omega$, 
\begin{equation} \label{transport}
\we_t+\nabla\pe(\ek)=0,
\end{equation}
where $\pe\circ\ek=\qe$ in $\Omega$. Since $\nabla\pe(\ek)\in L^2(0,T;H^{-1}(\Omega))$, this equality is true in $L^2(0,T;H^{-1}(\Omega))$ as well.

\subsection {Step 3. $w_\epsilon$ is bounded in $L^2(0,T;H^1(\Omega))$ 
independently of $\epsilon$}

Denoting $\ue=\we\circ\ek^{-1}$, by integrating (\ref{transport}) in time
from $0$ to $t$, we obtain the important formula
\begin{equation}
\label{g5}
\d \curl\ue(\ek)=\curl u_0+\int_0^t B(\uk,\ue)(\ek)  \ \ \ \text{ in }
L^2(0,T;H^{-1}(\Omega)) \,,
\end{equation}
with
\begin{align*}
B(\uk,\ue)&=
-(\uk^i,_2 \ue_3,_i -\uk^i,_3\ue_2,_i,\ \uk^i,_3 \ue_1,_i -\uk^i,_1\ue_3,_i,\ \uk^i,_1 \ue_2,_i -\uk^i,_2\ue_1,_i).
\end{align*}

\begin{remark}
Note well that our approximated and penalized $\kappa$-problem 
preserves the structure of the original Euler equations as can be seen by
(\ref{transport}).
As a  result, (\ref{g5}) contains only first-order derivatives of the 
velocity.
\end{remark}

Our next task is to prove that $w_\epsilon$ in $L^2(0,T;H^1(\Omega))$.  For
suppose that this was the case; then, (\ref{g5}) together with bounds
on the divergence of $w_\epsilon$ and $w_\epsilon \cdot N$ on $\Gamma$, provide
bounds for $w_\epsilon$ in $L^2(0,T;H^1(\Omega))$ 
(by the Hodge elliptic estimate (\ref{divcurl0})) which are independent of
$\epsilon >0$.

We proceed by showing that appropriately convolved velocity fields are bounded 
independently of the parameter of convolution in $L^2(0,T;H^1(\Omega))$.
This is the first instance that our horizontal convolution by layers is crucially 
required. 

\subsubsection{ For any subdomain $\omega \subset\subset \Omega$, 
$w_\epsilon \in L^2(0,T;H^1(\omega))$}\label{6.3.1}

We analyze the third component of (\ref{g5}), the other components being 
treated similarly. This leads us to the following equality in 
$L^2(0,T;H^{-1}(\Omega))$:
\begin{align*}
(\bak)_2^j \we,_j^1-(\bak)_1^j \we,_j^2=-\curl u_0^3+\int_0^t [\ -\bar v,_j^i (\bak)_2^j\we,_l^1(\bak)_i^l+\bar v,_j^i (\bak)_1^j\we,_l^2(\bak)_i^l] .
\end{align*}

Our goal is to prove that $w_\epsilon \in H^1(\Omega)$.  To proceed,
we let $\sigma_p$ denote a standard sequence of Friederich's mollifier 
in ${\mathbb R}^3$ with support $B(0, 1/p)$,
and establish that  $\sigma_p \star w_\epsilon$ is bounded in
$H^1(\omega)$ for any $\omega \subset\subset \Omega$.
For this purpose, we choose 
$\psi\in\mathcal{D}(\Omega)$, and find that
\begin{align}
(\bak)_2^j (\psi\we),_j^1-(\bak)_1^j (\psi\we),_j^2=&-\psi\curl u_0^3+(\bak)_2^j \psi,_j\we^1-(\bak)_1^j \psi,_j\we^2\n\\
&+\int_0^t [\ -\bar v,_j^i (\bak)_2^j(\psi\we),_l^1(\bak)_i^l+\bar v,_j^i (\bak)_1^j(\psi\we),_l^2(\bak)_i^l] \n\\
&-\int_0^t [\ -\bar v,_j^i (\bak)_2^j\psi,_l\we^1(\bak)_i^l+\bar v,_j^i (\bak)_1^j\psi,_l\we^2(\bak)_i^l].
\label{g6}
\end{align}

In order to proceed, we shall need to identify curl-type structures
(in Lagrangian variables) for $\sigma_p \star w_\epsilon$; this requires
the following:
for $\d\frac{1}{p}\le \text{dist}(\text{supp}\psi,\Omega^c)$, 
and $f \in C^\infty(\Omega)$, we have the equality in $H^{-1} (\Omega)$
\begin{align*}
\sigma_p\star [f(\psi\we),_j]-f \sigma_p\star [(\psi\we),_j]=\int_{\R^3} (\sigma_p),_j (x-y) (f(y)-f(x)) \psi \we(y) dy\\
-\int_{\R^3} \sigma_p (x-y) f,_j(y) \psi \we(y) dy\,,
\end{align*}
showing that $\sigma_p\star [f(\psi\we),_j]-f\ 
\sigma_p\star [(\psi\we),_j]\in L^2(\Omega)$, with 
\begin{align}
\bigl\|\sigma_p\star [f(\psi\we),_j]-f\ \sigma_p\star [(\psi\we),_j]\bigr\|_0\le  C [\ \|\sigma,_j\|_{0,\R^3}+\|\sigma\|_{0,\R^3}] \|\nabla f\|_{L^\infty(\Omega)} \|\we\|_0.
\label{g7}
\end{align}
We thus infer from (\ref{g6}) and (\ref{g7}) that the vorticity structure 
satisfies
\begin{align}
(\bak)_2^j \sigma_p\star(\psi\we),_j^1-(\bak)_1^j \sigma_p\star(\psi\we),_j^2=&
\int_0^t [\ -\bar v,_j^i (\bak)_2^j\sigma_p\star(\psi\we),_l^1(\bak)_i^l\n\\
&+\bar v,_j^i (\bak)_1^j\sigma_p\star(\psi\we),_l^2(\bak)_i^l] +R_1,
\label{g8}
\end{align}
with $\|R_1\|_{L^2(0,T;L^2(\Omega))}\le  N(u_0)$, where $N(u_0)$ is defined
in (\ref{Nu0}).
Next, we infer from (\ref{g3}) and (\ref{g7}) that the divergence structure satisfies
\begin{align}
(\bak)_i^j \sigma_p\star(\psi\we),_j^i=R_2,
\label{g9}
\end{align}
with $\|R_2\|_{L^2(0,T;L^2(\Omega))}\le  N(u_0).$ Since we also have 
$\psi\we=0$ on $\Gamma$,  so that with
(\ref{divcurl0}), we have  a.e. in $(0,T)$
\begin{align*}
\|\sigma_p\star(\psi\we)(t)\|_1
\le  \|R_1(t)\|_0+\|R_2(t)\|_0+N(u_0)\int_0^t \|\sigma_p\star(\psi\we)\|_1 \,,
\end{align*}
and thus
\begin{align}
\int_0^T \|\sigma_p\star(\psi\we)\|_1^2\le N(u_0).
\label{g10}
\end{align} Since this inequality does not depend on $p$, this implies that
$\psi \we\in L^2(0,T;H^1(\Omega))$, and therefore $\we\in L^2(0,T;H^1(\omega))$,
with an estimate depending {\it a priori} on $\omega\subset\subset \Omega$.

\subsubsection{ The horizontal convolved-by-layers velocity fields 
are in $L^2(0,T;H^1(\Omega))$}\hf
Fix $l\in \{1,...,K\}$, and set 
$$W(l)=w_\epsilon\circ\theta_l \ \text{ and } \
\bbk=[\nabla({\bar\eta}_\kappa\circ\theta_l)]^{-1}.$$ 
Hence, in $(0,1)^2\times(\frac{1}{p},1)$ for $p>1$, 
the Lagrangian ``divergence-free'' 
constraint is given by
\begin{align}
(\bbk)_i^j (\alpha_l W(l)),_j^i=-(\bbk)_i^j \alpha_l,_j\ W(l)^i - 
\alpha_l\epsilon\,q_\epsilon(\theta_l)\,, \label{divcond}
\end{align}
where the crucial observation is that the right-hand side of (\ref{divcond})
is in $L^2(0,T; L^2([0,1]^3))$.

Now for $\d\frac{1}{m}\le \text{dist}(\text{supp}\alpha_l,\partial (0,1)
\times(0,1)^2)$, and $f$ smooth in $[0,1]^3$, we have by 
Lemma \ref{commutator} that for
$\beta=1,2$,
$\rho_m\star_h [f(\alpha_l W(l)),_\beta]-f \rho_m\star_h 
[(\alpha_l W(l)),_\beta]\in L^2((0,1)^2)$, with the estimate a.e in $\d(\frac{1}{p},1)$:
\begin{align*}
\bigl|\rho_m\star_h [f(\alpha_l W(l)),_\beta]- f\ \rho_m\star_h &[(\alpha_l W(l)),_\beta]\bigr|_{0,(0,1)^2\times\{y\}}\n\\
&\le  C_\rho |\nabla f|_{L^\infty((0,1)^2\times\{y\})} \ |W(l)|_{0,(0,1)^2\times\{y\}}.
\end{align*}
This leads to
\begin{align}
\bigl\|\rho_m\star_h [f(\alpha_l W(l)),_\beta]- f\ \rho_m\star_h [(\alpha_l W(l)),_\beta]&\bigr\|_{0,(0,1)^3}\n\\
&\le  C_\rho \|\nabla f\|_{L^\infty((0,1)^3)} \|W(l)\|_{0,(0,1)^3}.
\label{g12}
\end{align}
Now, for the case of the vertical derivative, we will need to express 
$W(l),_3$ in terms of $W(l),_1$, $W(l),_2$, $\curl_{\bar\eta_\kappa\circ\theta_l} W(l)$ and $\div_{\bar\eta_\kappa\circ\theta_l} W(l)$, where
\begin{align*}
\div_{\bar\eta_\kappa\circ\theta_l} W(l)&=(\bbk)_i^j W(l),_j^i,\\
\curl_{\bar\eta_\kappa\circ\theta_l}^1 W(l)&=(\bbk)^i_2 W(l),_i^3-(\bbk)^i_3 W(l),_i^2,\\
\curl_{\bar\eta_\kappa\circ\theta_l}^2 W(l)&=(\bbk)^i_3 W(l),_i^1-(\bbk)^i_1 W(l),_i^3,\\
\curl_{\bar\eta_\kappa\circ\theta_l}^3 W(l)&=(\bbk)^i_1 W(l),_i^2-(\bbk)^i_2 W(l),_i^1.
\end{align*}
Notice that the first three lines above can be written as the following
vector field:
\begin{align*}
(\div_{\bar\eta_\kappa\circ\theta_l} W(l), \curl_{\bar\eta_\kappa\circ\theta_l}^1 W(l), \curl_{\bar\eta_\kappa\circ\theta_l}^2 W(l))=\sum_{i=1}^3 M_i^\kappa W(l),_i,
\end{align*}
where the $M_i^\kappa$ are smooth matrix fields depending on $\bbk$. 
From condition (\ref{matrix}), since 
$$\det M_3^\kappa= (\bbk)_3^3\ \sum_{i=1}^3 [(\bbk)_i^3]^2\ge \ud,$$ 
we see that $M_3^\kappa$ is invertible on $[0,T]$ 
(regardless of the choice of $\bar v\in C_T$). Therefore,
\begin{align}
W(l),_3=\div_{\bar\eta_\kappa\circ\theta_l} W(l)\ V^\kappa+M^\kappa \curl_{\bar\eta_\kappa\circ\theta_l} W(l)+\sum_{i=1}^2 A_i^\kappa W(l),_i,
\label{g13}
\end{align}
where $M^\kappa$ and the $A_i^\kappa$ are smooth matrix fields depending on 
$\bbk$, and $V^\kappa$ is a vector field depending on $\bbk$. From (\ref{g6}), 
we have that
\begin{align*}
\curl_{\bar\eta_\kappa\circ\theta_l} W(l)=\curl u_0(\theta_l) +\sum_{i=1}^3 \int_0^t N_i^\kappa W(l),_i.
\end{align*}
where the $N_i^\kappa$ are smooth matrix fields depending on $\bbk$.
By using (\ref{g13}) and the fact that 
$\div_{\bar\eta_\kappa\circ\theta_l} W(l)\in L^2(0,T; L^2(\Omega))$ from
(\ref{g3}), we obtain after time differentiating that
\begin{align*}
[\curl_{\bar\eta_\kappa\circ\theta_l} W(l)]_t-N_3^\kappa M^\kappa \curl_{\bar\eta_\kappa\circ\theta_l} W(l)=\sum_{\beta=1}^2 P_\beta^\kappa W(l),_\beta,
+ N^\kappa_3 V^\kappa
\div_{\bar\eta_\kappa\circ\theta_l} W(l)
\end{align*}
where $P_\beta^\kappa$, $\beta=1,2$, are smooth matrix fields depending on $\bbk$. Therefore,
\begin{align}
\curl_{\bar\eta_\kappa\circ\theta_l} W(l)=A^\kappa\curl u_0(\theta_l) 
+A^\kappa \int_0^t( B_\beta^\kappa W(l),_\beta, +
\div_{\bar\eta_\kappa\circ\theta_l} W(l))
\label{g14}
\end{align}
where $A^\kappa$ and $B_\beta^\kappa$, $\beta=1,2$, are smooth matrix fields depending on $\bbk$. With (\ref{g12}) and (\ref{g14}), we infer in a similar way as for (\ref{g9}) that on $(0,1)^2\times(\frac{1}{p},1)$ we have
\begin{align}
\curl_{\bar\eta_\kappa\circ\theta_l} \rho_m\star_h [\alpha_l W(l)]= A^\kappa\sum_{\beta=1}^2 \int_0^t B_\beta^\kappa \rho_m\star_h[\alpha_l\circ\theta_l W(l)],_\beta+R_3,
\label{g15}
\end{align}
with $\|R_3\|_{L^2(0,T;L^2((0,1)^3))}\le N(u_0)$. Therefore, with (\ref{g13}) and (\ref{g14}), we have that
\begin{align}
 [\alpha_l(\theta_l) W(l)],_3= M^\kappa A^\kappa\sum_{\beta=1}^2 \int_0^t B_\beta^\kappa [\alpha_l(\theta_l) W(l)],_\beta+\sum_{\beta=1}^2 A_\beta^\kappa [\alpha_l(\theta_l)W(l)],_\beta+R_4,
\label{g16}
\end{align}
with $\|R_4\|_{L^2(0,T;L^2((0,1)^3))}\le N(u_0)$. Thus, for any test function
$\varphi\in H^1((0,1)^3)$,
\begin{align*}
\int_{(0,1)^2\times\frac{1}{p}} \alpha_l(\theta_l) W(l)\cdot\varphi=
\int_{(0,1)^2\times(\frac{1}{p},0)} [\alpha_l(\theta_l) W(l)],_3 \cdot\varphi
+\int_{(0,1)^2\times(\frac{1}{p},0)} \alpha_l(\theta_l) W(l) \cdot\varphi,_3.
\end{align*}
Now, since for $\beta=1,2$, we have
\begin{align*}
\int_{(0,1)^2\times(\frac{1}{p},0)} [\alpha_l(\theta_l) W(l)],_\beta \cdot\varphi=
-\int_{(0,1)^2\times(\frac{1}{p},0)} [\alpha_l(\theta_l) W(l)] \cdot\varphi,_\beta,
\end{align*}
using (\ref{g16}), we  infer that
\begin{align*}
\bigr|\int_{(0,1)^2\times\frac{1}{p}} \alpha_l(\theta_l) W(l)\cdot\varphi\bigr|\le C\ (\|W(l)\|_{0,(0,1)^3}+\|R_4\|_{0,(0,1)^3})\|\varphi\|_{1,(0,1)^3},
\end{align*}
implying (independently of $p>1$) the following trace estimate for $W(l)$ (not just its normal component):
\begin{equation}
\label{g17}
\int_0^T |\alpha_l(\theta_l)\ W(l)|^2_{-\ud,(0,1)^2\times\frac{1}{p}}\le N(u_0).
\end{equation}
Similarly as (\ref{g15}), we also have the divergence relation
\begin{align}
\div_{\bar\eta_\kappa\circ\theta_l} \rho_m\star_h [\alpha_l W(l)]= C^\kappa\sum_{\beta=1}^2 \int_0^t D_\beta^\kappa \rho_m\star_h[\alpha_l\circ\theta_l W(l)],_\beta+R_5,
\label{g18}
\end{align}
with $\|R_5\|_{L^2(0,T;L^2((0,1)^3))}\le N(u_0)$, and $C^\kappa$ and 
$D_\beta^\kappa$, $\beta=1,2$, are smooth matrix fields in terms of $\bbk$. 
From (\ref{g15}) and (\ref{g18}), we then infer, just as in (\ref{g10}), that
\begin{align*}
\int_0^T \|\rho_m\star_h [\alpha_l W(l)]\circ(\bar\eta_\kappa\circ\theta_l)^{-1}
\|^2_{1,\Omega_p^l}\le&\ N(u_0)
+\int_0^T  | \rho_m\star_h [\alpha_l W(l)]\circ(\bar\eta_\kappa\circ\theta_l)^{-1}\cdot \bar n_\kappa |^2_{\ud,\partial\Omega_p^l},
\end{align*}
where $\Omega_p^l=\theta_l((0,1)^2\times(\frac{1}{p},1))$. Thus,
\begin{align*}
\int_0^T \| \rho_m\star_h [\alpha_l(\theta_l) W(l)]
\|^2_{1,(0,1)^2\times(\frac{1}{p},1)}\le&\ N(u_0)
+\int_0^T  |\rho_m\star_h [\alpha_l(\theta_l) W(l)]|^2_{\ud,(0,1)^2\times\frac{1}{p}}.
\end{align*}
Now, from the properties of the convolution,
\begin{align*}
\frac{1}{m} \left|\rho_m\star_h [\alpha_l(\theta_l) W(l)] \,
\right|_{\ud,(0,1)^2\times\frac{1}{p}}
\le C \left|\rho_m\star_h [\alpha_l(\theta_l) W(l)] \,
\right|_{-\ud,(0,1)^2\times\frac{1}{p}},
\end{align*}
which, with (\ref{g17}), leads us (independently of $p>1$) to
\begin{align*}
\frac{1}{m^2} \int_0^T \| \rho_m\star_h [\alpha_l(\theta_l) W(l)]
\|^2_{1,(0,1)^2\times(\frac{1}{p},1)}\le N(u_0),
\end{align*}
for any 
$0<\d\frac{1}{m}\le \text{dist}(\text{supp}\alpha_l(\theta_l),
\partial (0,1)^2\times(0,1))$. Since this estimate holds for any $p>1$, 
we then infer that
\begin{align}
\frac{1}{m^2} \int_0^T \| \rho_m\star_h [\alpha_l(\theta_l) W(l)]
\|^2_{1,(0,1)^3}\le N(u_0),
\label{g19}
\end{align}
for any 
$0<\d\frac{1}{m}\le\text{dist}(\text{supp}\alpha_l(\theta_l),\partial(0,1)^2\times(0,1))$.
Therefore, $\rho_m\star_h [\alpha_l(\theta_l) W(l)]\in H^1((0,1)^3)$ 
(which was not a priori known since our convolution smooths only in the
horizontal directions), with a bound depending {\it a priori} on $m$.

\subsubsection{Control of the horizontal convolved-by-layers velocity fields 
independently of $m$}\hf

From (\ref{g15}) and (\ref{g18}), we infer that
\begin{align*}
&\int_0^T \|\rho_m\star_h [\alpha_l W(l)]\circ(\bar\eta_\kappa\circ\theta_l)^{-1}
\|^2_{1,\bar\eta_\kappa(\Omega)}
\\ & \qquad\qquad\qquad\qquad
\le\ N(u_0)
+
\int_0^T  | \rho_m\star_h [\alpha_l(\theta_l) W(l)]
\circ(\bar\eta_\kappa\circ\theta_l)^{-1}\cdot 
\nk |^2_{\ud,\partial\bar\eta_\kappa(\Omega)},
\end{align*}
and thus,
\begin{align}
&\int_0^T 
\| \rho_m \star_h [\alpha_l(\theta_l) W(l)] \|^2_{1,(0,1)^3} 
\label{g20}
\\
& \qquad\qquad
\le\ N(u_0)\n
+\int_0^T | \rho_m\star_h [\alpha_l(\theta_l) W(l)]
\cdot \nk(\bar\eta_\kappa\circ\theta_l)|^2_{\ud,(0,1)^2\times\{0\}}.
\nonumber
\end{align}

Next, we have for any $x\in (0,1)^2\times\{0\}$:
\begin{align*}
\rho_m\star_h [\alpha_l(\theta_l) W(l)]\cdot \nk(\bar\eta_\kappa\circ\theta_l)(x)&=
\rho_m\star_h [\alpha_l(\theta_l) W(l)\cdot \nk(\bar\eta_\kappa\circ\theta_l)](x)
+ f(x),
\end{align*}
with
\begin{align*}
f(x)=\int_{\R^2} \rho_m(x_H-y_H) \alpha_l(\theta_l) W(l)(y_H,x_3)\cdot[\nk(\bar\eta_\kappa\circ\theta_l)(x_H,x_3)-\nk(\bar\eta_\kappa\circ\theta_l)(y_H,x_3)] dy_H.
\end{align*}
Therefore, with (\ref{g20}), we obtain
\begin{align}
&\int_0^T \|\rho_m\star_h [\alpha_l(\theta_l) W(l)]
\|^2_{1,(0,1)^3} \n\\
&\qquad \le N(u_0) +|f|_{\ud,(0,1)^2\times\{0\}}
+\int_0^T  \left| \rho_m\star_h [\alpha_l(\theta_l) W(l)\cdot \nk(\bar\eta_\kappa\circ\theta_l)]\right|^2_{\ud,(0,1)^2\times\{0\}}\n\\
&\qquad \le N(u_0) +|f|_{\ud,(0,1)^2\times\{0\}}
+\int_0^T  | \alpha_l(\theta_l) W(l)\cdot \nk(\bar\eta_\kappa\circ\theta_l)|
^2_{\ud,(0,1)^2\times\{0\}}\n\\
&\qquad\le N(u_0) +|f|_{\ud,(0,1)^2\times\{0\}}+\int_0^T  
| \alpha_l w_\epsilon\cdot \nk(\bar\eta_\kappa)|^2_{\ud,\Gamma}\n\\
&\qquad\le N(u_0) +|f|_{\ud,(0,1)^2\times\{0\}},
\label{g21}
\end{align}
where we have used the trace control (\ref{g4}) 
for the last inequality in (\ref{g21}).
We now turn our attention to $|f|_{\ud,(0,1)^2\times\{0\}}$. We first have that
\begin{align}
\|f\|_{0,(0,1)^3}&\le  \frac{C}{m} \|\nk(\bar\eta_\kappa)\|_{H^3(\Omega)}
\|\rho_m\star_h\alpha_l(\theta_l)|W(l)|\|_{0,(0,1)^3}
\le  \frac{C}{m} N(u_0),
\label{g22}
\end{align}
where we have used the definition of $C_T$ to bound 
$\|\nk(\bar\eta_\kappa)\|_{H^3(\Omega)}$. Next, we have for $\beta=1,2$
that
\begin{align*}
f,_\beta(x)&=\int_{\R^2} \rho_m,_\beta(x_H-y_H) \alpha_l(\theta_l) W(l)(y_H,x_3)\cdot[\nk(\bar\eta_\kappa\circ\theta_l)(x_H,x_3)-\nk(\bar\eta_\kappa\circ\theta_l)(y_H,x_3)] dy_H\\
&\ \ \ + \int_{\R^2} \rho_m(x_H-y_H) \alpha_l(\theta_l) W(l)(y_H,x_3) dy \cdot \nk(\bar\eta_\kappa\circ\theta_l),_\beta(x),
\end{align*}
showing that
\begin{align}
\|f,_\beta\|_{0,(0,1)^3}
&\le  \frac{C}{m} \|\nk(\bar\eta_\kappa)\|_{H^3(\Omega)}\sum_{i=1}^3
\|\, |\rho_m,_\beta|\, \star_h\alpha_l(\theta_l)|W(l)^i|\|_{0,(0,1)^3} \n \\
& \qquad\qquad\qquad\qquad
+ \|\nk(\bar\eta_\kappa)\|_{H^3(\Omega)}\sum_{i=1}^3\|\rho_m\star_h\alpha_l(\theta_l)|W(l)^i|\|_{0,(0,1)^3}\n\\
&\le  C \|\nk(\bar\eta_\kappa)\|_{H^3(\Omega)}\sum_{i=1}^3\||(\rho,_\beta)_m|
\star_h\alpha_l(\theta_l)|W(l)^i|\|_{0,(0,1)^3} \n\\
&\qquad\qquad\qquad\qquad
+ \|\nk(\bar\eta_\kappa)\|_{H^3(\Omega)}\sum_{i=1}^3\|\rho_m
\star_h\alpha_l(\theta_l)|W(l)^i|\|_{0,(0,1)^3}\n\\
&\le  C \sum_{i=1}^3\||(\rho,_\beta)_m|
\star_h\alpha_l(\theta_l)|W(l)^i|\|_{0,(0,1)^3}
+ C \sum_{i=1}^3\|\rho_m\star_h\alpha_l(\theta_l)|W(l)^i|\|_{0,(0,1)^3}\n\\
&\le C \|\alpha_l(\theta_l)W(l)\|_{0,(0,1)^3} \le N(u_0).
\label{g23}
\end{align}
Next, for the vertical derivative,
\begin{align*}
f,_3(x)&=\int_{\R^2} \rho_m(x_H-y_H) [\alpha_l(\theta_l) W(l)],_3(y_H,x_3)\cdot[\nk(\bar\eta_\kappa\circ\theta_l)]_{(y_H,x_3)}^{(x_H,x_3)} dy_H\\
&\ \ \ + \int_{\R^2} \rho_m(x_H-y_H) [\alpha_l(\theta_l) W(l)](y_H,x_3)\cdot[\nk(\bar\eta_\kappa\circ\theta_l),_3]_{(y_H,x_3)}^{(x_H,x_3)} dy_H\,,
\end{align*}
where $[\cdot ]_{(y_H,x_3)}^{(x_H,x_3)} = [\cdot](x_H,x_3)-[\cdot](y_H,x_3)$.
Notice that for a smooth matrix field $A$ in $(0,1)^3$ and for $\beta=1,2$,
\begin{align*}
G(x)&=\int_{\R^2} \rho_m(x_H-y_H) A(y_h,x_3) [\alpha_l(\theta_l) W(l)],_\beta(y_H,x_3) \cdot[\nk(\bar\eta_\kappa\circ\theta_l)]_{(y_H,x_3)}^{(x_H,x_3)} dy_H,
\end{align*}
satisfies
\begin{align*}
G(x)=&-m\int_{\R^2} (\rho,_\beta)_m(x_H-y_H) A(y_h,x_3) [\alpha_l(\theta_l) W(l)](y_H,x_3) \cdot[\nk(\bar\eta_\kappa\circ\theta_l)]_{(y_H,x_3)}^{(x_H,x_3)} dy_H\\
&-\int_{\R^2} \rho_m(x_H-y_H) A,_\beta(y_h,x_3) [\alpha_l(\theta_l) W(l)](y_H,x_3) \cdot[\nk(\bar\eta_\kappa\circ\theta_l)]_{(y_H,x_3)}^{(x_H,x_3)} dy_H\\
&-\int_{\R^2} \rho_m(x_H-y_H) A(y_h,x_3) [\alpha_l(\theta_l) W(l)](y_H,x_3) \cdot[\nk(\bar\eta_\kappa\circ\theta_l),_\beta]_{(y_H,x_3)}^{(x_H,x_3)} dy_H,
\end{align*}
showing, just as for (\ref{g23}), 
that $\|G\|_{0,(0,1)^3}\le N(u_0)$. Therefore, with
 (\ref{g16}), we see that the first integral term appearing in the expression of $f,_3$ is bounded in a similar way, implying that
\begin{align}
\|f,_3\|_{0,(0,1)^3}\le N(u_0).
\label{g24}
\end{align}
Consequently, with (\ref{g22}), (\ref{g23}), (\ref{g24}), we obtain that
\begin{equation}
\label{g25}
\|f\|_{1,(0,1)^3}\le N(u_0).
\end{equation}
Therefore, (\ref{g21}) implies that
\begin{align}
\label{g25'}
\int_0^T \|\rho_m\star_h [\alpha_l(\theta_l) W(l)]
\|^2_{1,(0,1)^3}&\le N(u_0) .
\end{align}
\subsubsection{Control of $w_\epsilon$ in $L^2(0,T;H^1(\Omega))$}\hf

Since (\ref{g25'}) holds independently of $m$ sufficiently large, 
this implies that
\begin{align*}
\int_0^T \|\alpha_l(\theta_l) W(l)
\|^2_{1,(0,1)^3}&\le N(u_0) .
\end{align*}
Since we proved in subsection \ref{6.3.1} that $w_\epsilon$ is bounded in 
$L^2(0,T;H^1(\omega))$ independently of $\epsilon$ for each domain 
$\omega\subset\subset \Omega$, this provides us with the estimate
\begin{align}
\int_0^T \|w_\epsilon
\|^2_{1}&\le N(u_0),
\end{align}
independently of $\epsilon>0$.

\begin{remark} In the two-dimensional case, a simpler proof of Step 3 is 
possible, founded upon a scalar potential function for the velocity field.
For conciseness, we consider a simply-connected domain, the non-simply 
connected case being treated similarly by local charts. Once again, 
we let $\ue=\we\circ\ek^{-1}$.  From (\ref{g3.b}) and (\ref{g4}), 
let $w^\epsilon_{\tau}\in L^2(0,T;H^1(\Omega))$ such that
\begin{align*}
\operatorname{div}(w^\epsilon_{\tau}(\ek^{-1}))(\ek)&= \operatorname{div}(\we(\ek^{-1}))(\ek)  \ \ \text{in}\ L^2(0,T; L^2(\Omega)),
\\
w^\epsilon_{\tau}\cdot \nk(\ek)&= \we\cdot\nk(\ek)\ \ \text{in}\ L^2(0,T; H^1(\Gamma))\ .
\end{align*}
We infer the existence of $\psi^\epsilon\in L^2(0,T;H^1_0(\ek(\Omega)))$ such that $\ue=w^\epsilon_{\tau}(\ek^{-1})+(-\psi^\epsilon,_2,\psi^\epsilon,_1)$.
Now, from (\ref{g6}), we see that in $L^2(0,T;H^{-1}(\Omega))$, we have for $\bar\psi^\epsilon=\psi^\epsilon\circ\ek$,
\begin{equation}
\label{g6bis}
\d -(\bar{a}_\kappa)_i^k ((\bar{a}_\kappa)_i^j \bar\psi^\epsilon,_j),_k=f^\epsilon_{\tau}-\int_0^t A^\kappa_{ij} \bar\psi^\epsilon,_{ij},
\end{equation}
where $f^\epsilon_\tau$ is bounded in $L^2(0,T;L^2(\Omega))$. 
It is readily seen that $\bar\psi^\epsilon$ is the unique solution of this equation in
$L^2(0,T;H^1_0(\Omega))$.
We now establish that this uniqueness provides extra regularity for 
$\bar\psi^\epsilon$. By defining the mapping $\Theta$ from 
$L^2(0,T;H^2(\Omega)\cap H^1_0(\Omega))$ into itself by associating to 
any $\xi$ in this space, the solution $\Theta\xi$ (for
almost all $t\in[0,T]$) of
\begin{equation*}
\d -(\bak)_i^k ((\bak)_i^j \Theta\xi,_j),_k=f_{\tau}-\int_0^t A^\kappa_{ij} \xi,_{ij} ,
\end{equation*}
we see that for $t_1$ small enough (depending on Sobolev constants and on
$\|\uk\|_{L^\infty(0,T;H^3(\Omega))}$) $\Theta$ is contractive from 
$L^2(0,t_1;H^2(\Omega)\cap H^1_0(\Omega))$ into itself, which provides a 
fixed-point for $\Theta$ in this space. It is thus a solution of (\ref{g6bis}) 
on $[0,t_1]$. By uniqueness of such a solution, we have that 
$\bar\psi^\epsilon\in L^2(0,t_1;H^2(\Omega))$ and thus that 
$\we\in L^2(0,t_1;H^1(\Omega))$. By defining a mapping similar to $\Theta$, but
this time starting from $\d t_2\in [\frac{t_1}{2},t_1]$ such that 
$\we(t_2)\in H^1(\Omega)$ instead of $u_0$ (which ensures that the new 
$f_\tau$ is still in $L^2(0,t_2;L^2(\Omega))$), we obtain the same conclusion
on $[t_2,t_2+t_1]$, leading us to $\d \we\in L^2(0,\td t_1;H^1(\Omega))$. 
By induction, we then find $\we\in L^2(0,T;H^1(\Omega))$.
\end{remark}

\begin{remark}
Whereas Hodge decompositions with vector potentials $\psi$ are possible in 
higher dimension, it turns out that a Dirichlet condition $\psi=0$ for the
associated elliptic problem is not possible. This in turn is problematic for 
any uniqueness argument in $L^2(0,T;H^1(\Omega))$ for $\psi$, since it does not
seem possible to find a boundary condition that would be naturally associated 
to the second order operators appearing on both sides of the three-dimensional 
analogous of (\ref{g6bis}).
\end{remark}

\section{Pressure as a Lagrange multiplier}
\label{5}
We will need two Lagrange multiplier lemmas for our pressure function
in our analysis as the penalization parameter $\epsilon\rightarrow 0$. We
begin with a lemma that is necessary for a new Hodge-type decomposition of
the velocity field.

\begin{lemma}\label{lemma_lagrange}
For all $l \in H^\ud(\Omega)$, $t\in [0,T]$, there exists a constant
$C>0$ and
$\phi(l) \in H^\td(\Omega)$ such that $(\bak)_i^j (t) \phi^i,_j =l$ in $\Omega$ and
\begin{equation}\label{v-p} 
\|\phi(l)\|^2_{\td} \le  C\|l\|^2_{\ud}. 
\end{equation}
\end{lemma}
\begin{proof}
Let $\psi(l)$ be the solution of
\begin{subequations}
\label{laplacien}
\begin{align}
(\bak)_i^j[ (\bak)_i^k \psi(l),_k],_j&=l\ \text{in}\ \Omega\\
\psi(l)&=0\ \text{on}\ \Gamma.
\end{align}
\end{subequations}
We then see that $\phi^i(l)=(\bak)_i^j \psi(l),_j$ satisfies the statement of 
the lemma.
The inequality (\ref{v-p}) is a simple consequence of the properties of
$l$ and of the condition $\bar v\in C_T$.
\end{proof}

We can now follow \cite{SolSca1973}. For $p\in H^\ud(\Omega)'$, 
 define the linear functional on $H^\td(\Omega)$
by $\langle p,(\bak)_i^j (t) \varphi^i,_j\rangle_ {\ud}$, where $\varphi\in H^\td(\Omega)$.
By the Riesz representation theorem, there is a bounded linear operator
$Q (t): (H^\ud(\Omega))'\rightarrow  H^\td(\Omega)$ such that
$$ 
\forall \varphi\in H^\td(\Omega),\ \langle p,\ (\bak)_i^j (t) \varphi^i,_j\rangle_ {\ud}=
(Q(t)p,\ \varphi)_{\td}. 
$$ 
Letting $\varphi=Q(t)p$ shows that 
\begin{equation}\label{Qp0}
\|Q(t)p\|_{\td} \le C
\|p \|_ {H^\ud(\Omega)'}
\end{equation}
for some constant $C>0$. 
Using Lemma \ref{lemma_lagrange},  we see that
\begin{equation*}
\forall l\in H^\ud(\Omega),\ \langle p,\ l\rangle_ {\ud}=
(Q(t)p,\ \phi(l))_{\td}, 
\end{equation*}
and thus
\begin{equation}\label{Qp}
\|p\|_{H^\ud(\Omega)'}\le C \|Q(t)p\|_{\td},
\end{equation}
which shows that $R(Q(t))$ is closed in $H^\td(\Omega)$.
Let ${\mathcal V}_{\bar v}(t)= \{ v\in L^2(\Omega) \ | \ (\bak)^j_i(t) v^i,_j(t)=0\}$.
Since ${\mathcal V}_{\bar v}(t)\cap H^\td(\Omega) = R(Q(t))^\perp$, it follows that 
\begin{equation}\label{hodge}
H^\td(\Omega) = R(Q(t)) \oplus_ {H^\td(\Omega)} {\mathcal V}_{\bar v}(t)\cap H^\td(\Omega).
\end{equation}

We can now introduce our first Lagrange multiplier
\begin{lemma} \label{Lagrange}
Let ${\mathfrak L}(t)\in H^\td (\Omega)'$ be such that ${\mathfrak L}(t) \varphi=0$ for any
$\varphi\in {\mathcal V}_{\bar v}(t)\cap H^\td(\Omega)$. Then there exists a unique $q(t)\in H^\ud (\Omega)'$, which is termed the pressure function, satisfying 
$$\forall \varphi\in H^\td (\Omega),\ \ 
{\mathfrak L}(t) (\varphi)=\langle q(t),\ (\bak)_i^j \varphi^i,_j\rangle_{\ud}.$$ 
Moreover, there is a $C>0$ (which does not depend on $t\in [0,T]$ and on the choice of $\bar v\in C_T$) such that$$\|q(t)\|_{H^\ud (\Omega)'}\le C\ \|{\mathfrak L}(t)\|_{H^\td(\Omega)'}.$$
\end{lemma}

\begin{proof}
By the decomposition (\ref{hodge}), for $\varphi\in{H^\td(\Omega, {\mathbb R}^3)}$, we let
$\varphi=v_1+v_2$, where $v_1 \in {\mathcal V}_{\bar v} (t)\cap H^\td(\Omega)$ and $v_2 \in R(Q(t))$.  From our assumption, it follows that
$$ 
{\mathfrak L}(t)(\varphi) = {\mathfrak L}(t)(v_2) = ( \psi(t), v_2)_{H^\td(\Omega)} 
= ( \psi(t), \varphi)_{H^\td(\Omega)}, $$ \ for a unique  \ $\psi(t) \in R(Q(t))$.

 From the definition of $Q(t)$ we then get the existence of a unique 
$q(t)\in H^\ud (\Omega)'$ such that 
$$\forall \varphi\in H^\td (\Omega),\ \ {\mathfrak L}(t) (\varphi)
=\langle q(t),\ (\bak)_i^j \varphi^i,_j\rangle_{\ud}.$$
The estimate stated in the lemma is then a simple consequence of (\ref{Qp}).
\end{proof}

We also need the case where the pressure function is in $H^\ud(\Omega)$. We 
start, as above, with a simple elliptic result:

\begin{lemma}\label{lemma_lagrangebis}
For all $l \in H^\ud(\Omega)'$, $t\in [0,T]$, there exists a constant
$C>0$ and
$\phi(l) \in H^\ud(\Omega)$ such that $(\bak)_i^j (t) \phi^i,_j =l$ in $\Omega$ and
\begin{equation}\label{v-pbis} 
\|\phi(l)\|^2_{\ud} \le  C\|l\|^2_{H^\ud(\Omega)'}. 
\end{equation}
\end{lemma}
\begin{proof}Let $\psi(l)$ be the solution of (\ref{laplacien}). Since $\psi$ is linear and continuous from $H^1(\Omega)'$ into $H^1(\Omega)$ and from $L^2(\Omega)$ into $H^2(\Omega)$, by interpolation, we have that
$\psi$ is linear and continuous from $H^\ud(\Omega)'$ into $H^\td(\Omega)$. We then see that $\phi^i(l)=(\bak)_i^j \psi(l),_j$ satisfies the statement of the
lemma.
\end{proof}

 For $p\in H^\ud(\Omega)$, we define the linear functional on $X(t)$
by $\langle (\bak)_i^j (t) \varphi^i,_j,p\rangle_ {\ud}$, where $\varphi\in X(t)=\{\psi\in H^\ud(\Omega)|\ (\bak)_i^j \psi^i,_j\in H^\ud(\Omega)'\}$.
By the Riesz representation theorem, there is a bounded linear operator
$Q (t): H^\ud(\Omega)\rightarrow  X(t)$ such that
$$ 
\forall \varphi\in X(t),\ \langle (\bak)_i^j (t) \varphi^i,_j,p\rangle_ {\ud}=
(Q(t)p,\ \varphi)_{X(t)}. 
$$ 
Letting $\varphi=Q(t)p$ shows that 
\begin{equation}\label{Qp0bis}
\|Q(t)p\|_{X(t)} \le C
\|p \|_ {H^\ud(\Omega)}
\end{equation}
for some constant $C>0$. 
Using Lemma \ref{lemma_lagrangebis},  we see that
\begin{equation*}
\forall l\in H^\ud(\Omega)',\ \langle l,\ p\rangle_ {\ud}=
(Q(t)p,\ \phi(l))_{X(t)}, 
\end{equation*}
and thus
\begin{equation}\label{Qpbis}
 \|p\|_{H^\ud(\Omega)}\le C \|Q(t)p\|_{X_t},
\end{equation}
which shows that $R(Q(t))$ is closed in $X(t)$.
Since ${\mathcal V}_{\bar v}(t)\cap X(t) = R(Q(t))^\perp$, it follows that 
\begin{equation}\label{hodgebis}
X(t) = R(Q(t)) \oplus_ {X(t)} {\mathcal V}_{\bar v}(t)\cap X(t).
\end{equation}

Our second Lagrange multiplier lemma can now be stated.
\begin{lemma} \label{Lagrangebis}
Let ${\mathfrak L}(t)\in X(t)'$ be such that ${\mathfrak L}(t) \varphi=0$ for any
$\varphi\in {\mathcal V}_{\bar v}(t)\cap H^\ud(\Omega)$. Then there exists a unique $q(t)\in H^\ud (\Omega)$, which is termed the pressure function, satisfying 
$$\forall \varphi\in X(t),\ \ 
{\mathfrak L}(t) (\varphi)=\langle (\bak)_i^j \varphi^i,_j, q(t)\rangle_{\ud}.$$ 
Moreover, there is a $C>0$ (which does not depend on $t\in [0,T]$ and on the choice of $\bar v\in C_T$) such that$$\|q(t)\|_{H^\ud (\Omega)}\le C\ \|{\mathfrak L}(t)\|_{X(t)'}.$$
\end{lemma}

\begin{proof}
By the decomposition (\ref{hodgebis}), for $\varphi\in X(t)$, we let
$\varphi=v_1+v_2$, where $v_1 \in {\mathcal V}_{\bar v} (t)\cap H^\ud(\Omega)$ and $v_2 \in R(Q(t))$.  From our assumption, it follows that
$$ 
{\mathfrak L}(t)(\varphi) = {\mathfrak L}(t)(v_2) = ( \psi(t), v_2)_{X(t)} 
= ( \psi(t), \varphi)_{X(t)}, $$ \ for a unique  \ $\psi(t) \in R(Q(t))$.

 From the definition of $Q(t)$ we then get the existence of a unique 
$q(t)\in H^\ud (\Omega)$ such that 
$$\forall \varphi\in X(t),\ \ {\mathfrak L}(t) (\varphi)
=\langle (\bak)_i^j \varphi^i,_j,\ q(t)\rangle_{\ud}.$$
The estimate stated in the lemma is then a simple consequence of (\ref{Qpbis}).
\end{proof}

\section{Existence of a solution to the linearized smoothed 
$\kappa$-problem (\ref{smoothlinear})}
\label{6}
 
In this section, we prove the existence of a solution $w$ to 
the linear problem (\ref{smoothlinear}), constructed as the limit 
$\epsilon\rightarrow 0$. 

The analysis requires establishing the regularity of the weak solution.
Note that the extra regularity on $u_0$ is needed in order to ensure the 
regularity property for $w$, $q$, and their time derivatives as stated in the 
next theorem, without having to consider the variational limits of the time differentiated penalized problems. 

\begin{theorem}
\label{uniqueweak}
Suppose  that $u_0 \in H^{13.5}(\Omega)$ and $\Omega$ is of class $C^\infty$.
Then, there exists a unique weak solution $w$
to the linear problem (\ref{smoothlinear}), which is moreover in 
$L^2 (0,T;H^{13.5}(\Omega))$. Furthermore, 
\begin{gather*}
\partial_t^i w\in L^2 (0,T;H^{13.5-3i}(\Omega))\cap 
L^\infty (0,T;H^{12.5-3i}(\Omega)), \ \ \ i=1,2,3,4\,,\\
\partial_t^i q\in L^2 (0,T;H^{11.5-3i}(\Omega))\cap 
L^\infty (0,T;H^{10.5-3i}(\Omega)), \ \ \ i=0,1,2,3\,.
\end{gather*}
\end{theorem}

\begin{proof} 
\noindent{\bf Step 1. The limit as $\epsilon \rightarrow 0$.}

Let $\epsilon=\frac{1}{m}$; we first pass to the weak limit as 
$m\rightarrow \infty$.  
The inequality (\ref{g1}) provides the following bound, independent of $\epsilon$: $$\int_0^T \frac{1}{\epsilon}\|(\bak)_i^j \we^i,_j\|^2_{0}+|\we\cdot\nk(\ek)|^2_1+\|\we\|^2_0\ dt\le  N(u_0)$$ which provides a subsequence 
$\{w_{\frac{1}{m_l}}\}$ such that
\begin{subequations}
\label{weakconvergence}
\begin{align}
w_{\frac{1}{m_l}} \rightharpoonup  w \ \ \text{ in } \ \  L^2(0,T; L^2(\Omega))\ ,\\
(\bak)_i^j w_{\frac{1}{m_l}}^i,_j \rightharpoonup  (\bak)_i^j w^i,_j \ \ \text{ in } \ \  L^2(0,T; L^2(\Omega))\ ,\\
w_{\frac{1}{m_l}}\cdot\nk(\ek) \rightharpoonup  w\cdot\nk(\ek) \ \ \text{ in } \ \  L^2(0,T; H^1(\Gamma))\ .
\end{align}
\end{subequations}
The justification for $w\cdot\nk(\ek)$ being the third weak limit in (\ref{weakconvergence}) comes from the identity $(\bak)_i^j \we^i,_j=\operatorname{div}(\we\circ\ek^{-1})(\ek)$ and the fact that $\nk$ is the normal to $\ek(\Omega)$.

Moreover, since (\ref{g1}) also shows that $\|(\bak)_i^j w_{\frac{1}{m}}^i,_j\|_{L^2(0,T;L^2(\Omega))}\rightarrow 0$ as $m\rightarrow\infty$, we then have
$\|(\bar a_\kappa)_i^j  w^i,_j\|_{L^2(0,T;L^2(\Omega))}=0$, {\it i.e.}
\begin{equation}
\label{divfree}
(\bar a_\kappa)_i^j  w^i,_j=0\ \text{in}\ L^2(0,T;L^2(\Omega)).
\end{equation}

Now, let us denote $u=w\circ\ek^{-1}$, so that thanks to (\ref{divfree}) and (\ref{g5}) we have
\begin{subequations}
\label{divcurl}
\begin{align}
\operatorname{div}u&=0\ \text{in}\ \ek(\Omega)\ ,\\
\operatorname{curl}u(\ek)&=\curl u_0+\int_0^t B(\nabla\uk,\nabla u)\ \text{in}\ H^{-1}(\Omega).
\end{align}
\end{subequations}
By proceeding as in Step 3 of Section \ref{4}, the trace regularity $(u\cdot\nk)(\ek)\in L^2(0,T;H^1(\Gamma))$ and the system (\ref{divcurl}) then yield
\begin{equation*}
\|w\|_{L^2(0,T;H^\td(\Omega))}\le N(u_0),
\end{equation*}
where $N(u_0)$ is defined in (\ref{Nu0}).

\noindent{\bf Step 2. The equation for $w$ and the pressure.}

Now, for any $y\in L^2(0,T;H^\td(\Omega))$ and $l= (\bar a_\kappa)_i^j y^i,_j$, 
we see that for a solution $\varphi$
almost everywhere on $(0,T)$ of the elliptic problem
\begin{align*}
(\bak)_i^j [\bar J_\kappa^{-1}(\bak)_i^k \varphi,_k],_j&=l\ \text{in}\ \Omega\\
\varphi&=0\ \text{on}\ \Gamma,
\end{align*}
if we let $e^i=\bar J_\kappa^{-1}(\bak)_i^k \varphi,_k$, and set $v=y-e$, we have that $e$ and $v$ 
are both in 
$L^2(0,T;H^\td(\Omega))$, with 
\begin{align*}
\int_0^T [ \|e\|^2_\td+\|v\|^2_\td] \le & C \int_0^T \|y\|^2_\td,\\
(\bak)_i^j v^i,_j =& 0.
\end{align*}
Since $(\bar a_\kappa)_i^j  w^i,_j=0\ \text{in}\ L^2(0,T;L^2(\Omega))$, we infer that $(\bar a_\kappa)_i^j  w_t^i,_j=-[(\bar a_\kappa)_i^j]_t  w^i,_j\in L^2(0,T;H^\ud(\Omega))$, and that
\begin{equation*}
\langle\bar J_\kappa\ w_t,e\rangle_\td=([(\bar a_\kappa)_i^j]_t  w^i,_j,\varphi)_0.
\end{equation*}

But $v$ also satisfies the variational equation
\begin{align*}
  \int_0^T \langle\bar J_\kappa\ \we_{t},v\rangle_{\td}  
+\kappa \int_0^T&  [{\we}\cdot\nk(\ek)),v\cdot\nk(\ek)]_1
=  \sigma \int_0^T  [L_{\bar g}\en\cdot\nk(\ek) ,v\cdot\nk(\ek)]_0 \  ,
\end{align*}
leading to 
\begin{align*}
\lim_{\epsilon\rightarrow 0}  \int_0^T \langle\bar J_\kappa\ \we_{t},y\rangle_{\td}=&   
\int_0^T  ([(\bar a_\kappa)_i^j]_t  w^i,_j,\varphi)_0
+  \sigma \int_0^T  [L_{\bar g}\en\cdot\nk(\ek) ,v\cdot\nk(\ek)]_0\\&
-\kappa \int_0^T [ {w}\cdot\nk(\ek),v\cdot\nk(\ek)]_1. 
 \end{align*}
We then see that as $\epsilon\rightarrow 0$,
\begin{equation}
\label{wt}
\int_0^T \|\we_t\|^2_{H^\td(\Omega)'}\le  N(u_0).
\end{equation}
By standard arguments, we infer that $\we_t\rightharpoonup w_t$ in $L^2(0,T;H^\td(\Omega)')$. This ensures that $w\in \mathcal{C}^0([0,T];L^2(\Omega))$, and the condition 
$\we(0)=u_0$ provides $w(0)=u_0$. Furthermore, we also have for 
any $\phi\in L^2(0,T;H^\td(\Omega))$ such that 
$(\bak)_i^j \phi^i,_j=0$ in $(0,T)\times\Omega$, the variational equation
\begin{align*}
  \int_0^T \langle\bar J_\kappa\ w_{t},\phi\rangle_{\td}  
&+\kappa \int_0^T [{w}\cdot\nk(\ek),\phi\cdot\nk(\ek)]_1 
=  \sigma \int_0^T [L_{\bar g}\en\cdot\nk(\ek) ,\phi\cdot\nk(\ek)]_0 \ .
\end{align*}
Next, since $w_t\in L^2(0,T;H^\td(\Omega)')$, the Lagrange multiplier lemma \ref{Lagrange} shows that there exists $q\in L^2(0,T;H^\ud(\Omega)')$ such that for
any $\phi\in L^2(0,T;H^\td(\Omega))$,
\begin{align}
  \int_0^T \langle\bar J_\kappa\ w_{t},\phi\rangle_{\td}  
&+\kappa \int_0^T [{w}\cdot\nk(\ek),\phi\cdot\nk(\ek)]_1 \n
\\ &
-\int_0^T \langle q,(\bak)_i^j \phi^i,_j\rangle_\ud=  \sigma \int_0^T [L_{\bar g}
\en\cdot\nk(\ek) ,\phi\cdot\nk(\ek)]_0 \ .
\label{wvar}
\end{align}
Now, if we have another solution $\tilde w\in L^2(0,T;H^\td(\Omega))$ such that
$\tilde w(0)=u_0$ and $\tilde w_t\in L^2(0,T;H^\td(\Omega)')$, we then see,
by using $w-\tilde w$ as a test function in the difference between (\ref{wvar})
and its counterpart with $\tilde w$, that we get $w-\tilde w=0$, ensuring 
uniqueness to the weak solution of (\ref{smoothlinear}).

\noindent {\bf Step 3. Regularity of $w$.}
We can now study the regularity of $w$ via difference quotient techniques.
 We will denote $\R^3_+=\{x\in\R^3|\ x_3> 0\}$, $S_0=B(0,1)\cap \{x\in\R^3|\ x_3= 0\}$ and $B_+(0,r)=B(0,r)\cap\R^3_+$ .
 We denote by $\theta$ a $C^\infty$ diffeomorphism from $B(0,1)$ into a neighborhood $V$ of a point $x_0\in\Gamma$ such that
$\theta(B(0,1)\cap\R^3_+)=V\cap\Omega$, with $\det\nabla\theta=1$. We consider the smooth cut-off function $\psi(x)=e^{\frac{1}{|x|^2-\ud}}$ if $x\in B(0,\ud)$, and $\psi(x)=0$ elsewhere, and with the use of the test function $[D_{-h}[\psi D_h (w\circ\theta)]]\circ\theta^{-1}\in L^2(0,T;H^\td(\Omega))$ in (\ref{wvar}), with $h=|h| e_\alpha (\alpha=1,2)$, we obtain:
\begin{align*}
I_1+\kappa I_2+I_3=  \sigma \int_0^T [L_{\bar g}\en\cdot\nk(\ek) ,[D_{-h}[\psi D_h (w\circ\theta)]]\circ\theta^{-1}\cdot\nk(\ek)]_{0} \ ,
\end{align*}
with
\begin{align*}
I_1&=\int_0^T \langle\bar J_\kappa\ w_{t},[D_{-h}[\psi D_h (w\circ\theta)]]\circ\theta^{-1}\rangle_{\td},\\
I_2&=\int_0^T [\partial ({w}\cdot\nk(\ek)),\partial ([D_{-h}[\psi D_h (w\circ\theta)]]\circ\theta^{-1}\cdot\nk(\ek))]_0,\\
I_3&=-\int_0^T \langle q,(\bak)_i^j[D_{-h}[\psi D_h (w\circ\theta),_j^i]]\circ\theta^{-1}\rangle_\ud.
\end{align*}
For $I_1$, we simply have
\begin{align}
I_1=&\|\sqrt{\psi}\ w\circ\theta(t)\|^2_{L^2(B_+(0,1))}-\|\sqrt{\psi}\ u_0\circ\theta\|^2_{L^2(B_+(0,1))}\n\\
&\ + \int_0^T \langle D_h[\bar J_\kappa(\theta)]\ w_{t}\circ\theta,\psi D_h (w\circ\theta)\rangle_{\td}\n\\
\ge & \|\sqrt{\psi}\ w\circ\theta(t)\|^2_{L^2(B_+(0,1))} -N(u_0)-\int_0^T \|w_t\|_{H^\td(\Omega)'} \|D_h[\bar J_\kappa(\theta)]\psi D_h (w\circ\theta)\|_{\td}\n\\
\ge & \|\sqrt{\psi}\ w\circ\theta(t)\|^2_{L^2(B_+(0,1))} -C_\delta N(u_0)-\delta \int_0^T  \|\sqrt{\psi} D_h (w\circ\theta)\|^2_{\td},
\label{I1}
\end{align}
where we have used (\ref{wt}) for the last inequality, 
and where the choice of $\delta>0$ will be made precise later.

For $I_2$, we have, if we define in $B_+(0,1)$, $W=w\circ\theta$ and $\NK=\nk(\ek)(\theta)$,
\begin{align}
I_2 &=\int_0^T \int_{S_0}\frac{G_{\alpha\beta}}{\sqrt{ a_0}} [W\cdot\NK],_\alpha [D_{-h}[\psi D_h W]\cdot\NK],_\beta\n\\
&=\int_0^T \int_{S_0}\frac{G_{\alpha\beta}}{\sqrt{ a_0}} [D_h W\cdot\NK],_\alpha [\psi D_h W\cdot\NK],_\beta\n\\
&\ \ + \int_0^T \int_{S_0} [D_h[ \frac{G_{\alpha\beta}}{\sqrt{ a_0}} [W\cdot\NK],_\alpha]- \frac{G_{\alpha\beta}}{\sqrt{ a_0}} [D_h W\cdot\NK],_\alpha ] [\psi D_h W\cdot\NK],_\beta\n\\
&\ \ + \int_0^T \int_{S_0}[\frac{G_{\alpha\beta}}{\sqrt{ a_0}} [ W\cdot\NK],_\alpha](\cdot+h) [\psi D_h W\cdot D_h\NK],_\beta,
\label{I21}
\end{align}
where 
$G_{\alpha\beta}=\theta,_\alpha\cdot\theta,_\beta$ and $a_0=\operatorname{det} 
G$.

In this section, we will denote by $\|\cdot\|_{s,\Theta}$ and
$|\cdot|_{s,\partial\Theta}$ the standard norms of $H^s(\Theta)$ and 
$H^s(\partial\Theta)$.

For the first term appearing in the right-hand side of the second inequality, we have
\begin{align*}
 \frac{g_{0\alpha\beta}}{\sqrt{ a_0}} [D_h W\cdot\NK],_\alpha [\psi D_h W\cdot\NK],_\beta &=\frac{g_{0\alpha\beta}}{\sqrt{ a_0}} \psi [D_h W\cdot\NK],_\alpha [D_h W\cdot\NK],_\beta\\
&\ \ +\frac{g_{0\alpha\beta}}{\sqrt{ a_0}} \sqrt{\psi}[D_h W\cdot\NK],_\alpha D_h W\cdot\NK \frac{\psi,_\beta}{\sqrt{\psi}},
\end{align*}
and thus, since $\psi$, $\sqrt{\psi}$ and $\d\frac{\nabla\psi}{\sqrt{\psi}}$ 
are chosen smooth, we infer that
\begin{align*}
\int_0^T\int_{S_0} \frac{g_{0\alpha\beta}}{\sqrt{ a}} [D_h W\cdot\NK],_\alpha [\psi D_h W\cdot\NK],_\beta &\ge C \int_0^T \| \sqrt{\psi} D_h W\cdot\NK \|^2_{1,S_0}\\
&\ \ - N(u_0) .
\end{align*}
The other terms in (\ref{I21}) are easily estimated leading to the estimate:
 \begin{equation}
\label{I2}
I_2\ge C \int_0^T \| \sqrt{\psi} D_h W\cdot\NK \|^2_{1,S_0} - N(u_0).
\end{equation}
Concerning $I_3$, we have
\begin{align*}
I_3&=-\int_0^T \langle q,(\bar b_\kappa )_i^j[D_{-h}[\psi D_h (W)]^i,_j\circ\theta^{-1}]\rangle_\ud,
\end{align*}
with $\bar b_\kappa=[\nabla(\bar\eta_\kappa\circ\theta)]^{-1}$. 
Now since $(\bar b_\kappa )_i^j W^i,_j=0$, we obtain
\begin{align*}
(\bar b_\kappa )_i^j D_{-h} D_h W^i,_j =&-D_{-h}D_h (\bar b_\kappa )_i^j W^i,_j
- D_h [(\bar b_\kappa )_i^j] (\cdot-h) D_{-h} W^i,_j\\
&- D_{-h} [(\bar b_\kappa )_i^j(\cdot+h)]  D_{h} W^i,_j,
\end{align*}
and thus
\begin{align}
|I_3|&\le C \int_0^T \|q\|_{H^\ud(\Omega)'} [\ \|\psi D_h W^i,_j\|_{\ud,B_+(0,1)}+
\|\frac{D_h \psi}{\sqrt{\psi}}\  \sqrt{\psi} D_h W^i,_j\|_{\ud, B_+(0,1)}\n\\
&\hskip 4cm
+\|\sqrt{\psi}D_h W\|_{\td,B_+(0,1)}]\n\\
&\le C_\delta N(u_0)+\delta \int_0^T |\sqrt{\psi} D_h W^i,_j|^2_{\ud,B_+(0,1)},
\label{I3}
\end{align}
where $\delta>0$ is arbitrary. Now, let $\Theta$ be a smooth domain included in $B_+(0,1)$ and containing $B_+(0,\ud)$. The inequalities (\ref{I1}), (\ref{I2}) and (\ref{I3}) yield
\begin{equation}
\label{traceh1}
\int_0^T |\sqrt{\psi} D_h W \cdot N^\kappa|^2_{1,\partial\Theta}\le C_\kappa\ N(u_0) + \delta \int_0^T \|\sqrt{\psi} D_h W\|^2_{\td,\Theta}.
\end{equation}
We now define in $B_+(0,1)$
\begin{align*}
\operatorname{div}_{\ek\circ\theta}W&=\operatorname{div}(W\circ\theta^{-1}\circ\ek^{-1})(\ek\circ\theta)=\operatorname{div}(u)(\ek\circ\theta),\\
\operatorname{curl}_{\ek\circ\theta}W&=\operatorname{curl}(u)(\ek\circ\theta).
\end{align*}
Thus, (\ref{divcurl}) translates in $B_+(0,1)$ to
\begin{align*}
\operatorname{div}_{\ek\circ\theta}W&=0,\\
[\operatorname{curl}_{\ek\circ\theta}W](t)&=[\operatorname{curl} u_0]\circ\theta+\int_0^t B(\nabla\uk,\nabla u)(\ek\circ\theta),\n\\
&=[\operatorname{curl} u_0]\circ\theta+\int_0^t B(\nabla\uk,\nabla W(\ek\circ\theta)^{-1} \nabla (\ek\circ\theta)^{-1})(\ek\circ\theta),
\end{align*}
and thus 
\begin{subequations}
\label{W1}
\begin{align}
\operatorname{div}_{\ek\circ\theta}(\sqrt{\psi} D_h W)&=-\sqrt{\psi} D_h (\bar b_\kappa)_i^j\  W^i,_j(\cdot+h) + \ud \frac{\psi,_j}{\sqrt{\psi}}\ (\bar b_\kappa)_i^j D_h W^i,\\
[\operatorname{curl}_{\ek\circ\theta}(\sqrt{\psi} D_h W)](t)&=R(W)+\int_0^t  B(\nabla\uk(\ek\circ\theta), \nabla [\sqrt{\psi} D_h W] [\nabla(\ek\circ\theta)]^{-1}),
\end{align}
\end{subequations}
with
\begin{align*}
\int_0^T \|R(W)\|^2_{\ud,\Theta}\le C N(u_0).
\end{align*}
With the trace estimate (\ref{traceh1}) and the control of $W$ in $L^2(0,T;H^\td(\Theta))$, we can then infer as we did in Step 3 of Section \ref{4} that
\begin{equation*}
\int_0^T \|\sqrt{\psi} D_h W\|^2_{\td,\Theta}\le C_\kappa\ N(u_0) + C_\kappa \delta \int_0^T \|\sqrt{\psi} D_h W\|^2_{\td,\Theta},
\end{equation*}
and thus with a choice of $\delta$ small enough,
\begin{equation*}
\int_0^T \|\sqrt{\psi} D_h W\|^2_{\td,\Theta}\le C_\kappa\ N(u_0) ,
\end{equation*}
yielding
\begin{equation*}
\int_0^T |\sqrt{\psi} D_h W|^2_{1,\partial\Theta}\le C_\kappa\ N(u_0).
\end{equation*}
Since this estimate is independent of $h$, we get the trace estimate
\begin{equation*}
\int_0^T |\sqrt{\psi}  W|^2_{2,\partial\Theta}\le C_\kappa\ N(u_0),
\end{equation*}
and thus with this trace estimate and the div and curl system (\ref{W1}), still with arguments similar as in Step 2 of Section \ref{4},
\begin{equation*}
\int_0^T |\sqrt{\psi} W|^2_{\cd,\Theta}\le C_\kappa\ N(u_0).
\end{equation*}
By patching together all the estimates obtained on each chart defining $\Omega$, we thus deduce that
\begin{equation}
\label{wcd}
\int_0^T \|w\|^2_{\cd}\le C_\kappa\ N(u_0).
\end{equation}
Now, for the pressure, we see that for any $y\in X^\ud(t)=\{\phi\in H^\ud(\Omega)|\ (\bak)_i^k(t) \phi^i,_k\in H^\ud(\Omega)'\}$, for $\varphi$ a solution of
the elliptic problem
\begin{align*}
(\bak)_i^j[ (\bak)_i^k \varphi^i,_k],_j&=(\bak)_i^k(t) y^i,_k\ \text{in}\ (H^\ud)'(\Omega)\\
\varphi&=0\ \text{on}\ \Gamma,
\end{align*}
we have by interpolation that $\varphi\in H^\td(\Omega)$. 
If we once again let $e=(\bak)_i^k \varphi,_k$, and set $v:=y-e$, we have that 
$e\in H^\ud(\Omega)$, $v\in V(t)=\{\phi\in H^\ud(\Omega)|\ (\bak)_i^k(t) \phi^i,_k=0\}$, with $\|e\|_{\ud}+\|v\|_{\ud}\le C \|y\|_{X^\ud(t)}$. 
Now, by proceeding in the same fashion as in Step 2 above, we see that thanks to our 
decomposition and the regularity (\ref{wcd}), 
$w_t\in L^2(0,T;X^\ud(t)')$ with
\begin{equation*}
\int_0^T \|w_t\|^2_{X^\ud(t)'}\le N(u_0).
\end{equation*}
By the Lagrange multiplier Lemma \ref{Lagrangebis}, we then infer
\begin{equation}
\label{qud}
\int_0^T \|q\|^2_{\ud}\le  N(u_0).
\end{equation}
Next, by using $D_{-h}D_h[\psi D_{-h}D_h w]$ as a test function in (\ref{wvar}),
we infer, similarly to how we obtained (\ref{wcd}),  that the estimates
 (\ref{wcd}) and (\ref{qud}) imply that
\begin{equation}
\label{wsd}
\int_0^T \|w\|^2_{\sd}\le  N(u_0).
\end{equation}
We now explain
the additional estimates employed for this higher-order differencing.
We need the fact that independently of any horizontal vector $h$, there
exists a constant $C>0$ such that for
$\text{Supp}\psi+h\subset\Theta$, we have that
\begin{align}
\forall f\in H^\td(\Theta),\ \|\sqrt{\psi}D_h f\|_{\ud,\Theta}\le& C\ \|f\|_{\td,\Theta},\n\\
\forall f\in H^\ud(\Theta),\ \|\sqrt{\psi} D_h f\|_{H^\ud(\Theta)'}\le& C\ \| f\|_{\ud,\Theta}.
\label{int}
\end{align}
The first inequality easily follows by interpolation. For the second one, if $f\in L^2(\Theta)$ we notice that for any $\phi\in H^1(\Theta)$, since the difference quotients are in an horizontal direction, 
\begin{align*}
\int_{\Theta} \sqrt{\psi} D_h f\ \phi &= \int_{\Theta} \sqrt{\psi} f D_{-h}\phi+ 
\int_{\Theta} D_{-h}\sqrt{\psi} f \phi(\cdot-h)\\
&\le C \|f\|_{0,{\Theta}}\|\phi\|_{1,{\Theta}},
\end{align*}
which shows that there exists $C>0$ such that 
\begin{equation*}
\forall f\in L^2(\Theta),\ \|\sqrt{\psi} D_h f\|_{H^1(\Theta)'}\le C\ \| f\|_{0,{\Theta}}.
\end{equation*}
By interpolating with the obvious inequality (for some $C>0$)
\begin{equation*}
\forall f\in H^1(\Theta),\ \|\sqrt{\psi} D_h f\|_{0,\Theta}\le C\ \| f\|_{1,\Theta},
\end{equation*}
we then get (\ref{int}).

Now, the pressure solves the elliptic equation
\begin{subequations}
\label{press}
\begin{align}
\Delta p&=-(\uk)^i,_j u^j,_i\ \text{in}\ \ek(\Omega),\\
p&=-[\sigma \Delta_{\bar g}\en\cdot\nk(\ek)\ \nk(\ek)+\kappa \Delta_{\bar 0}
(w\cdot\nk(\ek))\ \nk(\ek)](\ek^{-1})
\end{align}
\end{subequations}
Using the same change of variables that provides the pressure estimate 
(\ref{ex2}), and using the elliptic estimates for coefficients with
Sobolev class regularity as in \cite{Eb2002}, we find that
\begin{equation*}
\int_0^t \|q\|^2_{\td}\le  N(u_0, \sup_{[0,t]}|\bar w_\kappa|_4, 
\int_0^t\| w\|_\sd^2),
\end{equation*}
where the right-hand side is defined in (\ref{Nu0}).
 Therefore with (\ref{wsd}),
\begin{equation*}
\int_0^t \|q\|^2_{\td}\le  N(u_0, \sup_{[0,t]}\|\bar w_\kappa\|_4).
\end{equation*}

Higher-order regularity results follow successively by appropriate higher-order
difference quotients, leading to, for $n\ge 1$,
\begin{equation}
\label{wqnd}
\int_0^t \|w\|^2_{n+\td}+\int_0^t \|q\|^2_{n-\ud}\le C 
N(u_0,\sup_{[0,t]}|\bar w_\kappa|_{n+2}).
\end{equation}
Now, since $w_t=-(\an)_i^j q,_j$ in $\Omega$, we then infer that for $n\ge 2$,
\begin{equation}
\label{wtqnd}
\int_0^t \|w_t\|^2_{n-\td}\le  N(u_0,\sup_{[0,t]}|\bar w_\kappa|_{n+2}),
\end{equation}
and thus in $[0,t]$,
\begin{equation*}
\|w(t)\|_{13.5}\le \|u_0\|_{13.5}+\sqrt{t}\ N(u_0,\sup_{[0,t]}
|\bar w_\kappa|_{17}).
\end{equation*}
By Lemma \ref{smoothvbisk} (for the smoothing operation given in 
Definition \ref{smoothv} on $C_T$), we have that
\begin{equation}
\label{w12}
\|w(t)\|_{13.5}\le \|u_0\|_{13.5}+\sqrt{t}\ N_0(u_0,C^0_\kappa),
\end{equation}
where we use $C_\kappa^0$ to denote a fixed (nongeneric) constant which
depends on $\kappa$.
\section{Existence of a fixed-point solution of the smoothed $\kappa$-problem
with surface tension}
\label{7}
Let $A: (\bar w \in B_0) \mapsto w$, with $w$ a solution of
(\ref{smoothlinear}). By the relation (\ref{w12}), we see that if 
we take $T_\kappa\in (0,T)$ such that 
$$\sqrt{T_\kappa}N_0(u_0,C^0_\kappa)\le 1,$$  
then
\begin{equation}
\label{stable}
A(C_{T_\kappa})\subset C_{T_\kappa}.
\end{equation}

We now prove that $A$ is weakly lower semi-continuous in $C_{ T_\kappa}$. 
To this end, let $(\bar w^n)_{n=0}^\infty$ be a weakly convergent
sequence (in $L^2(0,T_\kappa;H^{13.5}(\Omega))$) toward a weak limit $\bar w$. Necessarily, $\bar w\in C_{T_\kappa}$. 

By the usual compactness theorems, we have the successive strong convergent
sequences
\begin{align*}
\bar\eta^n&\rightarrow \bar\eta\ \ \text{in}\ L^2(0,T_\kappa;H^{12.5}(\Omega)),\\
 (\bar\eta^n)_\kappa&\rightarrow \bar\eta_\kappa\ \ \text{in}\ L^2(0,T_\kappa;H^{12.5}(\Omega)).
\end{align*}
Now, if we let $w^n=A(\bar w^n)$, we obtain from the stability of 
$C_{T_\kappa}$ by $A$ and (\ref{wtqnd}) the following bounds:
\begin{align*}
\int_0^T \|w^n_t\|^2_{10.5}&\le C N(u_0),\\
\sup_{[0,T]}\| w_n\|_{13.5}&\le 2\|u_0\|_{13.5}+1.
\end{align*}
We thus have the existence of a weakly convergent subsequence $(w^{\sigma(n)})$ in the space $L^2(0,T_\kappa;H^{13.5}(\Omega))$, to a limit 
$l\in C_{T_\kappa} $. By compactness, from our bound on $w^n_t$, 
\begin{align*}
 w^{\sigma(n)}&\rightarrow l\ \ \text{in}\ L^2(0,{T}_\kappa;H^{12.5}(\Omega)).
\end{align*}

From the strong convergence of $(\bar\eta^n)^\kappa$, we then infer from the relations $(\bar a^n_\kappa),_j^i w^n,_j^i=0$ in $\Omega$, that 
\begin{equation}
\label{ldiv}
(\bar a_\kappa),_j^i l,_j^i=0\ \text{ in} \Omega.
\end{equation}

Moreover, we see that 
\begin{equation*}
p^{\sigma(n)}\rightharpoonup p\ \text{in}\ L^2(0,T_\kappa;H^{11.5}(\Omega)),
\end{equation*}
with $p$ the solution of
\begin{align*}
\triangle p&=-(\uk)^i,_j l^j,_i\ \text{in}\ \ek(\Omega),\\
p&=-[\sigma \Delta_{\bar g}\en\cdot\nk(\ek)+\kappa \Delta_0(w\cdot\nk(\ek)](\ek^{-1}).
\end{align*}

  From the relations (\ref{wvar}) for each $n$, we see from the previous 
weak and strong convergence that
\begin{align}
  \int_0^T \langle \bar J_\kappa l_{t},\phi\rangle_{\td}  
&+\kappa \int_0^T [{l}\cdot\nk(\ek),\phi\cdot\nk(\ek)]_1 
\n\\
&-\int_0^T \langle q,(\bak)_i^j \phi^i,_j\rangle_\ud=  \sigma \int_0^T 
[L_{\bar g}\en \cdot \nk(\ek) ,\phi\cdot \nk(\ek)]_0 \ ,
\label{lvar}
\end{align}
which together with (\ref{ldiv}), and the fact that 
$l\in C_{T_\kappa}$ implies that $l=A(\bar w)$. By uniqueness of the limit, we 
then infer that 
\begin{equation*}
w^n\rightharpoonup w\ \text{in}\ L^2(0,T_\kappa;H^{13.5}(\Omega)),
\end{equation*}

By the Tychonoff fixed-point theorem, we then conclude the existence of a 
fixed-point $\bar w=w$ in the closed convex set $C_{T_\kappa}$  of the 
separable Banach space $L^2(0,T_\kappa;H^{13.5}(\Omega))$. This fixed-point 
satisfies the smoothed system (\ref{smooth}), if we denote $\d \eta=\text{Id}+\int_0^\cdot w$ and $u=w\circ{\eta^\kappa}^{-1}$. It is also readily seen that $w$, $q$ and their time derivatives have the regularity stated in Theorem \ref{uniqueweak}.
\end{proof}

\section{Estimates for the divergence and curl}
\label{sec_divcurl}

\begin{definition}[Energy function for the smoothed $\kappa$-problem]
\label{kenergy}
We set
$$
E^{2D}_\kappa(t) = \sum_{k=0}^3 \|\partial_t^k \eta(t)\|_{4.5-k}^2 
+ \|v_{ttt}(t)\|_0^2  \ \text{ and } \ 
E^{3D}_\kappa(t) = \sum_{k=0}^4 \|\partial_t^k \eta(t)\|_{5.5-k}^2 
+ \|v_{tttt}(t)\|_0^2 \,. \ \ 
$$
We use $E_\kappa(t)$ to denote the energy function when the dimension is clear.
\end{definition}
We use these energy functions  
to construct solutions for the Euler equations.  
The increase in the derivative count from the 2D case to the 3D case is
necessitated by the Sobolev embedding theorem.
We will show that solutions of the $\kappa$-problem (\ref{smooth}) have
bounded energy $E_\kappa(t)$ for $t\in [0,T]$ when $T$ is taken sufficiently
small, and that the bound is, in fact, independent of $\kappa$; as such,
we will prove that the limit as $\kappa \rightarrow 0$ of the sequence of 
solutions
to the $\kappa$-problem converges to a solution of the Euler equations.

Our estimates begin with the following 
\begin{lemma}[Divergence and curl estimates] \label{lemma1}
Let ${\mathfrak n}:=\text{dim}(\Omega)=2$ or $3$.
Letting $L_1=\curl$ and $L_2=\div$, and let $\eta_0:=\eta(0)$ and 
$$
M_0:= P(\|u_0\|_{2.5 + {\mathfrak n}}, |\Gamma| _{4 + {\mathfrak n}}, 
\sk\|u_0\|_{1.5+3{\mathfrak n}}, \sk
|\Gamma|_{1+3{\mathfrak n}})
$$
 denote a polynomial function of its arguments.  
Then for $j=1,2$, 
\begin{align*}
&
\sup_{t\in[0,T]}
\|\sk L_j \eta(t)\|^2_{2.5+{\mathfrak n}} 
+ \sum_{k=0}^{{\mathfrak n}+1} \left(
\sup_{t\in[0,T]} \|L_j \partial_t^k\eta(t)\|^2_{1.5+{\mathfrak n}-k} 
 + \int_0^T \|\sk L_j \partial_t^k v\|^2_{2.5+{\mathfrak n}-k} \right) \\
&\qquad\qquad\qquad\qquad\qquad\qquad
\le M_0 + C\,T\, \pskE \,.
\end{align*}
\end{lemma}
\begin{proof}
In Eulerian variables, equation (\ref{smooth.b})  is written as
$u^i_t + u^i,_l (u_\kappa)^l + p,_i=0$, where the transport velocity
is the horizontally smoothed vector $u_\kappa$.
Taking the curl of this equation and using the formula
$(\curl u)^i = \epsilon_{ijk} u^k,_j$ with $\varepsilon_{ijk}$
denoting the permutation symbol, we see that
$\varepsilon_{ijk} [ \partial_t u^k,_j + u^k,_{jl}u_\kappa^l + u^k,_lu_\kappa^l,_j]
=0$.
Thus, defining the bilinear form $B^i(\nabla u, \nabla u_\kappa)= \epsilon_{ijk}
u^k,_l (u_\kappa)^l,_j$, we can write the vorticity equation as
$\frac{D}{Dt}\curl u = B(\nabla u, \nabla u_\kappa)$.   (When the transport velocity
is divergence-free, then $B$ is the familiar vortex-stretching term.)
Composing this equation with $\eta_\kappa$, switching to Lagrangian variables via
the chain rule, and integrating this from $0$ to $t$, we have
\begin{equation} \label{basic0}
\varepsilon_{ijk} v^k,_r {a_\kappa}^r_j 
= \curl u_0 + \int_0^t B_{a_\kappa}(\tau) d\tau, \ \ \ \
B^i_{a_\kappa} := \varepsilon_{ijk} J_\kappa^{-2} v^k,_r {a_\kappa}^r_l (v_\kappa)^l,_m {a_\kappa}^m_j. 
\end{equation}
This is the time-integrated Lagrangian form of the vorticity equation.  We
will need to space-differentiate this equation once more for the estimate
on $\curl\eta$.  Hence, 
\begin{equation} \label{basic}
\varepsilon_{ijk} \nabla v^k,_{r} {a_\kappa}^r_j = \nabla \curl u_0^i 
+\varepsilon_{ikj} v^k,_r \nabla {a_\kappa}^r_j + \int_0^t \nabla B_{a_\kappa}
(\tau) d\tau \,.
\end{equation}

We begin with the estimates for the case that ${\mathfrak n}=2$; we set
$E_\kappa(t)= E_\kappa^{2D}(t)$ and
proceed with the estimate for $\curl \eta$.  Using that
$\nabla v^k,_{r} {a_\kappa}^r_j  = 
\partial_t(\nabla \eta^k,_r {a_\kappa}^r_j)
- \nabla \eta^k,_{r} \partial_t{a_\kappa}^r_j$,
we see that
$$
\varepsilon_{ijk}
\partial_t(\nabla \eta^k,_r {a_\kappa}^r_j)
= \nabla \curl u_0^i +
\varepsilon_{ijk}
\nabla \eta^k,_{r} \partial_t{a_\kappa}^r_j
+\varepsilon_{ikj} v^k,_r \nabla {a_\kappa}^r_j + \int_0^t \nabla B_{a_\kappa}
(\tau) d\tau \,.
$$
Integrating once again in time from $0$ to $t$ yields 
$$
\varepsilon_{ijk} \nabla \eta^k,_r {a_\kappa}^r_j
= t\nabla \curl u_0^i +
\varepsilon_{ijk}
\int_0^t (\nabla \eta^k,_{r} \partial_t{a_\kappa}^r_j
+v^k,_r \nabla {a_\kappa}^r_j) + \int_0^t\int_0^{t'} 
\nabla B_{a_\kappa}  \,,
$$
where
\begin{align}
\nabla B_{a_\kappa} &= \varepsilon_{ijk}\left[ J_\kappa^{-2}
(\nabla v^k,_r {v_\kappa}^l,_m + v^k,_r \nabla{v_\kappa}^l,_m) 
{a_\kappa}^r_l{a_\kappa}^m_j
\right. \nonumber \\
& \qquad \left. 
+ J_\kappa^{-2}
v^k,_r {v_\kappa}^l,_m (\nabla {a_\kappa}^r_l{a_\kappa}^m_j
+{a_\kappa}^r_l\nabla {a_\kappa}^m_j) 
+ (\nabla J_\kappa^{-2}) \,
v^k,_r {v_\kappa}^l,_m {a_\kappa}^r_l{a_\kappa}^m_j
\right]
\label{Ba}
\end{align}
and 
\begin{equation}\label{derivative_a}
\begin{array}{c}
\nabla {a_\kappa}^m_j = J_\kappa^{-1} ( {a_\kappa}^s_r{a_\kappa}^m_j -
{a_\kappa}^m_r{a_\kappa}^s_j) \nabla {\eta_\kappa}^r,_s\,, \\
\partial_t {a_\kappa}^m_j = 
J_\kappa^{-1} ( {a_\kappa}^s_r{a_\kappa}^m_j -
{a_\kappa}^m_r{a_\kappa}^s_j) {v_\kappa}^r,_s \,,  \\ 
\nabla J_\kappa = {a_\kappa}^r_s \nabla {\eta_\kappa}^r,_s \,.
\end{array}
\end{equation}

Since $\|v_\kappa\|_s \le C \|v\|_s$ (and similarly for $\eta_\kappa$), we will
write (\ref{Ba}) and (\ref{derivative_a}) in the following way:
\begin{align*}
\nabla B_{a_\kappa} &\sim J_\kappa^{-2} a_\kappa^2 \nabla v \nabla^2 v
+ J_\kappa^{-3} a_\kappa^3 (\nabla v)^2 \nabla^2 \eta \,,\\
\nabla a_\kappa &\sim J_\kappa^{-1} a_\kappa^2 \nabla^2 \eta\,, \\
\partial_t a_\kappa &\sim J_\kappa^{-1} a_\kappa^2 \nabla v \,,\\
\nabla J_\kappa &\sim a_\kappa \nabla^2 \eta \,,
\end{align*}
where note that we are not distinguishing between $\eta_\kappa$ and $\eta$ or
between $v_\kappa$ and $v$ in the highest-order terms.  The point is that the
precise structure of these equations is not important for our estimates; we need
only be careful with the derivative count appearing in these expressions.  The
power on each expression is merely to indicate the number of times such a term
appears.

Next, with the fundamental theorem of calculus,
\begin{align*}
\varepsilon_{ijk} \nabla \eta^k,_r {a_\kappa}^r_j
= 
\nabla \curl \eta^i +
\varepsilon_{ijk} 
\nabla \eta^k,_r \int_0^t \partial_t {a_\kappa}^r_j\,,
\end{align*}
so that
\begin{equation}\label{zs3}
\nabla (\curl \eta^i - t \curl u_0^i)
= 
\varepsilon_{ijk}\left[
\nabla \eta^k,_r \int_0^t \partial_t {a_\kappa}^r_j +
\int_0^t (\nabla \eta^k,_{r} \partial_t{a_\kappa}^r_j +v^k,_r \nabla {a_\kappa}^r_j) 
\right]
+ \int_0^t\int_0^{t'} \nabla B_{a_\kappa}  \,.
\end{equation}

Let
$$
F:=P(J_\kappa^{-1}, a_\kappa, \nabla v) \text{ and }
F_1:=P(F, \nabla^2 \eta, \nabla^2 v)
$$
denote polynomial functions of their arguments.  We then express (\ref{zs3}) as
$$
\nabla \curl \eta^i 
\sim 
 t \nabla\curl u_0^i+
\nabla^2 \eta \int_0^t F + \int_0^t F\, \nabla^2 \eta + \int_0^t\int_0^{t'} F\,
(\nabla^2 \eta + \nabla^2 v) \,,
$$
and taking two more spatial derivatives yields
\begin{align*}
&\nabla^3 \curl \eta^i 
\sim 
 t \nabla^3\curl u_0^i+
\nabla^4 \eta \int_0^t F + \nabla^3 \eta \int_0^t F_1 
+ \nabla^2 \eta \int_0^t F_1  \\
&\qquad\qquad
+\left(\int_0^t\int_0^{t'} + \int_0^t\right) \left[
F_1 \, (\nabla^3 \eta + \nabla^3 v) + F\, \nabla^4 \eta \right] 
+ \int_0^t\int_0^{t'}\left[ F_1\, \nabla^3 v + F \, \nabla^4 v \right] \,.
\end{align*}

Since $\int_0^t\int_0^{t'} F \, \nabla^4 v = - \int_0^t\int_0^{t'} F_t\,
\nabla^2\eta + \int_0^t F_t\, \nabla^4 \eta$,
\begin{align}
&\nabla^3 \curl \eta^i 
\sim 
 t \nabla^3\curl u_0^i+
\nabla^4 \eta \int_0^t F + (\nabla^3\eta+\nabla^2 \eta) \int_0^t F_1  
\nonumber\\
&\qquad\qquad
+\left(\int_0^t\int_0^{t'} + \int_0^t\right) \left[
F_1 \, (\nabla^3 \eta + \nabla^3 v) + (F+F_t)\, \nabla^4 \eta \right] 
+ \int_0^t\int_0^{t'} F_1\, \nabla^3 v  \,.
\label{zzss0}
\end{align}
We use interpolation to compute $\|\nabla^3\curl \eta\|_{0.5}=\|\curl \eta\|_{3.5}$. 
We begin with the highest-order term: 
\begin{align*}
\left\|\int_0^t (F+F_t)\, \nabla^4 \eta \right\|_0 
&\le \sup_{t\in[0,T]}
\left\|F+F_t\right\|_{L^\infty} \, \left\|\int_0^t\eta\right\|_4 
\le C \sup_{t\in[0,T]}\left\|F+F_t\right\|_2 \, \left\|\int_0^t\eta\right\|_4 \,
\\
\left\|\int_0^t (F+F_t)\, \nabla^4 \eta \right\|_1 
&\le \sup_{t\in[0,T]}
(\|F+F_t\|_{L^\infty} \, \left\|\int_0^t\eta\right\|_5 
+\sup_{t\in[0,T]} \|F+F_t\|_{L^4} \, \left\|\int_0^t\nabla^4\eta\right\|_{L^4}\\
&\qquad \qquad \qquad 
\le C 
\sup_{t\in[0,T]}\|F+F_t\|_2 \, \left\|\int_0^t\eta\right\|_5 \,.
\end{align*}
Since 
$\|F+F_t\|_2 \le C \|F\|_{L^\infty} \,\|v_t\|_{2.5}$, by the 
interpolation theorem 7.17  in \cite{Adams1978}, 
$$
\left\|\int_0^t (F+F_t)\, \nabla^4 \eta \right\|_{0.5} 
\le C \sup_{t\in[0,T]}\|F\|\, \|v_t\|_{2.5} \left\|\int_0^t\eta\right\|_{4.5} 
\le C\,T\, \sup_{t\in[0,T]}\|F\|\, \|v_t\|_{2.5} \, \|\eta\|_{4.5}
$$
The other terms have similar estimates in the $H^{0.5}(\Omega)$-norm, so that
\begin{equation}\label{curl_eta_2D}
\sup_{t\in[0,T]}\|\curl \eta\|_{3.5}^2 \le T \|u_0\|_{4.5} ^2
+ t\sup_{t'\in[0,t]} \|\int_0^t v\|_{3.5}^2
+ T\sup_{t\in[0,T]} \|v_t\|_{2.5} \|\eta\|_{4.5} ^2\,.
\end{equation}

By differentiating (\ref{zzss0}) once more in space, the same interpolation
estimates show that
$$
\sup_{t\in[0,T]}\|\sk\curl \eta\|_{4.5}^2 \le T \|\sk u_0\|_{5.5} ^2
+ T \int_0^T \| v\|_{4.5}^2
+ T\sup_{t\in[0,T]} \|v_t\|_{2.5} \|\sk\eta\|_{5.5} ^2\,.
$$

Next,
we rewrite (\ref{basic0}) as
\begin{equation}\label{curlvss}
\curl v
= \curl u_0 +
\varepsilon_{ijk} v^k,_r \int_0^t \partial_t {a_\kappa}^r_j 
+ \int_0^t B_{a_\kappa} \,.
\end{equation}
Using the fact that $H^{s}(\Omega)$ is a multiplicative algebra for $s>1$,
It follows from (\ref{Ba}) and (\ref{derivative_a}) that
$\sup_{t\in[0,T]}\|\curl v(t)\|_{2.5} \le \|u_0\|_{4.5} +  C\, T\, \pskE$.
Differentiating the above expression for $\curl v$ yields
$$
\curl v_t = \varepsilon_{ijk} v^k,_r \partial_t {a_\kappa}^r_j + B_{a_\kappa}
+\varepsilon_{ijk} \partial_t v^k,_r \int_0^t \partial_t {a_\kappa}^r_j + B_{a_\kappa} \,,
$$
so with the fundamental theorem of calculus and our generic polynomial function $F$,
\begin{equation}\label{zs3b}
\curl v_t \sim P(\nabla u_0) + \nabla v_t \int_0^t F + \int_0^t F_t, \ \ \ F_t \sim
F\, \nabla v_t \,.
\end{equation}
Again using the properties of the multiplicative algebra, we see that
$\sup_{t\in[0,T]}\|\curl v_t(t)\|_{1.5}^{2} \le P(\|u_0\|_{4.5}) +  C\, T\, \pskE$.  From 
the time differentiation of (\ref{zs3b}), 
\begin{align}
\curl v_{tt} &\sim \nabla v_t  F + \nabla v_{tt}\int_0^t F \nonumber \\
& \sim \nabla v_t(0) P(\nabla u_0) + \nabla v_{tt} \int_0^t F
+ \int_0^t [F \nabla v_{tt}  + F \nabla v_t \nabla v_t ] \,.`
\label{zs4}
\end{align}
We must estimate the $H^{0.5}(\Omega)$-norm of the three terms on 
the right-hand side of (\ref{zs4}) by using interpolation.  Let $L$
denote the linear form given by $L(w) = \int_0^t F w$.  Then
$$\|L(w)\|_0 \le C_0 \|\int_0^t w\|_0, \ \ \ 
C_0= \sup_{[0,t]} \|F\|_{L^\infty} .$$
Letting $F_1:=P(J_\kappa^{-1},a_\kappa, \nabla v, \nabla^2\eta, \nabla^2 v)$, it
is easy to check that
$$\|L(w)\|_1 \le C_1 \|\int_0^t w\|_1, \ \ \ 
C_1= \sup_{[0,t]} \|F_1\|_{L^\infty} .$$
By the interpolation theorem 7.17  in \cite{Adams1978},
$$
\|\int_0^t F \nabla v_{tt} \|_{0.5} \le \sqrt{C_0}\sqrt{C_1}  \, \,
\|\int_0^t v_{tt}\|_{1.5} \,,
$$
so that by Jensen's inequality and Sobolev embedding,
$\|\int_0^t F \nabla v_{tt} \|_{0.5}^2 \le C\, T\, \pskE$. All of
the other time-dependent terms in (\ref{zs4}) have the same bound by
the same interpolation procedure.  For the time $t=0$ term, interpolation
provides the estimate 
$$\|P(\nabla u_0) \nabla v_t(0)\|_{0.5} \le
C P(\|u_0\|_{4.5}) \|v_t(0)\|_{1.5} \le 
C P(\|u_0\|_{4.5}) \|q(0)\|_{2.5} \le M_0.$$
The initial pressure $q(0)$ solves the Dirichlet problem
\begin{align*}
\Delta q(0) & = {\mathfrak i}_0:=(u_0)^i,_j({u_0})_\kappa)^j,_i \text{ in } \Omega \,, \\
q(0) &={\mathfrak b}_0:= \frac{\sqrt{g_0}}{\sqrt{{g_0}_\kappa}} {\Pi_0}^i_j g_0^{\alpha\beta}
{\eta_0}^i,_{\alpha \beta} N^i_\kappa + \kappa \Delta_0(u_0\cdot N_\kappa)
\text{ on } \Gamma \,.
\end{align*}
Since $\|q(0)\|_{2.5}\le C(\|{\mathfrak i}_0\|_{0.5} + | {\mathfrak b}_0 |_{2})\le
M_0$, we see that
$\sup_{t\in[0,T]}\|\curl v_{tt}(t)\|_{0.5}^{2} \le M_0 +  C\, T\, \pskE$. 

Differentiating (\ref{zs4}) with respect to time, we see that
$\curl v_{ttt} \sim \nabla v_{ttt} \int_0^t F 
+  F\, \nabla v_{tt} + F \, \nabla v_t \, \nabla v_t$ so that by the
fundamental theorem of calculus 
\begin{align*}
 F\, \nabla v_{tt} + F \, \nabla v_t \, \nabla v_t & = 
F(_0) [ \nabla v_{tt}(0) + \nabla v_t(0)\,\nabla v_t(0)] \\
&\qquad +
\int_0^t [ F\, \nabla v_{ttt} + F \, \nabla v_t \, \nabla v_{tt}
+ F \, \nabla v_t \, \nabla v_{t}\, \nabla v_{t}]\,,
\end{align*}
so that 
\begin{align}
&\int_0^T \| \sk \curl v_{ttt}\|_{0.5}^2  \le
\int_0^T \left\|\sk F(0) [ \nabla v_{tt}(0) + \nabla v_t(0)\,\nabla v_t(0)] 
\right\|_{0.5}^2
\nonumber \\
&\qquad
\int_0^T \left\| \sk \nabla v_{ttt}\int_0^t F\right\|_{0.5}^2 +
\int_0^T \left\| 
\int_0^t [ F\, \sk\nabla v_{ttt} + F \, \nabla v_t \, \sk\nabla v_{tt}
+ F \, \nabla v_t \, \nabla v_{t}\, \sk\nabla v_{t}] \right\|_{0.5}^2
\,.
\label{zs6}
\end{align}
We repeat the
interpolation estimates between $L^2(\Omega)$ and $H^1(\Omega)$ just as
for the estimates for $\curl v_{tt}$;  for example,
$$
\left\|\int_0^t F \sk\nabla v_{ttt} \right\|_{0.5} \le 
\sqrt{C_0}\sqrt{C_1}  \, \,
\left\|\int_0^t \sk v_{ttt}\right\|_{1.5} \,,
$$
so that by Jensen's inequality and Sobolev embedding,
$$
\|\int_0^t F \sk \nabla v_{ttt} \|_{0.5}^2 \le C\, t\, \sup_{[0,t]}E_\kappa \,
\|\sk v_{ttt}\|_{0.5}^2\,.$$
Thus, integrating from $0$ to $T$ gives the estimate
$$
\int_0^T\|\int_0^t F \sk \nabla v_{ttt} \|_{0.5}^2 \le 
C\, T\, \pskE \,
\int_0^T \|\sk v_{ttt}\|_{0.5}^2 \le  C\, T\, \pskE \,.
$$
The other time-dependent terms in (\ref{zs6}) have the same bound by the same 
argument. The time $t=0$ terms require analysis of the elliptic problem for
$q_t(0)$:
\begin{align*}
\Delta q_t(0) & \sim  {\mathfrak i}_1 := \nabla [ P(\nabla u_0)\, \nabla q(0)] +
F\, \nabla v_t(0) \text{ in } \Omega \,, \\
q_t(0) &\sim {\mathfrak b}_1:= 
Q(\partial \eta_0) \partial^2 u_0 + Q(\partial \eta_0) \partial u_0 
\partial ^2 \eta_0 \\
& \qquad \qquad\qquad
+ \kappa \Delta_0 (v_t(0) \cdot N_\kappa) + \kappa \Delta_0( u_0 \cdot
Q(\partial {\eta_0}_\kappa)  \partial {u_0}_\kappa)
\text{ on } \Gamma \,.
\end{align*}
By interpolation estimates (as above),
$$
\int_0^T \left\|\sk F(0) [ \nabla v_{tt}(0) + \nabla v_t(0)\,\nabla v_t(0)] 
\right\|_{0.5}^2
\le \kappa T\,\|P (\nabla u_0)\|_{L^\infty}^2 \|v_{tt}(0)\|_{1.5}^2,
$$
and since time-differentiation of the Euler equations shows that
$$
v_{tt}(0) = - \nabla u_0\, \nabla q(0) - \nabla q_t(0),
$$
interpolation provides the estimate
\begin{align*}
\sk\|v_{tt}(0)\|_{1.5} & \le \sk\| \nabla^2 u_0 \, \nabla q(0)\|_{0.5} +
\sk\| \nabla u_0 \, \nabla^2 q(0)\|_{0.5} + \sk\|q_t(0)\|_{2.5} \\
& \le M_0 + \sk\|{\mathfrak i}_1\|_{0.5} + \sk|{\mathfrak b}_1|_2 \\
& \le M_0 + \sk|{\mathfrak b}_1|_2 \\
\end{align*}
where we have used the elliptic estimate $\|q(0)\|_{2.5} \le M_0$ (from above)
for both the second and third inequalities.  (The remaining estimate  for
$|{\mathfrak b}_1|_2 $ places the regularity constraints on the polynomial
function $M_0$ in the hypothesis of the lemma.)  Because $H^2(\Gamma)$ is
a multiplicative algebra, the bound for
$\sk|{\mathfrak b}_1|_2 $ is controlled by the highest-order terms
$\sk |v_t(0)|_4$ and $\sk |{u_0}_\kappa|_5\le \sk C \|u_0\|_{5.5}$.
 Now, 
$$
\sk |v_t(0)|_4 \le \sk \|q(0)\|_{5.5} \le  \sk 
\| {\mathfrak i}_0\|_{3.5} + \sk | {\mathfrak b}_0|_5\,,
$$
and
$\| {\mathfrak i}_0\|_{3.5}$ is bounded by $P(\|u_0\|_{4.5})$ while
the highest-order terms in $\sk | {\mathfrak b}_0|_5$ require bounds
on $\sk \|\eta_0\|_{7.5}$ and $\sk \|u_0\|_{7.5}$  With our
definition of $M_0$, we see that 
$$\kappa T\,\|P (\nabla u_0)\|_{L^\infty}^2 \|v_{tt}(0)\|_{1.5}^2
\le M_0$$ and hence $\int_0^T \|\sk \curl v_{ttt}\|_{0.5}^2 \le M_0 + 
C\, T\, \pskE$.  

The proof that 
$\int_0^T \|\sk \curl v_{tt}\|_{1.5}^2 \le M_0 + 
C\, T\, \pskE$ is essentially identical.

The divergence estimates begin with the fundamental equation
${a_\kappa}^j_i v^i,_j =0$.  By taking one derivative of this equation and
integrating-by-parts in time, we find that
$$
\nabla \div \eta = \nabla \eta^i,_j \int_0^t \partial_t {a_\kappa}^j_i
+ \int_0^t( \partial_t {a_\kappa}^j_i \nabla \eta^i,_j -
\nabla {a_\kappa}^j_i v^i,_j) \,.
$$
Computing the $H^{2.5}(\Omega)$-norm of this equation yields the estimate
$\sup_{t\in[0,T]} \| \div \eta(t)\|_{3.5}^2 \le M_0 + C\, T\, \pskE$.
The divergence estimates for $v$, $v_t$, $v_{tt}$, $\sk v_{tt}$,
and $\sk v_{ttt}$ follow the same argument as the corresponding curl estimates.

\vspace{.2 in}

In the case that ${\mathfrak n}=3$, the estimates are found in the same way,
with one minor change.  Set $E_\kappa(t)=E_\kappa^{3D}(t)$.  The estimates for $\curl \eta$, which rely on Sobolev
embedding, require greater regularity on $v_t$.
The estimate (\ref{curl_eta_2D} becomes
\begin{equation}\nonumber
\sup_{t\in[0,T]}\|\curl \eta\|_{3.5}^2 \le T \|u_0\|_{4.5} ^2
+ t\sup_{t'\in[0,t]} \|\int_0^t v\|_{3.5}^2
+ T\sup_{t\in[0,T]} \|v_t\|_{3} \|\eta\|_{4.5} ^2\,,
\end{equation}
and similarly,
$$
\sup_{t\in[0,T]}\|\sk\curl \eta\|_{4.5}^2 \le T \|\sk u_0\|_{5.5} ^2
+ T \int_0^T \| v\|_{4.5}^2
+ T\sup_{t\in[0,T]} \|v_t\|_{3} \|\sk\eta\|_{5.5} ^2\,.
$$
\end{proof}

\section{Some geometric identities}
\label{Section_apriori}
\def\nk{{\tilde n_\kappa}}
\def\vk{ {\tilde{v}_\kappa} }
\def\tk{{\tau_\kappa}}
\def\ek{{\tilde \eta_\kappa}}
\def\gk{{\tilde g_\kappa}}
\def\pk{{\Pi^\kappa}}
\def\ak{{{\tilde a}^k}}
\def\intT{\int_{T'-\epsilon}^{T'}}
\def\intT{\int_{0}^{T}}

We will usually omit writing $dS_0$ in our surface integrals, and for
convenience we set $\sigma=1$.
Let $\Pi^i_j$ denote the projection operator  onto the  direction
normal to $\eta(\Gamma)$, defined as
$\Pi^i_j = \delta^i_j - g^{\alpha \beta} \eta^i_{,\alpha} \delta_{jl} 
\eta^l_{,\beta}$, where $g$ is the induced metric on $\eta(\Gamma)$ defined
in (\ref{gstuff}).  The mean curvature vector motivates us to introduce
the projection operator $\Pi$.  In particular, we have the important formula 
\begin{align}
-\sqrt{g}Hn\circ\eta& = \sqrt{g}\Delta_g(\eta)\nonumber\\
& = \sqrt{g}\left[
g^{\alpha\mu}( \delta^{ij} - g^{\nu\beta}\eta^i,_\beta\eta^j,_\nu) 
\eta^j,_{\mu \alpha} 
+( g^{\alpha\beta} g^{\mu\nu} - g^{\alpha\nu} g^{\mu\beta})
\eta^i,_\beta \eta^j,_\nu \eta^j,_{\mu\alpha})
\right] \nonumber\\
&= \sqrt{g} g^{\mu\alpha} \Pi^i_j \eta^j,_{\mu\alpha}\,. \label{importantF}
\end{align}
where the last equality follows since
$( g^{\alpha\beta} g^{\mu\nu} - g^{\alpha\nu} g^{\mu\beta})
\eta^i,_\beta \eta^j,_\nu \eta^j,_{\mu\alpha} =0$.
For a vector field $F$ on $\Gamma$, $\Pi \, F = 
[n \cdot F] \, n$, i.e., $\Pi = n\otimes n$.

We let
\begin{equation}\label{Q}
{Q}(\partial \eta) = f_1(\partial \eta)/f_2(\sqrt{g})
\,,
\end{equation}
denote a generic rational function where $f_1$ and $f_2$ are smooth functions.
We record for later use that
that $n= \frac{a^T N}{|a^TN|}=
\frac{\eta,_1 \times \, \eta,_2}{|\eta,_1 \times \, \eta,_2|}$ and
that $|a^TN|= \sqrt{\det g}$ on $\Gamma$, as
\begin{align*}
|\eta,_1 \times \, \eta,_2|^2
& = \varepsilon_{ijk}\eta^j,_1 \eta^k,_2 \,\varepsilon_{irs}\eta^r,_1 \eta^s,_2\\ 
&= (\delta_{jr}\delta_{ks}-\delta_{js}\delta_{kr}) \eta^k,_1 
\eta^j,_s \eta^r,_1 \eta^s,_2 = |\eta,_1|^2|\eta,_2|^2 -  
[\eta,_1 \cdot \eta,_2]^2 =  \det g \,,
\end{align*}
where $\varepsilon_{ijk}$ denotes the permutation symbol of $(1,2,3)$.
We will use the symbol 
${Q}$ to denote any smooth (tensor) function that can be represented as
(\ref{Q}). 

\begin{remark}
The $L^\infty$-norm of the numerator of ${Q}$ is bounded by 
a polynomial of the energy function, while the 
$L^\infty$-norm of the denominator of ${Q}$ is uniformly controlled
by (\ref{deteta.a}).   Thus, the generic constant $C$ which appears in the
following inequalities may depend on a polynomial of $\det g_0$.
In particular, $\|{Q}(\partial \eta)\|_{\L^\infty} 
\le C(\det g_0) \|P (\partial \eta)\|_{L^\infty}$.
\end{remark}

For a vector field $F$ on $\Gamma$, 
$F \cdot N =  F \cdot n + F \cdot (N-n)$ and 
$$|N-n|_{L^\infty}
 \le \int_0^t |n_t|_{L^\infty} = \int_0^t |{Q}(\partial \eta) 
\partial v|_{L^\infty} \le C\, t\, P(E_\kappa(t))\,,$$
the last inequality following from (\ref{deteta.a}). 
If $|\Pi F |_s \le M_0 + C P(E_\kappa(t))$, then $|F\cdot N|_s$ satisfies the
same bound.

\section{$\kappa$-independent estimates for the smoothed problem and
existence of solutions in 2D}\label{kapriori}
All of the variables in the smoothed $\kappa$-problem (\ref{smooth}) implicitly
depend on the parameter $\kappa$.  In this section,
where we study the asymptotic behavior of the solutions to (\ref{smooth}) 
as $\kappa \rightarrow 0$, we will make this dependence explicit by placing 
a $\sim$  over each of the variables.  We set $E_\kappa(t)= E_\kappa^{2D}(t)$.
\begin{remark}\label{2Dvs3D}
The only difference between the 2D and 3D cases arise from the embedding of
$\tilde v_t \in L^\infty(\Omega)$.  In 2D, $\tilde v_t \in H^{2.5}(\Omega)$ is 
sufficient,
while in 3D, we need $\tilde v_t \in H^{3}(\Omega)$.
\end{remark}

The pressure function $\tilde q$  can be formulated to solve either a Dirichlet
problem with boundary condition (\ref{kbc}) or a Neumann problem found
by taking the inner-product of the Euler equations with ${\tilde a_\kappa}^T N$.
We use the latter.

\begin{lemma}[Pressure estimates]
With $(\tilde v, \tilde q)$ a solution of the $\kappa$-problem (\ref{smooth})
\begin{align}
& \|\tilde q(t)\|^2_{3.5}
+\|\tilde q_{t}(t)\|^2_{2.5}
+\|\tilde q_{tt}(t)\|^2_{1}
 \le C \, P( E_\kappa(t)) \,. \label{qestimate} 
\end{align}
\end{lemma}
\begin{proof}
Denoting $\tilde a_\kappa$ by $A$,
we define the divergence-form elliptic operator $L_A$ and corresponding
Neumann boundary operator $B_A$ as
$$
L_A = \partial_j ( \tilde J_\kappa^{-1} A^j_i A^l_i \partial_l)\,, \ \ \ 
B_A = \tilde J_\kappa^{-1} A^j_i A^l_i  N_j \partial_l \,.
$$

For $k=0,1,2$, we analyze the Neumann problems
\begin{align}
L_a (\partial^k_t \tilde q) = f_k \ \text{ in } \ \Omega \ \ & \ \ 
B_a(\partial^k_t \tilde q) = g_k \ \text{ on } \ \Gamma \label{Neumann}
\end{align}
with
\begin{alignat*}{2}
f_0&=  \partial_t A^j_i \, \tilde v^i,_j \qquad &&g_0= - \tilde v_t\cdot 
\sqrt{\tilde g_\kappa} \tilde n_\kappa\\
f_1&= 
-{L_A}_t(\tilde q) - \partial_t^2A^j_i\, \tilde v^i,_j- \partial_tA^j_i\, \tilde v_t^i,_j
\qquad  &&g_1= {B_A}_t(\tilde q) - \tilde v_t \cdot (\sqrt{\tilde g_\kappa} \tilde n_\kappa)_t
- \tilde v_{tt} \cdot \sqrt{\tilde g_\kappa} \tilde n_\kappa \\
f_2&= 
-2{L_A}_{t}(\tilde q_t) -{L_A}_{tt}(\tilde q) - \partial_t^3A^j_i\, \tilde v^i,_j
  \qquad &&g_2= 2{B_A}_t(\tilde q_t) + {B_A}_{tt}(\tilde q) - \tilde v_{ttt}\cdot \sqrt{\tilde g_\kappa}
                                                            \tilde   n_\kappa \\
&\qquad 
- 2\partial_t^2A^j_i\, \tilde v_t^i,_j-\partial_tA^j_i\, \tilde v_{tt}^i,_j
  \qquad &&
 \qquad -2\tilde v_{tt}\cdot (\sqrt{\tilde g_\kappa} \tilde n_\kappa)_t 
 -\tilde v_{t}\cdot (\sqrt{\tilde g_\kappa} \tilde n_\kappa)_{tt}
\,.
\end{alignat*}

For $s\ge 1$, elliptic estimates 
provide the inequality
\begin{equation}\label{elliptic}
\| \partial^k_t \tilde q(t) \|_s \le C_s [P( \|\eta\|_{4.5})\|f_k\|_{s-2} + 
|g_k|_{s-3/2} + \| \tilde q \|_0] \,,
\end{equation}
where $\|\cdot \|_{-1}$ denotes the norm on $[H^1(\Omega)]'$.  We
remark that the usual $H^s$ elliptic estimates
require that coefficient have the regularity $\partial^{s-1}(A^l_i
A^j_i) \in L^\infty(\Omega)$, however  
$\partial^{s-1}(A^l_i A^j_i) \in L^2(\Omega)$ is sufficient.  See
see \cite{Eb2002} or the quasilinear estimates in
\cite{Taylor1996}.

As we cannot guarantee that solutions $\tilde q$ to the $\kappa$-problem (\ref{smooth})
have zero average, we use $\|\tilde q\|_0 \le C\|\tilde q\|_1$ and the $H^1$ elliptic estimate
for the Dirichlet problem $L_A(q)= f_0$ in $\Omega$ with 
$-q= \Delta_{\tilde g}\tilde \eta\cdot \tilde n_\kappa + \kappa \Delta_0( \tilde v\cdot \tilde n_\kappa)$ on
$\Gamma$.  Thus, $\|\tilde q\|_1 \le C (\|f_0\|_0 +
|\Delta_{\tilde g}\tilde \eta\cdot \tilde n_\kappa + \kappa \Delta_0( \tilde v\cdot \tilde n_\kappa)|_{0.5})
\le C\, P(E_\kappa(t))$.

From (\ref{derivative_a}), it is clear that
$\|f_0\|_{1.5}^2 + |g_0|_{2}^2 \le C P(E_\kappa(t))$; thus, from the elliptic
estimate, 
\begin{equation}\label{q0}
\|\tilde q\|_{3.5}^2  \le C P(E_\kappa(t)).
\end{equation}
Next, we must show
that $\|f_1\|_{0.5}^2 + |g_1|_1^2 \le C \, P(E_\kappa(t))$.
But $f_1 \sim P(\tilde J_\kappa^{-1},A,\nabla \tilde v_\kappa) ( [\nabla \tilde v_t]^2 + \nabla^2 \tilde q)$
so that with (\ref{q0}), $\|f_1\|_{0.5}^2\le C\, P(E_\kappa(t))$, with the same bound
for $|g_1|_1^2$, so that $ \|\tilde q_t\|_{2.5}^2  \le C P(E_\kappa(t))$. Using this,
we find,
in the same fashion, that $\|f_2\|_0 \le C\, P(E_\kappa(t))$.  The
normal trace theorem, read in Lagrangian variables, states that
if $\tilde v_{ttt} \in L^2(\Omega)$ with $\|A^j_i \tilde v_{ttt}^i,_j\|_0 \in  L^2(\Omega)$,
then $ \tilde v_{ttt} \cdot \sqrt{\tilde g_\kappa} \tilde n_\kappa \in H^{-0.5}(\Gamma)$ with 
the estimate
$| \tilde v_{ttt}\cdot \sqrt{\tilde g_\kappa} \tilde n_\kappa|_{-0.5}^2 \le C \, P(E_\kappa(t))$.
Since $\|\operatorname{Tr}(A\, \nabla \tilde v_{ttt})\|_0^2 = \|\operatorname{Tr}( 3
A_{t}\, \nabla \tilde v_{tt} + 3A_{tt}\, \nabla \tilde v_{t} + A_{ttt}\, \nabla \tilde v)\|_0^2 
\le C\, P(E_\kappa(t))$, and using the above estimates for $\tilde q$ and $\tilde q_t$, we find
that $|g_2|_{-0.5} \le C\, P(E_\kappa(t))$, thus completing the proof.
\end{proof}

Our smoothed $\kappa$-problem (\ref{smooth}) uses the boundary condition
(\ref{smooth.e}) which we write as
\begin{equation}\label{kbc}
\tilde q\, \nk  =  \frac{\sqrt{\tilde g}}{\sqrt{\tilde g_\kappa}}
\tilde H\, \tilde n \cdot \nk\, \nk- \kappa \Delta_{ 0}
(\tilde v\cdot \nk)\, \nk\,,
\end{equation}
where (we remind the reader) $\kappa>0$ is the artificial viscosity, 
$\Delta_{0}=
\sqrt{\tilde g_\kappa}^{-1}
\partial_\alpha (\sqrt{g_0} g_0^{\alpha \beta} \partial_\beta)$, 
$\tilde n$ is the unit normal along the boundary 
$\tilde \eta(t)(\Gamma)$ and $\nk$ is
the unit normal along the smoothed $\kappa$-boundary 
$\tilde \eta_\kappa(t)(\Gamma)$.


We begin with an energy estimate for the third time-differentiated problem.
Although we are doing the estimates for the 2D domain $\Omega$, we keep the
notation of the 3D problem as well as terms that only arise in 3D when
differentiating the mean curvature vector. Thus, when we turn to the
3D problem in Section \ref{kapriori3}, the modifications will be trivial.

\begin{lemma}[Energy estimates for the third time-differentiated 
$\kappa$-problem]  
\label{lemma10.2}
For $M_0$ taken as in Lemma \ref{lemma1} and $\delta >0$,  solutions of
the $\kappa$-problem (\ref{smooth}) satisfy:
\begin{align}
&\sup_{t\in[0,T]}\left[
\|\tilde v_{ttt}\|_0^2 + |\tilde v_{tt}\cdot\tilde  n|_1^2
\right]
+\int_0^T |\sk\partial_t^3\tilde v \cdot \nk|_1^2 
\le M_0 + T \, P(\sup_{t\in [0,T]} E_\kappa(t))
+ \delta  \sup_{t\in [0,T]} E_\kappa(t)  \nonumber \\
&\qquad
+ C \sup_{t\in [0,T]} [
P(\|\tilde v_t\|_{2.5}^2) + P(\|\tilde v\|_{3.5}^2)+ P(\|\tilde\eta\|_{4.5}^2)]
+ C P(\|\sk \tilde v_{tt}\|_{L^2(0,T;H^{2.5}(\Omega))}^2) \,.
\label{ss_kttt}
\end{align}
\end{lemma}

\begin{proof}
Letting $A= \tilde a_\kappa$, and
testing $\partial_t^3 (\tilde J_\kappa \tilde v_t^i) +
\partial_t^3({A}^k_i \tilde q,_k) =0$ with  with 
$\partial_t^3 \tilde v^i$ shows that
\begin{equation}\label{eulerttt}
\int_0^T{\frac{1}{2}} \int_{\Omega} \partial_t^3(\tilde J_\kappa \tilde v_t^i)
\partial_t^3 \tilde v^i
- \int_0^T\int_{\Omega} \partial^3_t(  A^k_i q) \ \partial_t^3 \tilde 
v^i,_k = -
\int_0^T\int_{\Gamma} \partial_t^3 (\sqrt{\tilde g_\kappa}
\tilde q \nk (\tilde \eta_\kappa)) \cdot \partial_t^3 \tilde v \ 
dS_0  \,.
\end{equation}

\noindent
{\bf Step 1. Boundary integral term.}
We rewrite the modified boundary condition (\ref{kbc}) as
\begin{equation}\label{kbca}
\tilde q\, \nk  = \frac{\sqrt{g}}{\sqrt{g_\kappa}}
\left[\tilde H\tilde n + \tilde H\tilde n\cdot(\nk -\tilde n)\, \tilde n 
+ \tilde H\tilde n\cdot\nk\, (\nk -\tilde n)\right]
- \kappa \Delta_{0}(\tilde v\cdot \nk)\, \nk\,.
\end{equation}
We first consider the boundary integral on the right-hand side of 
(\ref{eulerttt}) with only the first term on the right-hand side of 
(\ref{kbca}):
\begin{align}
&-\int_0^T\int_{\Gamma} \partial^3_{t} (\sqrt{\tilde g}\tilde H\tilde n^i \circ \eta) \partial^3_t \tilde v^i \, 
dS_0  \nonumber \\
&=
-\int_0^T\int_{\Gamma}\sqrt{\tilde g} \tilde g^{\alpha\beta} \Pi^i_j \partial^2_t \tilde v^j_{,\beta}
\partial^3_t \tilde v^i_{,\alpha} 
- \int_0^T\int_{\Gamma} 
\sqrt{\tilde g}[ \tilde g^{\mu\nu}\tilde g^{\alpha\beta} -
\tilde g^{\alpha\nu}\tilde g^{\mu\beta}] \tilde \eta^j_{,\nu} \partial^2_t \tilde v^j_{,\mu} \,
\tilde \eta^i_{,\beta} \partial^3_t \tilde v^i_{,\alpha} \nonumber \\
&\qquad  \qquad
+\int_0^T \int_{\Gamma} Q_{ij}^{\alpha\beta}(\partial \tilde \eta, \partial \tilde v) \,
\partial_t \tilde v^j_{,\beta}\, \partial^3_t \tilde v^i_{,\alpha}
+\int_0^T \int_{\Gamma} Q_{i}^{\alpha}(\partial \tilde \eta, \partial \tilde v) \,
 \partial^3_t \tilde v^i_{,\alpha}  \label{Hn_ttt} \\
& =:  I + II + III + IV \,. \nonumber
\end{align}
The first term $I$ on the right-hand side of  (\ref{Hn_ttt}) is given by
\begin{align*}
I &= \left.
-\frac{1}{2}\int_{\Gamma} \sqrt{\tilde g} (\Pi^i_j \partial^2_t\tilde v^j_{,\beta})
\tilde g^{\alpha\beta}(\Pi^i_k \partial^2_t\tilde v^k_{,\alpha})\right]^T_0
+\int_0^T\int_{\Gamma} Q^{\alpha\beta}_{jk}(\partial \tilde \eta,\partial \tilde v)
\partial^2_t \tilde v^k_{,\alpha}\, \partial^2_t \tilde v^j_{,\beta} \,, \\
\end{align*}
where we use the notation $f]^T_0 = f(T) -f(0)$.
Since $\Pi^i_j \tilde v^j_{tt},_\beta =(\Pi^i_j \tilde v^j_{tt}),_\beta  -{\Pi^i_j},_\beta
\tilde v^j_{tt}$ and ${\Pi^i_j},_\beta = {Q}^{i\gamma}_{jl}(\partial\tilde \eta)
\tilde  \eta^l,_{\gamma \beta}$ with ${Q}(\partial\tilde \eta)$ defined by (\ref{Q}),
for $\delta >0$,
\begin{align*}
-\frac{1}{2}\int_{\Gamma} \sqrt{\tilde g} (\Pi^i_j \partial^2_t\tilde v^j_{,\beta})
\tilde g^{\alpha\beta}(\Pi^i_k \partial^2_t\tilde v^k_{,\alpha})
& \le -{\frac{1}{2}} | \Pi \tilde v_{tt} |_1^2  + \delta | \Pi \tilde v_{tt} |_1^2
+ (1+C_\delta) 
\left| {Q}^{i\alpha\beta}_{jl}(\partial \tilde \eta)\, \tilde \eta^l,_{\alpha
\beta}\, \tilde v_{tt}^j\right|_0^2 \,,
\end{align*}
where the constant $C_\delta$ depends inversely on $\delta$.  Since
for any $t \in [0,T]$
$$
\left| \left( {Q}^{i\alpha\beta}_{jl}(\partial \tilde \eta)\, \tilde \eta^l,_{\alpha
\beta}\, \tilde v_{tt}^j\right)_t \right|_0^2(t)  \le C P(E_\kappa(t)),
$$
it follows that
$$
I\le -{\frac{1}{2}} \sup_{t\in[0,T]}| \Pi \tilde v_{tt}|_1^2 + M_0(\delta)+ \delta \sE + C \, T \psE \,.
$$

The second term $II$ requires some care (in the way in which the terms are
grouped together).  Letting
\begin{align}
&
{\mathcal A}^1=
\left[
\begin{array}{cc}
\tilde \eta_{,1}\cdot \partial_t^2\tilde v_{,1} &\tilde \eta_{,1}\cdot \partial_t^2\tilde v_{,2} \\
\tilde \eta_{,2}\cdot \partial_t^2\tilde v_{,1} &\tilde \eta_{,2}\cdot \partial_t^2\tilde v_{,2} 
\end{array}
\right], \
{\mathcal A}^2=
\left[
\begin{array}{cc}
\tilde v_{,1}\cdot \partial_t^2\tilde v_{,1} &\tilde \eta_{,1}\cdot \partial_t^2\tilde v_{,2} \\
\tilde v_{,2}\cdot \partial_t^2\tilde v_{,1} &\tilde \eta_{,2}\cdot \partial_t^2\tilde v_{,2} 
\end{array}
\right], \nonumber \\
& \qquad\qquad\qquad\qquad
{\mathcal A}^3=
\left[
\begin{array}{cc}
\tilde \eta_{,1}\cdot \partial_t^2\tilde v_{,1} &\tilde v_{,1}\cdot \partial_t^2\tilde v_{,2} \\
\tilde \eta_{,2}\cdot \partial_t^2\tilde v_{,1} &\tilde v_{,2}\cdot \partial_t^2\tilde v_{,2} 
\end{array}
\right] \,,
\label{detA}
\end{align}
we find that
\begin{align}
II &= \int_0^T\int_{\Gamma} \det{\tilde g^{-{\frac{1}{2}}}} \, ( \partial_t \det{{\mathcal A}^1}
- \det{{\mathcal A}^2}- \det{{\mathcal A}^3}) \nonumber \\
   &= \int_0^T\int_{\Gamma} -(\det{g^{-{\frac{1}{2}}}})_t \,  \det{{\mathcal A}^1}
- \det \tilde g^{-{\frac{1}{2}}}(\det{{\mathcal A}^2}+ \det{{\mathcal A}^3}) + \left.\int_{\Gamma}
\det{\tilde g^{-{\frac{1}{2}}}} \det{{\mathcal A}^1} \right]_0^T\,.
\label{term_II}
\end{align}
For $\alpha=1,2$, let $V_\alpha = \tilde \eta_{,\alpha} \cdot \partial_t^2 \tilde v$;
thus ${V_\alpha},_\beta = {\mathcal A}^1_{\alpha\beta} 
+ \tilde \eta^i,_{\alpha\beta} \tilde v^i_{tt}$
so that
\begin{align*}
\det{{\mathcal A}^1_{\alpha\beta}} & = \det( {V_\alpha},_\beta 
- \tilde \eta^i,_{\alpha\beta} \tilde v^i_{tt}) \\
& =\det{{V_\alpha},_\beta} - \det (\tilde \eta^i,_{\alpha\beta} \tilde v^i_{tt}) 
+ P_{ij}(\partial^2 \tilde \eta) \tilde v^i_{tt} \tilde v^j_{tt} 
+ P^\alpha_{ij}(\partial^2 \tilde \eta) \tilde v^i_{tt} \tilde v_{tt}^j,_\alpha
\,.
\end{align*}
With $A = \cf (\partial V)$, $\det \partial V = A^{\beta\alpha}
{V_\alpha},_\beta$.  It follows
that 
$$\int_\Gamma \det{\tilde g^{-{\frac{1}{2}}}} \det \partial V = - \int_\Gamma 
(\det{\tilde g^{-{\frac{1}{2}}}})_{,\beta} A^{\beta\alpha} V,_\alpha,$$ 
as 
$A^{\beta\alpha},_\beta =0$ since $A$ is the cofactor matrix.  Hence,
\begin{equation}\label{ssss7}
\int_\Gamma \det {\tilde g^{-{\frac{1}{2}}}} \det {\mathcal A}_1 
= 
\int_\Gamma 
 P_{ij}(\partial^2 \tilde \eta) \tilde v^i_{tt} \tilde v^j_{tt} 
+ P^\alpha_{ij}(\partial^2 \tilde \eta) \tilde v^i_{tt} \tilde v_{tt}^j,_\alpha
\,,
\end{equation}
so that
\begin{align*}
II & \le \int_0^T \int_{\Gamma} Q^{\alpha\beta}_{ij} 
(\partial \tilde \eta, \partial \tilde v) \,
\partial^2_t \tilde v^i_{,\alpha}\, \partial^2_t \tilde v^j_{,\beta}
+
\left.\
\int_\Gamma 
[ P_{ij}(\partial^2 \tilde \eta) \tilde v^i_{tt} \tilde v^j_{tt} 
+ P^\alpha_{ij}(\partial^2 \tilde \eta) \tilde v^i_{tt} \tilde v_{tt}^j,_\alpha]
\right]^T_0 \,.
\end{align*}
By the fundamental theorem of calculus and Young's inequality, for
$\delta >0$,
\begin{align*}
&\int_\Gamma[P^\alpha_{ij}(\partial^2\tilde \eta) \tilde v^i_{tt} \tilde v^j_{tt},_\alpha](T)  =
\int_\Gamma[P^\alpha_{ij}(\partial^2\tilde \eta) \tilde v^i_{tt}](0)\, \tilde v^j_{tt},_\alpha(T)
+ \int_\Gamma\int_0^T
{[P^\alpha_{ij}(\partial^2\tilde \eta) \tilde v^i_{tt}]}_t \, dt\, \tilde v^j_{tt},_\alpha(T) \\
&\qquad \le M_0(\delta) + \delta \|\tilde v_{tt}(T)\|_{1.5}^2
+ T\left ( \int_\Gamma \sup_{t\in [0,T]} | {(P^\alpha_{ij} (\partial^2\tilde \eta)
\tilde v^i_{tt})}_t |^2 dx \right)^{\frac{1}{2}} \|\tilde v^j_{tt}(T)\|_{1.5}^{\frac{1}{2}}\,.
\end{align*}
Since ${[P^\alpha_{ij}(\partial^2
\tilde \eta) \tilde v^i_{tt}]}_t \in L^\infty(0,T; L^2(\Gamma))$, we conclude that
$$
II\le M_0(\delta)+ \delta \sE + C \, T \psE \,.
$$
A temporal integration by parts in 
the third and fourth terms on the right-hand side of (\ref{Hn_ttt}) yields
\begin{align*}
III+ IV & =
-\int_0^T \int_{\Gamma}  { [ Q^{\alpha\beta}_{ij} 
(\partial \tilde \eta, \partial \tilde v) \tilde v^j_t,_\beta
+ Q^\alpha_i(\partial\tilde \eta, \partial \tilde v)] }_t
 \, \tilde v^j_{tt},_\beta \\
& \qquad\qquad\qquad\qquad +
\left.\int_{\Gamma} [Q^{\alpha\beta}_{ij} (\partial \tilde \eta,\partial \tilde v) 
\tilde v^j_{t},_\beta
+ Q^\alpha_i(\partial\tilde \eta,\partial \tilde v) ] \tilde v^i_{tt},_\alpha
\right]^T_0  \,,
\end{align*}
which has the same bound as term $II$; it follows that
\begin{align}
&
\int_0^T\int_{\Gamma} \partial_t^3 (\sqrt{\tilde g}\tilde H\tilde n \circ \tilde \eta) \cdot \partial_t^3 \tilde v \nonumber \\
& \qquad\qquad\qquad\qquad
= -{\frac{1}{2}}\sup_{t\in[0,T]} |\Pi \tilde v_{tt}|_1^2 +
M_0(\delta)+ \delta \sE + C \, T \psE \,.
\label{ttt1}
\end{align}
\begin{remark} \label{remark_det}
The determinant structure which appears in (\ref{term_II}) is 
crucial in order to obtain the desired estimate.  In particular, the
term $\det {\mathcal A}_1$ is linear in the highest-order derivative
$\partial \partial_t^2 v$ rather than quadratic (as it a priori appears).
\end{remark}

There are three remaining boundary integral terms 
appearing on the right 
hand side of (\ref{eulerttt}) arising from (\ref{kbca}); 
the terms involving $\kappa$ are
\begin{align}
&-\kappa \int_0^T
\left( 
\left[\partial_t^3 (\tilde v \cdot \nk), \partial_t^3\tilde v \cdot \nk\right]_1
+3\left[\partial_t^2 (\tilde v \cdot \nk), \partial_t^3\tilde v \cdot \partial_t \nk\right]_1
+3\left[\partial_t (\tilde v \cdot \nk), \partial_t^3\tilde v \cdot \partial_t^2 \nk\right]_1 
\right.
\nonumber
\\
& 
\qquad\qquad\qquad \qquad\qquad\qquad\qquad \qquad
\left.
+\left[\tilde v \cdot \nk, \partial_t^3\tilde v \cdot \partial_t^3 \nk\right]_1
\right)\,. \label{kextra}
\end{align}
The first term in (\ref{kextra}) provides both the energy
contribution $\int_0^T |\sk\partial_t^3 \tilde v \cdot \nk|_1^2$ as
well as error terms.  We start the analysis with the most difficult error term,
\begin{equation}\label{difficult}
\kappa  \int_0^T\{[ \tilde v \cdot \partial \partial_t^3 \nk, 
\partial(  \nk \cdot \partial_t^3\tilde v)]_0,
\end{equation}
whose highest-order contribution has
an integrand (modulo $L^\infty$ terms) of the form  $\partial^2 
\sk\tilde {v_\kappa}_{tt}$ $\sk\partial \tilde v_{ttt}$.


With $\nk= (\partial_1\ek \times \partial_2\ek)/\sqrt{g}=: {Q}(\partial \tilde \eta_\kappa)$, ${Q}$ given
by (\ref{Q}), the highest-order term in
$\partial\partial_t^3 \nk$ is 
${Q}(\partial\tilde\eta_\kappa) \partial^2 \partial_t^2 
\tilde v_\kappa$, so that with ${\mathcal R}_1$ denoting a lower-order 
remainder term, and using (\ref{deteta.b}), we have that

\begin{align*}
&-\kappa \int_0^T[ \tilde v \cdot \partial \partial_t^3 \nk, \partial( \nk \cdot
\partial_t^3 \tilde v)]_0
 \le C
\sup_{t\in[0,T]}|P(\tilde v,\partial\tilde \eta_\kappa)|_{L^\infty} 
\int_0^T|\sk\partial_t^3\tilde v \cdot \nk|_1 \, 
| \sk\partial_t^2 \tilde v_\kappa|_{2} + {\mathcal R}_1\\
&\qquad
\le C \sup_{t\in[0,T]}|P(\tilde v,\partial\tilde \eta_\kappa)|_{L^\infty} 
|\sk\partial_t^3\tilde v \cdot \nk|_{L^2(0,T; H^1(\Gamma))}
|\sk\partial_t^2 \tilde v_\kappa|_{L^2(0,T; H^2(\Gamma))}
+ {\mathcal R}_1\\
&\qquad
\le C_\delta \left[\sup_{t\in[0,T]}|P(\tilde v,\partial\tilde \eta_\kappa)|_{L^\infty} 
\|\sk\tilde {v_\kappa}_{tt}\|_{L^2(0,T; H^{2.5}(\Omega))}\right]^2
+ \delta |\sk \tilde v_{ttt} \cdot \nk|_{L^2(0,T; H^1(\Gamma))}^2
+ {\mathcal R}_1\\
&\qquad
 \le M_0+ C\,T\,\pskE + \| \sk \tilde {v_\kappa}_{tt} \|^4_{L^2(0,T; H^{2.5}(\Omega))}
  + \delta \skE \,,
\end{align*}
where ${\mathcal R}_1$ also satisfies
${\mathcal R}_1 \le C\,T\, \pskE + \delta \skE$.
The second term in (\ref{kextra}) has a highest-order contribution with the
same type of integrand, and its analysis (and bound) is identical.  The third
and fourth terms in (\ref{kextra}) are effectively lower-order by one derivative
with respect to the worst case analyzed above.

Next, we estimate 
$\int_0^T\int_{\Gamma} \partial_t^3\{ \tilde H\tilde n\cdot \nk\, 
(\nk -\tilde n)\}\, \partial_t^3 \tilde v$. 
Since $\nk={Q}(\tilde\eta_\kappa)$ and since
$|\nk -\tilde n| \le \sup_\kappa |\partial Q(\partial\eta_\kappa)| \cdot
|\partial \eta_\kappa - \partial \eta|$, then  our assumed bounds
(\ref{deteta}) together with (\ref{Linf_est}) imply that
\begin{align}
|\nk - \tilde n|_{L^\infty} 
\le C\,\sk \, |P(\partial\tilde \eta, \partial^2 \tilde \eta )|_{L^\infty} |
\tilde\eta |_{2.5} 
\le C\, \kappa\, P(E_\kappa(t))\,.
\label{n0}
\end{align}
Similarly,
\begin{align}
|\partial\nk - \partial\tilde n|_{L^\infty} 
\le C\,\sk \, |P(\partial\tilde \eta, \partial^2 \tilde \eta )|_{L^\infty} 
|\tilde\eta |_{3.5} \,.
\label{n}
\end{align}
Also by (\ref{Linf_est}), for $k=1,2,3$,
\begin{align}
|\partial_t^k\nk -\partial_t^k \tilde n|_{L^\infty} 
\le C\sk \,  | P(\partial\tilde \eta, \partial^2\tilde\eta)|_{L^\infty}\,
| \partial_t^{k-1}\tilde v|_{2.5} 
\,,
\label{ntt}
\end{align}
and
\begin{align}
|\partial_t^2\vk -\partial_t^2 \tilde v|_0
&\le C\, \sk\, |\tilde v_{tt}|_{1.5}\,.
\label{vtt}
\end{align}
Taking three time-derivatives of formula (\ref{importantF}), we see that
the highest order term in 
$\partial_t^3(\sqrt{\tilde g}\tilde H\tilde n)$ is $Q(\partial \tilde \eta) \partial^2
\tilde v_{tt}$.  Thus, the highest-order term in the integral
$$
\int_0^T\int_{\Gamma} \partial_t^3(\sqrt{\tilde g}\tilde H\tilde n)^s
\nk^s \, \, (\nk^r -\tilde n^r) \partial_t^3\tilde v^r
$$
is estimated using an integration by parts in space.  The highest derivative
count
occurs when the tangential derivative is moved onto the $\tilde v_{ttt}$ term
giving us
\begin{align*}
\int_0^T\int_{\Gamma} Q^s_i(\partial \tilde \eta) \partial \tilde v^i_{tt}
\nk^s \, \, (\nk^r -\tilde n^r) \partial\tilde v^r_{ttt}
&\le
C\int_0^T 
|P(\partial \tilde \eta)|_{L^\infty}\, |\tilde v_{tt}|_1 \,
|\nk -\tilde n |_{L^\infty}\, |\partial_t^3 \tilde v|_1 \\
& 
\le
C\int_0^T 
|P(\partial \tilde \eta, \partial\tilde \eta)|_{L^\infty}\, 
|\tilde \eta|_{2.5}\, |\tilde v_{tt}|_1 \,
|\sk \tilde v_{ttt}|_1 \\
& 
\le  C_\delta\, T\, \pskE + \delta \int_0^T |\sk \tilde v_{ttt}|_1^2\,,
\end{align*}
where (\ref{n}) is used for the second inequality.  If instead, integration
by parts places the tangential derivative on $\nk -\tilde n$, then (\ref{n})
provides the same estimate for this term.  The other terms are clearly
lower-order.

Thanks to (\ref{ntt}),
\begin{align*}
\int_0^T\int_{\Gamma} \partial_t(\sqrt{\tilde g}\tilde H\tilde n^s )
\partial_t^2\{\nk^s\, (\nk^r -\tilde n^r)\} \, \partial_t^3\tilde v^r
&\le
C\int_0^T |P(\partial \tilde \eta, \partial^2 \tilde \eta)|_{L^\infty}\, 
|\partial^2 \tilde v|_0\,
|\tilde v_t|_{2.5} \, |\sk \tilde v_{ttt}|_0 \\
& 
\le  C_\delta\, T\, \pskE + \delta \sE\,.
\end{align*}


We next consider the integral 
\begin{align}
&\int_0^T\int_{\Gamma} \sqrt{\tilde g}\tilde H\tilde n^i  
\partial_t^3\{ \nk^i\, (\nk^r -\tilde n^r)\}\,
\partial_t^3 \tilde v^r 
\nonumber
\\
&\qquad
=
\int_0^T\int_{\Gamma} \sqrt{\tilde g}\tilde H\tilde n \cdot 
\partial_t^3\nk\, \, (\nk -\tilde n)\cdot  \tilde v_{ttt}
+
\int_0^T\int_{\Gamma} \sqrt{\tilde g}\tilde H\tilde n \cdot 
\nk\, \, \partial_t^3(\nk -\tilde n)\cdot  \tilde v_{ttt} + {\mathcal R}_2
\nonumber
\\
&\qquad
=: I + II + {\mathcal R}_2 \,, \label{improve}
\end{align}
where ${\mathcal R}_2$ is a lower-order term.
For term $I$, we use the estimate $| \nk -\tilde n|_{L^\infty} \le C\, \kappa
|\tilde \eta|_{3.5}$ One $\sk$ goes with $\partial_t^3\nk$ and the other
$\sk$ goes with $\tilde v_{ttt}$.  Thus,
$|I| \le  C\, T\, \pskE + \delta \skE$.  

To study $II$, 
we set $f=\sqrt{\tilde g} \tilde H \tilde n \cdot \nk$, and consider the
term $\int_0^T\int_\Gamma f\, \tilde n_{ttt}\cdot \tilde v_{ttt}$.
We expand $v$ into its normal and tangential components:  set
$\tau_\alpha = \tilde \eta,_\alpha$, so that
$$
\tilde v = v^\tau \,\tau + v^n\, \tilde n, \text{ where } v^\tau \, \tau = (\tilde v\cdot \tau_\alpha) \, \tau_\alpha\, \text{ and } v^n =\tilde v \cdot \tilde n \,.
$$
Then
\begin{align*}
\tilde v_{ttt} &= v^\tau_{ttt} \tau + 3 v^\tau_{tt} \tau_{t}
+ 3v^\tau_{t} \tau_{tt}+ v^\tau \tau_{ttt} 
+ v^n_{ttt} \tilde n + 3 v^n_{tt} \tilde n_{t}
+ 3v^n_{t} \tilde n_{tt}+ v^n \tilde n_{ttt}\,.
\end{align*}
The most difficult term to estimate comes from the term $v^\tau_{ttt} \tau$,
which gives the integral 
$\int_0^T\int_\Gamma f\, \tilde n_{ttt}\cdot \tau \, v^\tau_{ttt}$.

First,
notice that $\tilde n_{ttt} \cdot \tau$ is equal to 
$-\tilde n \cdot \tau_{ttt}$, plus lower order terms that have at most two
time derivatives on either $\tilde n$ or $\tau$, and 
$\tilde n \cdot \tau_{ttt} = \tilde n \cdot \partial_\beta \tilde v_{tt}$ for 
$\beta=1$ or $2$.
Next, the $\kappa$ problem states that 
$\tilde v^i_{ttt} = (A^k_i \tilde q,_k)_{tt}$, where recall that
$$
A = \tilde a_\kappa\,.
$$
We have the formula
$$
\tilde \eta^i,_\beta A^k_i \partial_k \tilde q_{tt} = \mathcal{J}^\beta
\partial_\beta \tilde q_{tt} \ (\text{no sum on } \beta) \ \ \ \text{ for } \ \beta=1,2\,,
$$
where
$$
\mathcal{J}^1= \tilde \eta^i,_1 \, A^1_i\,, \ \ \ \ 
\mathcal{J}^2= \tilde \eta^i,_2 \, A^2_i\,.
$$
(In the case that $\kappa=0$, $\mathcal {J}^\beta = J =1$.)
Using this,
we see that the highest-order term in our integral is given by
\begin{equation}\label{S0}
\int_0^T\int_\Gamma \mathcal{J} f\, 
(\tilde n \cdot \partial_\beta \tilde v_{tt}) \partial_\beta \tilde q_{tt} \,.
\end{equation}

Second, write $\tilde q_{tt}$ as
\begin{equation}\label{S1}
\tilde q_{tt} = -\left[ \frac{\sqrt{\tilde g}}{\sqrt{\gk}}
\left[
\Delta_{\tilde g}(\tilde \eta) \cdot \tilde n
+ \Delta_{\tilde g}(\tilde \eta) \cdot (\nk -\tilde n)
\right]
+\kappa \Delta_0 (\tilde v \cdot \nk)
\right]_{tt} \,.
\end{equation}
We begin by substituting the first term on the right-hand side of (\ref{S1})
into (\ref{S0}); the highest-order contribution comes from 
$\partial_\beta\partial_t^2 
\Delta_{\tilde g}(\tilde \eta) = Q(\partial \tilde \eta, \partial \tilde
\eta_\kappa) \tilde g^{\mu \nu} \tilde n \cdot \partial_t \tilde v,_
{\mu\nu \beta}$.
Integrating by parts with respect to $\partial_\nu$, the highest-order term
in our integral is given by
$$
\int_0^T \int_\Gamma 
 Q(\partial \tilde \eta, \partial \tilde \eta_\kappa)\, f\, 
(\tilde n \cdot \tilde v_{t},_ {\mu \beta}) \, 
\tilde g^{\mu \nu} \,
(\tilde n \cdot \tilde v_{tt},_ {\nu \beta}) \,.
$$
Letting $G_{ij}^{\mu\nu}:=
Q(\partial \tilde \eta, \partial \tilde \eta_\kappa)\, f\, 
\tilde n_i \, \tilde n_j
=
Q(\partial \tilde \eta, \partial \tilde \eta_\kappa) \partial^2\eta$,
integration by parts in time yields 
\begin{align} 
&-\int_0^T \int_\Gamma 
\partial_t G_{ij}^{\mu\nu}\, \,\,
\tilde v^j_{t},_ {\nu \beta}\, 
\tilde v^i_{t},_ {\mu \beta}
+\left. \int_\Gamma 
\partial_t G_{ij}^{\mu\nu}\, \,\,
\tilde v^j_{t},_ {\nu \beta}\, 
\tilde v^i_{t},_ {\mu \beta}\right]_0^T \nonumber \\
& \qquad\qquad
\le C\, T\pskE + M_0 + C \sup_{t\in [0,T]} |G_t|_{L^\infty} \,\| \tilde v_t\|_{2.5}^2
\nonumber \\
& \qquad\qquad
\le M_0 + C\, T\pskE + C \sup_{t\in [0,T]} [
P(\|\tilde v_t\|_{2.5}^2) + P(\|\tilde v\|_{3.5}^2)+ P(\|\tilde\eta\|_{4.5}^2)]
\,.\label{G2D}
\end{align}
For the second term on the RHS of (\ref{S1}), the highest-order term gives
the integral
\begin{align*}
&\int_0^T \int_\Gamma
f\,Q(\partial \tilde \eta, \partial \tilde \eta_\kappa) \, 
(\tilde n \cdot \tilde v_{tt},_\beta) \,
\tilde g^{\mu \nu}  \Pi \tilde v_t,_ {\mu\nu \beta}
\cdot (\nk -\tilde n) 
\le C\, T\, \pskE + \delta \int_0^T \|\sk \tilde v_t\|_{3.5}^2,
\end{align*}
where we used $|\nk -\tilde n|_{L^\infty} < C\, \kappa |\tilde \eta|_{3.5}$
again.

For the third term on the RHS of (\ref{S1}), the highest-order term gives
the integral
\begin{align*}
&\kappa\int_0^T \int_\Gamma
f\,Q(\partial \tilde \eta, \partial \tilde \eta_\kappa) \, 
(\tilde n \cdot \tilde \partial^2 v_{tt}) \,
(\nk \cdot \partial^2 \tilde v_{tt}) \\
&\qquad \qquad
\le M_0 + C\, T\, \pskE +  \|\sk \tilde v_{tt}\|_{L^2(0,T;H^{2.5}(\Omega)}^4 \,.
\end{align*}
We have thus estimated the integral
$\int_0^T\int_{\Gamma} \partial_t^3\{ \tilde H\tilde n\cdot \nk \, 
(\nk- \tilde n) \}\, \partial_t^3 \tilde v$.
The remaining integral
$\int_0^T\int_{\Gamma} \partial_t^3\{ \tilde H\tilde n\cdot (\nk -\tilde n)\, 
\nk \}\, \partial_t^3 \tilde v$ has the same bound.

\noindent
{\bf Step 2. The pressure term.}
We next consider the pressure term  in
(\ref{eulerttt}):
\begin{align}
-\int_0^T\int_{\Omega} \partial_t^3(A^k_i \tilde q) \partial_t^3  \tilde v^i,_k   
& = -\int_0^T\int_{\Omega} \partial_t^3  \tilde v^i,_k
\left[ \partial_t^3 A^k_i  \tilde q + 3 \partial_t^2{A^k_i}\,  \tilde q_t 
+ 3 \partial_t{A^k_i} \, \tilde q_{tt} + A^k_i \partial_t^3 \tilde q \right] \nonumber \\
&=: I + II + III + IV \,. \label{qttt_term}
\end{align}
We record the following identities:
\begin{subequations}
\label{a123t}
\begin{align}
\partial_t A^k_i  & = \tilde J_\kappa^{-1}(A^s_r A^k_i - A^k_r A^s_i) {\tilde v_\kappa}^r,_s
\,, \label{a123t.a} \\
\partial^2_t A^k_i &   = \tilde J_\kappa^{-1}(A^s_r A^k_i - A^k_r A^s_i) \partial_t{\tilde v_\kappa}^r,_s
+P^k_i(\tilde J_\kappa^{-1},A,\nabla \tilde v_\kappa)
\,, \label{a123t.b} \\
\partial^3_t A^k_i &   = \tilde J_\kappa^{-1}(A^s_r A^k_i - A^k_r A^s_i) \partial^2_t{\tilde v_\kappa}^r,_s
+ P^k_j (\tilde J_\kappa^{-1},A,\nabla {\tilde v_\kappa}) \partial_t{\tilde v_\kappa}^j,_i
\,. \label{a123t.c}
\end{align}
\end{subequations}
With (\ref{a123t.c}), and $f^{sk}_{ri}:=\tilde J_\kappa^{-1}(A^s_r A^k_i - A^k_r A^s_i)\tilde q$
term $I$ is written as
\begin{align*}
I  &= \int_0^T \int_\Omega [ f^{sk}_{ri} \partial_t^3 \tilde v^i,_k \partial_t^2
\vk^r,_s + \partial_t^3 \tilde v^i,_k \partial_t \vk^j,_i
P^k_j(\tilde J_\kappa^{-1},A,\nabla \vk)] \\
&=: I_a + I_b \,.
\end{align*}
We fix a chart $\theta_l$ in a neighborhood of the boundary $\Gamma$ and
let $\xi= \sqrt{\alpha_l}$, where once again, we remind the reader that 
$\{\alpha_i\}_{i=1}^L$ denotes the partition
of unity associated to the charts $\{\theta\}_{l=1}^L$. With ${\mathcal I}_a$
denoting the restriction $I_a|_{U_l}$, where $U_l\cap \Omega = \theta_l
((0,1)^3)$, and letting $\rho:= \rho_{\frac{1}{\kappa}}$ and
$\theta:=\theta_l$, we have that
\begin{align*}
{\mathcal I}_a & = \int_0^T \int_{(0,1)^3} 
f^{sk}_{ri}(\theta) \partial_t^3 \tilde v^i,_k (\theta)\, \xi\, \rho\star_h \rho\star_h
\xi \partial_t^2 \tilde v^r,_s(\theta) + f(\theta) \nabla\tilde v_{ttt}
G(\xi, \nabla \xi) \tilde v_{tt} \\
&=: {{\mathcal I}_a}_1+ {{\mathcal I}_a}_2 \,,
\end{align*}
where $G(\xi,\nabla \xi)$ is a bilinear function which arises
when the gradient acts on $\xi$ rather than $\tilde v_{tt}$.
The term ${{\mathcal I}_a}_1$ is the difficult term which requires
forming an exact derivative, and this, in turn, requires commuting the
convolution with $f^{sk}_{ri}$.   We let
$$
V(\theta) = \rho \star_h \xi \tilde v(\theta)
$$
so that using the symmetry property (\ref{selfadjoint}), we see that
\begin{align*}
{{\mathcal I}_a}_1 &= 
\int_0^T\int_{(0,1)^3} \rho\star_h [f^{sk}_{ri} \xi \partial_t^3
\tilde v^i,_k(\theta)] \partial_t^2 V^r,_s(\theta)  \\
&
 = 
\int_0^T\int_{(0,1)^3}\left[ f^{sk}_{ri} \partial_t^3
V^i,_k(\theta)] \partial_t^2 V^r,_s(\theta) 
+ (\rho \star_h [f^{sk}_{ri} \xi\partial_t^3\tilde v^i,_k]
- f^{sk}_{ri} \rho\star_h [\xi\partial_t^3\tilde v^i,_k])\,
V_{tt}^r,_s \right]\\
&=:
{{{\mathcal I}_a}_1}_i+{{{\mathcal I}_a}_1}_{ii}\,.
\end{align*}
Since $f^{sk}_{ri}$ is symmetric with respect to $i$ and $r$ and $k$ and $s$,
we see that
\begin{align}
{{{\mathcal I}_a}_1}_i 
 &= \frac{1}{2}\int_0^T\int_{(0,1)^3} f^{sk}_{ri}(\theta) 
\partial_t\left[ 
\partial_t^2 V^i,_k(\theta)\partial_t^2 V^r,_s(\theta) \right] 
\nonumber
\\
&
=
 -\frac{1}{2}\int_0^T\int_{(0,1)^3} \partial_t f^{sk}_{ri}(\theta) 
\partial_t^2 V^i,_k(\theta)\partial_t^2 V^r,_s(\theta) 
+
\frac{1}{2}\left.\int_{(0,1)^3} f^{sk}_{ri}(\theta) 
\partial_t^2 V^i,_k(\theta)\partial_t^2 V^r,_s(\theta)\right]^T_0 \,.
\nonumber
\end{align}
We sum over $l=1,...,L$.
The spacetime integral is bounded by $C T \psE$. For the
the space integral at time $t=T$, we employ the fundamental theorem of calculus:
\begin{align*}
\int_\Omega  [f^{sk}_{ri}V_{tt}^i,_k V_{tt}^r,_s ] (T) & =
\int_\Omega  V_{tt}^i,_k(T) V_{tt}^r,_s(T) f^{sk}_{ri}(0)
+ \int_\Omega  V_{tt}^i,_k(T) V_{tt}^r,_s(T) \int_0^T \partial_tf^{sk}_{ri} \\
& \le \|V_{tt}(T)\|_1^2 \, \|\tilde q_0\|_2 + \|V_{tt}(T)\|_1^2 \,\, \|\int_0^T |f_t\|
_{L^\infty} \\
& \le \|v_{tt}(T)\|_{1.5}^{2/3} \, \|v_{tt}(T)\|_{0}^{1/3} 
\, \|\tilde q_0\|_2 +  C T \psE \\
& \le {\frac{\delta}{2}} 
\|v_{tt}(T)\|_{1.5}^{2} + C(\delta) \|v_{tt}(T)\|_{0}^{1/2} 
\, \|\tilde q_0\|_2^{3/2} +  C T \psE \\
& \le \delta \|v_{tt}(T)\|_{1.5}^{2} + C(\delta) 
\, \|\tilde q_0\|_2^{2} +  C T \psE \\
& \le \delta \skE + M_0(\delta) + C T\psE \,,
\end{align*}
where Young's inequality has been used.   
For
${{{\mathcal I}_a}_1}_{ii}$, the commutation lemma \ref{commutator} shows
that
\begin{align*}
{{{\mathcal I}_a}_1}_{ii}
&\le C \kappa^{\frac{3}{2}}\int_0^T \|f\|_3\, \|\sqrt{\kappa} \tilde v_{ttt}\|_1 \,
\| \tilde v_{tt}\|_1 \\
&\le \delta \int_0^T \|\sqrt{\kappa} \tilde v_{ttt}\|_1^2 + C_\delta \int_0^T \|f\|_3^2\,
\|\tilde  v_{tt}\|_1^2 \,.
\end{align*}
Summing over $l=1,...,L$, we integrate-by-parts in time and write the term
${{{\mathcal I}_a}_2}$ as
\begin{align*}
{{{\mathcal I}_a}_2} & = -\int_0^T\int_\Omega f\nabla \tilde v_{tt} G(\xi,
\nabla \xi) (f \tilde v_{tt})_t + \left.\int_\Omega f\nabla \tilde v_{tt} G(\xi,
\nabla \xi) f \tilde v_{tt}\right]^T_0\,.
\end{align*}
This is estimated in the same way as term
${{{\mathcal I}_a}_1}_{i}$.  The term $I_b$ is handled in the identical
fashion with the same bound.  Thus, we have shown that
$$
I \le \delta \skE + M_0(\delta) + C T\psE \,.
$$

Using (\ref{a123t.b}) for term $II$, integration by parts in time gives
the identical bound as for term $I$.
For term $III$, a different approach is employed;
we use (\ref{a123t.a}) and integration by parts in space rather than time, 
and let 
$F^{sk}_{ri}:=3\tilde J_\kappa^{-1}(A^s_r A^k_i - A^k_r A^s_i)$ to find that
\begin{align*}
III  = - \int_0^T \int_\Omega \tilde v_{ttt}^i[ \vk^r,_s F^{sk}_{ri} q_{tt}],_k
+ \int_0^T \int_\Gamma \tilde v_{ttt}^i   F^{sk}_{ri} N_k q_{tt} \vk^r,_s 
\end{align*}
The Cauchy-Schwarz inequality together with the pressure estimate 
(\ref{qestimate}) give the bound $CT \psE$ for the first term on the right-hand
side.  The boundary integral term requires integration by parts in time:
\begin{align*}
\int_0^T \int_\Gamma \tilde v_{ttt}^i   F^{sk}_{ri} N_k \tilde q_{tt} \vk^r,_s 
&= 
\left. \int_\Gamma \tilde v_{tt}^i   F^{sk}_{ri} N_k \tilde q_{tt} \vk^r,_s 
 \right]^T_0
+\int_0^T \int_\Gamma v_{tt}^i [  F^{sk}_{ri} N_k \tilde q_{tt} \vk^r,_s ]_t\\
&=:  III_a + III_b\,.
\end{align*}
First, note that
\begin{align}
F^{sk}_{ri} N_k &= 3\tilde J_\kappa^{-1} A^s_r
\left[ (\tilde q A^k_i N_k)_{tt} - 2 \tilde q_t \partial_t A^k_i N_k -\tilde q
\partial_t^2 A^k_i N_k \right] \nonumber\\
& \qquad - 3\tilde J_\kappa^{-1} A^s_i
\left[ (\tilde q A^k_r N_k)_{tt} - 2 \tilde q_t \partial_t A^k_r N_k -\tilde q
\partial_t^2 A^k_r N_k \right]  \label{F} \,.
\end{align}
Next, substitute the boundary condition (\ref{kbca}), written as
\begin{equation}
\label{kbca2}
-\tilde q A^k_i N_k = \sqrt{\tilde g}\Delta_{\tilde g}(\tilde \eta^j)\left[
\delta_{ji} + ((\nk)_j - \tilde n_j) \tilde n_i + \tilde n_j((\nk)_i -
\tilde n_i)\right] + \kappa (\sqrt{g_0} g_0^{\alpha\beta} 
[\tilde v\cdot \tilde n],_\beta),_\alpha \tilde n_i \,,
\end{equation}
into (\ref{F}).
The two bracketed terms in (\ref{F}) are essentially the same, so it suffices to analyze
just the first term.  We begin by considering the term
$\sqrt{\tilde g}\Delta_{\tilde g}(\tilde \eta^i)$ in (\ref{kbca2}).

Then $III_a$ can be written as
\begin{align}
III_a &=
\left.
3 \int_\Gamma \tilde J_\kappa^{-1}\left\{ 
\partial_t^2(\sqrt{\tilde g} \tilde g^{\alpha\beta} \tilde \eta^i,_\beta)
,_\alpha \vk^r,_s A^s_r
-2 q_t \partial_t A^k_i N_i \vk^r,_s - q \partial_t^2 A^k_i N_i \vk^r,_s 
\right\}
\tilde v_{tt}^i 
\right]^T_0 \label{IIIa}
\\
&=
\left.
3 \int_\Gamma 
(\sqrt{\tilde g} \tilde g^{\alpha\beta} \tilde \eta^i,_\beta)
_{tt} ( \tilde J_\kappa^{-1}\vk^r,_s A^s_r \tilde v^i_{tt}),_\alpha
-2 \tilde J_\kappa^{-1}q_t \partial_t A^k_i N_i \vk^r,_s \tilde v^i_{tt}
- \tilde J_\kappa^{-1}q \partial_t^2 A^k_i N_i \vk^r,_s \tilde v_{tt}^i 
\right]^T_0 \nonumber \\
&\le \delta \sE + M_0(\delta) + CT \psE \,, \nonumber
\end{align}
the last inequality following from  the fundamental theorem of calculus
and the same argument we have  used above.

In order
to estimate $III_b$, because we do not have a trace estimate for
$\partial_t^3 A$, we let $Q_r(\partial \ek):= \sqrt{\gk} (\nk)_r$ and
compute
\begin{align}
\partial_t Q_r &= Q^\alpha_{ri} \vk^i,_\alpha \,,  \ \ \
\partial_t^2 Q_r = Q^{\alpha\beta}_{rij} \vk^i,_\alpha \vk^j,_\beta 
+ Q^\alpha_{ri} \partial_t\vk^i,_\alpha \,, \nonumber\\
\partial_t^3 Q_r &= Q^{\alpha\beta\gamma}_{rijk} \vk^i,_\alpha \vk^j,_\beta 
\vk^k,_j + 3 Q^{\alpha\beta}_{rij} \vk^i,_\alpha \partial_t \vk^j,_\beta
+ Q^\alpha_{ri} \partial_t^2 \vk^i,_\alpha \,. \label{Qttt}
\end{align}
Since $\sqrt{\gk} (\nk)_r \tilde q_{tt} = (\sqrt{\gk} (\nk)_r \tilde q)_{tt} -
2 \partial_t Q_r\, \tilde q_t - \partial_t^2 Q_r \, \tilde q$, it follows that
\begin{align}
III_b &= -3 \int_0^T \int_\Gamma \tilde J_\kappa^{-1}\left[ 
\tilde v_{tt}^r (\sqrt{\gk} (\nk)_r \tilde q)_{ttt} \vk^i,_s A^s_i
+ \tilde v_{tt}^r (Q_r \tilde q)_{tt} (\vk^i,_s A^s_i)_t \right. \nonumber \\
&\qquad\left.
- \tilde v_{tt}^r (2 Q^\alpha_{rj}\vk^j,_\alpha \vk^i,_s A^s_i \tilde q_t +
Q^{\alpha \beta}_{rlj} \vk^l,_\alpha \vk^j,_\beta \vk^i,_s A^s_i \tilde q +
Q^\alpha_{rl}\partial_t\vk^l,_\alpha \vk^i,_s A^s_i \tilde q)_t
\right] \,. \label{IIIb}
\end{align}
Using the pressure estimates and
by definition of our energy function, for $t \in (0,T)$,
\begin{gather}
\tilde q_{tt}(t) \in H^{0.5}(\Gamma), \ \ \  Q_r(t), \partial_t Q_r(t) \in 
L^\infty(\Gamma), \ \ \ \partial_t^3 Q_r(t) \in L^2(\Gamma), \nonumber \\
\partial \tilde v_t(t) \in H^{1.5}(\Gamma), \ \ \ \partial \tilde v_{tt}(t) \in L^2(\Gamma) \,;
\nonumber
\end{gather}
thus, all of the terms, except the first, on right-hand side of (\ref{IIIb}) 
can be easily bounded by $CT \psE$.  Integrating by parts in space, the
first term in (\ref{IIIb}) has the following estimate:
\begin{align}
&
3 \int_0^T \int_\Gamma \sqrt{\gk} \gk^{\alpha\beta} \Pi^r_k
\tilde v_{tt}^k,_\beta (\tilde v_{tt}^r \vk^i,_s J_\kappa^{-1} A^s_i),_\alpha 
+
\sqrt{\gk} (\gk^{\mu\nu} g^{\alpha\beta} - g^{\alpha\nu} g^{\mu\beta})
\tilde\eta^j,_\nu \tilde v_{tt}^j,_\mu \tilde\eta^r,_\beta
(\tilde v_{tt}^r \vk^i,_s \tilde J_\kappa^{-1} A^s_i),_\alpha  \nonumber\\
& \qquad\qquad 
+ [P^{\alpha\beta}_{ij}(\partial\tilde \eta,\partial \tilde v)
 \tilde v_{tt}^j,_\beta + P^\alpha_i(\partial \tilde \eta, \partial \tilde v)]
(\tilde v_{tt}^r \vk^i,_s \tilde J_\kappa^{-1} A^s_i),_\alpha 
 \le C T \psE \,. \label{IIIb_estimate}
\end{align}

The remaining three terms in the boundary condition (\ref{kbca2}) are now
considered.
The additional integrals which arise in the (\ref{IIIb}) are given by
\begin{align*}
&-3\int_0^T\int_{\Gamma} \partial_t^3 \{
\sqrt{\tilde g}\tilde H\tilde n\cdot(\nk -\tilde n)\, \tilde n_r
+ \sqrt{\tilde g} \tilde H\tilde n\cdot \nk\, ((\nk)_r -\tilde n_r) 
+ \kappa\sqrt{\tilde g_\kappa} \Delta_{\bar 0}(\tilde v\cdot \nk)\nk
\}\, 
\tilde v^r,_s \tilde a^s_i \tilde v^i_{tt} \\
&\qquad =: J_1+J_2+J_3 \,.
\end{align*}
Term $J_3$ with the artificial viscosity provides the integral
$\kappa\int_0^T [\partial_t^3(\tilde v\cdot \nk), 
P(\nabla \tilde \eta,\nabla \tilde v) \cdot \partial_t^2\tilde v]_1$. The 
highest-order terms in this integral are estimated as
\begin{align*}
&\kappa\int_0^T \left\{[ P^{\alpha\beta}_{ij}(\partial \nk, \nabla \tilde \eta,
\nabla \tilde v) \,
\partial_t^3 \tilde v^i_{,\beta}\,, \partial_t^2\tilde v^j_{,\alpha}]_0 + 
[ P^{\alpha\beta\gamma}_{ij}(\partial \nk, \nabla \tilde \eta,\tilde v,\nabla \tilde v) \,
\partial_t^2 \tilde v_\kappa^i,_{\beta\gamma}\,, \partial_t^2\tilde v^j_{,\alpha}]_0 
\right\} \\
& \qquad \le C\, \sk  \int_0^T \|P(\nabla\tilde \eta, \nabla \tilde v)\|_{L^\infty}\,
| \partial_t^2 \tilde v|_1  \,
\{ |\sk\partial_t^3 \tilde v|_1 + |\sk\partial_t^2 \tilde v_\kappa|_2 \}\, 
\\
& \qquad \le C(\delta)\, 
\int_0^T \|P(\nabla \tilde \eta,\nabla \tilde v) \|_{L^\infty}\, | \partial_t^2 \tilde v|^2_1 
+ \kappa\delta\int_0^T \left(|\sk\partial_t^3 \tilde v|^2_1+|\sk\partial_t^2 \tilde v_\kappa|^2_2\right)\\
&\qquad\le  C(\delta)\, T\, \pskE +   \delta \skE \,.
\end{align*}
The lower-order terms in $J_3$ also have the same bound.  As to
terms $J_1$ and $J_2$, the estimates for the
terms with $(\sqrt{\tilde g}\tilde H\tilde n)_{ttt}$ are obtained exactly as
in (\ref{IIIb_estimate}).  For the terms that contain $\tilde n_{ttt}$, we
use the formula (\ref{Qttt}) for the third time derivative of the unit normal;
it immediately follows that terms $J_1$ and $J_2$ are also bounded by
$\delta \skE + CT \pskE$.  


We need only consider now  the additional terms in (\ref{IIIa}) from
the remaining three terms in the boundary condition (\ref{kbca2}).  The only
novelty is in the highest-order 
integral coming from integration by parts in space in the $\kappa$ term:
$$\kappa\int_\gamma P^{\alpha\beta}_{ij}(T) v^i_{tt},_\alpha(T)
v^j_{tt},_\beta(T)\,$$ 
where $P^{\alpha\beta}_{ij}(t)$ and 
$\partial_t P^{\alpha\beta}_{ij}(t)$ are both in $L^\infty(\Gamma)$ for
each $t\in[0,T]$.  Using the fundamental theorem of calculus and
the fact that $\sqrt\kappa v_{ttt} \in L^2(0,T; H^{1.5}(\Omega))$ together
with Jensen's inequality shows that this term is bounded by $M_0(\delta) +
\delta \skE + CT\pskE$.

To study term $IV$, we use (\ref{a123t}) together with the
incompressibility condition $(\tilde v^i,_k A^k_i)_{ttt}=0$
to find that
\begin{align*}
IV= -\int_0^T \int_\Omega 
[(3 \tilde v_{tt}^i,_k \partial_t A^k_i + \tilde v^i,_k \partial^3_t A^k_i)
\tilde q_{ttt} +
 3 \tilde v_t^i,_k \partial_t^2 A^k_i \tilde q_{ttt}]
=:  IV_a + IV_b \,.
\end{align*}
For $IV_b$, we integrate by parts in time:
\begin{align*}
IV_b= \int_0^T \int_\Omega 
 3 (\tilde v_t^i,_k \partial_t^2 A^k_i)_t \tilde q_{tt}
- \left.
\int_\Omega
 3 \tilde v_t^i,_k \partial_t^2 A^k_i \tilde q_{tt}
\right]^T_0 \,.
\end{align*}
Since $\partial_t^3 A$ is bounded in $H^{0.5}(\Omega)$,
the spacetime integral is easily bounded by $CT\psE$; 
meanwhile, the remaining space integral  satisfies
\begin{align*}
{IV_b}_2 & \le 3\int_\Omega [\tilde v_t^i,_k(0)\partial_t^2A^k_i(0)] q_{tt}(T)
+\int_\Omega \int_0^T [ \tilde v_t^i,_k\partial_t^2A^k_i]_t \, dt \
q_{tt}(T) + M_0 \\
& \le  \delta \|q_{tt}(T)\|_0^2 + M_0(\delta)  +
T\sup_{t \in [0,T]} \| (\tilde v_t^i,_k \partial_t^2A^k_i)_t\|_0 \,
\|q_{tt}(T)\|_0 \\
&\le \delta \sE + M_0(\delta) + CT \psE \,.
\end{align*}
With $F^{sk}_{ri}:= \tilde J_\kappa^{-1} (A^s_rA^k_i - A^k_rA^s_i)$, $IV_a$
is written as
\begin{align}
IV_a &= -\int_0^T \int_\Omega  
(3 \partial_t^2 \tilde v^r,_s F^{sk}_{ri} \vk ^i,_k
+\partial_t^2 \vk^r,_s F^{sk}_{ri} \tilde v ^i,_k) \tilde q_{ttt} 
+ \tilde v^i,_k P^k_j(J_\kappa^{-1}, A, \nabla v_\kappa) \partial_t \vk^j,_i
\nonumber \\
&=:{IV_a}_1+{IV_a}_2+{IV_a}_3 \,.
\label{IVa} 
\end{align}
Term  ${IV_a}_3$ is estimated in the same way as term $IV_b$.
For term ${IV_a}_2$,  we integrate by parts in space to find that
\begin{align*}
{IV_a}_2 &= 
-\int_0^T \int_\Gamma
\partial_t^2 \vk^r F^{sk}_{ri} \tilde v ^i,_k \tilde q_{ttt} N_s
+
\int_0^T \int_\Omega  
\partial_t^2 \vk^r (F^{sk}_{ri} \tilde v ^i,_k \tilde q_{ttt}),_s 
:= {{IV_a}_2}_i + {{IV_a}_2}_{ii}
\end{align*}

The first integral ${{IV_a}_2}_i$
is handled identically to term $III_b$ to give the
bound $CT \psE$. We write the second integral as
\begin{align*}
{{IV_a}_2 }_{ii} =
\int_0^T \int_\Omega  
[\partial_t^2 \vk^r F^{sk}_{ri} \tilde v ^i,_k \tilde q_{ttt},_s 
+\partial_t^2 \vk^r (F^{sk}_{ri} \tilde v ^i,_k),_s \tilde q_{ttt}],
\end{align*}
integrate by parts in time, and obtain
\begin{align*}
{{IV_a}_2 }_{ii} & =
-\int_0^T \int_\Omega  
[(\partial_t^2 \vk^r F^{sk}_{ri} \tilde v ^i,_k))_t \tilde q_{tt},_s 
+(\partial_t^2 \vk^r (F^{sk}_{ri} \tilde v ^i,_k)_t,_s \tilde q_{tt}] 
\\
& \qquad\qquad\qquad+
\left. \int_\Omega  
[(\partial_t^2 \vk^r F^{sk}_{ri} \tilde v ^i,_k) \tilde q_{tt},_s 
+(\partial_t^2 \vk^r (F^{sk}_{ri} \tilde v ^i,_k),_s \tilde q_{tt}]\right]_0^T
\,.
\end{align*}
Since $\partial_t^3\vk$ is bounded in $L^2(\Omega)$, and $(F\, \nabla v)_t$
is bounded in $L^\infty(\Omega)$,
the spacetime integral is bounded by $CT\psE$.
Next, we analyze the highest-order term in the
 remaining temporal boundary integral in ${IV_a}_2$:
\begin{align*}
\int_\Omega 
&[(\partial_t^2 \vk^r F^{sk}_{ri} \tilde v ^i,_k) \tilde q_{tt},_s ](T) \\
&\qquad\qquad\qquad
=
\int_\Omega 
(\partial_t^2 \vk^r(0) F^{sk}_{ri}(0) \tilde v ^i,_k(0)) \tilde q_{tt},_s (T)
+ \int_\Omega \int_0^T 
(\partial_t^2 \vk^r F^{sk}_{ri} \tilde v ^i,_k)_t \tilde q_{tt},_s (T)\\
&\qquad\qquad\qquad
\le M_0(\delta) + \delta \|q_{tt}(T)\|_1^2 + T \sup_{t\in [0,T] }
\| (\partial_t^2 \vk^r F^{sk}_{ri} \tilde v ^i,_k )_t \|_0 \, \|q_{tt}(T)\|_1 \\
&\qquad\qquad\qquad
\le  M_0(\delta)+ \delta \sE  + CT \psE \,.
\end{align*}
The remaining term is analyzed in the same fashion.

\noindent
{\bf Step 3. The inertia term.}
Finally, the inertia term in (\ref{eulerttt}) satisfies
\begin{align*}
&\int_0^T \int_\Omega \partial_t^3 (\tilde J_\kappa \tilde
v_t^i) \tilde v_{ttt}^i
= \frac{1}{2}\|\tilde J_\kappa^{\frac{1}{2}} \tilde v_{ttt}(T)\|_0^2
- \frac{1}{2}\|\tilde v_{ttt}(0)\|_0^2 \\
&\qquad\qquad\qquad
+\int_0^T \int_\Omega \left[ 2.5\partial_t \tilde J_\kappa |\tilde v_{ttt}|^2
+ 3\partial_t^2 \tilde J_\kappa \tilde v_{tt}^i\tilde v_{ttt}^i
+ \partial_t^3 \tilde J_\kappa \tilde v_{t}^i\tilde v_{ttt}^i
\right]\,.
\end{align*}
Since $\partial_t \tilde J_\kappa 
= \text{Trace} (\tilde a_\kappa\, \nabla \tilde v_\kappa)$, 
$\partial_t^2 \tilde J_\kappa
= \text{Trace} (\tilde a_\kappa\, \nabla \partial_t\tilde v_\kappa)
+ P(\tilde J_\kappa^{-1}, \tilde a_\kappa, \nabla \tilde v_\kappa)$, 
and 
$\partial_t^2 \tilde J_\kappa
= \text{Trace} (\tilde a_\kappa\, \nabla \partial_t^2\tilde v_\kappa)
+ P(\tilde J_\kappa^{-1}, \tilde a_\kappa, \nabla \tilde v_\kappa)\,
\nabla \partial_t \tilde v_\kappa$, then
using condition (\ref{deteta.c}), we see
that 
$$\sup_{t\in[0,T]}\|v_{ttt}\|_0^2 \le \|\tilde v_{ttt}(0)\|_0^2 
+ C\, T\, \pskE.$$

From (\ref{Neumann}) evaluated at $t=0$, we see that 
$\|\tilde v_{ttt}(0)\|_0^2  \le M_0$, so the lemma is proved.
\end{proof}

\begin{lemma}[Energy estimates for the second time-differentiated 
$\kappa$-problem]  \label{lemma10.3}
For $M_0$ taken as in Lemma \ref{lemma1} and $\delta >0$,  solutions of
the $\kappa$-problem (\ref{smooth}) satisfy:
\begin{align}
&\sup_{t\in[0,T]} |\partial_2\tilde v_t \cdot \tilde n|_0^2
+ \int_0^T |\sk \partial^2\tilde v_{tt} \cdot \nk|_0^2 \nonumber \\
&\qquad 
\le M_0 + T \, P(\sup_{t\in [0,T]} E_\kappa(t))
+ \delta  \sup_{t\in [0,T]} E_\kappa(t) 
+C\,P(\|\sk \tilde v_t\|^2_{L^2(0,T;H^{3.5}(\Omega))}) \,.
\label{ss_kttx}
\end{align}
\end{lemma}

\begin{proof}
We let $\partial \partial_t^2$ act on (\ref{smooth.b}) and test with
$\zeta_i^2 \partial \tilde v_{tt}$ where $\zeta_i^2 = \alpha_i$, and
$\alpha_i$ is an element of our partition of unity.  This localizes the
analysis to a neighborhood of the boundary $\Gamma$ where the tangential
derivative is well-defined. In this neighborhood, we use a normal coordinate
system spanned by $(\partial_1 \eta_0, \partial_2 \eta_0, N)$.

We follow the proof of Lemma \ref{lemma10.2}
and replace $\partial_t^3$ with $\partial \partial_t^2$.   There are only
two differences between the analysis of the second and third time-differentiated
problems.  The first difference can be found in the analogue of term
$III$ in (\ref{Hn_ttt}), which now reads $\int_0^T\int_\Gamma P(\partial 
\tilde \eta, \partial \tilde v)\, \partial \tilde v_t\, \partial^2 \tilde
v_{tt}$.   After integration by parts in space, this term is bounded by
$C\, T\, \pskE$; however, this term requires a bound on $|\tilde v_{tt}|_1$,
which requires us to study the third time-differentiated problem.  (Compare
this with the third time-differentiated problem wherein integration by parts
in time forms an exact derivative which closes the estimate.)

The second difference is significant. Because  the
energy function places
$$\tilde v_{ttt} \in L^\infty(0,T; L^2(\Omega)) \text{ and }
\tilde v_{tt} \in L^\infty(0,T; H^{1.5}(\Omega)),$$ 
there is a one-half
derivative improvement that accounts for (\ref{ss_kttx}) being better than
(\ref{ss_kttt}).

In particular, the analogue of term $II$ in (\ref{improve}) is
$\int_0^T\int_\Gamma \sqrt{\tilde g} \tilde H \tilde n\cdot\nk\,
\partial \partial_t^2 (\nk -\tilde n) \cdot \partial \tilde v_{tt}$,
and since $|\partial \partial_t^2(\nk -\tilde n)|_0 \le C\,
|P(\partial \eta)|_{L^\infty}\, |\tilde v_t|_2$, then this integral is
easily seen to be bounded by $C\, T\, \pskE$.  (This is in sharp contrast
to the difficult analysis required at the level of the third time-differentiated
problem which follows equation (\ref{improve}).)  

All of the other estimates follow identically to the proof of Lemma \ref{lemma10.2} with $\partial_t^3$ being replaced by $\partial \partial_t^2$.
\end{proof}

\begin{lemma}[Energy estimates for the time-differentiated 
$\kappa$-problem]  
For $M_0$ taken as in Lemma \ref{lemma1} and $\delta >0$,  solutions of
the $\kappa$-problem (\ref{smooth}) satisfy:
\begin{align}
&\sup_{t\in[0,T]} |\partial^3 \tilde v \cdot \tilde n|_0^2
+ \int_0^T |\sk \partial_3\tilde v_{t} \cdot \nk|_0^2 \nonumber \\
&\qquad 
\le M_0 + T \, P(\sup_{t\in [0,T]} E_\kappa(t))
+ \delta  \sup_{t\in [0,T]} E_\kappa(t) 
+C\,P(\|\sk \tilde v\|^2_{L^2(0,T;H^{4.5}(\Omega))}) \,.
\label{ss_ktxx}
\end{align}
\end{lemma}
\begin{proof}
After replacing $\partial\partial_t^2$ with $\partial^2\partial_t$,
the proof is the same as the proof of Lemma \ref{lemma10.3}.
\end{proof}
\begin{lemma}[Energy estimates for the $\kappa$-problem]  
For $M_0$ taken as in Lemma \ref{lemma1} and $\delta >0$,  solutions of
the $\kappa$-problem (\ref{smooth}) satisfy:
\begin{align}
&\sup_{t\in[0,T]} |\partial^4 \tilde \eta \cdot \tilde n|_0^2
+ \int_0^T |\sk \partial_4\tilde v \cdot \nk|_0^2 
\le M_0 + T \, P(\sup_{t\in [0,T]} E_\kappa(t))
+ \delta  \sup_{t\in [0,T]} E_\kappa(t)  \,.
\label{ss_kxxx}
\end{align}
\end{lemma}
\begin{proof}
Let $\partial^3$ act on (\ref{smooth.b}) and test with $\partial^3 \tilde v$.
All of the terms are estimated as in Lemma \ref{lemma10.3}, except the
analogue of (\ref{difficult}) which reads, after replacing $\partial_t^3$ with
$\partial^3$, as
\begin{equation}\nonumber
\kappa  \int_0^T[ \tilde v \cdot \partial^4 \nk, 
\nk \cdot \partial^4\tilde v)]_0 \,.
\end{equation}
Since the energy function places 
$\sk \tilde \eta \in L^\infty(0,T; H^{5.5}(\Omega))$, we see that this
integral is bounded by $\delta \skE + C\, T\, \pskE$.
\end{proof}

To the above energy estimates, we add one elliptic estimate
arising from the modified boundary condition (\ref{kbca}).  
We will make use of the following identity:
\begin{align}\label{form1}
[\sqrt{g}\Delta_g(\eta^i)],_\gamma  
= [\sqrt{g} g^{\alpha \beta} \Pi ^i_j \eta^j,_{\beta\gamma} 
+ \sqrt{g}(g^{\nu\mu}g^{\alpha\beta} - g^{\alpha\nu}g^{\beta\mu})
\eta^i,_\beta \eta^j,_\nu \eta^j,_{\mu \gamma}],_\alpha \,.
\end{align}

\begin{lemma}[Elliptic estimate for $\sk\tilde \eta$]
Let $M_0$ be given as in Lemma \ref{lemma1}.  Then for $\delta >0$,
\label{lemma_elliptic}
\begin{equation}\label{ketaestimate}
\sup_{t \in [0,T]}|\sk\tilde \eta|^2_{5}(t)
 \le M_0(\delta) +\delta \skE + C\,T\,  \pskE \,,
\end{equation}
\end{lemma}
\begin{proof}
Letting $Q^{\alpha i}_j:= Q^{\alpha \beta i}_j(\partial \tilde \eta)$ denote a 
smooth function of
$\partial \tilde \eta$, from (\ref{form1}), we see that
\begin{align*}
\partial^3[\sqrt{\tilde g}\Delta_{\tilde g}(\tilde \eta^i)],_\gamma  &= 
 [\sqrt{\tilde g} \tilde g^{\alpha \beta} \Pi ^i_j \partial^3
\tilde \eta^j,_{\beta \gamma} + 
\sqrt{\tilde g}(\tilde g^{\nu\mu}
\tilde g^{\alpha\beta} -\tilde  g^{\alpha\nu}\tilde g^{\beta\mu})
\tilde \eta^i,_\beta \tilde \eta^j,_\nu \partial^3\tilde \eta^j,_{\mu \gamma}],_\alpha \\
& \qquad\qquad 
 + [\partial Q^{\alpha \beta i}_j \partial^2 \tilde \eta^j ,_{\beta \gamma}
 + \partial^2 Q^{\alpha \beta i}_j \partial \tilde \eta^j ,_{\beta\gamma}
 + \partial^3 Q^{\alpha \beta i}_j \tilde \eta^j ,_{\beta\gamma}],_\alpha
\end{align*}
The estimate (\ref{ketaestimate}) is obtained by letting 
$\partial^3\partial_\gamma$ act on the
modified boundary condition (\ref{kbc}) and then testing this  with 
$\zeta^2  \tilde g^{\gamma\delta}\Pi^i_l \partial^3 \tilde \eta^l,_\delta$
where $\zeta^2 = \alpha_i$.

For convenience we drop the subscript $i$ from $\Omega_i$ and $\Gamma_i$. 
(Recall that $\alpha_i$ denotes the partition of unity introduced in Section 2.)
For the surface tension term,
integration by parts with respect to $\partial_\alpha$ yields
\begin{align}
&\intT[\partial^3 (\sqrt{\tilde g}\tilde H\tilde n 
\circ \tilde\eta),_\gamma \,, \zeta^2  \tilde g^{\delta\gamma} 
\Pi \partial^3 \tilde \eta,_\delta ]_0 
 \le \intT\left[ -|\zeta\sqrt{g}^{\frac{1}{2}}  
\Pi \partial^5\eta|_0^2
+ C
|F|_{L^\infty} \, |\eta|_{4}\, |\Pi\partial_5\eta|_{0} \right] \,,
\nonumber
\end{align}
where $F:= P(\zeta, \partial \eta, \partial v, \partial^2\eta)$ is a 
polynomial of its arguments. 
To get this estimate we have used the fact that
$$
\sqrt{\tilde g}(\tilde g^{\nu\mu}
\tilde g^{\alpha\beta} -\tilde  g^{\alpha\nu}\tilde g^{\beta\mu})
\tilde \eta^i,_\beta \tilde \eta^j,_\nu \partial^3\tilde \eta^j,_{\mu \gamma}]
\tilde g^{\gamma\delta}\Pi^i_l \partial^3 \tilde \eta^l,_{\delta\alpha}=0,
$$
since $ \tilde \eta^i,_\beta \Pi^i_l =0$. (This ensures that the error term
is linear in $|\Pi \partial^5 \tilde \eta|_0$ rather than quadratic.)

We next analyze the artificial viscosity term.  The testing procedure gives us
the integral
$$-\intT\kappa [\partial^3\partial_\gamma 
\Delta_0(\tilde v\cdot \nk) \nk, \zeta^2 \tilde g^{\gamma\delta}
\Pi\partial^3 \tilde \eta,_\delta]_0.$$
The positive term comes from $\partial^3\partial_\gamma$ acting on $\tilde v$.
This gives, after integration by parts in space, the highest-order
integrand
$\kappa(\partial^3 \tilde v,_{\alpha\gamma} \cdot \nk)
g^{\gamma\delta} g_0^{\alpha\beta} (\nk \cdot \Pi \partial^3 \tilde \eta
,_{\beta\delta})$, where $\Pi = \tilde n\otimes \tilde n$.
We can write this term as
$$\kappa(\partial^3 \tilde v,_{\alpha\gamma} \cdot \nk)
g^{\gamma\delta} g_0^{\alpha\beta} (\nk \cdot \partial^3 \tilde \eta
,_{\beta\delta}) 
+
\kappa(\partial^3 \tilde v,_{\alpha\gamma} \cdot \nk)
g^{\gamma\delta} g_0^{\alpha\beta} (\nk \cdot(\Pi - \Pi_\kappa) \partial^3 \tilde \eta
,_{\beta\delta}),$$
where $\Pi_\kappa = \nk\otimes \nk $.
The first term is an exact derivative in time, and yields 
$$
\frac{\kappa}{2} \frac{d}{dt} |\partial^5 \tilde \eta \cdot \nk|^2 
- \kappa \partial^5 \tilde \eta^i \,\partial^5\tilde \eta^j\,\, \nk^i 
\partial_t \nk^j \,.
$$
The space integral of the second term is estimated by
$C\,  |F|_{L^\infty} |\sk v|_5\, |\sk \tilde \eta|_5\, |\Pi_\kappa
-\Pi|_{L^\infty}$ and $|\Pi_\kappa - \Pi|_{L^\infty}
\le C\, \kappa\, |\tilde \eta|_{3.5}$. From (\ref{kbc})
$$
\kappa \partial^3\Delta_0 (v\cdot \nk) = - \partial^3(
\sqrt{\gk}^{-1} \sqrt{\tilde g} \Delta
_{\tilde g} (\tilde \eta)\cdot \nk) + \partial^3 q.
$$
Thus,
$$
C\,\sk |\kappa \tilde v|_5 \le \kappa |F|_{L^\infty}( |\tilde v|_3 + 
|\sk \tilde v|_4) + 
(\kappa^{\frac{3}{2}} +\kappa^{\frac{1}{2}} ) |F|_{L^\infty} | \tilde \eta|_4
+ |F|_{L^\infty} |\sk\tilde \eta|_5  
+ \sk|\tilde q|_3 \,,  
$$
so that
\begin{align*}
&\int_0^T\int_\Gamma\kappa(\partial^3 \tilde v,_{\alpha\gamma} \cdot \nk)
g^{\gamma\delta} g_0^{\alpha\beta} 
(\nk \cdot(\Pi - \Pi_\kappa) \partial^3 \tilde \eta,_{\beta\delta}) \\
& \qquad
\le C\, \int_0^T \left[ \sk
|F|_{L^\infty}( |\tilde v|_3 + |\sk \tilde v|_4 + | \tilde \eta|_4
+ |\tilde q|_3) \,|\sk\tilde \eta|_5  
+ \,|F|_{L^\infty} |\sk\tilde \eta|_5^2\right] \,.
\end{align*}
Having finished the estimates for the terms leading to the positive energy 
contribution, we next consider the most difficult of the error terms.  
This occurs when 
$\partial^3 \partial_\gamma$ acts on $\nk$, producing the integral
$\intT \kappa [\zeta^2\,  \tilde v \cdot \partial^4 \nk, 
\nk \cdot \partial^4 \tilde \eta]_0.$

To analyze this term, we let 
$\tilde v =\tilde v^\gamma \partial_\gamma \tilde\eta_\kappa + \tilde v_n \nk;$
we remind the reader
that for $\gamma=1,2$, $\partial_\gamma \tilde \eta_\kappa$ spans the 
tangent space of $\tilde\eta_\kappa(t)(\Gamma)$ at $\tilde\eta_\kappa(x,t)$,
so that $\partial_\gamma \tilde\eta_\kappa$ is orthogonal to $\nk$.  
It follows that 
$v^\gamma \partial_\gamma \tilde \eta_\kappa \cdot \partial^5\nk$ is equal to
$-v^\gamma \nk \cdot \partial^5 \partial_\gamma \tilde \eta_\kappa$ plus lower
order terms $v^\gamma{\mathcal R}^\gamma_4(\tilde \eta)$, which have at most only five 
tangential derivative of $\tilde\eta_\kappa$.
Note also 
that since $\nk\cdot \nk=1$, $\nk\cdot \partial^5\nk$ is also a sum of lower
order terms which have at most only five 
tangential derivative of $\tilde\eta_\kappa$.

Thus,
$$
\intT\kappa [\zeta^2\,  \tilde v \cdot \partial^5 \nk, 
\nk \cdot \partial^5 \tilde \eta]_0
= 
\intT\kappa [\zeta^2\,  \tilde v^\gamma \nk \cdot \partial^5 \partial_\gamma
 \ek + \tilde v^\gamma{\mathcal R}^\gamma_4(\tilde \eta)
, \nk \cdot \partial^5 \tilde \eta]_0  \,,
$$
where the remainder term satisfies 
$$\left|\intT[ \kappa\tilde v^\gamma{\mathcal R}^\gamma_4(\tilde \eta)\, ,
\,  \nk\cdot \partial^5 \eta ]_0\right| 
\le C\, \intT | F|_{L^\infty}\, [ |\sk\tilde v|_4 \,|\sk\tilde \eta|_5 
+|\sk\tilde \eta|_5^2]\,.
$$

We must form an exact derivative from the remaining highest order term
\begin{equation}\label{hot}
\intT\kappa [\zeta^2\,  \tilde v^\gamma \nk \cdot \partial^5 \partial_\gamma
\ek , \nk \cdot \partial^5 \tilde \eta]_0 \,,
\end{equation}
and this will require commuting the horizontal convolution operator, so
that the $\tilde \eta$ on the right side of the $L^2(\Gamma)$ inner-product
also has a convolution operator, and is hence converted to a
$\nk\cdot \partial^4 \rho_{\frac{1}{\kappa}}\star_h\tilde \eta$ term.  
With this accomplished, we
will be able to pull-out the $\partial_\gamma$ operator and form an exact
derivative, which can be bounded by our energy function.

Noting that on $\Gamma$ the horizontal convolution $\star_h$
restricts to the usual convolution $*$ on ${\mathbb R}^2$, we have that
$$
\tilde \eta_\kappa= \sum_{i=1}^K  
\sqrt{\alpha_i}\left[ \rho_{\frac{1}{\kappa}}  * [\rho_{\frac{1}{\kappa}} * (\sqrt{\alpha_i}
\tilde\eta \circ \theta_i)] \right] \circ \theta_i^{-1} 
$$
For notational convenience, we set $\rho= \rho_{1/\kappa}$, $\zeta= 
\sqrt{\alpha_i}$, and $R=[0,1]^2 = \theta_i^{-1} (\Gamma\cap U_i)$.  
It follows that (\ref{hot}) can be expressed as
\begin{equation}\label{hot1}
\intT \kappa \sum_{i=1}^K \int_R 
(\tilde v^\gamma \nk )\circ \theta_i \cdot  \partial_\gamma \partial^5 
\left[ \zeta(\theta_i) \rho  * \rho * (\zeta
\tilde\eta) \circ \theta_i \right]  \, 
\left(\nk \cdot
\partial^5 \tilde \eta\right) \circ \theta_i
\end{equation}
With $\gk := \partial^5 \rho * \zeta(\theta_i)\tilde\eta(\theta_i)$, 
we see that
\begin{equation}\label{hot2}
\kappa 
\nk \cdot  \partial_\gamma\partial^5 \left[ \zeta(\theta_i) \rho  * \rho * 
(\zeta
\tilde\eta \circ \theta_i) \right]
= \kappa \nk \cdot (\partial_\gamma \partial^5 \zeta(\theta_i)) \rho *\rho*(\zeta
\tilde\eta(\theta_i) + 
\kappa \zeta \nk \cdot \rho * \partial_\gamma \gk
+ \kappa {\mathcal R}_5(\tilde\eta)\,,
\end{equation}
where the remainder ${\mathcal R}_5(\tilde\eta)$ has at most five tangential
derivatives on $\tilde \eta$.  Substitution of (\ref{hot2}) into (\ref{hot1})
yields three terms, corresponding to the three terms on the right-hand side
of (\ref{hot2}).  
For the first term, we see that
\begin{align}
&\intT \kappa \sum_{i=1}^K \int_R 
(\partial_\gamma \partial^5 \zeta(\theta_i)) \tilde v^\gamma\nk \cdot \, \rho *\rho*(\zeta
\tilde\eta(\theta_i)  \, \left(\nk \cdot \partial^5 \tilde\eta(\theta_i)
\right)  \nonumber 
\le C \intT  \|\sk\theta_i\|_{6.5} 
|\tilde\eta(\theta_i)|_{L^\infty}
|\sk\tilde\eta|_5 \,.
\end{align}

The  second term on the right-hand side of (\ref{hot2}) gives the integral
$$
\intT \kappa \sum_{i=1}^K \int_R 
(\nk  \cdot \rho * \partial_\gamma \gk) \, \, (\nk \cdot \gk) + {\mathcal R}_6
(\tilde\eta),
$$
where the remainder ${\mathcal R}_6(\tilde\eta)$ is lower-order containing
terms which have at most four tangential derivatives on $\tilde \eta$ and
five on $\zeta(\theta_i)$.

We fix 
$i \in \{1,...,K\}$, drop the explicit composition with $\theta_i$, we set
\begin{align*}
\triangle_{\nk, g_\kappa} &=  \nk \cdot \rho *\gk - \rho *(\nk \cdot \gk)\,, \\
\triangle_{\nk, \zeta\partial^5\tilde \eta} 
&=  \nk \cdot \rho *\zeta\partial^5\tilde \eta - \rho *(\nk \cdot \zeta\partial^5\tilde \eta)\,,
\end{align*}
and analyze the following integral:
\begin{align*}
&\int_R \{\zeta\nk \cdot \rho *\gk\} \, \{ \nk \cdot \partial^5 \tilde \eta\}
=
\int_R \{\rho *(\nk \cdot \gk)\} \, \{ \nk \cdot \zeta\partial^5 \tilde \eta\}
+ \int_R \triangle_{\nk,g_\kappa} \, \{ \nk \cdot \zeta\partial^4 \tilde \eta\} \,, \\
&\qquad =
\int_R \{\nk \cdot \gk\} \, \{ \nk \cdot \partial^4\rho* \zeta\tilde \eta\}
-\int_R \{\nk \cdot \gk\} \, \triangle_{\nk, \zeta\partial^5\tilde \eta} 
+ \int_R \triangle_{\nk,g_\kappa} \, \{ \nk \cdot \zeta\partial^5 \tilde \eta\} 
+{\mathcal R}_7(\tilde \eta)\,, \\
\end{align*}
where the remainder ${\mathcal R}_7(\tilde \eta)$ comes from commuting $\partial^5$ with
the cut-off function $\zeta$ and has the same bound as ${\mathcal R}_4(\tilde \eta)$.
The first term on the right-hand side is a perfect derivative and for the
remaining terms we use Lemma \ref{commutator} together with the estimate
$\kappa |\tilde g_\kappa|_{0,R} \le C \|\tilde \eta\|_{5.5}$ to find that
$$
\kappa\int_R \{\zeta\nk \cdot \rho *\gk\} \, \{ \nk \cdot \partial^5 \tilde \eta\}
\le C\, |F|_{L^\infty}\,|\sk\tilde \eta|_5^2 .
$$

Thus, summing over $i \in \{1,...,K\}$,  
\begin{equation}\nonumber
\kappa
\intT[\zeta^2\, \tilde v \cdot \partial^5 \nk, \nk \cdot \partial^5 \tilde\eta]_0
\le  C\,\intT |F|_{L^\infty}\, (|\sk\Gamma|_6 + |\sk\tilde v|_4 
+ |\sk \tilde \eta|_5)\, |\sk\tilde \eta|_5 \,,
\end{equation}
where $|\sk\Gamma|_6:= \max_{i\in\{1,...,K\}} |\sk\theta_i|_6$.

It is easy to see that
$$
\intT [\partial^3\partial_\gamma ( \sqrt{\tilde g} \Delta_{\tilde g}(\tilde
\eta) \cdot (\nk -\tilde n) \nk), \zeta^2 \tilde g^{\delta\gamma} \Pi 
\partial^3 \tilde \eta,_\delta] 
\le  C\,\intT |F|_{L^\infty}\,  |\sk\tilde \eta|_5^2\,
$$
with the same bound for
$ \intT [\partial^2\partial_\gamma ( \sqrt{\tilde g} \Delta_{\tilde g}(\tilde
\eta) \cdot \nk (\nk -\tilde n) , \zeta^4 \tilde g^{\delta\gamma} \Pi 
\partial^2 \tilde \eta,_\delta] $.
With (\ref{deteta}),  we infer that
\begin{align*}
|\sk \partial^5 \tilde \eta\cdot \nk|_0^2 (t) 
\le M_0+ C\,\intT|\tilde q|_{3}^2 
+  C\,\intT |F|_{L^\infty}\, (|\sk\Gamma|_6 + |\sk\tilde v|_4 + |\sk \tilde \eta|_5)\, |\sk\tilde \eta|_5 \,.
  \end{align*}
Adding to this inequality the curl estimate (\ref{curl_eta_2D}) for $\sk\tilde \eta$ and
the divergence estimate (which has the same bound as the curl estimate), 
and using Young's inequality to get
$ \intT |F|_{L^\infty}\, |\sk\tilde v|_4 \, |\sk\tilde \eta|_5 \le 
\delta \intT|\sk\tilde v|_4^2 + C_\delta
\intT |F|_{L^\infty}^2\, |\sk\tilde \eta|_5^2$,
we see
that
\begin{align*}
&\sup_{t\in[0,T]}| \sk \tilde \eta|_5^2 
\le M_0+  C\, T\, \sup_{t\in[0,T]}\left(\|v_t\|_{2.5}^2 + \|v\|_{3.5}^2+
\|\eta\|_{4.5}^4 \right) \\
&\qquad
+  C\,T\, \sup_{t\in[0,T]}\left(|F|_{L^\infty}^2\, (|\sk\Gamma|_6^2 
+ |\sk \tilde \eta|_5^2 \right)
+ \delta \intT|\sk\tilde v|_4^2 \,,
\end{align*}
from which the lemma follows.
\end{proof}

\subsection{Removing the additional regularity assumptions on the initial
data} \label{initdatacondition}
At this stage, we explain 
how we can remove the extra regularity assumptions on the initial
data, $u_0$ and $\Omega$, so that the constant $M_0$ depends on 
$\|u_0\|_{4.5}$ and $|\Gamma|_{5.5}$ rather than
$\sk\|u_0\|_{10.5}$  and $\sk|\Gamma|_{7}$ as stated in Lemma \ref{lemma1}.  
The modification requires the following regularization of the initial data:  set
$u_0 = \rho_{e^{-\kappa}}\star E_{\Omega_\kappa}(u(0))$, where
$\Omega_\kappa$ is obtained by smoothing $\Omega$ via convolution with
$\rho_{e^{-\kappa}}$, i.e., we use
$\rho_{e^{-\kappa}}* \bar\theta_i$ as our family of charts.
We make use of the fact that
$$
P(\|\sk u_0\|_{10.5}) \le C P\|u_0\|_{4.5}) \,, \ \ \ 
P(|\sk \Gamma|_{10}) \le C P(|\Gamma|_{5.5}) \,, \ \ \ 
$$
which follows by integration by parts of six tangential derivatives onto the 
mollifier  $\rho_{e^{-\kappa}}$; this
and results in the constant $C>0$ being independent of $\kappa$.

\subsection{The limit as $\kappa \rightarrow 0$}
\begin{proposition}\label{prop3}
With $M_0 = P(\|u_0\|_{4.5}, |\Gamma|_{5.5})$ a polynomial of its arguments
and for
$\tilde M_0 > M_0$,
\begin{align}
\sup_{t\in[0,T]} E_\kappa(t) \le \tilde M_0\,,
\label{Ekbound}
\end{align}
where $T$ depends on the data, but not on $\kappa$.
\end{proposition}
\begin{proof}
Summing the inequalities (\ref{ss_kttt}), (\ref{ss_kttx}), (\ref{ss_ktxx}), 
(\ref{ss_kxxx}), and (\ref{ketaestimate}),  
and using Lemma \ref{lemma1} and Proposition \ref{prop1},
we find that
\begin{align*}
\sup_{t\in[0,T]} E_\kappa(t) \le M_0 + C\, T\, P(\sup_{t\in[0,T]} 
E_\kappa(t)) + \delta \sup_{t\in[0,T]} E_\kappa(t)
\,,
\end{align*}
where the polynomial $P$ and the constant $M_0$ do not depend on
$\kappa$.  Choose $\delta< 1$.  Then,
from the continuity of the left-hand side in $T$, we may choose 
$T$ sufficiently small and independent of $\kappa$, to ensure that 
(\ref{Ekbound}) holds.  (See \cite{CoSh2005b} for a detailed account
of such polynomial  inequalities.)
\end{proof}

Proposition \ref{prop3} provides the weak convergence as $\kappa \rightarrow 0$
of subsequences
of $(\v,\q)$ toward a limit which we denote by $(v,q)$ in the
same space.  We then set $\eta= \operatorname{Id} + \int_0^t v$, and
$u = v \circ \eta^{-1}$. 
It is obvious that $\v_\kappa$, arising from the double horizontal convolution
by layers of $\v$, satisfies 
$\v_\kappa \rightarrow v$ in $L^2(0,T; H^{3.5}(\Omega))$, and therefore
$\tilde\eta^\kappa \rightarrow \eta$ in $L^2(0,T; H^{4}(\Omega)$.
It follows that ${\operatorname{div}}u=0$ in 
$\eta(\Omega)$ in the limit as $\kappa \rightarrow 0$ in (\ref{leuler.c}).
Thus, the limit
$(v,q)$ is a solution to the problem (\ref{leuler}), and satisfies
$E_0(t) \le \tilde M_0$.   We then take $T$ even small, if necessary, to 
ensure that (\ref{deteta}) holds, which follows from the fundamental theorem
of calculus.

\section{A posteriori elliptic estimates}
Solutions of the Euler equations gain regularity with respect to the
$E_\kappa(t)$ from elliptic estimates of the boundary condition (\ref{leuler.d}),
which we write as $\sqrt{g}Hn(\eta) = \sqrt{g} q n$.

Replacing $\partial_\gamma$ with $\partial_t$ in (\ref{form1}), we have the
identities
\begin{align} \label{form2}
\partial_t(\sqrt{g}Hn \circ \eta)^i
=- [\sqrt{g} g^{\alpha \beta} \Pi ^i_j v^j,_{\beta} 
+ \sqrt{g}(g^{\nu\mu}g^{\alpha\beta} - g^{\alpha\nu}g^{\beta\mu})
\eta^i,_\beta \eta^j,_\nu v^j,_{\mu }],_\alpha \,.
\end{align}
and
\begin{align} \label{form3}
\partial_t^2(\sqrt{g}Hn \circ \eta)^i
=- [\sqrt{g} g^{\alpha \beta} \Pi ^i_j v_t^j,_{\beta} 
+ \sqrt{g}(g^{\nu\mu}g^{\alpha\beta} - g^{\alpha\nu}g^{\beta\mu})
\eta^i,_\beta \eta^j,_\nu v_t^j,_{\mu } 
+ Q^{i \alpha\beta}_{j} v^j,_\beta],_\alpha \,,
\end{align}
where $Q^{i \alpha\beta}_{j} = Q(\partial \eta)$ is a rational function
of $\partial \eta$.

\begin{lemma} \label{euler_elliptic}
Taking $\tilde M_0$ as in Proposition \ref{prop3}, and letting
${\mathcal M}_0$ denote a polynomial function of $\tilde M_0$,  for $T$
taken sufficiently small,
$$
\sup_{t\in[0,T]}
\left[
| \Gamma(t)|_{5.5} + \|v(t)\|_{4.5} + \|v_t(t)\|_3\right] \le {\mathcal M}_0 \,.
$$
\end{lemma}
\begin{proof}
We begin with the estimate for $v_t$.  
Following the proof Lemma \ref{lemma_elliptic},
we let $\partial_\gamma \partial_t^2$
act on the
boundary condition (\ref{leuler.d}) and test with  
$-\zeta^2 g^{\gamma\delta}\Pi^i_k v^k_t,_\delta$, where $\zeta^2 = \alpha_i$,
an element of our partition of unity.
Using (\ref{form3}), we see that
$-\int_\Gamma \partial_\gamma \partial_t^2[ \sqrt{g} Hn(\eta)] \cdot 
\zeta^2 g^{\gamma\delta} v_t,_\delta = \int_\Gamma (\sqrt{g} q n)_{tt} \cdot
[\zeta^2 g^{\gamma\delta} v_t,_{\delta}],_\gamma$.  Using (\ref{Ekbound}),
letting $\tilde C$ denote a constant that depends on $\tilde M_0$, and summing
over the partition of unity, we find that
\begin{align}
|\partial^2 v_t\cdot n|_0^2 \le  \tilde C\, [|v|_2 + |\eta|_2\, 
(|v_t|_1+ |v|_1) + |(q\sqrt{g}n)_{tt}|_0 ] \, |\partial^2v_t\cdot n|_{0} 
+ \tilde C |v_t|_1^2 |\eta|_3 \,.
\label{j1}
\end{align}

This follow since
$$
\left[\sqrt{g}(g^{\nu\mu}g^{\alpha\beta} - g^{\alpha\nu}g^{\beta\mu})
\eta^i,_\beta \eta^j,_\nu v_t^j,_{\mu\gamma }  \right]
\left[
 g^{\gamma\delta}\Pi^i_k v^k_t,_{\delta \alpha}\right] =0 \,,
$$
while
$$
\int_\Gamma
\left[\sqrt{g}(g^{\nu\mu}g^{\alpha\beta} - g^{\alpha\nu}g^{\beta\mu})
\eta^i,_\beta \eta^j,_\nu v_t^j,_{\mu\gamma }  \right]
\left[
 (g^{\gamma\delta}\Pi^i_k),_\alpha v^k_t,_{\delta }\right] \le
\tilde C\, |v_t|_1\, (|\partial^2 v_t \cdot n|_0 + |v_t|_1\,|\eta|_3)\,.
$$
Applying Young's inequality  to (\ref{j1}) yields, after adjusting the constant,
\begin{align*}
|\partial^2 v_t\cdot n|_0^2 \le  \tilde C\, [|\eta|_3^2 + |v|_2 ^2
+ |v_t|_1^2 + |q_{tt}|_0^2]  \,.
\end{align*}
A similar computation shows that
\begin{align*}
|\partial^3 v_t\cdot n|_0^2 \le  \tilde C\, [|\eta|_4^2 + |v|_3 ^2
+ |v_t|_2^2 + |q_{tt}|_1^2]  \,.
\end{align*}
Thus, by interpolation
\begin{align*}
\sup_{t\in[0,T]}
|\partial^2 v_t\cdot n|_{0.5}^2 \le  \tilde C\, 
\sup_{t\in[0,T]} [|\eta|_{3.5}^2 + |v|_{2.5} ^2 + |v_t|_{1.5}^2 + |q_{tt}|_{0.5}^2] 
\le {\mathcal M}_0 \,. 
 \end{align*}

Computing the $H^2(\Omega)$-norm of (\ref{zs3b}), we find that
$$\sup_{t\in[0,T]}\|\curl v_t \|_2 \le {\mathcal M}_0 + 
\tilde C\, T\, \sup_{t\in[0,T]}\|v_t \|_3\,, $$
with the same estimate for $\sup_{t\in[0,T]}\|\div v_t \|_2$.
Hence, for $T$ taken sufficiently small, we infer from Proposition \ref{prop1}
that  
\begin{equation}\label{vth3}
\sup_{t\in[0,T]}\|v_t\|_3 \le {\mathcal M}_0.
\end{equation}

Next, we let $\partial_\gamma \partial^2 \partial_t$
act on the boundary condition (\ref{leuler.d}) and test with  
$-\zeta^2 g^{\gamma\delta}\Pi^i_k \partial^2v^k,_\delta$.  Using (\ref{form1}),
we find that
\begin{align*}
\sup_{t\in[0,T]}
|\partial^4 v\cdot n|_0^2 \le  \tilde C\, 
\sup_{t\in[0,T]} [|\eta|_4^2 + |v|_3 ^2 + |q_t|_2^2]   \le {\mathcal M}_0\,.
\end{align*}
Computing the $H^{3.5}(\Omega)$-norm of (\ref{curlvss}), and again taking
$T$ sufficiently small, we see that $\sup_{t\in[0,T]}\|v\|_{4.5} \le
{\mathcal M}_0$.

In order to prove our remaining estimate, we need a convenient
reparameterization of $\Gamma(t)$  via a height function $h$ in
the normal bundle over $\Gamma$.

Consider the isometric immersion $\eta_0:(\Gamma,g_0) \to
({\mathbb R}^3,\text{Id})$. Let ${\mathcal B}=\Gamma\times(-\epsilon,\epsilon)$ where 
$\epsilon$ is chosen sufficiently small so that the
map
$B:{\mathcal B}\to\bbR^3: (y,z) \mapsto y+zN(y)$
is itself an immersion, defining a tubular neighborhood of $\eta_0(\Gamma)$ 
in $\bbR^3$. We can choose a coordinate system
$\frac{\partial}{\partial y^\alpha}$, $\alpha=1,2$ and 
$\frac{\partial}{\partial z}$.

Let $G=B^*(\id)$, denote the induced metric on ${\mathcal B}$, and note that
$G(y,z)=G_z(y)+dz\otimes dz,$
where $G_z$ is the metric on the surface $\Gamma\times\{z\}$, and that 
$G_0=g_0$.
Let $h:\Gamma\to(-\epsilon,\epsilon)$ be a smooth function and consider the 
graph of $h$ in ${\mathcal B}$, parameterized by
$\phi:\Gamma\to{\mathcal B}:y\mapsto (y,h(y))$. 
The tangent space to graph(h), considered as a submanifold of ${\mathcal
B}$, is spanned at a point $\phi(x)$ by the vectors
$$\phi_*(\frac{\partial}{\partial y^\alpha}) = \frac{\partial\phi}{\partial y^\alpha}=\frac{\partial}{\partial y^\alpha} +
\frac{\partial h}{\partial y^\alpha}\frac{\partial}{\partial z},$$
and the normal to graph(h) is given by
\begin{equation}\label{nn}
n(y)=J_h^{-1}(y)\Big(-G^{\alpha\beta}_{h(y)}\frac{\partial h}{\partial y^\alpha}\frac{\partial}{\partial y^\beta} +
\frac{\partial}{\partial z}\Big)
\end{equation}
where $J_h=(1+h_{,\alpha}G^{\alpha\beta}_{h(y)} h_{,\beta})^{1/2}$. Therefore, 
twice the mean curvature $H$ is defined to be the
trace of $\nabla n$ while
\begin{align*}
(\nabla n)_{ij}=G(\nabla^{\mathcal B}_{\frac{\partial}{\partial w^i}} n,\frac{\partial}{\partial w^j})
\end{align*}
where $\frac{\partial}{\partial w^\alpha} = \frac{\partial}{\partial y^\alpha}$ for $\alpha=1$, $2$ and $\frac{\partial}{\partial
w^3}=\frac{\partial}{\partial z}$. Substituting the formula (\ref{nn}) for $n$,
we see that
\begin{align*}
(\nabla n)_{\alpha\beta} =&\ G\Big(\nabla^{\mathcal B}_{\frac{\partial}{\partial y^\alpha}}\Big[-J_h^{-1}
G_h^{\gamma\delta}h_{,\gamma}\frac{\partial}{\partial y^\delta}+J_h^{-1}\frac{\partial}{\partial z}\Big],\frac{\partial}
{\partial y^\beta}\Big)\\
=&\ -(G_h)_{\delta\beta}(J_h^{-1}G_h^{\gamma\delta}h_{,\gamma})_{,\alpha} + F^1_{\alpha\beta}(y,h,\partial h); \\
(\nabla n)_{33} =&\ G\Big(\nabla^{\mathcal B}_{\frac{\partial}{\partial z}}\Big[-J_h^{-1}
G_h^{\gamma\delta}h_{,\gamma}\frac{\partial}{\partial y^\delta}+J_h^{-1}\frac{\partial}{\partial z}\Big],\frac{\partial}
{\partial z}\Big) \\
=&\ F^2_{\alpha\beta}(y,h,\partial h)
\end{align*}
for some functions $F_{\alpha\beta}^1$ and $F_{\alpha\beta}^2$, $\alpha,\beta=1,2$. Letting $\gamma_0$ denote the Christoffel
symbols associated to the metric $g_0$ on $\Gamma$, we find that the curvature 
of graph($h$) is given by
\begin{align}
L_h(h):= H = -(J_h^{-1}G_h^{\gamma\delta}h_{,\gamma})_{,\delta} + 
J_h^{-1}([\gamma_0]^j_{j3} - G_h^{\gamma\delta}h_{,\gamma}[\gamma_0]^j_{j\delta})\,.
\label{divform}
\end{align}
Note that the metric $G_h = P(h)$,
and that the highest-order term is in divergence form, while the lower-order
term is a polynomial in $\partial h$.   The function $h$ determines
the {\it height}, and hence shape, of the surface $\Gamma(t)$ 
above $\Gamma$.  

Given a signed height function $h: \Gamma_0 \times [0,T) \rightarrow
{\mathbb R}$, for each $t\in [0,T)$, define the {\it normal} map
$$
\eta^\nu : \Gamma_0 \times [0,T) \rightarrow \Gamma(t), \ \ \ (y,t) \mapsto y+ h(y,t) N(y)\,.
$$
Then, there exists a unique {\it tangential} map $\eta^\tau: \Gamma_0 \times
[0,T)\rightarrow \Gamma_0 $ (a diffeomorphism as long as $h$ remains
a graph) such that $\eta|_\Gamma(t)$ has the decomposition
$$\eta|_\Gamma(\cdot, t) = \eta^\nu(\cdot,t) \circ \eta^\tau(\cdot, t), \ \ \
\eta|_\Gamma(y,t) = \eta^\tau(y,t) + h(\eta^\tau(y,t),t) N(\eta^\tau(y,t))\,.$$
The boundary condition (\ref{euler.c}) 
can be written as $\sigma L_h(h)= q\circ (\eta^\tau)^{-1}$.

The operator $L_h$ is a quasilinear elliptic operator; 
from the standard regularity theory for quasilinear elliptic operators with 
$H^3$ coefficients on a compact manifold, we have the elliptic estimate
$$
|h|_{5.5} \le \tilde C\, |q \circ (\eta^\tau)^{-1}|_{3.5} \le \tilde C\,
\|q\|_{4} \,.
$$
By (\ref{Neumann}), we see that for all $t\in[0,T]$,
$$
\|q\|_4 \le \tilde C\, \|a\|_2^2\, \|v\|_3^2 + \tilde C \,|\eta|_3 \, |v_t|_2
\le {\mathcal M}_0\,,
$$
the last inequality following from (\ref{vth3}).

Since $\Gamma(t)=$graph$h(t)$,
this estimate shows that $\Gamma(t)$ maintains its $H^{5.5}$-class regularity
for $t\le T$.
\end{proof}

\section{$\kappa$-independent estimates for the smoothed problem and
existence of solutions in 3D} \label{kapriori3}

The 3D analysis of the $\kappa$-problem requires assuming that 
the initial data
$ u_0 \in H^{5.5}(\Omega)$
and $ \Gamma$ of class $H^{6.5}$.  This is necessitated by Sobolev embedding
$\|\tilde v_t\|_{L^\infty} \le C \|\tilde v_t\|_3$.

By replacing the third-time differentiated problem with the fourth 
time-differentiated problem the identical analysis as in Section \ref{kapriori}
yields
$$
E^{3D}_\kappa(t) \le \tilde M_0\,,
$$
where $\tilde M_0$ is a polynomial of $\| u_0\|_{5.5}$ and $| \Gamma|_{6.5}$.
(In fact, our analysis in Section \ref{kapriori} used all of the 3D terms
and notation, so no changes are required other than raising the regularity
by one derivative.)

We let $(v,q)$ again denote the limit of $(\tilde v, \tilde q)$ as $\kappa
\rightarrow 0$.  The identical limit process as in 2D, shows that 
$(v,q)$ is a solution of the Euler equations.

Having a solution $(v,q)$ to the Euler equation, we can use the 
a posteriori estimates (\ref{euler_elliptic}) as a priori estimates for
solutions of the Euler equations.  We see that
$$
\sup_{t\in[0,T]}
\left[
E_\kappa^{2D}(t) +
| \Gamma(t)|_{5.5} + \|v(t)\|_{4.5} + \|v_t(t)\|_3\right] \le {\mathcal M}_0 \,,
$$
where ${\mathcal M}_0$ is polynomial function of 
$\|u_0\|_{4.5}$ and $|\Gamma|_{5.5}$.   
This key point here, is that the elliptic estimate for $v_t \in H^3(\Omega)$
improves the regularity given by $E_\kappa^{2D}(t)$ and allows for the
required Sobolev embedding theorem to hold.

Since our initial data is a priori assumed regularized as in 
Subsection \ref{initdatacondition}, we see that solutions of the Euler
equations in 3D only depend on ${\mathcal M}_0$.

\section{Uniqueness of solutions to (\ref{leuler})}\label{uniqueness}

Suppose that
 $(\eta^1,v^1,q^1)$ and $(\eta^2,v^2,q^2)$ are both solutions of 
(\ref{leuler})  with initial data $u_0 \in H^{5.5}(\Omega$ and $\Gamma \in H^{6.5}$.  
Setting
$$
{\mathcal E}_\eta(t) = \sum_{k=0}^4 \|\partial_t^k \eta(t)\|_{5.5-k}^2 ,
$$
by the method of Section \ref{kapriori} with $\kappa=0$, we infer that
both ${\mathcal E}_{\eta^1}(t)$ and ${\mathcal E}_{\eta^2}(t)$ are bounded by 
a constant ${\mathcal M}_0$ depending on the data $u_0$ and $\Gamma$ on a time
interval $0\le t\le T$ for $T$ small enough.

Let
$$
w: = v^1-v^2, \ \  r: = q^1 -  q^2, \  \text{ and } \xi := \eta^1 - \eta^2 \,.
$$
Then $(\xi,w,r)$ satisfies
\begin{subequations}
  \label{lunique}
\begin{alignat}{2}
\xi&=\int_0^t w\ \ \   &&\text{in} \ \Omega \times (0,T]\,, \label{lunique.a}\\
\partial_t w^i+ (a^1)^k_i\,r,_k&=(a^2-a^1)^k_i\, q^2,_k  &&\text{in} \ \Omega \times (0,T]\,, \label{lunique.b}\\
   (a^1)^j_i w^i,_j  &= (a^2-a^1)^j_i {v^2}^i,_j     &&\text{in} \   \Omega \times (0,T]
\,, \label{lunique.c}\\
   r n_1 &=  -\sigma\Pi^1 {g^1}^{\alpha\beta} \xi,_{\alpha\beta}
-\sigma\sqrt{g^1} \Delta_{g^1-g^2}(\eta^2) \quad\quad   &&\text{on} \ \Gamma\times(0,T]
                                                 \,, \label{lunique.d} \\
   (\xi,w) &= (0, 0)  &&\text{on} \ \Omega\times\{t=0\} 
                                                 \,. \label{lunique.e} 
\end{alignat}
\end{subequations}
Set
$$
E(t) = \sum_{k=0}^3 \|\partial_t^k \xi(t)\|_{4.5-k}^2 .
$$
We will show that $E(t)=0$, which shows that $w=0$.
To do so, 
we analyze the  forcing terms on the right-hand side of (\ref{lunique.b}) 
and (\ref{lunique.c}),
as well as the  term $\sigma\Delta_{g^1-g^2}(\eta^2)$ in (\ref{lunique.d}).
We begin with the third time-differentiated problem, and study the integral
$\int_0^T \int_\Omega \partial_t^3 [ (a^1-a^1) \, \nabla q^2]\, w_{ttt}$.
The highest-order term is 
\begin{align*}
\int_0^t\int_\Omega (a^1-a^2) \, \nabla q^2_{ttt} w_{ttt} &
\le {\frac{1}{2}}{\mathcal M}_0\int_0^t\|a^1-a^2\|_{L^\infty}^2 + 
{\frac{1}{2}}\int_0^T \|w_{ttt}\|_0^2 
\le C\, t\, P(E(t))\,.
\end{align*}
The third space differentiated and mixed-derivative problems have forcing
terms that can be similarly bounded.

The difference in pressure $r$ satisfies, using the notation of (\ref{Neumann}),
satisfies the following Neumann problem:
\begin{align*}
L_{a_1}(r) & = -\partial_t {a^1}^j_i w^i,_j + [{a^1}^j_i (a^2-a^1)^k_i q^2,_k],_j
\text{ in } \Omega \,, \\
B_{a_1}(r) & = -w_t\cdot \sqrt{g^1} n^1 + {a^1}^j_i (a^2-a^1)^k_i q^2,_k N_j
\text{ on } \Gamma \,. 
\end{align*}
Since $P(\|\eta^1\|_{4.5})$ is bounded by some constant $C=P({\mathcal M}_0)$,
(\ref{elliptic}) provides the estimate
$$
\|r\|_{3.5} \le C [ \|a^1_t\|_{1.5} \, \|w\|_{2.5} + \|a^1\|_{2.5}\, \|q^2\|_{3.5}\, \|a^1-a^2\|_{2.5}
+ |\sqrt{g^1}n^1|_{2}\, \|w_t\|_{2.5}] \,.
$$
Since $\|a^1-a^2\|_{2.5} \le C \|\xi\|_{3.5}$, and $\|a^1_t\|_{1.5}$, 
$\|a^1\|_{2.5}$, 
$\|q^2\|_{3.5}$, and $|\sqrt{g^1}n^1|_{2}$ are all bounded ${\mathcal M}_0$,
we see that $\|r(t)\|_{3.5} \le C P(E(t))$.  Similar estimates for the 
time derivatives of $r$ show that 
$\|r(t)\|_{3.5} + \|r_t(t)\|_{2.5}+ \|r_{tt}(t)\|_{1}
 \le C P(E(t))$.  

This shows that the energy estimates of Section \ref{kapriori} 
go through unchanged for equation (\ref{lunique});  therefore, using 
(\ref{lunique.e}), we see that
 $\sup_{t\in[0,T]} E(t) \le 
C\,T\,P(\sup_{t\in[0,T]} E(t) )$.


\def\R{\mathbb R}
\def\n{\nonumber}
\def\d{\displaystyle}
\def\w{\tilde w}
\def\a{\bar a}
\def\b{\tilde b}
\def\v{\tilde v}
\def\q{\tilde q}
\def\e{\tilde\eta}
\def\u{\tilde u}
\def\p{\tilde p}
\def\E{\tilde{E}}
\def\o{\operatorname}
\def\hf{\hfill\break}
\def\ek{{\eta^{\kappa}}}
\def\vk{{v^\kappa}}
\def\uk{{u^\kappa}}
\def\pk{{p^\kappa}}
\def\qk{{q^\kappa}}
\def\nk{{n^\kappa}}
\def\NK{{N^\kappa}}
\def\s{\starh}
\def\ur{\tilde u^{\kappa}}
\def\ud{{\frac{1}{2}}}
\def\td{{\frac{3}{2}}}
\def\cd{{\frac{5}{2}}}
\def\ci{{\frac{5}{2}}}
\def\sd{{\frac{7}{2}}}
\def\se{{\frac{7}{2}}}
\def\ud{{\frac{1}{2}}}
\def\o{\Omega}
\def\s{\starh}
\def\en{\bar{\eta}}
\def\uk{{\bar{u}^{\kappa}}}
\def\nk{{\bar{n}^{\kappa}}}
\def\ek{{{{\bar{\eta}^{\kappa}}}}}
\def\an{{\bar{a}^\kappa}}
\def\tak{{{\tilde a}^\kappa}}
\def\tvk{{{\tilde v}^\kappa}}
\def\un{{\bar{u}^n}}
\def\div{\operatorname{div}}
\def\curl{\operatorname{curl}}
\def\tek{{{\tilde\eta}^\kappa}}
\def\tuk{{{\tilde u}^\kappa}}
\def\tbk{{{\tilde b_l}^\kappa}}
\def\ttek{{{\tilde\eta}_{l\kappa}}}
\def\ttvk{{{\tilde v}_\kappa}}
\def\ttuk{{{\tilde u}_{l\kappa}}}
\def\O{{(0,1)^2}}
\def\starh{\star_h}

\section{The zero surface tension case $\sigma=0$} \label{L1}

In this, the second part of the paper, we use our 
methodology to prove well-posedness of the
free-surface Euler equations with $\sigma=0$ and the Taylor sign condition
(\ref{lindblad}) imposed, 
previously established by Lindblad  in \cite{Li2004}. The main advantages
of our method over the Nash-Moser approach of \cite{Li2004} are the 
significantly shorter proof and the fact that we 
provide directly the optimal space in which the problem is set, 
instead of having to separately perform an optimal energy study once a solution
is known as in \cite{ChLi2000}. If one uses a Nash-Moser approach without 
performing the analysis of \cite{ChLi2000}, then one obtains results with 
much higher regularity requirements than necessary, as for instance in 
\cite{La2005} for the irrotational water-wave problem without surface tension. 
We also obtain lower regularity results than those given by the functional
framework of \cite{ChLi2000} for the 3D case.

We will extensively make use of the horizontal convolution by layers 
defined in Section \ref{2}, and just as in the first part of the paper, 
for $v\in L^2(\Omega)$ and $\kappa\in (0,\kappa_0)$, we define the smoothed
velocity $v^\kappa$ by
\begin{equation*}
\d v^{\kappa}=\sum_{i=1}^K \sqrt{\alpha_i} \bigl[\rho_{\frac{1}{\kappa}}\starh [\ \rho_{\frac{1}{\kappa}}\starh ((\sqrt{\alpha_i} v)\circ\theta_i)]\bigr]\circ\theta_i^{-1}+\sum_{i=K+1}^{L} \alpha_i  v \,.
\end{equation*}
The horizontal convolution by layers is of crucial importance for defining an 
approximate problem whose asymptotic behavior will be compatible with the formal 
energy laws for smooth solutions of the original (unsmoothed) problem
(\ref{euler}), since the regularity 
of the moving domain will appear as a surface integral term. 
In this second part of the paper, the properties 
of these horizontal convolutions will be featured in a more 
extensive way than in the surface tension case of the first part of the paper.

We remind the reader that this type of smoothing satisfies the usual properties
of the standard convolution; in particular, independently of $\kappa$, we have 
the existence of $C>0$ such that for any $v\in H^s(\Omega)$:
\begin{equation*}
\|v^{\kappa}\|_s\le C\ \|v\|_s,\ \text{and}\  \ |v^\kappa|_{s-\ud+p}\le {C} \kappa^{-p} |v|_{s-\ud}\ \text{for}\ p\ge 0.
\end{equation*}
We will denote for any $l\in \{1,...,K\}$, the following transformed functions from $v$ and $\eta$ that will naturally arise at the variational level:
\begin{definition}
\begin{equation*}
\d v_{l\kappa}=\rho_{\frac{1}{\kappa}}\starh (\sqrt{\alpha_l} v\circ\theta_l) \ \text{in}\ (0,1)^3,
\end{equation*}
\begin{align*}
\eta^\kappa&=\text{Id}+\int_0^t v^\kappa\ \text{in}\ \Omega,\\
\eta_{l\kappa}&=\theta_l+\int_0^t v_{l\kappa}\ \text{in}\ (0,1)^3.
\end{align*}
\end{definition}

\begin{remark}
The regularity of the moving free surface will be provided by  control of 
each $\eta_{l\kappa}$ in a suitable norm independently of the parameter 
$\kappa$.
\end{remark}

\section{The smoothed $\kappa$-problem and its linear fixed-point formulation}
\label{L3}

As it turns out, the smoothed problem associated to the zero surface tension 
Euler equations can be found quite simply and naturally, and involves only 
transport-type arguments in an Eulerian framework. Also, the construction of a 
solution is
easier if we assume more regularity on the domain and initial velocity than in 
Theorem \ref{ltheorem2}. We shall therefore assume until Section \ref{L12} that
$\Omega$ is of class $H^{\frac{9}{2}}$ and $u_0\in H^{\frac{9}{2}}(\Omega)$. 
In Section \ref{L12}, we will show how this restriction can be  removed.

Letting $u=v\circ{\eta^\kappa}^{-1}$, we consider the  following 
sequence of approximate problems in which the transport velocity $u^\kappa$
is smoothed:
\begin{subequations}
\label{smoothl}
\begin{align}
 u_t+\nabla_{u^\kappa}u+\nabla p&=0\ \text{in}\ \eta^\kappa(t,\Omega),\\
\operatorname{div}u&=0\ \text{in}\ \eta^\kappa(t,\Omega),\\
p&=0\  \text{on}\ \eta^\kappa(t,\Gamma),\\
u(0)&=u_0\ \text{in}\ \Omega.
\end{align}
\end{subequations}

In order to solve this smoothed problem, we will use a linear problem whose 
fixed point will provide the desired solution. If we denote by $\bar v$ an 
arbitrary element of $C_T$ defined in Section \ref{L4}, and $\ek$ the 
corresponding Lagrangian flow defined above, then we
search for $w$ such that if 
$u=w\circ{(\bar\eta^\kappa)^{-1}}$ and $\uk=\bar v^\kappa\circ(\ek) ^{-1}$, we have that
\begin{subequations}
\label{smoothlinearl}
\begin{align}
 u_t+\nabla_{{\bar u}^\kappa}u +\nabla p&=0\ \text{in}\ \ek(t,\Omega),\\
\operatorname{div}u&=0\ \text{in}\ \ek(t,\Omega),\\
p&=0 \  \text{on}\ \ek(t,\Gamma),\\
u(0)&=u_0\ \text{in}\ \Omega.
\end{align}
\end{subequations}
A fixed point $w=\bar v$ to this problem then provides a solution to
(\ref{smoothl}). In the following sections, $\bar v\in C_T$ is assumed given, 
and $\kappa$ is in $(0,\kappa_0)$. $\kappa$ is fixed until Section \ref{L6} 
where we study the asymptotic behavior of the problem (\ref{smoothl}) as 
$\kappa\rightarrow 0$.
\begin{remark}
Note that, for this problem, we do not add any parabolic artificial viscosity, in order to keep the transport-type structure of the Euler equations and to preserve the condition $p=0$ on the free boundary.
\end{remark}

\section{Existence of a solution to (\ref{smoothl})}
\label{L4}

\subsection{A Closed convex set}
\begin{definition}
For $T>0$, we define the following closed convex set of the Hilbert space $L^2(0,T;H^{\se}(\Omega))$:  
\begin{align*}
C_T=\{ v \in L^2(0,T;H^{\se}(\Omega))| \  \sup_{[0,T]} \|v\|_{\se}\le 2 \|u_0\|_{\se}+1 \},
\end{align*}
\end{definition}
 It is clear that $C_T$ is  non-empty (since it contains the constant in time function $u_0$), and is
a convex, bounded and closed set of the separable Hilbert space $L^2(0,T;H^{\se}(\Omega))$.  

By choosing $T(\|\nabla u_0\|_{{\sd}}+1)\le C_\Omega \epsilon_0$,
condition (\ref{eta}) holds for
$\eta = \text{Id}+\int_0^t v$ and any $v \in C_T$ and thus  (\ref{aequation}) 
is well-defined.

We then see that, by taking $T$ smaller if necessary, we have the existence of 
$\kappa_1>0$ such that for any $\kappa\in (0,\kappa_1)$, we have the 
injectivity of $\eta^\kappa(t)$ on $\Omega$ for any $t\in [0,T]$, and 
$\nabla\eta^\kappa$ satisfies condition (\ref{eta}). We then 
denote $a^\kappa=[\nabla\eta^\kappa]^{-1}$, and we let $n^\kappa(\eta^\kappa(x))$ denote
the exterior unit normal to $\eta^\kappa(\Omega)$ at $\eta^\kappa(x)$, with $x\in\Gamma$. 
We now set $\kappa_2=\min(\kappa_0,\kappa_1)$, and assume in the following 
that $\kappa\in (0,\kappa_2)$.

\subsection{Existence and uniqueness for the smoothed problems (\ref{smoothlinearl}) and (\ref{smoothl})}
Suppose that $\bar v\in C_T$ is given. Now, for $v\in C_T$ given, 
we define $p$ on $\ek(t, \Omega)$ by
\begin{subequations}
\label{ex1}
\begin{align}
\triangle p&=-\uk_i,_j u_j,_i\ \text{in}\ \ek(t, \Omega),\\
p&=0\ \text{on}\ \ek(t, \Gamma),
\end{align}
\end{subequations}
where $u=v\circ(\ek)^{-1}$. We next define $\tilde v$ in $\Omega$ by
\begin{equation}
\tilde v(t)=u_0+\int_0^t [\nabla p] (t',
\ek(t',\cdot) ) dt',
\end{equation}
and we now explain why the mapping $v\rightarrow \tilde v$ has a fixed point in $C_T$ for $T>0$ small enough. 
For each $t\in [0,T]$, let $\Psi(t)$ denote the solution of
$\Delta \Psi(t) =0$ in $\Omega$ with $\Psi(t) = \bar \eta^\kappa(t)$ on 
$\Gamma$.  For $\kappa$ and $T$ taken sufficiently small $|\bar\eta^\kappa - \text{Id}|_4
<<1 $ so that $\Psi(t)$ is an embedding, and satisfies
\begin{equation} \label{ex1a}
\| \Psi(t)\|_{4.5} \le C |\bar \eta^\kappa(t)|_4\,.
\end{equation}
Letting $Q(x,t) = p(\Psi(x,t), t)$ and $A(x,t) = [ \nabla \Psi(x,t)]^{-1}$, 
(\ref{ex1}) can be written as
\begin{align*}
[A^k_i A^j_i q,_k],_j &= -[\uk_i,_j u_j,_i](\Psi(x,t),t) \ \ \text{ in } \ \ 
\Omega \,, \\
Q &=0 \ \ \text{ on } \ \ \Gamma\,.
\end{align*}
By elliptic regularity (with Sobolev
class regularity on the coefficients \cite{Eb2002}) almost everywhere in 
$(0,T)$ and using (\ref{ex1a}), 
\begin{align}
\|p\|_{\frac{9}{2},\ek(t,\cdot)}\le C P(|\ek|_4) \|\bar v\|_\sd\|v\|_\sd P(\|\ek\|_\sd),
\label{ex2}
\end{align}
where $P$ denotes a generic polynomial.
Now, with the definition of $T$ and $C_T$, along with the properties of the convolution that allow us to state that $$|\ek|_4\le \frac{1}{\kappa} |\ek|_3,$$ (since the derivatives involved are along the boundary, allowing our convolution by layers to smooth in these directions), this provides the following estimate:
\begin{align*}
\|p\|_{\frac{9}{2},\ek(t,\cdot)}\le C_\kappa P(\|\ek\|_\sd)\|v\|_\sd\le C_\kappa\|v\|_\sd,
\end{align*}
where we have used the definition of $C_T$, 
and where $C_\kappa$ denotes a generic constant depending on $\kappa$.  
Consequently, we get in $[0,T]$,
\begin{align}
\|\tilde v(t)\|_\sd &\le \|u_0\|_\sd+\int_0^t \|p\|_{\frac{9}{2},\ek(t',\cdot)}\|\ek\|_{\frac{7}{2}}\n\\
&\le \|u_0\|_\sd+C_\kappa \int_0^t \|v\|_\sd\|\ek\|_{\frac{7}{2}}.
\label{ex3}
\end{align} 
With the definition of $C_T$, this yields:
\begin{equation*}
\sup_{[0,T]}\|\tilde v(t)\|_\sd \le \|u_0\|_\sd +C_\kappa T.
\end{equation*}
Now, for $T_\kappa\in(0,T)$ such that $T_\kappa C_\kappa\le 1$, we see that $\tilde v\in C_{T_\kappa}$, which ensures that the closed convex set $C_{T_\kappa}$ is stable under the mapping $v\rightarrow \tilde v$. We could also show that this mapping is also sequentially weakly continuous in $L^2(0,{T_\kappa}; H^\se(\Omega))$. Therefore, by the Tychonoff fixed point theorem, there exists a fixed point $v=\tilde v$ in $C_{T_\kappa}$. Now, to see the uniqueness of this fixed point, we see that if another fixed point $\check v$ existed, we would have by the linearity of the mapping $v\rightarrow p$ and the estimates (\ref{ex2}) and (\ref{ex3}) an inequality of the type:
\begin{equation*}
\|(v-\check v)(t)\|_\se \le C \int_0^t \|v-\check v\|_\se,
\end{equation*}
which establishes the uniqueness of the fixed point.  By construction, if we denote $u=v\circ(\ek)^{-1}$, this fixed point satisfies the equation on $(0,{T_\kappa})$:
\begin{equation*}
 u_t+\uk_i  u,_i+\nabla p=0\ \text{in}\ \ek(t,\Omega).
\end{equation*}
 Besides from the definition of $p$ in (\ref{ex1}), we have 
\begin{equation*}
 \div u_t+\uk_i  \div u,_i=0\ \text{in}\ \ek(t,\Omega),
\end{equation*}
{\it i.e.}
\begin{equation*}
 \div u (t,\ek(t,x))=\div u_0(x)=0\ \text{in}\ \Omega.
\end{equation*}
This precisely shows that $u=v\circ(\ek)^{-1}$ is the unique solution of the linear system (\ref{smoothlinearl}) on $(0,{T_\kappa})$.

Now, we see that we again have a mapping $\bar v\rightarrow v$ from $C_{T_\kappa}$ into itself, which is also sequentially weakly lower semi-continuous. It therefore also has a fixed point $v_\kappa$ in $C_{T_\kappa}$, which is a solution of (\ref{smoothl}). 

In the following we study the limit as $\kappa\rightarrow 0$ of the time of existence ${T_\kappa}$ and of $v_\kappa$. We will also denote for the sake of conciseness $v_\kappa$, $u_\kappa=v_\kappa\circ{\eta^{\kappa}}^{-1}$  and $(u_\kappa)^\kappa$ respectively by $\tilde v$, $\tilde u$ and $\tilde u^\kappa$.

\section{Conventions about constants, the time of existence $T_\kappa$, and the dimension of the space}
\label{L5}
From now on, until Section \ref{L12}, we shall stay in $\R^2$ for the sake of notational convenience. In Section \ref{L12}, we shall explain the differences for
the three dimensional case.
In the remainder of the paper, we will denote any constant depending on $\|u_0\|_{\frac{9}{2}}$ as $N(u_0)$. So, for instance, with $q_0$ solution of
\begin{subequations}
\begin{align}
\triangle q_0&=-u_0,_j^iu_0,_i^j\ \text{in}\ \Omega,\\
q_0&=0\ \text{on}\ \Gamma,
\end{align}
\end{subequations}
we have by elliptic regularity $\|q_0\|_{\frac{9}{2}}\le N(u_0)$ (since $\Omega$ is assumed in $H^{\frac{9}{2}}$ until Section \ref{L12}). 

We will also denote generic constants by the letter $C$. Moreover, we will denote
$$\|\Omega\|_s=\sum_{i=1}^K \|\theta_i\|_{s,\O}.$$
Furthermore, the time $T_\kappa>0$ will be chosen small enough so that on $[0,T_\kappa]$, we have for our solution $\v$ given by Section \ref{L4}:
\begin{subequations}
\label{assume}
\begin{align}
\ud\le \text{det} \nabla\tek &\le \td\ \text{in}\ \Omega,
\label{assume.a}\\
\|\e\|_3&\le |\Omega|+1,
\label{assume.b}\\
\|\q\|_3&\le \|q_0\|_3+1
\label{assume.c}\\
\|\v\|_\cd&\le \|u_0\|_\cd+1,
\label{assume.d}
\end{align}
\end{subequations} 
The right-hand sides appearing in the three last inequalities shall be denoted by a generic
constant $C$ in the estimates that we will perform. In what follows, we will prove that this can be achieved on a time independent of $\kappa$.

\section{A continuous in time space energy appropriate for the asymptotic process}
\label{L6}
\begin{definition}
We choose  $0\le \xi\in\mathcal{C}^\infty(\overline{\Omega})$ such that $\text{Supp}\xi\subset\cap_{i=K+1}^{L}[\text{Supp}\alpha_i]^c$ and $\xi=1$ in a neighborhood of $\Gamma$. We then pick $0\le \beta\in\mathcal{D}(\Omega)$ such that $\beta=1$ on $[\text{Supp}\xi]^c$. We then define on $[0,T_\kappa]$:
\begin{align}
\tilde E(t)=&\sup_{[0,t]}\bigl[\sum_{l=1}^K\|\sqrt{\alpha_l}(\theta_l)\ttek\|_{\se,(0,1)^2}+\|\tek\|_\se+\|\beta \tilde\eta\|_\se+\|\v\|_3+\|\q\|_\se \bigr]\nonumber\\
&+\sup_{[0,t]}\sum_{l=1}^K\bigl[\kappa\|\sqrt{\alpha_l}\v\circ\theta_l\|_{\se,(0,1)^2}+\kappa^\td\|\sqrt{\alpha_l}\v\circ\theta_l\|_{4,(0,1)^2}\bigr]+1\,. \label{Elin}
\end{align}
\end{definition}
\begin{remark}
Note the presence of $\kappa$-dependent coefficients in $\tilde E(t)$, that 
indeed arise as a necessity for our asymptotic study. The corresponding terms, 
without the $\kappa$, would of course not be asymptotically controlled.
\end{remark}
\begin{remark}
The 1 is added to the norm to ensure that $\tilde E\ge 1$, which will sometimes be convenient, whereas not necessary.
\end{remark}

Now, since from Section \ref{L4}, $\v\in \mathcal{C} ^0 (0,T_\kappa;H^\se(\Omega))$ (in a way not controlled asymptotically, which does not matter for our purpose), we have $\tilde\eta\in\mathcal{C}^0([0,T_\kappa];H^\se(\Omega))$. Next with the definition (\ref{ex1}), and the definition of our fixed point $\v$, we have for $\tilde p=\q\circ(\tek)^{-1}$:
\begin{align*}
\triangle \p&=-\tuk_i,_j \u_j,_i\ \text{in}\ \tek(t, \Omega),\\
\p&=0\ \text{on}\ \tek(t, \Gamma),
\end{align*}
which shows that $\q\in\mathcal{C}^0([0,T_\kappa];H^\se(\Omega))$. Consequently, $\tilde E$ is a continuous function on $[0,T_\kappa]$.

 We will then prove that this continuous in time energy is controlled by the same type of
polynomial law as (41) of \cite{CoSh2005b}, which will provide a control independent of $\kappa$ on a time independent of $\kappa$.

\section{A commutation-type lemma.}
\label{L7}

 We will need the following Lemma in order to later on identify exact in time energy laws from terms arising from our convolution by horizontal layers:
\begin{lemma}
\label{convolutionbis}
Let $\delta_0>0$ be given. Independently of $\kappa\in (0,\delta_0)$, there exists $C>0$ such that for any $g\in H^\ud((0,1)^2)$ and $f\in H^\cd((0,1)^2)$ such that 
$$\delta_0<\min(\hbox{dist}({supp}\ fg, \{1\}\times [0,1]), \text{dist}(\text{supp}\ fg, \{0\}\times [0,1])),$$
 we have,
\begin{align*}
\bigl\|\rho_{\frac{1}{\kappa}}\starh[f g]-
f  \rho_{\frac{1}{\kappa}}\starh g\bigr\|_{\ud,(0,1)^2}\le C\ \|\kappa g\|_{\ud,(0,1)^2}\|f\|_{\cd,(0,1)^2}+ C\kappa^\ud\ \| g\|_{0,(0,1)^2}\|f\|_{\cd,(0,1)^2}.
\end{align*}
\end{lemma}
\begin{proof}
Let 
$\Delta=\rho_{\frac{1}{\kappa}}\starh[f  g]-
f\ \rho_{\frac{1}{\kappa}}\starh g$. Then, we have:
\begin{align*}
\Delta(x)=\int_{x_1-\kappa}^{x_1+\kappa}\ \rho_{\frac{1}{\kappa}}({x_1-y_1})[f(y_1,x_2)-f(x_1,x_2)]\ g(y_1,x_2)\ dy_1,
\end{align*}
this integral being well-defined because of our condition on the support of $fg$.
We then have, since $H^\td$ is embedded in $L^\infty$ in $2d$,
\begin{align*}
|\Delta(x)|\le C \kappa \|f\|_{\cd,(0,1)^2} \int_{x_1-\kappa}^{x_1+\kappa}\ \rho_{\frac{1}{\kappa}}({x_1-y_1})\ |g(y_1,x_2)|\ dy_1,
\end{align*}
showing that
\begin{align}
\|\Delta\|_{0,\O}&\le C \kappa \|f\|_{\cd,(0,1)^2} \| \rho_{\frac{1}{\kappa}}\starh |g|\|_{0,\O}\n\\
&\le C \kappa \|f\|_{\cd,(0,1)^2} \| g\|_{0,\O}.
\label{si1}
\end{align}
Now, let $p\in\{1,2\}$. Then,
$$\Delta,_p=\rho_{\frac{1}{\kappa}}\starh[f  g,_p]-
f\ \rho_{\frac{1}{\kappa}}\starh g,_p+\rho_{\frac{1}{\kappa}}\starh[f,_p  g]-
f,_p\ \rho_{\frac{1}{\kappa}}\starh g.$$
The difference between the two first terms of the right-hand side of this identity can be treated in a similar fashion as (\ref{si1}), leading us to:
\begin{align}
\|\Delta,_p\|_{0,\O}&\le C \kappa \|f\|_{\cd,(0,1)^2} \|g\|_{1,\O}+\|\rho_{\frac{1}{\kappa}}\starh[f,_p  g]\|_{0,\O}+\|f,_p\ \rho_{\frac{1}{\kappa}}\starh g\|_{0,\O}\n\\
&\le C \kappa \|f\|_{\cd,(0,1)^2} \|g\|_{1,\O}+\|f,_p  g\|_{0,\O}+\|f,_p\|_{L^\infty(\O)} \| \rho_{\frac{1}{\kappa}}\starh g\|_{0,\O}\n\\
&\le C \kappa \|f\|_{\cd,(0,1)^2} \|g\|_{1,\O}+2 \|f,_p\|_{L^\infty(\O)} \| g\|_{0,\O}\n\\
&\le C \kappa \|f\|_{\cd,(0,1)^2} \|g\|_{1,\O}+C \|f\|_{\cd,\O} \| g\|_{0,\O}.
\label{si2}
\end{align}
Consequently, we obtain by interpolation from (\ref{si1}) and (\ref{si2}):
\begin{align*}
\|\Delta\|_{\ud,\O}
&\le C \kappa \|f\|_{\cd,(0,1)^2} \|g\|_{\ud,\O}+C\kappa^\ud \|f\|_{\cd,\O} \| g\|_{0,\O}.
\end{align*}
\end{proof}
We then infer the following result, whose proof follows the same patterns as the previous one:
\begin{lemma}
\label{convolution}
Let $\delta_0>0$ be given. Independently of $\kappa\in (0,\delta_0)$, there exists $C>0$ such that for any $g\in H^s((0,1)^2)$ ($s=\td,\cd$) and for any $f\in H^\cd((0,1)^2)$ such that 
$$\delta_0<\min(\hbox{dist}({supp}\ fg, \{1\}\times [0,1]), \text{dist}(\text{supp}\ fg, \{0\}\times [0,1])),$$
we have
\begin{align*}
\bigl\|\rho_{\frac{1}{\kappa}}\starh[f g]-
f  \rho_{\frac{1}{\kappa}}\starh g\bigr\|_{s,(0,1)^2}\le C\ \|\kappa g\|_{s,(0,1)^2}\|f\|_{\cd,(0,1)^2}+ C\kappa^\ud \| g\|_{s-\ud,(0,1)^2}\|f\|_{\cd,(0,1)^2}.
\end{align*}
\end{lemma}

\section{Asymptotic regularity of the divergence and curl of $\tilde\eta_{l\kappa}$.}
\label{L8}
In this Section, we state the necessary a priori controls that we have on the divergence and curl of various transformations of $\v$ and $\e$. This process has to be justified again, since the functional framework substantially differs from the case with surface tension, with in this case one time derivative on the velocity corresponding to half a  space derivative.

We will base our argument on the fact that the divergence and curl of $\tilde u$ satisfy the following transport type equations:
\begin{subequations}
\begin{align}
D_t\text{div}\tilde u&=0,
\label{si11.a}\\
D_t\text{curl}\tilde u+\tuk_i,_1 \tilde u^2,_i-\tuk_i,_2\tilde u^1,_i&=0.
\label{si11.b}
\end{align}
\end{subequations}
We now study the consequences of these relations on the divergence and curl of $\e$ in the interior of $\Omega$, and of each $\tilde\eta_{l\kappa}, (1\le l\le N)$.

\subsection{Estimate for $\div (\beta\tilde\eta),_s$}
From (\ref{si11.a}), we then infer in $\Omega$ that:
$(\tak)_i^j \v,_j^i=0.
$ Thus, for $s=1, 2$
\begin{equation*}
(\tak)_i^j (\beta\v),_{sj}^i=-\beta(\tak)_i^j,_s \v,_{j}^i+[(\tak)_i^j (\beta\v),_{sj}^i-\beta (\tak)_i^j \v,_{sj}^i],
\end{equation*}
and by integration in time,
\begin{align}
(\tak)_i^j (\beta\tilde\eta),_{sj}^i(t)&=(\tak)_i^j (\beta\tilde\eta),_{sj}^i(0)+\int_0^t [-\beta(\tak)_i^j,_s \v,_{j}^i+((\tak)_i^j (\beta\v),_{sj}^i-\beta (\tak)_i^j \v,_{sj}^i)]\n\\
&\ \ \ +\int_0^t {(\tak)_i^j}_t (\beta\e),_{sj}^i.
\label{diveta1}
\end{align}
Consequently,
\begin{align}
\div (\beta\tilde\eta),_s(t)=&[-(\tak)_i^j+\delta_i^j] (\beta\tilde\eta),_{sj}^i(t)+(\tak)_i^j (\beta\tilde\eta),_{sj}^i(0)+\int_0^t {(\tak)_i^j}_t (\beta\e),_{sj}^i\n\\
&+\int_0^t [-\beta(\tak)_i^j,_s \v,_{j}^i+((\tak)_i^j (\beta\v),_{sj}^i-\beta (\tak)_i^j \v,_{sj}^i)]\n\\
=&[-\int_0^t{(\tak)_i^j}_t] (\beta\tilde\eta),_{sj}^i(t)+(\tak)_i^j (\beta\tilde\eta),_{sj}^i(0)+\int_0^t {(\tak)_i^j}_t (\beta\e),_{sj}^i\n\\
&+\int_0^t [-\beta(\tak)_i^j,_s \v,_{j}^i+((\tak)_i^j (\beta\v),_{sj}^i-\beta (\tak)_i^j \v,_{sj}^i)],
\label{diveta1bis}
\end{align}
showing that
\begin{align}
\|\div (\beta\tilde\eta),_s(t)\|_\td&\le C t  \sup_{ [0,t]}[\ \|\tak_t\|_\td\ \|\beta\tilde\eta\|_\se] +C+Ct\sup_{[0,t]}[\ \|\tak\|_\cd\|\v\|_\cd]\le C t \tilde E(t)+C,
\label{divetabeta}
\end{align}
where we have used our convention stated in Section \ref{L5}.

\subsection{Estimate for $\div[\ttek,_s\circ\theta_l^{-1}\circ(\tek)^{-1}]$}
Since $\text{det}\nabla\theta_l=1$, we then infer in $(0,1)^2$, with $\tbk=[\nabla\tek\circ\theta_l]^{-1}$ that
\begin{equation*}
(\tbk)_i^j (\v\circ\theta_l),_j^i=0.
\end{equation*}
Therefore, as for (\ref{diveta1}), 
\begin{align}
(\tbk)_i^j ((\sqrt{\alpha_l}\tilde\eta)\circ\theta_l),_{sj}^i(t)=&\ (\tbk)_i^j (\sqrt{\alpha_l}\tilde\eta\circ\theta_l),_{sj}^i(0)\n\\
& +\int_0^t {(\tbk)_i^j}_t (\sqrt{\alpha_l}\e\circ\theta_l),_{sj}^i
-\int_0^t \sqrt{\alpha_l}(\theta_l)(\tbk)_i^j,_s (\v\circ\theta_l),_{j}^i\n\\
 &+\int_0^t [(\tbk)_i^j (\sqrt{\alpha_l}\v\circ\theta_l),_{sj}^i-\sqrt{\alpha_l}(\theta_l) (\tbk)_i^j (\v\circ\theta_l),_{sj}^i].
\label{diveta3}
\end{align}
Consequently,
\begin{align*}
\rho_{\frac{1}{\kappa}}\starh[(\tbk)_i^j ((\sqrt{\alpha_l}\tilde\eta)\circ\theta_l),_{sj}^i](t)=
&\rho_{\frac{1}{\kappa}}\starh\bigl[(\tbk)_i^j (\sqrt{\alpha_l}\tilde\eta\circ\theta_l),_{sj}^i(0)\bigr]\n\\
& +\int_0^t \rho_{\frac{1}{\kappa}}\starh\bigl[{(\tbk)_i^j}_t (\sqrt{\alpha_l}\e\circ\theta_l),_{sj}^i- \sqrt{\alpha_l}(\theta_l)(\tbk)_i^j,_s (\v\circ\theta_l),_{j}^i\bigr]\n\\
 &+\int_0^t \rho_{\frac{1}{\kappa}}\starh\bigl[(\tbk)_i^j (\sqrt{\alpha_l}\v\circ\theta_l),_{sj}^i-\sqrt{\alpha_l}(\theta_l) (\tbk)_i^j (\v\circ\theta_l),_{sj}^i\bigr]\n\\
=& \int_0^t \rho_{\frac{1}{\kappa}}\starh\bigl[{(\tbk)_i^j}_t (\sqrt{\alpha_l}\e\circ\theta_l),_{sj}^i\bigl] +R,
\end{align*}
with $\|R\|_\td\le Ct \tilde E(t)+C.$
Next, thanks to Lemma \ref{convolution},
\begin{align}
&\bigl\|\rho_{\frac{1}{\kappa}}\starh\bigl[{(\tbk)_i^j} (\sqrt{\alpha_l}\tilde\eta\circ\theta_l),_{sj}^i(t)\bigr]-{(\tbk)_i^j}\ttek,_{sj}^i(t)\bigr\|_{\td,\O}\n\\
&\le  C  \|(\tbk)_i^j \|_{\cd,(0,1)^2}\ \|\kappa (\sqrt{\alpha_l}\tilde\eta\circ\theta_l),_{sj}^i(0)+\kappa \int_0^t (\sqrt{\alpha_l}\tilde v\circ\theta_l),_{sj}^i\|_{\td,\O}\n\\
&\ \ \ + C  \|(\tbk)_i^j \|_{\cd,(0,1)^2}\kappa^\ud \|(\sqrt{\alpha_l}\tilde\eta\circ\theta_l),_{sj}^i(0)+\int_0^t (\sqrt{\alpha_l}\tilde v\circ\theta_l),_{sj}^i\|_{1,\O}]\n\\
&\le C \kappa^\ud \tilde E(t)+C t \tilde E(t)^2+C.
\label{diveta6}
\end{align}
By successively integrating by parts in time and using Lemma \ref{convolution},
\begin{align}
&\bigl\|\int_0^t\rho_{\frac{1}{\kappa}}\starh\bigl[{(\tbk)_i^j}_t (\sqrt{\alpha_l}\tilde \eta\circ\theta_l),_{sj}^i(t)\bigr]-\int_0^t{(\tbk)_i^j}_t\rho_{\frac{1}{\kappa}}\starh\bigl[ (\sqrt{\alpha_l}\tilde \eta\circ\theta_l),_{sj}^i\bigr]\bigr\|_{\td,\O}\n\\
&\le \bigl\|\int_0^t\rho_{\frac{1}{\kappa}}\starh\bigl[{(\tbk)_i^j} (\sqrt{\alpha_l}\tilde v\circ\theta_l),_{sj}^i(t)\bigr]-\int_0^t{(\tbk)_i^j}\rho_{\frac{1}{\kappa}}\starh\bigl[ (\sqrt{\alpha_l}\tilde v\circ\theta_l),_{sj}^i\bigr]\bigr\|_{\td,\O}\n\\
 &+\bigl\|\bigl[\rho_{\frac{1}{\kappa}}\starh\bigl[{(\tbk)_i^j} (\sqrt{\alpha_l}\tilde \eta\circ\theta_l),_{sj}^i(t)\bigr]-{(\tbk)_i^j}\rho_{\frac{1}{\kappa}}\starh\bigl[ (\sqrt{\alpha_l}\tilde \eta\circ\theta_l),_{sj}^i\bigr]\bigr]_0^t\bigr\|_{\td,\O}\n\\
&\le C \kappa^\ud \tilde E(t)+C t \tilde E(t)^2+C.
\label{diveta5}
\end{align}

Consequently, with (\ref{diveta3}), (\ref{diveta5}) and (\ref{diveta6}), we infer
\begin{align*}
\bigl\|\div[\ttek,_s\circ\theta_l^{-1}\circ(\tek)^{-1}](\tek\circ\theta_l)(t)&-\int_0^t {(\tbk)_i^j}_t\rho_{\frac{1}{\kappa}}\starh\bigl[ (\sqrt{\alpha_l}\tilde \eta\circ\theta_l),_{sj}^i]\bigr\|_{\td,\O}\\
& \le  C \kappa^\ud \tilde E(t)+C t \tilde E(t)^2+C,
\end{align*}
showing that
\begin{align}
\bigl\|\div[\ttek,_s\circ\theta_l^{-1}\circ(\tek)^{-1}]\bigr\|_{\td,\tek(\theta_l(\O))} \le C \kappa^\ud \tilde E(t)+C t \tilde E(t)^2+C.
\label{diveta}
\end{align}

We now study the curl of the same vector fields as in the two previous subsections.

\subsection{Estimate for $\curl (\beta\tilde\eta),_s$}
From (\ref{si11.b}), we obtain:
\begin{equation*}
(\tak)_1^j \v,_j^2-(\tak)_2^j \v,_j^1=\curl\u(0)+\int_0^t [\ -\tvk,_j^i (\tak)_1^j\v,_k^2(\tak)_i^k+\tvk,_j^i (\tak)_2^j\v,_k^1(\tak)_i^k] .
\end{equation*}
 Therefore, for $s=1, 2,$
\begin{align*}
(\tak)_1^j (\beta\v),_{sj}^2-(\tak)_2^j (\beta\v),_{sj}^1=&\beta\curl\u(0),_s
-\beta[(\tak)_1^j,_s (\v),_{j}^2-(\tak)_2^j,_s (\v),_{j}^1]\\
&+
(\tak)_1^j [(\beta\v),_{sj}^2-\beta\v,_{sj}^2]-(\tak)_2^j [(\beta\v),_{sj}^1
-\beta\v,_{sj}^1]\\
&+\int_0^t \beta[\ -\tvk,_j^i (\tak)_1^j\v,_k^2(\tak)_i^k+\tvk,_j^i (\tak)_2^j\v,_k^1(\tak)_i^k],_s,
\end{align*}
which implies by integration in time,
\begin{align*}
(\tak)_1^j (\beta\e),_{sj}^2-(\tak)_2^j (\beta\e),_{sj}^1=&\int_0^t [(\tak_t)_1^j (\beta\e),_{sj}^2-(\tak_t)_2^j (\beta\e),_{sj}^1] +t \beta\curl\u(0),_s\\
&-\beta\int_0^t [(\tak)_1^j,_s (\v),_{j}^2-(\tak)_2^j,_s (\v),_{j}^1]+\int_0^t [f+g]\\
&+\int_0^t(\tak)_1^j [(\beta\v),_{sj}^2-\beta\v,_{sj}^2]-\int_0^t(\tak)_2^j [(\beta\v),_{sj}^1
-\beta\v,_{sj}^1],
\end{align*}
with
$\d f(t')=\int_0^{t'} \beta[\ -\tvk,_j^i (\tak)_1^j\v,_k^2(\tak)_i^k],_s$ and $\d g(t')=\int_0^{t'}\beta[ \tvk,_j^i (\tak)_2^j\v,_k^1(\tak)_i^k],_s$.
Now, since $H^\td$ is a Banach algebra in 2d,
\begin{align}
\bigl\|(\tak)_2^j (\beta\e),_{sj}^1-(\tak)_1^j (\beta\e),_{sj}^2\bigr\|_\td&
\le C \int_0^t \|\tak_t\|_\td \|\beta\e\|_\se  + t\|u_0\|_\se+\int_0^t \|\tak\|_\cd\|\v\|_\cd\n\\
&\ \ \  +\int_0^t\|f+g\|_\td\n\\
&\ \ \ \le  N(u_0)+ Ct \tilde E(t)+\int_0^t\|f+g\|_\td.
\label{curleta1}
\end{align}
We now notice that
\begin{align*}
f(t)&=-\int_0^{t} \beta[\ \tvk,_{sj}^i (\tak)_1^j\v,_k^2(\tak)_i^k+\tvk,_{j}^i (\tak)_1^j\v,_{sk}^2(\tak)_i^k]-\int_0^t \beta \tvk,_{j}^i \v,_k^2 [(\tak)_1^j(\tak)_i^k],_s\\
&=\int_0^{t} \beta[\ \tek,_{sj}^i [(\tak)_1^j\v,_k^2(\tak)_i^k]_t +\e,_{sk}^2[\tvk,_{j}^i (\tak)_1^j(\tak)_i^k]_t]-\int_0^t \beta \tvk,_{j}^i \v,_k^2 [(\tak)_1^j(\tak)_i^k],_s\\
&\ \ \ +\bigl[\beta[\ \tek,_{sj}^i (\tak)_1^j\v,_k^2(\tak)_i^k +\e,_{sk}^2\tvk,_{j}^i (\tak)_1^j(\tak)_i^k]\bigr]_0^t,
\end{align*}
which allows us to infer that
\begin{align*}
\|f(t)\|_\td&\le \int_0^{t} (\|\tek\|_\se+\|\beta\e\|_\se+\|\e\|_\cd) [\|\tak_t\|_\td\|\tak\|_\td\|\v\|_\cd+ \|\tak\|_\td\|\tak\|_\td\|\v_t\|_\cd]\\
&\ \ \  +\int_0^t \|\v\|^2_\cd \|\tak\|_\cd^2+(\|\tek\|_\se+\|\beta\e\|_\se+\|\e\|_\cd)(t) \|\tak(t)\|^2_\td\|\v\|_\cd+  N(u_0)\\
&\ \ \ \le C t \tilde E(t)^2 +N(u_0)+C \tilde E(t).
\end{align*} 
 Since $g(t)$ can be estimated in a similar fashion, (\ref{curleta1}) provides us with
\begin{align*}
\bigl\|(\tak)_2^j (\beta\e),_{sj}^1-(\tak)_1^j (\beta\e),_{sj}^2\bigr\|_\td\le C t \tilde E(t)^2 +N(u_0).
\end{align*}
showing that,
\begin{align}
\bigl\|\curl(\beta\e),_{s}\bigr\|_\td&\le \bigl\|\bigl[\int_0^t (\tak_t)_2^j \bigr](\beta\e),_{sj}^1-\bigl[\int_0^t (\tak_t)_1^j \bigr](\beta\e),_{sj}^2\bigr\|_\td  + C t \tilde E(t)^2 +N(u_0)\n\\
&\le C t \tilde E(t)^2 +N(u_0).
\label{curletabeta}
\end{align}

\subsection{Estimate for $\curl[\ttek,_s\circ\theta_l^{-1}\circ(\tek)^{-1}]$}
In a similar fashion as we obtained (\ref{diveta}), we also have here
\begin{align}
\bigl\|\curl[\ttek,_s\circ\theta_l^{-1}\circ(\tek)^{-1}]\bigr\|_{\td,\tek(\Omega)}\le C t \tilde E(t)^2 + N(u_0).
\label{curleta}
\end{align}

\subsection{Estimate for $\kappa \div[(\sqrt{\alpha_l}\v_t\circ\theta_l),_s\circ\theta_l^{-1}\circ(\tek)^{-1}]$}
By time differentiating (\ref{diveta3}) twice in time, we find:
\begin{align}
(\tbk)_i^j ((\sqrt{\alpha_l}\tilde v_t)\circ\theta_l),_{sj}^i(t)=&\ -{(\tbk)_i^j}_{tt} ((\sqrt{\alpha_l}\tilde \eta)\circ\theta_l),_{sj}^i(t)-2{(\tbk)_i^j}_t ((\sqrt{\alpha_l}\tilde v)\circ\theta_l),_{sj}^i(t)\n\\
& +[{(\tbk)_i^j}_t (\sqrt{\alpha_l}\e\circ\theta_l),_{sj}^i]_t
- [\sqrt{\alpha_l}(\theta_l)(\tbk)_i^j,_s (\v\circ\theta_l),_{j}^i]_t\n\\
 &+ [(\tbk)_i^j (\sqrt{\alpha_l}\v\circ\theta_l),_{sj}^i-\sqrt{\alpha_l}(\theta_l) (\tbk)_i^j (\v\circ\theta_l),_{sj}^i]_t.
\label{divvt1}
\end{align}
Therefore,
\begin{align}
\kappa\|(\tbk)_i^j ((\sqrt{\alpha_l}\tilde v_t)\circ\theta_l),_{sj}^i(t)\|_{\td,\O}\le&\ C\kappa \tilde E(t)^2\n\\
&\  +\|{(\tbk)_i^j}_t\|_{\td,\O}\|\kappa ((\sqrt{\alpha_l}\tilde v)\circ\theta_l),_{sj}^i(t)\|_{\td,\O}\n\\
&\ +\|\kappa[(\tbk)_i^j,_s]_t\|_{\td,\O} \|(\sqrt{\alpha_l}\v\circ\theta_l),_{j}^i\|_{\td,\O}\n\\
&\le\ C \tilde E(t)^2+\|\kappa[(\tbk)_i^j,_s]_t\|_{\td,\O} \|(\sqrt{\alpha_l}\v\circ\theta_l),_{j}^i\|_{\td,\O}.
\label{diveta2bis}
\end{align}
Next, we for instance have:
\begin{align*}
\kappa{(\tbk,_s)_1^1}_t&=\kappa \frac{\tilde v_2^\kappa\circ\theta_l,_{2s}}{\text{det}\nabla(\tek\circ\theta_l)}-\kappa [\text{det}(\tek\circ\theta_l)]_t,_{s}\frac{\tilde \eta_2^\kappa\circ\theta_l,_{2}}{\text{det}^2\nabla(\tek\circ\theta_l)}\\
&\ \ \  -\kappa [\text{det}(\tek\circ\theta_l)]_t\frac{\tilde \eta_2^\kappa\circ\theta_l,_{s2}}{\text{det}^2\nabla(\tek\circ\theta_l)}-\kappa [\text{det}(\tek\circ\theta_l)],_s\frac{\tilde v_2^\kappa\circ\theta_l,_{2}}{\text{det}^2\nabla(\tek\circ\theta_l)},
\end{align*}
which shows that, 
\begin{align*}
\|\kappa{(\tbk,_s)}_t\|_{\td,\O}&\le C \tilde E(t),
\end{align*}
and with (\ref{diveta2bis}) this implies:
\begin{align*}
\kappa\|(\tbk)_i^j ((\sqrt{\alpha_l}\tilde v_t)\circ\theta_l),_{sj}^i(t)\|_{\td,\O}\le&\ C \tilde E(t)^2.
\end{align*}
and thus, still by writing $\d \tek(t)=\tek(0)+\int_0^t \tilde v^\kappa$, we finally have
\begin{align}
\kappa \bigl\|\div[(\sqrt{\alpha_l}\v_t\circ\theta_l),_s\circ\theta_l^{-1}\circ(\tek)^{-1}]\bigr\|_{\td,\tek(\Omega)}\le C \tilde E(t)^2.
\label{divkappavt}
\end{align}

With the same type of arguments, we also have the following asymptotic estimates:

\subsection{Estimate for $\kappa \curl[(\sqrt{\alpha_l}\v_t\circ\theta_l),_s\circ\theta_l^{-1}\circ(\tek)^{-1}]$}
\begin{align}
\kappa \bigl\|\curl[(\sqrt{\alpha_l}\v_t\circ\theta_l),_s\circ\theta_l^{-1}\circ(\tek)^{-1}]\bigr\|_{\td,\tek(\Omega)}\le C \tilde E(t)^2.\label{curlkappavt}
\end{align}

\subsection{Estimate for $\kappa^2 \div[(\sqrt{\alpha_l}\v_t\circ\theta_l),_s\circ\theta_l^{-1}\circ(\tek)^{-1}]$}
\begin{align}
\kappa^2 \bigl\|\div[(\sqrt{\alpha_l} \v_t\circ\theta_l),_s\circ\theta_l^{-1}\circ(\tek)^{-1}]\bigr\|_{\cd,\tek(\Omega)}\le C \tilde E(t)^2.
\label{divkappa2vt}
\end{align}

\subsection{Estimate for $\kappa^2 \curl[(\sqrt{\alpha_l}\v_t\circ\theta_l),_s\circ\theta_l^{-1}\circ(\tek)^{-1}]$}
\begin{align}
\kappa^2 \bigl\|\curl[ (\sqrt{\alpha_l} \v_t\circ\theta_l),_s\circ\theta_l^{-1}\circ(\tek)^{-1}]\bigr\|_{\cd,\tek(\Omega)}\le C \tilde E(t)^2.
\label{curlkappa2vt}
\end{align}

\begin{remark} Since we will time integrate the previous quantities, the absence of a small parameter in front of $\tilde E(t)^2$ is not problematic.
\end{remark}

\subsection{Asymptotic control of $\div\tuk (\tek)$} We have
\begin{align*}
\text{div}\tuk (\tek)&=(\tak)_k^j\tvk,_j^k\\
&=(\tak)_k^j[\sqrt{\alpha_i}[\rho_{\frac{1}{\kappa}}\starh (\v_{i\kappa})^k](\theta_i^{-1})],_j\\
&=[(\tbk)_k^j[\sqrt{\alpha_l}(\theta_l)\rho_{\frac{1}{\kappa}}\starh (\v_{l\kappa})],_j^k](\theta_l^{-1}).
\end{align*}
Now, thanks to Lemma \ref{convolution}, this leads us to
\begin{align*}
\text{div}\tuk (\tek)&=[\rho_{\frac{1}{\kappa}}\starh [(\tbk)_k^j\sqrt{\alpha_l}(\theta_l) (\v_{l\kappa}),_j^k]](\theta_l^{-1})+r_1,
\end{align*}
with
\begin{align*}
\|r_1\|_\cd\le C \|\nabla\tek\|_\ci ( \|\kappa\nabla \v\|_\cd +\kappa^\ud\|\nabla\v\|_2)\le C \tilde E(t)^2.
\end{align*}
Next, we notice that 
\begin{align*}
(\tbk)_k^l (E_t^{l\kappa}),_l^k &=(\tbk)_k^j \rho_{\frac{1}{\kappa}}\starh[(\sqrt{\alpha_l}\v)(\theta_l)],_j^k\\
&= \rho_{\frac{1}{\kappa}}\starh[\sqrt{\alpha_l}(\theta_l) (\tbk)_k^j \v(\theta_l),_j^k]+r_2,
\end{align*}
with in virtue of Lemma \ref{convolution},
\begin{align*}
\|r_2\|_\cd\le C \|\nabla\tek\|_\ci ( \|\kappa\nabla \v\|_\cd +\kappa^\ud\|\nabla\v\|_2)\le C \tilde E(t)^2.
\end{align*}
Now, since $(\tbk)_k^l \v(\theta_l),_l^k=0$ in $\O$, this finally provides us with
\begin{equation}
\label{divtuk}
\|\div\tuk(\tek)\|_{\cd}\le C \tilde E(t)^2.
\end{equation}

\section{Asymptotic regularity of $\kappa \v_t$ and of $\kappa^2 \v_t$}
\label{L9}
This Section is devoted to the asymptotic control of $\kappa\v_t$ and $\kappa^2 \v_t$ in spaces smoother than the natural regularity $H^\cd(\Omega)$ for $\v_t$,
 the idea still being that one degree in the power of $\kappa$ allows one more degree of space regularity.

\subsection{Asymptotic control of $\kappa \v_t$ in $H^\se(\Omega)$}
Our starting point will be the fact that since $\q=0$ on $\Gamma$, we have for
any $l\in \{1,...,K\}$ on $(0,1)\times\{0\}$:
\begin{equation*}
{\v_t\circ\theta_l}+\frac{{\q\circ\theta_l},_2}{\text{det}\nabla\tek(\theta_l)}{{\tek\circ\theta_l},_1^\perp}=0,
\end{equation*}
where $x^\perp=(-x_2,x_1)$. Therefore, we have on $(0,1)\times\{0\}$:
\begin{align*}
({\sqrt{\alpha_l}\v_t\circ\theta_l}),_{111}\cdot\tek\circ\theta_l,_1&=-[\frac{\sqrt{\alpha_l}(\theta_l){\q\circ\theta_l},_2}{\text{det}\nabla\tek(\theta_l)}],_{11}{{\tek\circ\theta_l},_{11}^\perp}\cdot\tek\circ\theta_l,_1\\
&\ \ \ -[\frac{\sqrt{\alpha_l}(\theta_l){\q\circ\theta_l},_2}{\text{det}\nabla\tek(\theta_l)}],_{1}{{\tek\circ\theta_l},_{111}^\perp}\cdot\tek\circ\theta_l,_1\\
&\ \ \ - \frac{\sqrt{\alpha_l}(\theta_l){\q\circ\theta_l},_2}{\text{det}\nabla\tek(\theta_l)}{{\tek\circ\theta_l},_{1111}^\perp}\cdot\tek\circ\theta_l,_1,
\end{align*}
showing that
\begin{align}
\kappa\|({\sqrt{\alpha_l}\v_t\circ\theta_l}),_{111}\cdot\tek\circ\theta_l,_1\|_{0,\partial\O}&\le \kappa \tilde E(t)^2 +C\kappa\|\sqrt{\alpha_l}(\theta_l){{\tek\circ\theta_l},_{1111}^\perp}\|_{0,\partial\O}.
\label{kappa0}
\end{align}
By definition,
\begin{align*}
\sqrt{\alpha_l}(\theta_l){\tek\circ\theta_l},_{1111}=\sqrt{\alpha_l}(\theta_l)
\sum_{i=1}^K[\sqrt{\alpha_i}(\theta_l) [\rho_{\frac{1}{\kappa}}\starh E^{i\kappa}](\theta_i^{-1}\circ\theta_l)],_{1111},
\end{align*}
the sum being restricted to the indices $i$ such that $\theta_i(\O)$ and $\theta_l(\O)$ have a non empty intersection. We then have
\begin{align}
\sqrt{\alpha_l}(\theta_l){\tek\circ\theta_l},_{1111}=\sqrt{\alpha_l}(\theta_l)
[\sum_{i=1}^K\sqrt{\alpha_i}(\theta_l) [\rho_{\frac{1}{\kappa}}\starh E^{i\kappa}],_{i_1i_2i_3i_4}(\theta_i^{-1}\circ\theta_l) a_{il,1111}^{i_1i_2i_3i_4}+\Delta],
\label{kappa1}
\end{align}
with
$a_{il,1111}^{i_1i_2i_3i_4}=(\theta_i^{-1}\circ\theta_l),_1^{i_1}(\theta_i^{-1}\circ\theta_l),_1^{i_2}(\theta_i^{-1}\circ\theta_l),_1^{i_3}(\theta_i^{-1}\circ\theta_l),_1^{i_4}$, and
\begin{align}
\kappa\|\sqrt{\alpha_l}(\theta_l)\Delta\|_{0,\partial\O}&\le C \kappa \|\Omega\|_{\frac{9}{2}} \sup_{i}\|\rho_{\frac{1}{\kappa}}\starh[\sqrt{\alpha_i}\e\circ\theta_i]\|_{3,\O}\n\\
&\ \ \ +C\kappa \sup_{i}\|\rho_{\frac{1}{\kappa}}\starh[\sqrt{\alpha_i}\e\circ\theta_i]\|_{\se,\O}\n\\
& \le C \kappa (\|\Omega\|_{\frac{9}{2}}+1) \tilde E(t)\n\\
&\le C \kappa \tilde E(t)^2.
\label{kappa2}
\end{align}
Now, we notice that for $x_1\in (0,1)$ such that $\theta_l(x_1,0)\in\theta_i(\O)$, we necessarily have, since for all $k\in \{1,...,K\}$, $\theta_k([0,1]\times\{0\})=\partial\Omega\cap\theta_k([0,1]^2)$, that $\theta_i^{-1}\circ\theta_l(x_1,0)=(f_{il}(x_1),0)$, showing that on $(0,1)\times\{0\}$, we have $\sqrt{\alpha_l}(\theta_l) a_{il,1111}^{i_1i_2i_3i_4}=0$ except when $i_1=i_2=i_3=i_4=1$. Therefore, (\ref{kappa1}) can be expressed as
\begin{align}
\sqrt{\alpha_l}(\theta_l){\tek\circ\theta_l},_{1111}=\sqrt{\alpha_l}(\theta_l)
[\sum_{i=1}^K\sqrt{\alpha_i}(\theta_l) [\rho_{\frac{1}{\kappa}}\starh E^{i\kappa}],_{1111}(\theta_i^{-1}\circ\theta_l) a_{il,1111}^{1111}+\Delta].
\label{kappa3}
\end{align}
Now, from the properties of our convolution by layers, we have (since the derivatives are horizontal) that
\begin{align}
\kappa\|[\rho_{\frac{1}{\kappa}}\starh E^{i\kappa}],_{1111}\|_{0,(0,1)\times\{0\}}\le C \|E^{i\kappa},_{111}\|_{0,(0,1)\times\{0\}}.
\label{kappa4}
\end{align}
Thus, with (\ref{kappa1}), (\ref{kappa2}) and (\ref{kappa4}), we infer
\begin{align*}
\kappa\|\sqrt{\alpha_l}(\theta_l)\tek\circ\theta_l,_{1111}\|_{0,\partial\O} \le C \kappa  \tilde E(t)+C \tilde E(t),
\end{align*}
which coupled with (\ref{kappa0}) provides us with
\begin{align*}
\kappa\|({\sqrt{\alpha_l}\v_t\circ\theta_l}),_{111}\cdot\tek\circ\theta_l,_1\|_{0,\partial\O}&\le C \tilde E(t)^2 + C.
\end{align*}
This provides us the trace estimate:
\begin{align}
\kappa\|({\sqrt{\alpha_l}\v_t\circ\theta_l}),_{1}(\theta_l^{-1}\circ(\tek)^{-1})\cdot&\tek\circ\theta_l,_1(\theta_l^{-1}\circ(\tek)^{-1})\|_{2,\partial\tek(\theta_l(\O))}\n\\
&\le C \tilde E(t)^2 + C.
\label{kappa5}
\end{align}
Consequently with the divergence and curl estimates (\ref{divkappavt}) and (\ref{curlkappavt}) and the trace estimate (\ref{kappa5}), we infer by elliptic regularity:
\begin{align}
&\kappa\|({\sqrt{\alpha_l}\v_t\circ\theta_l}),_{1}(\theta_l^{-1}\circ(\tek)^{-1})\|_{\cd,\tek(\theta_l(\O))}\le C \tilde E(t)^2 + C.
\label{kappavt}
\end{align}

\begin{remark}
It is the presence of $\|\Omega\|_{\frac{9}{2}}$ in the inequalities leading
to (\ref{kappa2}) which
explains the assumption of $\Omega$ in $H^{\frac{9}{2}}$. It is, however, not 
essential as will be shown in Section \ref{L12}. One way to see this, is to 
smooth the initial domain by a convolution with the parameter $\kappa$ to form
$\Omega^\kappa$. Then, 
by the properties of the convolution,
 $\kappa \|\Omega^\kappa\|_{\frac{9}{2}}\le C \|\Omega\|_\se$.
\end{remark}

\subsection{Asymptotic control of $\kappa^2 \v_t$ in $H^{\frac{9}{2}}(\Omega)$}

In a similar way as in the previous subsection, we would obtain the trace estimate:
\begin{align*}
\kappa^2\|({\sqrt{\alpha_l}\v_t\circ\theta_l}),_{1}(\theta_l^{-1}\circ(\tek)^{-1})\cdot&\tek\circ\theta_l,_1(\theta_l^{-1}\circ(\tek)^{-1})\|_{3,\partial\tek(\theta_l(\O))}\n\\
&\le C \tilde E(t)^2 + C,
\end{align*}
which coupled with (\ref{divkappa2vt}) and (\ref{curlkappa2vt}) provides
\begin{align}
&\kappa^2\|({\sqrt{\alpha_l}\v_t\circ\theta_l}),_{1}(\theta_l^{-1}\circ(\tek)^{-1})\|_{\se,\tek(\theta_l(\O))}\le C \tilde E(t)^2 + C.
\label{kappa2vt}
\end{align}

\section{Basic energy law for the control of $\v$ and $\ttek$ independently of $\kappa$.} 
\label{L10}
We will use a different type of energy than in \cite{ChLi2000}, namely:
\begin{definition}
$$\d H^\kappa(t)=\ud\sum_{l=1}^K \int_{(0,1)^2}\xi_l\circ\theta_l|(\tilde v\circ\theta_l),_{111}|^2,$$
 where $\xi_l=\xi\ \alpha_l$, $\xi$ being defined in Section \ref{L6}.
\end{definition}
\begin{remark}
 The main differences with respect to the energy of \cite{ChLi2000} are in the absence in our energy of any restriction to the tangent components, allowing a more convenient set of estimates, and in a setting in Lagrangian variables.
\end{remark}
We have:
\begin{align*}
H^\kappa_t(t)&=\sum_{l=1}^K \int_{(0,1)^2} \xi_l\circ\theta_l \tilde v_t\circ\theta_l,_{111}\tilde v\circ\theta_l,_{111}\n\\
&= -\sum_{l=1}^K \int_{(0,1)^2} \xi_l\circ\theta_l ((\tak)_j^k\q,_k)\circ\theta_l,_{111}{\tilde v}_j\circ\theta_l,_{111}\n\\
&= -\sum_{l=1}^K \int_{(0,1)^2} \xi_l\circ\theta_l [(\tbk)_j^p\q\circ\theta_l,_p],_{
111}{\tilde v}_j\circ\theta_l,_{111},
\end{align*}
where $\tbk=[\nabla(\tek\circ\theta_l)]^{-1}$ . Next, we see that $H^\kappa_t=-[H_1+H_2+H_3]$, with
\begin{align*}
H_1(t)&=\sum_{l=1}^K \int_{(0,1)^2} \xi_l(\theta_l) (\tbk)_j^p,_{111}\q\circ\theta_l,_p{\tilde v}_j\circ\theta_l,_{111},\\
H_2(t)&=\sum_{l=1}^K \int_{(0,1)^2} \xi_l\circ\theta_l [(\tbk)_j^p\q\circ\theta_l,_{p111}]{\tilde v}_j\circ\theta_l,_{111},\\
H_3(t)&=\sum_{l=1}^K \int_{(0,1)^2} \xi_l\circ\theta_l [[(\tbk)_j^p\q\circ\theta_l,_p],_{
111}-(\tbk)_j^p\q\circ\theta_l,_{l111}-(\tbk)_j^p,_{111}\q\circ\theta_l,_l]{\tilde v}_j\circ\theta_l,_{111}.
\end{align*}
We immediately have for the third term:
\begin{equation}
\label{ecl1}
|H_3(t)|\le C \tilde E(t)^2.
\end{equation}
Next, for $H_2$, since $(\xi_l \q)\circ\theta_l=0$ on $\partial (0,1)^2$, 
\begin{align*}
H_2(t)&=-\sum_{l=1}^K \int_{(0,1)^2} \xi_l\circ\theta_l\ (\tbk)_j^p\q\circ\theta_l,_{111}{\tilde v}_j\circ\theta_l,_{p111}\\
&\ \ \ -\sum_{l=1}^K \int_{(0,1)^2} [\xi_l\circ\theta_l (\tbk)_j^p],_p \q\circ\theta_l,_{111}{\tilde v}_j\circ\theta_l,_{111}.
\end{align*}
We then notice that from the divergence condition, we have $(\tbk)_j^p\ \v_j\circ\theta_l,_p=0$ in $(0,1)^2$, implying
\begin{align*}
H_2(t)&=\sum_{l=1}^K \int_{(0,1)^2} \xi_l\circ\theta_l\ \q\circ\theta_l,_{111}(\tbk)_j^p,_{111} {\tilde v}_j\circ\theta_l,_{p}\\
&\ \ \ +\sum_{l=1}^K \int_{(0,1)^2} \xi_l\circ\theta_l\ \q\circ\theta_l,_{111}(\tbk)_j^p,_{11} {\tilde v}_j\circ\theta_l,_{p1}\\
&\ \ \ +\sum_{l=1}^K \int_{(0,1)^2} \xi_l\circ\theta_l\ \q\circ\theta_l,_{111}(\tbk)_j^p,_{1} {\tilde v}_j\circ\theta_l,_{p11}\\
&\ \ \ -\sum_{l=1}^K \int_{(0,1)^2} [\xi_l\circ\theta_l (\tbk)_j^p],_p \q\circ\theta_l,_{111}{\tilde v}_j\circ\theta_l,_{111}.
\end{align*}
Now, in a way similar as (\ref{int}), we have for any $f\in H^\ud((0,1)^2)$
\begin{equation*}
\|\xi_l\circ\theta_l f,_1\|_{H^\ud ((0,1)^2)'}\le C \|f\|_{H^\ud ((0,1)^2)},
\end{equation*}
since the derivative is in the horizontal direction.
By applying this result to $f={(\tbk)_j^p},_{11}$ for the first integral appearing in the equality above, and by using the continuous embedding of $H^1$ into $L^6$ ($6\in (1,\infty)$), and of $H^\ud$ into $L^3$ ($ 3\in (1,4)$) for the other integrals, we then obtain:
\begin{align}
|H_2(t)|&\le C \| {\tak}\|_\ci \|\q\|_\se \|{\tilde v}\|_3\le C \tilde E(t)^3.
\label{ecl2}
\end{align}
We now come to $H_1$, which will require more care, and will provide us with the
regularity of $\ttek(\Omega)$ in $H^\se$ independently of $\kappa$. We have:
\begin{align}
H_1(t)&= \int_{(0,1)^2} \xi_l\circ\theta_l (\tbk)_j^p,_{111} [\q\circ\theta_l],_p\ {\tilde v}_j\circ\theta_l,_{111}\n\\
&=H_{11}(t)+H_{12}(t) -\int_{(0,1)^2} \xi_l\circ\theta_l \Delta_j^p \ [\q\circ\theta_l],_p\ {\tilde v}_j\circ\theta_l,_{111},
\label{ecl3}
\end{align}
with
\begin{align*}
H_{11}(t)&=\int_{(0,1)^2} \xi_l\circ\theta_l \frac{(\text{Cof}\nabla(\tek\circ\theta_l))_j^p,_{111}}{\text{det}\nabla(\tek\circ\theta_l)}[\q\circ\theta_l],_p\ {\tilde v}_j\circ\theta_l,_{111},\n\\
H_{12}(t)&=-\int_{(0,1)^2} \xi_l\circ\theta_l (\text{Cof}\nabla(\tek\circ\theta_l))_j^p\frac{[\text{det}\nabla(\tek\circ\theta_l)],_{111}}{[\text{det}\nabla(\tek\circ\theta_l)]^2}[\q\circ\theta_l],_p\ {\tilde v}_j\circ\theta_l,_{111},
\end{align*}
and
\begin{align*}
\Delta_j^p=\bigl[\frac{(\text{Cof}\nabla(\tek\circ\theta_l))_j^p}{\text{det}\nabla(\tek\circ\theta_l)}\bigr],_{111}-\frac{(\text{Cof}\nabla(\tek\circ\theta_l))_j^p,_{111}}{\text{det}\nabla(\tek\circ\theta_l)}+(\text{Cof}\nabla(\tek\circ\theta_l))_j^p\frac{[\text{det}\nabla(\tek\circ\theta_l)],_{111}}{[\text{det}\nabla(\tek\circ\theta_l)]^2},
\end{align*}
so that
\begin{align}
\bigl|\int_{(0,1)^2}  \Delta_j^p\ [(\xi_l\q)\circ\theta_l],_p\ {\tilde v}_j\circ\theta_l,_{111}\bigr|\le C \|\v\|_3.
\label{ecl4}
\end{align}
We now turn our attention to the other terms of (\ref{ecl3}), and to shorten notations, we will set:  $\tilde Q_l=\q\circ\theta_l$.
 We first study the perturbation $H_{12}$, which would not appear if the volume preserving condition was respected by our
smoothing by convolution. It turns out that we do need the double convolution by layers appearing in the definition of $v^\kappa$ in order to identify time derivatives of space energies.  
We first notice that since $\theta_l$ does not depend on $t$, we have:
\begin{equation*}
(\tek\circ\theta_l)_t=\tuk\circ(\tek\circ\theta_l),
\end{equation*}
from which we infer in $(0,1)^2$, since $\theta_l$ is volume preserving,
\begin{equation}
[\text{det}(\nabla\tek\circ\theta_l)]_t=\text{div}\tuk (\tek\circ\theta_l)\ \text{det}\nabla(\tek\circ\theta_l).
\label{ec5}
\end{equation}
\subsection{Study of $H_{12}$.}
We have after an integration by parts in time, and the use of (\ref{ec5}):
\begin{align*}
 \int_0^t H_{12} =\sum_{i=1}^3 H_{12}^i+R_{12},
\end{align*}
with
\begin{align*}
H_{12}^1&=\int_0^t\int_{(0,1)^2} \xi_l\circ\theta_l (\text{Cof}\nabla(\tek\circ\theta_l))_j^p\frac{\text{div}\tuk (\tek\circ\theta_l)\ [\text{det}\nabla(\tek\circ\theta_l)],_{111}}{[\text{det}\nabla(\tek\circ\theta_l)]^2}\tilde Q_l,_p\ {\tilde \eta}_j\circ\theta_l,_{111},\\
H_{12}^2&=\int_0^t\int_{(0,1)^2} \xi_l\circ\theta_l (\text{Cof}\nabla(\tek\circ\theta_l))_j^p\frac{[\text{div}\tuk (\tek\circ\theta_l)],_{111}\ \text{det}\nabla(\tek\circ\theta_l)}{[\text{det}\nabla(\tek\circ\theta_l)]^2}\tilde Q_l,_p\ {\tilde \eta}_j\circ\theta_l,_{111},\\
H_{12}^3&=-\int_{(0,1)^2} \xi_l\circ\theta_l (\text{Cof}\nabla\tilde (\tek\circ\theta_l))_j^p\frac{ [\text{det}\nabla(\tek\circ\theta_l)],_{111}}{[\text{det}\nabla(\tek\circ\theta_l)]^2}\tilde Q_l,_p\ {\tilde \eta}_j\circ\theta_l,_{111}(t),
\end{align*}
and 
\begin{align}
|R_{12}(t)|\le C t \tilde E(t)^3 +C.
\label{ecl6}
\end{align}
\subsubsection{\bf Study of $H_{12}^1$.}
 For the sake of conciseness, we denote $$\d A_{jl}=\xi_l\circ\theta_l(\text{Cof}\nabla(\tek\circ\theta_l))_j^p \frac{\text{div}\tuk (\tek\circ\theta_l)}{[\text{det}\nabla(\tek\circ\theta_l)]^2}[\q\circ\theta_l],_p.$$ 
We then see, by expanding the third space derivative of the determinant in the integrand of $H_{12}^1$, that $H_{12}^1=\sum_{i=1}^4 H_{12}^{1i}+R_{12}^1$ with the $H_{12}^{1i}$ being estimated as $H_{12}^{11}$ that we precise below and $R_{12}^1$ being a remainder estimated as (\ref{ecl6}). 
By definition of $\tek$ we have,
if we denote $$E^{i\kappa}=\rho_{\frac{1}{\kappa}}\starh  ((\sqrt{\alpha_i} \tilde\eta)\circ\theta_i),$$ 
\begin{align*}
H_{12}^{11}&=\int_0^t\int_{(0,1)^2} \bigl[ A_{jl} (\tek_2\circ\theta_l),_2\ \bigl[\bigl[\sqrt{\alpha_i}(\theta_i)\ \rho_{\frac{1}{\kappa}}\starh E_1^{i\kappa}\bigr](\theta_i^{-1}\circ\theta_l)\bigr],_{1111}\\
&\hskip 3cm [{\tilde\eta_j\circ\theta_i(\theta_i^{-1}\circ\theta_l)}],_{111}\bigr]+R_{12}^{1},
\end{align*}
with $|R_{12}^1|\le Ct$ and 
where, because of the term $\sqrt{\alpha_i}(\theta_l)$, the only indexes $i$ and $l$ appearing in this sum are the ones for which
$\theta_l((0,1)^2)\cap\theta_i((0,1)^2)\ne\emptyset$. Only such indexes will be considered later on when such terms arise. From our assumed regularity on $\Omega$ in $H^{\frac{9}{2}}$, we then have
\begin{align*}
H_{12}^{11}&=\int_0^t\int_{(0,1)^2} \bigl[ A_{jl} (\tek_2\circ\theta_l),_2\ \bigl[\bigl[\rho_{\frac{1}{\kappa}}\starh E_1^{i\kappa}\bigr](\theta_i^{-1}\circ\theta_l)\bigr],_{1111}\\
&\hskip 3cm [\sqrt{\alpha_i}{\tilde\eta_j\circ\theta_i(\theta_i^{-1}\circ\theta_l)}],_{111}\bigr]+R_{12}^{11},
\end{align*}
with $|R_{12}^{11}|\le C t \tilde E(t)^2.$
 We next have, since the charts $\theta_i$ are volume preserving,
\begin{align*}
H_{12}^{11}&=\int_0^t\int_{\theta_i^{-1}(\theta_l(0,1)^2)} \bigl[ [A_{jl}(\tek_2\circ\theta_l),_2](\theta_l^{-1}\circ\theta_i)\ [\rho_{\frac{1}{\kappa}}\starh E_1^{i\kappa}],_{j_1j_2j_3j_4}\\
&\hskip 4cm c^{j_1j_2j_3j_4}_{il,1111}c^{i_1i_2i_3}_{il,111} [\sqrt{\alpha_i}\tilde\eta_j\circ\theta_i],_{i_1i_2i_3}\bigl]+R,
\end{align*}
with $|R|\le C t \tilde E(t)^2$ and 
\begin{subequations}
\label{cij}
\begin{align}
c_{il,111}^{j_1j_2j_3j_4}&=[(\theta_i^{-1}\circ\theta_l)^{j_1},_1(\theta_i^{-1}\circ\theta_l)^{j_2},_1(\theta_i^{-1}\circ\theta_l)^{j_3},_1(\theta_i^{-1}\circ\theta_l)^{j_4},_1](\theta_l^{-1}\circ\theta_i),\\
c_{il,111}^{i_1i_2i_3}&=[(\theta_i^{-1}\circ\theta_l)^{i_1},_1(\theta_i^{-1}\circ\theta_l)^{i_2},_1(\theta_i^{-1}\circ\theta_l)^{i_3},_1](\theta_l^{-1}\circ\theta_i).
\end{align}
\end{subequations}
Next, we notice that the term $A_{jl}(\theta_l^{-1}\circ\theta_i)$ introduces a factor $\alpha_l\circ\theta_i$ which is non zero only if $x\in \theta_i^{-1}(\theta_l(0,1)^2)$, leading us to
\begin{align*}
H_{12}^{11}&=\int_0^t\int_{(0,1)^2} \bigl[ [A_{jl}(\tek_2\circ\theta_l),_2](\theta_l^{-1}\circ\theta_i)\ \rho_{\frac{1}{\kappa}}\starh E_1^{i\kappa},_{j_1j_2j_3j_4}\\
&\hskip 4cm c^{j_1j_2j_3j_4}_{il,1111}c^{i_1i_2i_3}_{il,111} [\sqrt{\alpha_i}\tilde\eta_j\circ\theta_i],_{i_1i_2i_3}\bigl]+R,
\end{align*}
where $\theta_l^{-1}\circ\theta_i$ is extended outside of $\theta_i^{-1}(\theta_l(0,1)^2)$ in any fashion. This argument of replacing an integral on
a subset of $(0,1)^2$ by an integral on $(0,1)^2$ will be implicitly repeated at other places later on. Now, since $\rho$ is even,
\begin{align}
H_{12}^{11}&=\int_0^t\int_{(0,1)^2} E_1^{i\kappa},_{j_1j_2j_3j_4} \rho_{\frac{1}{\kappa}}\star\bigl[ [A_{jl}(\tek_2\circ\theta_l),_2](\theta_l^{-1}\circ\theta_i)\ \n\\
&\hskip 5cm c^{j_1j_2j_3j_4}_{il,1111}c^{i_1i_2i_3}_{il,111} [\sqrt{\alpha_i}\tilde\eta_j\circ\theta_i],_{i_1i_2i_3}\bigl]+R.
\label{ecl16bis}
\end{align}
Now, let us call $f=[A_{jl}(\tek_2\circ\theta_l),_2](\theta_l^{-1}\circ\theta_i) c^{j_1j_2j_3j_4}_{il,1111}c^{i_1i_2i_3}_{il,111}$ and $g=[\sqrt{\alpha_j}\tilde\eta_1\circ\theta_i],_{i_1i_2i_3}$. We notice that $\|f\|_{\cd,\O}$ is the natural norm associated to $\tilde E(t)$. Here we cannot use directly Lemma \ref{convolution} for the case where all the $j_i=3$, since $E^{i\kappa}_{3333}$ is not necessarily in $H^\ud((0,1)^2)'$ a priori. Instead, we write
\begin{align*}
f(y_1,x_2)=f(x_1,x_2)+(y_1-x_1). f,_1(x_1,x_2)+\int_{x_1}^{y_1}[f,_1(x,x_2)-f,_1(x_1,x_2)]dx,
\end{align*}
which shows that on $\O$:
\begin{align*}
\rho_{\frac{1}{\kappa}}\starh[fg] (x_1,x_2)=&f(x_1,x_2) \rho_{\frac{1}{\kappa}}\starh g (x_1,x_2) \\
&+f,_1(x_1,x_2)\int_\R \rho_{\frac{1}{\kappa}}(y_1-x_1) (y_1-x_1)g(y_1,x_2) dy_1\\
&+\int_\R \rho_{\frac{1}{\kappa}}(y_1-x_1) \int_{x_1}^{y_1}[f,_1(x,x_2)-f,_1(x_1,x_2)]dx\ g(y_1,x_2) dy_1.
\end{align*}
This implies 
\begin{align}
H_{12}^{11}&=\int_0^t\int_{(0,1)^2} \bigl[ [A_{jl}(\tek_2\circ\theta_l),_2](\theta_l^{-1}\circ\theta_i)\  E_1^{i\kappa},_{j_1j_2j_3j_4}\n\\
&\hskip 3cm c^{j_1j_2j_3j_4}_{il,1111}c^{i_1i_2i_3}_{il,111} \rho_{\frac{1}{\kappa}}\starh [\sqrt{\alpha_i}(\theta_i)[\tilde\eta_1\circ\theta_i],_{i_1i_2i_3}]\bigl] +R-R_1-R_2,
\label{ecl16ter}
\end{align}
with
\begin{align*}
R_1=&\int_0^t\int_{(0,1)^2} \bigl[E_1^{i\kappa},_{j_1j_2j_3j_4}(x_1,x_2) f,_1(x_1,x_2)\n\\
&\hskip 3cm \int_\R \rho_{\frac{1}{\kappa}}(y_1-x_1) (y_1-x_1)g(y_1,x_2) dy_1\bigr] \ dx_1dx_2,\n\\
R_2=& \int_0^t\int_{(0,1)^2} \bigl[E_1^{i\kappa},_{j_1j_2j_3j_4}(x_1,x_2)\n\\
&\hskip 2cm \int_\R \rho_{\frac{1}{\kappa}}(y_1-x_1) \int_{x_1}^{y_1}[f,_1(x,x_2)-f,_1(x_1,x_2)]dx\ g(y_1,x_2) dy_1\bigr]\ dx_1dx_2.
\end{align*}
Now, for $R_2$, we notice that since $$|f,_1(x,x_2)-f,_1(x_1,x_2)|\le C |f|_{\cd,\O} |x-x_1|^\ud\le C |\tek|_{\se} |x-x_1|^\ud,$$ we have
\begin{align}
R_2\le & C\int_0^t\|\tek\|_\se \int_{(0,1)^2} \bigl[|E_1^{i\kappa},_{j_1j_2j_3j_4}(x_1,x_2)|\int_\R \rho_{\frac{1}{\kappa}}(y_1-x_1) \kappa^\td |g(y_1,x_2)| dy_1\bigr] dx_1dx_2\n\\
\le & C\int_0^t\|\tek\|_\se \int_{(0,1)^2} \bigl[|\kappa^\td E_1^{i\kappa},_{j_1j_2j_3j_4}|\ \rho_{\frac{1}{\kappa}}\starh  |g|\bigr]\n\\
\le & C\int_0^t\|\tek\|_\se \kappa^\td \bigl[\|E_1^{i\kappa},_{j_1j_2j_3j_4}(0)\|_{0,\O}+\int_0^t \|(E_1^{i\kappa})_t,_{j_1j_2j_3j_4} \|_{0,\O}\bigr]  \|g\|_{0,\O}\n\\
\le & Ct \kappa^\td  \tilde E(t) +Ct \tilde E(t)^3,
\label{ecl16qua}
\end{align}
where we have used the fact that $\|\kappa^\td \sqrt{\alpha_i}\tilde v_1\circ\theta_i],_{i_1i_2i_3i_4}\|_{0,\O}$ is contained in the definition of $\tilde E(t)$.
We now turn our attention to $R_1$. We first remark that 
\begin{align*}
f,_1=& ([A_{jl}(\tek_2\circ\theta_l),_2](\theta_l^{-1}\circ\theta_i) c^{i_1i_2i_3}_{il,111}),_1 c^{j_1j_2j_3j_4}_{il,1111}\\
&+\sum_{n=1}^4\bigl[ [A_{jl}(\tek_2\circ\theta_l),_2](\theta_l^{-1}\circ\theta_i) c^{i_1i_2i_3}_{il,111}[(\theta_i^{-1}\circ\theta_l),_1^{j_n}(\theta_l^{-1}\circ\theta_i)],_1\\
&\hskip 5cm \Pi_{p\ne n}(\theta_i^{-1}\circ\theta_l),_1^{j_p}(\theta_l^{-1}\circ\theta_i)\bigr],
\end{align*}
which implies that
\begin{align*}
R_1=R_1^1+\sum_{n=1}^4 R_1^{i_n},
\end{align*}
with
\begin{align*}
R_1^1&=\int_0^t\int_{(0,1)^2}  c^{j_1j_2j_3j_4}_{il,1111} \bigl[E_1^{i\kappa},_{j_1j_2j_3j_4}(x_1,x_2) ([A_{jl}(\tek_2\circ\theta_l),_2](\theta_l^{-1}\circ\theta_i) c^{i_1i_2i_3}_{il,111}),_1\n\\
&\hskip 3cm \int_\R \rho_{\frac{1}{\kappa}}(y_1-x_1) (y_1-x_1)g(y_1,x_2) dy_1\bigr] \ dx_1dx_2,\n\\
R_1^{i_n}&= \int_0^t\int_{(0,1)^2} \bigl[E_1^{i\kappa},_{j_1j_2j_3j_4}(x_1,x_2) A_{jl}(\tek_2\circ\theta_l),_2](\theta_l^{-1}\circ\theta_i) \n\\
&\hskip 3cm [(\theta_i^{-1}\circ\theta_l),_1^{j_n}(\theta_l^{-1}\circ\theta_i)],_1
\Pi_{p\ne n}(\theta_i^{-1}\circ\theta_l),_1^{j_p}(\theta_l^{-1}\circ\theta_i)\n\\
&\hskip 3cm c^{i_1i_2i_3}_{il,111}\int_\R \rho_{\frac{1}{\kappa}}(y_1-x_1) (y_1-x_1)g(y_1,x_2) dy_1\bigr] \ dx_1dx_2.
\end{align*}
Let us study $R_1^1$. If we denote $h(x_1,x_2)=\int_\R \rho_{\frac{1}{\kappa}}(y_1-x_1) (y_1-x_1)g(y_1,x_2) dy_1$, since $(E_1^{i\kappa},_{j_2j_3j_4}\circ\theta_i^{-1}\circ\theta_l),_1(\theta_l^{-1}\circ\theta_i)=(\theta_i^{-1}\circ\theta_l),_1^{j_1}(\theta_l^{-1}\circ\theta_i) E_1^{i\kappa},_{j_1j_2j_3j_4}$, 
\begin{align*}
R_1^1&=\int_0^t\int_{(0,1)^2}  c^{j_2j_3j_4}_{il,111} \bigl[(E_1^{i\kappa},_{j_2j_3j_4}\circ\theta_i^{-1}\circ\theta_l),_1(\theta_l^{-1}\circ\theta_i) \n\\
&\hskip 3cm ([A_{jl}(\tek_2\circ\theta_l),_2](\theta_l^{-1}\circ\theta_i) c^{i_1i_2i_3}_{il,111}),_1 h\bigr]\n\\
&=\int_0^t\int_{(0,1)^2}  \bigl[(E_1^{i\kappa},_{j_2j_3j_4}\circ\theta_i^{-1}\circ\theta_l),_1 \n\\
&\hskip 3cm [c^{j_2j_3j_4}_{il,111}([A_{jl}(\tek_2\circ\theta_l),_2](\theta_l^{-1}\circ\theta_i) c^{i_1i_2i_3}_{il,111}),_1 h](\theta_i^{-1}\circ\theta_l)\bigr].
\end{align*}
Since the derivative of $(E_1^{i\kappa},_{j_2j_3j_4}\circ\theta_i^{-1}\circ\theta_l)$ is in the horizontal direction, we infer similarly as in (\ref{int}) that
\begin{align*}
R_1^1
&\le \int_0^t\int_{(0,1)^2}  \bigl[\bigl\|E_1^{i\kappa},_{j_2j_3j_4}\circ\theta_i^{-1}\circ\theta_l\bigr\|_{\ud,\O} \n\\
&\hskip 3cm \bigl\|[c^{j_2j_3j_4}_{il,111}([A_{jl}(\tek_2\circ\theta_l),_2](\theta_l^{-1}\circ\theta_i) c^{i_1i_2i_3}_{il,111}),_1 h](\theta_i^{-1}\circ\theta_l)\bigr\|_{\ud,\O}\bigr].
\end{align*}
Since we have by interpolation $\|h\|_{\ud,\O}\le C \kappa \|g\|_{\ud,\O}$, we then infer:
\begin{align}
|R_1^1|\le C t \tilde E(t)^2.
\label{ecl16qui}
\end{align}
In a similar fashion, for $R_1^{i_1}$ we can identify an horizontal derivative: 
\begin{align*}
&E_1^{i\kappa},_{j_1j_2j_3j_4}[(\theta_i^{-1}\circ\theta_l),_1^{j_{1}}(\theta_l^{-1}\circ\theta_i)],_1 \Pi_{p=2}^4(\theta_i^{-1}\circ\theta_l),_1^{j_{p}}(\theta_l^{-1}\circ\theta_i)\\
&=[(\theta_i^{-1}\circ\theta_l),_1^{j_{1}}(\theta_l^{-1}\circ\theta_i)],_1 \bigl[\Pi_{p=2}^3(\theta_i^{-1}\circ\theta_l),_1^{j_{p}} (E_1^{i\kappa},_{j_1j_2j_3}\circ\theta_i^{-1}\circ\theta_l),_1\bigr](\theta_l^{-1}\circ\theta_i),
\end{align*}
which leads for the same reasons as for $R_1^1$ to $|R_1^{i_1}|\le C t \tilde E(t)^2$. Since the other $R_1^{i_n}$ are similar in structure, we have
\begin{align}
|R_1^{i_n}|\le C t \tilde E(t)^2.
\label{ecl16hex}
\end{align}
Consequently from (\ref{ecl16ter}), (\ref{ecl16qua}), (\ref{ecl16qui}) and
(\ref{ecl16hex}), we infer 
\begin{align*}
H_{12}^{11}&=\int_0^t\int_{(0,1)^2} \bigl[ [A_{jl}(\tek_2\circ\theta_l),_2](\theta_l^{-1}\circ\theta_i)\  E_1^{i\kappa},_{j_1j_2j_3j_4}\\
&\hskip 3cm c^{j_1j_2j_3j_4}_{il,1111}c^{i_1i_2i_3}_{il,111} \rho_{\frac{1}{\kappa}}\starh [\sqrt{\alpha_i}(\theta_i)[\tilde\eta_1\circ\theta_i],_{i_1i_2i_3}]\bigl] +r_{12}^{11},
\end{align*}
with
\begin{align}
|r_{12}^{11}(t)|\le C t \tilde E(t)^3+C t\kappa^\td  \tilde E(t).
\label{ecl16sept}
\end{align}
Since $E_1^{i\kappa},_{j_1j_2j_3j_4} c^{j_1j_2j_3j_4}_{il,1111}=c^{j_2j_3j_4}_{il,111} (E_1^{i\kappa},_{j_2j_3j_4}\circ\theta_i^{-1}\circ\theta_l),_1(\theta_l^{-1}\circ\theta_i)$, we infer as for $R_1^1$ that
\begin{align*}
|H_{12}^{11}(t)|&\le C t \sup_{j,l}\sup_{[0,t]}\|A_{jl}\|_{\cd,(0,1)^2}\|\ttek\|^2_\se+|r_{12}^{11}|\le  C t \tilde E(t)^3+C t\kappa^\td \tilde E(t).
\end{align*}
The other $H_{12}^{1i}$ are estimated in the same fashion, leading us to
\begin{align}
|H_{12}^1(t)|\le  C t \tilde E(t)^3+C t\kappa^\td  \tilde E(t).
\label{ecl7}
\end{align}

\subsubsection{\bf Study of $H_{12}^2$.}
Next, for $H_{12}^2$, we first notice from the asymptotic regularity result (\ref{divtuk}) on $\div\tuk(\tek)$, that $H_{12}^2$ can be treated in the same fashion as $H_{12}^1$, leading to
\begin{align}
|H_{12}^2(t)|\le  C t \tilde E(t)^3+C t\kappa^\td \tilde E(t).
\label{ecl8}
\end{align}

\subsubsection{\bf Study of $H_{12}^3$.}
We simply write
\begin{align*}
-H_{12}^3&=\int_{(0,1)^2}  \xi_l(\theta_l)(\text{Cof}\nabla(\tek\circ\theta_l))_j^p\frac{ [\text{det}\nabla\theta_l],_{111}}{[\text{det}\nabla(\tek\circ\theta_l)]^2}\tilde Q_l,_p\ {\tilde \eta}_j\circ\theta_l,_{111}(t)+R_{12}^3\\
&\ \ \  +\int_{(0,1)^2} \xi_l(\theta_l) (\text{Cof}\nabla(\tek\circ\theta_l))_j^p\frac{\int_0^t [\text{det}\nabla(\tek\circ\theta_l)],_{111}\text{div}(\tuk\circ\theta_l)}{[\text{det}\nabla(\tek\circ\theta_l)]^2}\tilde Q_l,_p\ {\tilde \eta}_j\circ\theta_l,_{111}(t)\\
&\ \ \  +\int_{(0,1)^2} \xi_l(\theta_l) (\text{Cof}\nabla(\tek\circ\theta_l))_j^p\frac{\int_0^t \text{det}\nabla(\tek\circ\theta_l)[\text{div}(\tuk\circ\theta_l)],_{111}}{[\text{det}\nabla(\tek\circ\theta_l)]^2}\tilde Q_l,_p\ {\tilde \eta}_j\circ\theta_l,_{111}(t),
\end{align*}
with $R_{12}^3$ being bounded by a term similar as the right-hand side of (\ref{ecl8}). We also see that the first term of this equality can be estimated by a bound similar as the right-hand side of (\ref{ecl8}). The third term is treated in a way similar as $H_{12}^1$, in order to put a convolution in front of $(\tilde\eta_j\circ\theta_l),_{111}$. There is no difference linked to the fact that the integral from $0$ to $t$ does not apply on all terms as for $H_{12}^1$ since $\rho_{\frac{1}{\kappa}}$ and the $\theta_l$ do not depend on time. The fourth term follows the same treatment as $H_{12}^2$, leading us to
\begin{align*}
|H_{12}^3(t)|\le C t \tilde E(t)^3+C t\kappa^\td  \tilde E(t),
\end{align*}
which with (\ref{ecl7}) and (\ref{ecl8}) implies
\begin{align}
|H_{12}(t)|\le C t \tilde E(t)^3+C t\kappa^\td  \tilde E(t).
\label{ecl9}
\end{align}

\subsection{Study of $H_{11}$.}
As for $H_{11}$, we have if we still denote $E^{i\kappa}=\rho_{\frac{1}{\kappa}}\starh  ((\sqrt{\alpha_i} \tilde\eta)\circ\theta_i)$ and $\epsilon^{mn}$ the sign of the permutation between $(m,n)$ and $(1,2)$:
\begin{align*}
H_{11}(t)&=\epsilon^{mn}\epsilon^{rs} \int_{(0,1)^2}  \xi_l(\theta_l) \frac{(\tek_m\circ\theta_l),_{r111}}{\text{det}\nabla(\tek\circ\theta_l)}[\q\circ\theta_l],_s\ {\v_n(\theta_l)},_{111}\\
&=\epsilon^{mn}\epsilon^{rs} \int_{(0,1)^2} \bigl[\xi_l(\theta_l) \frac{[\q\circ\theta_i(\theta_i^{-1}\circ\theta_l)],_s}{\text{det}\nabla(\tek\circ\theta_l)}\ \bigl[[\sqrt{\alpha_i}(\theta_i)\rho_{\frac{1}{\kappa}}\starh E_m^{i\kappa}](\theta_i^{-1}\circ\theta_l)\bigr],_{r111}\\
&\hskip 3cm [{\v_n\circ\theta_i(\theta_i^{-1}\circ\theta_l)}],_{111}\bigr]\\
&=\epsilon^{mn}\epsilon^{rs}\int_{(0,1)^2} \bigl[ \xi_l(\theta_i)\frac{[\q\circ\theta_i],_{i_1}}{\text{det}\nabla(\tek\circ\theta_l) (\theta_l^{-1}\circ\theta_i)}\ \bigl[\rho_{\frac{1}{\kappa}}\starh E_m^{i\kappa}\bigr],_{j_1j_2j_3j_4}c^{j_1j_2j_3j_4}_{il,r111}\\
&\hskip 3cm c^{i_1i_2i_3i_4}_{il,s111} [{\sqrt{\alpha_i}\v_n\circ\theta_i],_{i_2i_3i_4}}\bigl]+R_{11},
\end{align*}
with
\begin{align*}
|R_{11}(t)|\le C \tilde E(t)^2,
\end{align*}
and 
\begin{align*}
c^{j_1j_2j_3j_4}_{il,r111}&=\bigl[(\theta_i^{-1}\circ\theta_l),_r^{j_1}(\theta_i^{-1}\circ\theta_l),_1^{j_2}(\theta_i^{-1}\circ\theta_l),_1^{j_3}(\theta_i^{-1}\circ\theta_l),_1^{j_4}\bigr](\theta_l^{-1}\circ\theta_i).
\end{align*}
Therefore,
\begin{align*}
H_{11}&=\epsilon^{mn}\epsilon^{rs}\int_{(0,1)^2} \xi_l(\theta_i)[ \q(\theta_i)],_{i_1} \rho_{\frac{1}{\kappa}}\starh E_m^{i\kappa},_{j_1j_2j_3j_4} [\sqrt{\alpha_i}\v_n\circ\theta_i],_{i_2i_3i_4} h_{rs}^{(ji)_{1234} } +R_{11},
\end{align*}
with 
\begin{equation}
\label{hrs}
\d h_{rs}^{(ji)_{1234}}=\bigl[\frac{c^{j_1j_2j_3j_4}_{il,r111} c^{i_1i_2i_3i_4}_{il,s111}} {\text{det}\nabla(\tek\circ\theta_l)(\theta_l^{-1}\circ\theta_i)}\bigr].
\end{equation}
Similarly as in the study of $H_{12}^{11}$ (from equations (\ref{ecl16bis}) to (\ref{ecl16sept})) we have:
\begin{align*}
H_{11}&=\epsilon^{mn}\epsilon^{rs}\int_{(0,1)^2} \xi_l(\theta_i)[\q(\theta_i)],_{i_1} h_{rs}^{(ji)_{1234}} E_m^{i\kappa},_{j_1j_2j_3j_4} \rho_{\frac{1}{\kappa}}\starh [\sqrt{\alpha_i}\v_n\circ\theta_i],_{i_2i_3i_4}\\
&\ \ \  +S_{11} +R_{11},
\end{align*}
with, 
\begin{align*}
|S_{11}|\le C t \tilde E(t)^3+C t\kappa^\td  \tilde E(t).
\end{align*}
By integrating by parts in space (and using $\xi_l\q(\theta_i)=0$ on $\partial (0,1)^2$),
\begin{align}
H_{11}&=-\epsilon^{mn}\epsilon^{rs}\int_{(0,1)^2} (\xi_l \q)(\theta_i)h_{rs}^{(ji)_{1234}}  \rho_{\frac{1}{\kappa}}\starh[\sqrt{\alpha_i}\tilde\eta_m\circ\theta_i],_{j_1j_2j_3j_4} \rho_{\frac{1}{\kappa}}\starh[\sqrt{\alpha_i}\v_n\circ\theta_i],_{i_1i_2i_3i_4}\n\\
&\ \ \  +H_{11}^1+H_{11}^2+S_{11} +R_{11},
\label{ecl10}
\end{align}
with
\begin{subequations}
\begin{align}
H_{11}^1&=-\epsilon^{mn}\epsilon^{rs}\int_{(0,1)^2} (\xi_l \q)(\theta_i) h_{rs}^{(ji)_{1234}}  \rho_{\frac{1}{\kappa}}\starh[\sqrt{\alpha_i}\tilde\eta_m\circ\theta_i],_{i_1j_1j_2j_3j_4} \rho_{\frac{1}{\kappa}}\starh[\sqrt{\alpha_i}\v_n\circ\theta_i],_{i_2i_3i_4},\label{N}\\
H_{11}^2&= -\epsilon^{mn}\epsilon^{rs}\int_{(0,1)^2} \q(\theta_i) [\xi_l(\theta_i)h_{rs}^{(ji)_{1234}}],_{i_1}  \rho_{\frac{1}{\kappa}}\starh[\sqrt{\alpha_i}\tilde\eta_m\circ\theta_i],_{j_1j_2j_3j_4} \rho_{\frac{1}{\kappa}}\starh[\sqrt{\alpha_i}\v_n\circ\theta_i],_{i_2i_3i_4}.
\label{E1}
\end{align}
\end{subequations}
For $H_{11}^1$, by taking into account the symmetric role of $\{i_2,i_3,i_4\}$ and $\{j_2, j_3,j_4\}$, we obtain:
\begin{align*}
H_{11}^1=-\epsilon^{mn}\epsilon^{rs}\int_{(0,1)^2} (\xi_l \q)(\theta_i) h_{rs}^{(ji)_{1234}}  \rho_{\frac{1}{\kappa}}\starh[\sqrt{\alpha_i}\tilde\eta_m\circ\theta_i],_{i_1j_1i_2i_3i_4} \rho_{\frac{1}{\kappa}}\starh[\sqrt{\alpha_i}\v_n\circ\theta_i],_{j_2j_3j_4}.
\end{align*}
Next, since $h_{rs}^{(ji)_{1234}}=h_{sr}^{(ij)_{1234}}$, this implies
\begin{align*}
H_{11}^1&=-\epsilon^{mn}\epsilon^{rs}\int_{(0,1)^2} (\xi_l \q)(\theta_i) h_{sr}^{(ij)_{1234}}  \rho_{\frac{1}{\kappa}}\starh[\sqrt{\alpha_i}\tilde\eta_m\circ\theta_i],_{i_1j_1i_2i_3i_4} \rho_{\frac{1}{\kappa}}\starh[\sqrt{\alpha_i}\v_n\circ\theta_i],_{j_2j_3j_4}\\
&=\epsilon^{mn}\epsilon^{sr}\int_{(0,1)^2} (\xi_l \q)(\theta_i) h_{sr}^{(ij)_{1234}}  \rho_{\frac{1}{\kappa}}\starh[\sqrt{\alpha_i}\tilde\eta_m\circ\theta_i],_{j_1i_1i_2i_3i_4} \rho_{\frac{1}{\kappa}}\starh[\sqrt{\alpha_i}\v_n\circ\theta_i],_{j_2j_3j_4}.
\end{align*}
Therefore, by relabeling $sr$ as $rs$ and $ij$ as $ji$, we obtain by comparison to (\ref{N}):
\begin{align*}
H_{11}^1&=-H_{11}^1,
\end{align*}
and thus $H_{11}^1=0$. For $H_{11}^2$, we have by integrating by parts $H_{11}^2=H_{11}^{21}+H_{11}^{22}$, with
\begin{align*}
H_{11}^{21}&= \epsilon^{mn}\epsilon^{rs}\int_{(0,1)^2} [ \q(\theta_i) [\xi_l(\theta_i)h_{rs}^{(ji)_{1234}}],_{i_1}],_{j_1}  \rho_{\frac{1}{\kappa}}\starh[\sqrt{\alpha_i}\tilde\eta_m\circ\theta_i],_{j_2j_3j_4} \rho_{\frac{1}{\kappa}}\starh[\sqrt{\alpha_i}\v_n\circ\theta_i],_{i_2i_3i_4}\\
H_{11}^{22}&=\epsilon^{mn}\epsilon^{rs}\int_{(0,1)^2}  \q(\theta_i) [\xi_l(\theta_i)h_{rs}^{(ji)_{1234}}],_{i_1}  \rho_{\frac{1}{\kappa}}\starh[\sqrt{\alpha_i}\tilde\eta_m\circ\theta_i],_{j_2j_3j_4} \rho_{\frac{1}{\kappa}}\starh[\sqrt{\alpha_i}\v_n\circ\theta_i],_{j_1i_2i_3i_4}.
\end{align*}
First, for $H_{11}^{21}$ we have if we denote $E_{mn}=\epsilon^{mn}\rho_{\frac{1}{\kappa}}\starh[\sqrt{\alpha_i}\tilde\eta_m\circ\theta_i],_{j_2j_3j_4} \rho_{\frac{1}{\kappa}}\starh[\sqrt{\alpha_i}\tilde\eta_n\circ\theta_i],_{i_2i_3i_4}$:
\begin{align*}
\int_0^t H_{11}^{21}=&-\ud 
\int_0^t \epsilon^{rs}\int_{(0,1)^2} [[ \q(\theta_i) [\xi_l(\theta_i)h_{rs}^{(ji)_{1234}}],_{i_1}],_{j_1}]_t  E_{mn}\\
&+\ud\bigl[\epsilon^{rs}\int_{(0,1)^2} [ \q(\theta_i) [\xi_l(\theta_i)h_{rs}^{(ji)_{1234}}],_{i_1}],_{j_1} E_{mn} \bigr]_0^t.
\end{align*}
Therefore,
\begin{align*}
\bigl|\int_0^t H_{11}^{21}\bigr|\le & C \int_0^t |\epsilon^{rs}| \bigl\|[[ \q(\theta_i) [\xi_l(\theta_i)h_{rs}^{(ji)_{1234}}],_{i_1}],_{j_1}]_t\bigr\|_{\ud,\O}  \|\tilde\eta\|^2_\sd\\
&+C |\epsilon^{rs}|\bigl\|[ \q(\theta_i) [\xi_l(\theta_i)h_{rs}^{(ji)_{1234}}],_{i_1}],_{j_1}(0)\bigr\|_{L^\infty(\O)} 
\bigl[ \|\tilde\eta\|^2_\sd(t) + \|\tilde\eta\|^2_\sd(0) \bigr]
\\
&+C \int_0^{t}\bigl\|[[ \q(\theta_i) [\xi_l(\theta_i)h_{rs}^{(ji)_{1234}}],_{i_1}],_{j_1}]_t\bigr\|_{\ud,\O} 
\bigl[ \|\tilde\eta\|^2_\sd(t) + \|\tilde\eta\|^2_\sd(0) \bigr] \,.
\end{align*}
With the definition (\ref{hrs}) and (\ref{divtuk}) for the control of the time derivative of $\det(\nabla(\tilde\eta^\kappa))$ in $H^\cd(\Omega)$, we then infer:
\begin{align}
\bigl|\int_0^t H_{11}^{21}\bigr|\le & C t \tilde E(t)^4 + N(u_0).
\label{ecl11}
\end{align}
Next, for $H_{11}^{22}$, we have by relabeling $m$ and $n$
\begin{align*}
H_{11}^{22}&=\epsilon^{nm}\epsilon^{rs}\int_{(0,1)^2}  \q(\theta_i) [\xi_l(\theta_i)h_{rs}^{(ji)_{1234}}],_{i_1}  \rho_{\frac{1}{\kappa}}\starh[\sqrt{\alpha_i}\tilde\eta_n\circ\theta_i],_{j_2j_3j_4} \rho_{\frac{1}{\kappa}}\starh[\sqrt{\alpha_i}\v_m\circ\theta_i],_{j_1i_2i_3i_4}\\
&=-\epsilon^{mn}\epsilon^{rs}\int_{(0,1)^2}  \q(\theta_i) [\xi_l(\theta_i)h_{rs}^{(ji)_{1234}}],_{i_1}  \rho_{\frac{1}{\kappa}}\starh[\sqrt{\alpha_i}\tilde\eta_n\circ\theta_i],_{j_2j_3j_4} \rho_{\frac{1}{\kappa}}\starh[\sqrt{\alpha_i}\v_m\circ\theta_i],_{j_1i_2i_3i_4}.
\end{align*}
By taking into account the symmetric role of $\{i_2,i_3,i_4\}$ and $\{j_2, j_3,j_4\}$, we then obtain:
\begin{align}
H_{11}^{22}&=-\epsilon^{mn}\epsilon^{rs}\int_{(0,1)^2}  \q(\theta_i) [\xi_l(\theta_i)h_{rs}^{(ji)_{1234}}],_{i_1}  \rho_{\frac{1}{\kappa}}\starh[\sqrt{\alpha_i}\tilde\eta_n\circ\theta_i],_{i_2i_3i_4} \rho_{\frac{1}{\kappa}}\starh[\sqrt{\alpha_i}\v_m\circ\theta_i],_{j_1j_2j_3j_4}.
\label{E2}
\end{align}
Consequently, by (\ref{E1}) and (\ref{E2}),
\begin{align*}
2 H_{11}^{2}&=-\epsilon^{mn}\epsilon^{rs}\int_{(0,1)^2}  \q(\theta_i) [\xi_l(\theta_i)h_{rs}^{(ji)_{1234}}],_{i_1}  \rho_{\frac{1}{\kappa}}\starh[\sqrt{\alpha_i}\tilde\eta_m\circ\theta_i],_{j_1j_2j_3j_4} \rho_{\frac{1}{\kappa}}\starh[\sqrt{\alpha_i}\v_n\circ\theta_i],_{i_2i_3i_4}\\
&\ \ \ -\epsilon^{mn}\epsilon^{rs}\int_{(0,1)^2}  \q(\theta_i) [\xi_l(\theta_i)h_{rs}^{(ji)_{1234}}],_{i_1}  \rho_{\frac{1}{\kappa}}\starh[\sqrt{\alpha_i}\tilde\eta_n\circ\theta_i],_{i_2i_3i_4} \rho_{\frac{1}{\kappa}}\starh[\sqrt{\alpha_i}\v_m\circ\theta_i],_{j_1j_2j_3j_4}\\
&\ \ \ +H_{11}^{21}\\
&=-\epsilon^{mn}\epsilon^{rs}\int_{(0,1)^2}  \q(\theta_i) [\xi_l(\theta_i)h_{rs}^{(ji)_{1234}}],_{i_1}  [\rho_{\frac{1}{\kappa}}\starh[\sqrt{\alpha_i}\tilde\eta_m\circ\theta_i],_{j_1j_2j_3j_4} \rho_{\frac{1}{\kappa}}\starh[\sqrt{\alpha_i}\tilde\eta_n\circ\theta_i],_{i_2i_3i_4}]_t\\
&\ \ \ +H_{11}^{21}.
\end{align*}
Therefore,
\begin{align*}
\int_0^t H_{11}^2&=\ud\epsilon^{mn}\epsilon^{rs}\int_0^t \int_{(0,1)^2} [ \q(\theta_i)[\xi_l(\theta_i)h_{rs}^{(ji)_{1234}}],_{i_1}]_t  \tilde\eta_{i\kappa}^m,_{j_1j_2j_3j_4} \tilde\eta_{i\kappa}^n,_{i_2i_3i_4}\n\\
&\ \ \ -\ud\epsilon^{mn}\epsilon^{rs}\bigl[ \int_{(0,1)^2}  \q(\theta_i)[\xi_l(\theta_i)h_{rs}^{(ji)_{1234}}],_{i_1} \e_{i\kappa}^m,_{j_1j_2j_3j_4} \e_{i\kappa}^n,_{i_2i_3i_4}\bigr]_0^t\n +\ud\int_0^tH_{11}^{21}.
\end{align*}
Now, from (\ref{ecl10}) and $H_{11}^1=0$, we infer by integrating by parts in time:
\begin{align*}
\int_0^t H_{11}&=\ud\epsilon^{mn}\epsilon^{rs}\int_0^t \int_{(0,1)^2} [(\xi_l \q)(\theta_i)h_{rs}^{(ji)_{1234}}]_t  \tilde\eta_{i\kappa}^m,_{j_1j_2j_3j_4} \tilde\eta_{i\kappa}^n,_{i_1i_2i_3i_4}\n\\
&\ \ \ -\ud\epsilon^{mn}\epsilon^{rs}\bigl[ \int_{(0,1)^2} (\xi_l \q)(\theta_i)h_{rs}^{(ji)_{1234}} \e_{i\kappa}^m,_{j_1j_2j_3j_4} \e_{i\kappa}^n,_{i_1i_2i_3i_4}\bigr]_0^t\n\\
&\ \ \  +\int_0^t [S_{11} +R_{11}+\ud H_{11}^{2}].
\end{align*}
Now, we claim that the only couples $(i_1,j_1)$ contributing to the sum above are the ones with $i_1\ne j_1$. To see that, we notice that if $i_1=j_1$, then by simply relabeling $m$ and $n$, and using the symmetric role of $\{i_2,i_3,i_4\}$ and $\{j_2,j_3,j_4\}$
\begin{align*}
h_{rs}^{(ji)_{1234}}\epsilon^{mn} \e_{i\kappa}^m,_{i_1j_2j_3j_4} \e_{i\kappa}^n,_{i_1i_2i_3i_4}&=
h_{rs}^{(ji)_{1234}}\epsilon^{nm} \e_{i\kappa}^n,_{i_1j_2j_3j_4} \e_{i\kappa}^m,_{i_1i_2i_3i_4}\\
&=h_{rs}^{(ji)_{1234}}\epsilon^{nm} \e_{i\kappa}^n,_{i_1i_2i_3i_4} \e_{i\kappa}^m,_{j_1j_2j_3j_4}
\end{align*}
leading to $h_{rs}^{(ji)_{1234}}\epsilon^{mn} \e_{i\kappa}^m,_{i_1j_2j_3j_4} \e_{i\kappa}^n,_{i_1i_2i_3i_4}=0$, since $\epsilon^{mn}=-\epsilon^{nm}$.   Consequently if we denote
\begin{align*}
d^{(ji)_{234}}=\frac{h_{rs}^{(ji)_{1234}}}{[(\theta_i^{-1}\circ\theta_l),_r^{j_1}
(\theta_i^{-1}\circ\theta_l),_s^{i_1}](\theta_l^{-1}\circ\theta_i)},
\end{align*}
we have
\begin{align*}
\int_0^t H_{11}&=\ud\epsilon^{mn}\epsilon^{i_1j_1}\int_0^t \int_{(0,1)^2}\bigl[ [(\xi_l \q)(\theta_i)d^{(ji)_{234}}]_t \epsilon^{rs} \epsilon^{i_1j_1}[(\theta_i^{-1}\circ\theta_l),_r^{j_1}
(\theta_i^{-1}\circ\theta_l),_s^{i_1}](\theta_l^{-1}\circ\theta_i)\\
&\hskip 4cm \tilde\eta_{i\kappa}^m,_{j_1j_2j_3j_4} \tilde\eta_{i\kappa}^n,_{i_1i_2i_3i_4}\bigr]\n\\
&\ \ \ -\ud\epsilon^{mn}\epsilon^{i_1j_1}\bigl[ \int_{(0,1)^2} \bigl[(\xi_l \q)(\theta_i)d^{(ji)_{234}} \epsilon^{rs}\epsilon^{i_1j_1} [(\theta_i^{-1}\circ\theta_l),_r^{j_1}
(\theta_i^{-1}\circ\theta_l),_s^{i_1}](\theta_l^{-1}\circ\theta_i)\\
&\hskip 4cm \e_{i\kappa}^m,_{j_1j_2j_3j_4} \e_{i\kappa}^n,_{i_1i_2i_3i_4}\bigr]\bigr]_0^t\n\\
&\ \ \  +\int_0^t [S_{11} +R_{11}+\ud H_{11}^{2}].
\end{align*}
Now, for any fixed $(i_1,j_1)$, we have $\epsilon^{rs} \epsilon^{i_1j_1}[(\theta_i^{-1}\circ\theta_l),_r^{j_1}
(\theta_i^{-1}\circ\theta_l),_s^{i_1}]=-\text{det}(\nabla(\theta_i^{-1}\circ\theta_l))=-1$, leading us to
\begin{align*}
\int_0^t H_{11}&=-\ud\epsilon^{mn}\epsilon^{i_1j_1}\int_0^t \int_{(0,1)^2}[(\xi_l \q)(\theta_i)d^{(ji)_{234}}]_t  \tilde\eta_{i\kappa}^m,_{j_1j_2j_3j_4} \tilde\eta_{i\kappa}^n,_{i_1i_2i_3i_4}\n\\
&\ \ \ +\ud\epsilon^{mn}\epsilon^{i_1j_1}\bigl[ \int_{(0,1)^2} (\xi_l \q)(\theta_i)d^{(ji)_{234}} \e_{i\kappa}^m,_{j_1j_2j_3j_4} \e_{i\kappa}^n,_{i_1i_2i_3i_4}\bigr]_0^t\n\\
&\ \ \  +\int_0^t [S_{11} +R_{11}+\ud H_{11}^{2}].
\end{align*}
Next, by integrating by parts in space:
\begin{align}
\int_0^t H_{11}&=\ud\epsilon^{mn}\epsilon^{i_1j_1}\int_0^t \int_{(0,1)^2} [(\xi_l \q)(\theta_i)d^{(ji)_{234}}]_t,_{i_1}  \e_{i\kappa}^m,_{j_1j_2j_3j_4} \e_{i\kappa}^n,_{i_2i_3i_4}\n\\
&\ \ \ -\ud\epsilon^{mn}\epsilon^{i_1j_1}\bigl[ \int_{(0,1)^2} [(\xi_l \q)(\theta_i)d^{(ji)_{234}}],_{i_1}  \e_{i\kappa}^m,_{j_1j_2j_3j_4} \e_{i\kappa}^n,_{i_2i_3i_4}\bigr]_0^t\n\\
&\ \ \  +\int_0^t [S_{11} +R_{11}+\ud H_{11}^{2}],
\label{ecl12}
\end{align}
where we have used the fact that similarly as for $H_{11}^1$, we have
\begin{align*}
0&=\epsilon^{mn}\epsilon^{i_1j_1} d^{(ji)_{234}} \e_{i\kappa}^m,_{i_1j_1j_2j_3j_4} \e_{i\kappa}^n,_{i_2i_3i_4}.
\end{align*}
We now come to the study of the crucial term bringing the regularity of the surface.
\subsubsection{\bf Control of the trace of $\e_{i\kappa}$ on $\Gamma$}
Let us study the second term of the right-hand side of (\ref{ecl12}):
\begin{align*}
 H&=-\ud\epsilon^{mn}\epsilon^{i_1j_1} \int_{(0,1)^2} [(\xi_l \q)(\theta_i)d^{(ji)_{234}} ],_{i_1} \e_{i\kappa}^m,_{j_1j_2j_3j_4} \e_{i\kappa}^n,_{i_2i_3i_4},
\end{align*}
for which we have
\begin{align*}
 H&=-\ud\epsilon^{mn}\epsilon^{i_1j_1} \int_{(0,1)^2} \bigl[[[(\xi_l \q)(\theta_i)d^{(ji)_{234}}]\circ\theta_i^{-1}\circ(\tek)^{-1}\circ\tek\circ\theta_i ],_{i_1} \\
&\hskip 4cm [\e_{i\kappa}^m,_{j_2j_3j_4}\circ\theta_i^{-1}\circ(\tek)^{-1}\circ\tek\circ\theta_i],_{j_1} \e_{i\kappa}^n,_{i_2i_3i_4}\bigr]\\
&=-\ud\epsilon^{mn}\epsilon^{i_1j_1} \int_{(0,1)^2} \bigl[[[(\xi_l \q)(\theta_i)d^{(ji)_{234}}]\circ\theta_i^{-1}\circ(\tek)^{-1}],_{i'_1}\circ\tek\circ\theta_i\ [\tek\circ\theta_i ],_{i_1}^{i'_1} \\
&\hskip 3cm [\e_{i\kappa}^m,_{j_2j_3j_4}\circ\theta_i^{-1}\circ(\tek)^{-1}],_{j'_1}\circ\tek\circ\theta_i\ [\tek\circ\theta_i],_{j_1}^{j'_1} \e_{i\kappa}^n,_{i_2i_3i_4}\bigr].
\end{align*}
Now, for the same reason as before, the couples $(i'_1,j'_1)$ such $i'_1=j'_1$ will not contribute to the sum above, leading us to
\begin{align*}
 H&=-\ud\epsilon^{mn}\epsilon^{i'_1j'_1} \int_{(0,1)^2} \bigl[[[(\xi_l \q)(\theta_i)d^{(ji)_{234}}]\circ\theta_i^{-1}\circ(\tek)^{-1}],_{i'_1}\circ\tek\circ\theta_i\\
&\hskip 3cm \epsilon^{i_1j_1}\epsilon^{i'_1j'_1} [\tek\circ\theta_i ],_{i_1}^{i'_1}[\tek\circ\theta_i],_{j_1}^{j'_1} \\
&\hskip 3cm [\e_{i\kappa}^m,_{j_2j_3j_4}\circ\theta_i^{-1}\circ(\tek)^{-1}],_{j'_1}\circ\tek\circ\theta_i\  \e_{i\kappa}^n,_{i_2i_3i_4}\bigr]\\
&=-\ud\epsilon^{mn}\epsilon^{i'_1j'_1} \int_{(0,1)^2} \bigl[[[(\xi_l \q)(\theta_i)d^{(ji)_{234}}]\circ\theta_i^{-1}\circ(\tek)^{-1}],_{i'_1}\circ\tek\circ\theta_i\ \text{det}\nabla(\tek\circ\theta_i) \\
&\hskip 3cm [\e_{i\kappa}^m,_{j_2j_3j_4}\circ\theta_i^{-1}\circ(\tek)^{-1}],_{j'_1}\circ\tek\circ\theta_i\  \e_{i\kappa}^n,_{i_2i_3i_4}\bigr]\\
&=-\ud\epsilon^{mn}\epsilon^{i'_1j'_1} \int_{\tek(\theta_i((0,1)^2))} \bigl[[[(\xi_l \q)(\theta_i)d^{(ji)_{234}}]\circ\theta_i^{-1}\circ(\tek)^{-1}],_{i'_1} \\
&\hskip 4cm [\e_{i\kappa}^m,_{j_2j_3j_4}\circ\theta_i^{-1}\circ(\tek)^{-1}],_{j'_1}\  \e_{i\kappa}^n,_{i_2i_3i_4}\circ\theta_i^{-1}\circ(\tek)^{-1}\bigr]\\
&=I+J,
\end{align*}
with
\begin{align*}
I&=-\ud(\epsilon^{mn})^2 \int_{\tek(\theta_i((0,1)^2))} \bigl[[[(\xi_l \q)(\theta_i)d^{(ji)_{234}}]\circ\theta_i^{-1}\circ(\tek)^{-1}],_{m} \\
&\hskip 4cm [\e_{i\kappa}^m,_{j_2j_3j_4}\circ\theta_i^{-1}\circ(\tek)^{-1}],_{n}\  \e_{i\kappa}^n,_{i_2i_3i_4}\circ\theta_i^{-1}\circ(\tek)^{-1}\bigr]\\
J&=  \ud(\epsilon^{mn})^2 \int_{\tek(\theta_i((0,1)^2))} \bigl[[[(\xi_l \q)(\theta_i)d^{(ji)_{234}}]\circ\theta_i^{-1}\circ(\tek)^{-1}],_{n} \\
&\hskip 4cm [\e_{i\kappa}^m,_{j_2j_3j_4}\circ\theta_i^{-1}\circ(\tek)^{-1}],_{m}\  \e_{i\kappa}^n,_{i_2i_3i_4}\circ\theta_i^{-1}\circ(\tek)^{-1}\bigr].
\end{align*}
Next, we notice that
\begin{align*}
J&=  -\ud\sum_n \int_{\tek(\theta_i((0,1)^2))} \bigl[[[(\xi_l \q)(\theta_i)d^{(ji)_{234}}]\circ\theta_i^{-1}\circ(\tek)^{-1}],_{n} \\
&\hskip 4cm [\e_{i\kappa}^n,_{j_2j_3j_4}\circ\theta_i^{-1}\circ(\tek)^{-1}],_{n}\  \e_{i\kappa}^n,_{i_2i_3i_4}\circ\theta_i^{-1}\circ(\tek)^{-1}\bigr]+J_1,
\end{align*}
with the perturbation term
\begin{align*}
J_1&=  \ud\int_{\tek(\theta_i((0,1)^2))} \bigl[[[(\xi_l \q)(\theta_i)d^{(ji)_{234}}]\circ\theta_i^{-1}\circ(\tek)^{-1}],_{n} \\
&\hskip 4cm \div[\e_{i\kappa},_{j_2j_3j_4}\circ\theta_i^{-1}\circ(\tek)^{-1}]\  \e_{i\kappa}^n,_{i_2i_3i_4}\circ\theta_i^{-1}\circ(\tek)^{-1}\bigr]\\
&=J_1^1+J_1^2,
\end{align*}
where
\begin{align*}
J_1^1&=  \ud\int_{\tek(\theta_i((0,1)^2))} \bigl[[[(\xi_l \q)(\theta_i)c^{i_2i_3i_4}_{il,111}]\circ\theta_i^{-1}\circ(\tek)^{-1}],_{n} c^{j_2j_3j_4}_{il,111}\circ\theta_i^{-1}\circ(\tek)^{-1}\\
&\hskip 4cm \div[\e_{i\kappa},_{j_2j_3j_4}\circ\theta_i^{-1}\circ(\tek)^{-1}]\  \e_{i\kappa}^n,_{i_2i_3i_4}\circ\theta_i^{-1}\circ(\tek)^{-1}\bigr],\\
J_1^2&=  \ud\int_{\tek(\theta_i((0,1)^2))} \bigl[[(\xi_l \q)(\theta_i)c^{i_2i_3i_4}_{il,111}]\circ\theta_i^{-1}\circ(\tek)^{-1} [c^{j_2j_3j_4}_{il,111}\circ\theta_i^{-1}\circ(\tek)^{-1}],_n\\
&\hskip 4cm \div[\e_{i\kappa},_{j_2j_3j_4}\circ\theta_i^{-1}\circ(\tek)^{-1}]\  \e_{i\kappa}^n,_{i_2i_3i_4}\circ\theta_i^{-1}\circ(\tek)^{-1}\bigr].
\end{align*}
Now, for $J_1^1$, let us set
\begin{align*}
f_{il}&=[[(\xi_l \q)(\theta_i)c^{i_2i_3i_4}_{il,111}]\circ\theta_i^{-1}\circ(\tek)^{-1}],_{n} [\Pi_{p=2}^3 (\theta_i^{-1}\circ\theta_l),_1^{j_p}(\theta_l^{-1}\circ\theta_i)\e_{i\kappa}^n,_{i_2i_3i_4}](\theta_i^{-1}\circ(\tek)^{-1}),
\end{align*}
so that in order to identify an horizontal derivative on the highest order term we have:
\begin{align*}
J_1^1&=  \ud\int_{\tek(\theta_i((0,1)^2))} f_{il}\ (\theta_i^{-1}\circ\theta_l),_1^{j_4}(\theta_l^{-1}\circ(\tek)^{-1})
 \div[\e_{i\kappa},_{j_2j_3j_4}\circ\theta_i^{-1}\circ(\tek)^{-1}]\\
&=  \ud\int_{\tek(\theta_i((0,1)^2))} f_{il}\ (\theta_i^{-1}\circ\theta_l),_1^{j_4}(\theta_l^{-1}\circ(\tek)^{-1})
 (\theta_i^{-1}\circ(\tek)^{-1}),_k^l [\e_{i\kappa},_{j_2j_3j_4l}^k]\circ\theta_i^{-1}\circ(\tek)^{-1}\\
&=  \ud\int_{\tek(\theta_i((0,1)^2))} \bigl[ f_{il}\ (\theta_i^{-1}\circ\theta_l),_1^{j_4}(\theta_l^{-1}\circ(\tek)^{-1})
 (\theta_i^{-1}\circ(\tek)^{-1}),_k^l \\
&\hskip 4cm [\e_{i\kappa},_{j_2j_3l}^k\circ\theta_i^{-1}\circ(\tek)^{-1}],_q(\tek\circ\theta_i),_{j_4}^q(\theta_i^{-1}\circ(\tek)^{-1})\bigr]\\
&=  \ud\int_{\tek(\theta_i((0,1)^2))} \bigl[f_{il}\ (\theta_i^{-1}\circ\theta_l),_1^{j_4}(\theta_l^{-1}\circ(\tek)^{-1})
  \\
&\hskip 2cm [(\theta_i^{-1}\circ(\tek)^{-1}),_k^l\e_{i\kappa},_{j_2j_3l}^k\circ\theta_i^{-1}\circ(\tek)^{-1}],_q(\tek\circ\theta_i),_{j_4}^q(\theta_i^{-1}\circ(\tek)^{-1})\bigr]+r_1^1\\
&=  \ud\int_{\tek(\theta_i((0,1)^2))}\bigl[ f_{il}\ [(\tek\circ\theta_i),_{j_4}^q(\theta_i^{-1}\circ\theta_l) (\theta_i^{-1}\circ\theta_l),_1^{j_4}](\theta_l^{-1}\circ(\tek)^{-1})
  \\
&\hskip 4cm \div[\e_{i\kappa},_{j_2j_3}\circ\theta_i^{-1}\circ(\tek)^{-1}],_q \bigr]+r_1^1\\
&=  \ud\int_{\tek(\theta_i((0,1)^2))} f_{il}\ (\tek\circ\theta_l),_{1}^q (\theta_l^{-1}\circ(\tek)^{-1})\ \div[\e_{i\kappa},_{j_2j_3}\circ\theta_i^{-1}\circ(\tek)^{-1}],_q +r_1^1\\
&=  \ud\int_{\tek(\theta_i((0,1)^2))} f_{il}\  (\div[\e_{i\kappa},_{j_2j_3}\circ\theta_i^{-1}\circ(\tek)^{-1}]\circ\tek\circ\theta_l),_1(\theta_l^{-1}\circ(\tek)^{-1}) +r_1^1\\
\end{align*}
with $|r_1^1|\le C \|\tek\|^2_3$. 
Now, we notice that the presence of the factor $\xi_l\circ(\tek)^{-1}$ in $f_{il}$ implies that the integrand in the integral above is zero outside of
$\tek(\theta_l((0,1)^2))$. Similarly, the presence of $\rho_{\frac{1}{\kappa}}\starh[\sqrt{\alpha_i}\e\circ\theta_i],_{j_2j_3}\circ\theta_i^{-1}\circ(\tek)^{-1}$ implies that $x\in \tek(\theta_i(\O))$ in order for this integrand to be non-zero. Therefore,
\begin{align*}
J_1^1&= \ud\int_{\O} f_{il}(\tek\circ\theta_l)\  (\div[\e_{i\kappa},_{j_2j_3}\circ\theta_i^{-1}\circ(\tek)^{-1}]\circ\tek\circ\theta_l),_1 \text{det}\nabla(\tek\circ\theta_l) +r_1^1.
\end{align*}
Now, since the derivative of $\div[\e_{i\kappa},_{j_2j_3}\circ\theta_i^{-1}\circ(\tek)^{-1}]\circ\tek\circ\theta_l$ is taken in the horizontal direction, this implies:
\begin{align*}
|J_1^1|&\le C \|f_{il}\|_{\ud,\O} \|\div[\e_{i\kappa},_{j_2j_3}\circ\theta_i^{-1}\circ(\tek)^{-1}]\circ\tek\circ\theta_l\|_{\ud,\O}+|r_1^1|\\
&\le C \|f_{il}\|_{\ud,\O} \|\div[\e_{i\kappa},_{j_2j_3}\circ\theta_i^{-1}\circ(\tek)^{-1}]\|_{\ud,\tek(\theta_i(\O))}+C \|\text{Id}+\int_0^t \v\|_3^2.
\end{align*}
Now, since we have in the same fashion as (\ref{diveta}):
\begin{equation*}
\|\div[\e_{i\kappa},_{j_2j_3}\circ\theta_i^{-1}\circ(\tek)^{-1}]\|_{\ud,\tek(\theta_i(\O))}\le C t \tilde E(t)^2+C\kappa^\ud \tilde E(t)+C,
\end{equation*}
we then have
\begin{align}
|J_1^1|&\le Ct \tilde E(t)^3 +C+C\kappa^\ud \tilde E(t)^2.
\label{J11}
\end{align}
Next, for $J_1^2$, we notice that $[c^{j_1j_2j_3}_{il,111}\circ\theta_i^{-1}\circ(\tek)^{-1}],_n$ is a sum of product, each one containing a factor $(\theta_i^{-1}\circ\theta_l),_1^{j_p} (\theta_l^{-1}\circ(\tek)^{-1})$. This implies that $J_1^2$ can be treated in the same way as $J_1^1$, with the identification of an horizontal derivative on the highest order term of the integrand, leading to the same majorization. We can also treat $I$ in a similar fashion, due to the curl estimate (similar as (\ref{curleta})):
\begin{equation*}
\|\curl[\e_{i\kappa},_{j_2j_3}\circ\theta_i^{-1}\circ(\tek)^{-1}]\|_{\ud,\tek(\theta_i(\O))}\le C t \tilde E(t)^2 +N(u_0),
\end{equation*}
which finally provides us with:
\begin{align*}
 H&=-\ud \sum_{m\ne n}\int_{\tek(\theta_i((0,1)^2))} \bigl[[[(\xi_l \q)(\theta_i)d^{(ji)_{234}}]\circ\theta_i^{-1}\circ(\tek)^{-1}],_{m} \\
&\hskip 4cm [\e_{i\kappa}^n,_{j_2j_3j_4}\circ\theta_i^{-1}\circ(\tek)^{-1}],_{m}\  \e_{i\kappa}^n,_{i_2i_3i_4}\circ\theta_i^{-1}\circ(\tek)^{-1}\bigr]\\
&\ \ \  -\ud \sum_n\int_{\tek(\theta_i((0,1)^2))} \bigl[d^{(ji)_{234}}\circ\theta_i^{-1}\circ(\tek)^{-1}(\xi_l \q)((\tek)^{-1})],_{n} \\
&\hskip 4cm [\e_{i\kappa}^n,_{j_2j_3j_4}\circ\theta_i^{-1}\circ(\tek)^{-1}],_{n}\  \e_{i\kappa}^n,_{i_2i_3i_4}\circ\theta_i^{-1}\circ(\tek)^{-1}\bigr]+h^1\\
&=-\ud \sum_{m, n}\int_{\tek(\theta_i((0,1)^2))} \bigl[d^{(ji)_{234}}\circ\theta_i^{-1}\circ(\tek)^{-1}(\xi_l \q)((\tek)^{-1})],_{m} \\
&\hskip 4cm [\e_{i\kappa}^n,_{j_2j_3j_4}\circ\theta_i^{-1}\circ(\tek)^{-1}],_{m}\  \e_{i\kappa}^n,_{i_2i_3i_4}\circ\theta_i^{-1}\circ(\tek)^{-1}\bigr]
+h^1,
\end{align*}
with
\begin{align*}
|h^1(t)|&\le C t \tilde E(t)^3 +N(u_0)+(C\kappa^\ud+\delta) \tilde E(t)^2+C_\delta.
\end{align*}
Therefore, by integrating by parts,
\begin{align*}
 H&=-\frac{1}{4} \int_{\partial\tek(\theta_i((0,1)^2))} \bigl[d^{(ji)_{234}}\circ\theta_i^{-1}\circ(\tek)^{-1}(\xi_l \q)((\tek)^{-1})],_{m} \tilde n^\kappa_m \\
&\hskip 4cm [\e_{i\kappa}^n,_{j_2j_3j_4}\circ\theta_i^{-1}\circ(\tek)^{-1}]  \e_{i\kappa}^n,_{i_2i_3i_4}\circ\theta_i^{-1}\circ(\tek)^{-1}\bigr]+h^2,
\end{align*}
with
\begin{align*}
 h^2&=\frac{1}{4} \int_{\tek(\theta_i((0,1)^2))} \bigl[[d^{(ji)_{234}}\circ\theta_i^{-1}\circ(\tek)^{-1}(\xi \q)((\tek)^{-1})],_{mm} \\
&\hskip 4cm [\e_{i\kappa}^n,_{j_2j_3j_4}\circ\theta_i^{-1}\circ(\tek)^{-1}]  \e_{i\kappa}^n,_{i_2i_3i_4}\circ\theta_i^{-1}\circ(\tek)^{-1}\bigr]+h^1,
\end{align*}
so that, 
\begin{align*}
 |h^2(t)|&\le [\ \|\q\|_{\se}\|\e\|_3+\|\tek\|_{\se}\|q\|_3]\ \|\e\|_3^2
+|h^1(t)|\\
&\le (C\kappa^\ud+\delta) \tilde E(t)^2+C_\delta + Ct \tilde E(t)^4 +N(u_0).
\end{align*}
Now, since $\q=0$ and $\xi_l=\alpha_l$ on $\Gamma$, we infer 
\begin{align*}
 H&= -\frac{1}{4}\int_{\partial\tek(\theta_i((0,1)^2))}  \bigl[\alpha_l((\tek)^{-1})[p_0+\int_0^t p_t],_{m} [N_m+\int_0^t (\tilde n^\kappa_m)_t]\\
&\hskip 2cm [d_{il}^{j_{234}}\e_{i\kappa}^n,_{j_2j_3j_4}]\circ\theta_i^{-1}\circ(\tek)^{-1} [ d_{il}^{i_{234}}\e_{i\kappa}^n,_{i_2i_3i_4}]\circ\theta_i^{-1}\circ(\tek)^{-1}\bigr]+h^2,
\end{align*}
with $\d d_{il}^{k_{234}}=\frac{(\theta_i^{-1}\circ\theta_l),_1^{k_2}(\theta_i^{-1}\circ\theta_l),_1^{k_3}(\theta_i^{-1}\circ\theta_l),_1^{k_4}}{\sqrt{\text{det}\nabla(\tek\circ\theta_l)}}(\theta_l^{-1}\circ\theta_i)$. Therefore, 
 with the initial pressure condition $$p_0,_m N_m <-C<0\ \text{on}\ \Gamma,$$ we infer
\begin{align}
 -H&\le -C \int_{\partial\tek(\theta_i((0,1)^2))}  
 [\alpha_l(\theta_i) d_{il}^{j_{234}}\e_{i\kappa}^n,_{j_2j_3j_4} d_{il}^{i_{234}}\e_{i\kappa}^n,_{i_2i_3i_4}]\circ\theta_i^{-1}\circ(\tek)^{-1}\n\\
&\ \ \ +t \tilde E(t)^2+|h^2|\n\\
&\le -C \int_{\partial\tek(\theta_i((0,1)^2))}  
 [\alpha_i(\theta_i) d_{ii}^{j_{234}}\e_{i\kappa}^n,_{j_2j_3j_4} d_{ii}^{i_{234}}\e_{i\kappa}^n,_{i_2i_3i_4}]\circ\theta_i^{-1}\circ(\tek)^{-1}\n\\
&\ \ \ +(C\kappa^\ud+\delta) \tilde E(t)^2+C t \tilde E(t)^4 +C_\delta +N(u_0)\n\\
&\le -C \int_{\partial\tek(\theta_i((0,1)^2))}  
 [\alpha_i(\theta_i)\e_{i\kappa}^n,_{j_2j_3j_4}\e_{i\kappa}^n,_{i_2i_3i_4}]\circ\theta_i^{-1}\circ(\tek)^{-1}\n\\
&\ \ \ +(C\kappa^\ud+\delta) \tilde E(t)^2+C t \tilde E(t)^4 +C_\delta +N(u_0)\n\\
&\le -C \int_{\partial\tek(\theta_i((0,1)^2))}  
 [\alpha_i(\theta_i)\e_{i\kappa}^n,_{111} \e_{i\kappa}^n,_{111}]\circ\theta_i^{-1}\circ(\tek)^{-1}\n\\
&\ \ \ +(C\kappa^\ud+\delta) \tilde E(t)^2+C t \tilde E(t)^4 +C_\delta +N(u_0).
\label{ecl13}
\end{align}
Now, it is clear that the space integral in front of the first time integral in (\ref{ecl12}) can be treated in a similar way, except that we do not have a control on the sign of the boundary term as in (\ref{ecl13}), which does not matter since a time integral is applied to it. This therefore leads us to
\begin{align*}
-\int_0^t H_{11}&\le -C \int_{\partial\tek(\theta_i((0,1)^2)}  
 [\alpha_i(\theta_i)\e_{i\kappa}^n,_{111} \e_{i\kappa}^n,_{111}]\circ\theta_i^{-1}\circ(\tek)^{-1}\n\\
&\ \ \ +(C\kappa^\ud+\delta) \tilde E(t)^2+C t \tilde E(t)^4 +C_\delta +N(u_0),
\end{align*}
which , with (\ref{ecl1}) and (\ref{ecl2}), finally gives the trace control for each $\e_{i\kappa}$, as well as the control of $v$ around $\Gamma$:
\begin{align}
 H^\kappa (t) +  C \int_{\partial\tek(\theta_i((0,1)^2))}  |[\sqrt{\alpha_i}(\theta_i) \e_{i\kappa}^n,_{111}]&\circ\theta_i^{-1}\circ(\tek)^{-1}|^2 \n\\
& \le\delta \tilde E(t)^2+C_\delta t \tilde E(t)^4 +C_\delta N(u_0).
\label{ecl28}
\end{align}

\subsection{Asymptotic regularity of each $\ttek$}
Consequently, we infer that for each $l\in \{1,...,K\}$, we have the trace control 
\begin{align}
\|[\sqrt{\alpha_l}(\theta_l) \ttek,_1]\circ\theta_l^{-1}\circ(\tek)^{-1}\|^2_{2,\partial\tek\circ\theta_l(\O)}  \le \delta \tilde E(t)^2+C_\delta t \tilde E(t)^4 +C_\delta N(u_0).
\label{reg1}
\end{align}
Consequently, with the estimates (\ref{diveta}) and (\ref{curleta}) on the divergence and curl, we obtain by elliptic regularity:
\begin{align*}
\| [\sqrt{\alpha_l}(\theta_l)\ttek,_1]\circ\theta_l^{-1}\circ(\tek)^{-1}\|^2_{\cd,\tek\circ\theta_l(\O)}  &\le P(\|\tek\|_{\cd})[\delta \tilde E(t)^2+C_\delta t\tilde E(t)^4 +C_\delta N(u_0)]\n\\
& \le \delta \tilde E(t)^2+C_\delta t \tilde E(t)^4 +C_\delta N(u_0).
\end{align*}
Therefore,
\begin{align}
\| \sqrt{\alpha_l}(\theta_l)\ttek,_1\|^2_{\cd,\O}  &\le \delta \tilde E(t)^2+C_\delta t \tilde E(t)^4 +C_\delta N(u_0),
\label{reg2}
\end{align}
which implies that 
$\| \sqrt{\alpha_l}(\theta_l)\ttek,_{12}\|^2_{\td,\O}  $
and
$\| \sqrt{\alpha_l}(\theta_l)\ttek,_{2}\|^2_{2,\partial\O} $ are controlled by the same right-hand side as in (\ref{reg2}).
Consequently, with (\ref{diveta}) and (\ref{curleta}), we infer in the same way as we obtained (\ref{reg2}) from (\ref{reg1}) that
\begin{align*}
\|\sqrt{\alpha_l}(\theta_l) \ttek,_{2}\|^2_{\cd,\O}  &\le \delta \tilde E(t)^2+C_\delta t \tilde E(t)^4 +C_\delta N(u_0),
\end{align*}
and finally that
\begin{align}
\|\sqrt{\alpha_l}(\theta_l) \ttek\|^2_{\sd,\O}  &\le \delta \tilde E(t)^2+C_\delta t \tilde E(t)^4 +C_\delta N(u_0).
\label{reg3}
\end{align}

\subsection{Asymptotic regularity of $\tek$}

This also obviously implies that for the advected domain
\begin{align}
\| \tek\|^2_{\sd,\Omega^K}  &\le \delta \tilde E(t)^2+C_\delta t \tilde E(t)^4 +C_\delta N(u_0),
\label{reg4}
\end{align}
where $\d\Omega^K=\Omega\cap_{i=K+1}^{L} (\text{supp} \alpha_i)^c$.

From the divergence and curl estimates (\ref{divetabeta}) and (\ref{curletabeta}), we then infer that
\begin{equation}
\label{betatildeeta}
\|{\beta}\tilde\eta\circ\theta_l\|^2_{\sd,\O}\le \delta \tilde E(t)^2+C_\delta t \tilde E(t)^4 +C_\delta N(u_0),
\end{equation}
which with (\ref{reg3}) provides
\begin{equation}
\label{etakappa}
\|\tek\|^2_{\sd}\le \delta \tilde E(t)^2+C_\delta t \tilde E(t)^4 +C_\delta N(u_0,\kappa\|\Omega\|_{\frac{9}{2}}).
\end{equation}

\subsection{Asymptotic regularity of $\v$}
The relation (\ref{ecl28}) provides us with the asymptotic regularity of $\v$ near $\partial\Omega$. For the interior regularity, we notice that if we time-differentiate the analog of (\ref{diveta3}) for the cut-off $\beta\in\mathcal{D}(\Omega)$,
we obtain
\begin{equation}
\label{divv}
\bigl\|\text{div}((\beta\v\circ\theta_l),_{s}\circ\theta_l^{-1}\circ(\tek)^{-1})\bigr\|_{1,\tek(\Omega)}\le C.
\end{equation}
Similarly, we also have
\begin{equation}
\label{curlv}
\bigl\|\text{curl}(({\beta}\v\circ\theta_l),_{s}\circ\theta_l^{-1}\circ(\tek)^{-1})\bigr\|_{2,\tek(\Omega)}\le C.
\end{equation}
From (\ref{divv}) and (\ref{curlv}), elliptic regularity yields 
\begin{equation*}
\|{\beta}\v\circ\theta_l\|_{2,\O}\le C,
\end{equation*}
which, together with (\ref{ecl28}), provides
\begin{equation}
\label{v}
\|\v\|^2_{3}\le \delta \tilde E(t)^2+C_\delta t \tilde E(t)^4 +C_\delta N(u_0).
\end{equation}

\subsection{Asymptotic regularity of $\q$.}
From the elliptic system:
\begin{align*}
(\tak)_i^j[(\tak)_i^k \q,_k],_j&=-(\tuk,_i^j\u,_j^i) (\tek)\ \text{in}\ \Omega,\\
\q=0\ \text{on}\ \Gamma,
\end{align*}
we then infer on $(0,T_\kappa)$ (ensuring that (\ref{assume}) is satisfied):
\begin{align}
\|\q\|^2_\se\le C \|\tek\|^2_\se \le \delta \tilde E(t)^2+C_\delta t \tilde E(t)^4 +C_\delta N(u_0).
\label{q}
\end{align}

\subsection{Asymptotic regularity of $\kappa\sqrt{\alpha_l}\v\circ\theta_l$}
From (\ref{kappavt}), we have:
\begin{equation}
\label{10.34}
[\kappa\|\sqrt{\alpha_l}\v\circ\theta_l\|_{\se,\O}]^2\le [\kappa \|u_0\|_\se]^2+ \delta \tilde E(t)^2+C_\delta t \tilde E(t)^4 +C_\delta N(u_0).
\end{equation}

\subsection{Asymptotic regularity of $\kappa^\td\sqrt{\alpha_l}\v\circ\theta_l$}
From (\ref{kappa2vt}), we have:
\begin{equation*}
[\kappa^2\|\sqrt{\alpha_l}\v\circ\theta_l\|_{\frac{9}{2},\O}]^2\le [\kappa^2 \|u_0\|_{\frac{9}{2}}]^2+C_\delta t \tilde E(t)^4 +C_\delta N(u_0),
\end{equation*} which by interpolation leads to
\begin{equation}
\label{10.35}
[\kappa^\td\|\sqrt{\alpha_l}\v\circ\theta_l\|_{4,\O}]^2\le [\kappa^\td \|u_0\|_4]^2+ \delta \tilde E(t)^2+C_\delta t \tilde E(t)^4 +C_\delta N(u_0).
\end{equation}

\subsection{Asymptotic regularity of $\q_t$.}
From the elliptic system:
\begin{align*}
(\tak)_i^j[(\tak)_i^k \q_t,_k],_j&=-[\tuk,_i^j\u,_j^i (\tek)]_t - [(\tak)_i^j[(\tak)_i^k]_t \q,_k],_j\ \text{in}\ \Omega,\\
\q=0\ \text{on}\ \Gamma,
\end{align*}
we then infer on $(0,T_\kappa)$:
\begin{align}
\|\q_t\|^2_3 \le \delta \tilde E(t)^2+C_\delta t \tilde E(t)^4 +C_\delta N(u_0).
\label{qtassume}
\end{align}

\subsection{Asymptotic regularity of $\v_t$} Since $\v^i_t=-(\tak)_i^j q,_j$, we then infer on $(0,T_\kappa)$:
\begin{align}
\|\v_t\|^2_\cd\le C (\|\tek\|_\se +\|\q\|_\se)^2\le \delta \tilde E(t)^2+C_\delta t \tilde E(t)^4 +C_\delta N(u_0).
\label{vtassume}
\end{align}

\section{Time of existence independent of $\kappa$ and solution to the 
limit problem}
\label{L11}

By (\ref{etakappa}), (\ref{q}), (\ref{v}), (\ref{betatildeeta}), (\ref{reg3}), (\ref{10.34}), (\ref{10.35}) we then infer the control on $(0,T_\kappa)$:
\begin{equation*}
\tilde E(t)^2\le \delta \tilde E(t)^2+C_\delta t \tilde E(t)^4 +C_\delta N(u_0),
\end{equation*}
which for a choice of $\delta_0$ small enough provides us with
\begin{equation*}
\tilde E(t)^2\le C_{\delta_0} N(u_0)+C_{\delta_0} t \tilde E(t)^4.
\end{equation*}
Similarly as in Section 9 of \cite{CoSh2005b}, this provides us with a time of existence $T_\kappa=T_1$ independent of $\kappa$ and an estimate on $(0,T_1)$ independent of $\kappa$ of the type:
\begin{equation*}
\tilde E(t)^2\le N_0(u_0),
\end{equation*}
as long as the conditions (\ref{assume}) hold. Now, since $\d\|\e(t)\|_3\le \|\text{Id}\|_3+\int_0^t \|\v\|_3$, we see that condition (\ref{assume.b}) will be satisfied
for $\d t\le \frac{1}{N_0(u_0)}$. The other conditions in (\ref{assume}) are satisfied with similar arguments ((\ref{qtassume}) and (\ref{vtassume}) are used for (\ref{assume.c}) and (\ref{assume.d})). This leads us to a time of existence $T_2>0$ independent of $\kappa$ for which we have the estimate on $(0,T)$ 
\begin{equation*}
\tilde E(t)^2\le N_0(u_0),
\end{equation*}
which provides by weak convergence the existence of a solution $(v,q)$ of (\ref{euler}), with $\sigma=0$, on $(0,T)$.

\section{Optimal regularity}
\label{L12}

In this section, we assume that $\Omega$ is of class $H^\se$ in $\R^3$, that $u_0\in H^3(\Omega)$, and that the pressure condition is satisfied. We denote by $N(u_0)$ a generic constant depending on $\|u_0\|_3$. With these requirements, we will only get the $H^\se$ regularity of the moving domain $\eta(\Omega)$ and not of the mapping $\eta$.

Due to the fact that $H^\td$ is not continuously embedded in $L^\infty$ in 
the case that $\Omega$ is three-dimensional, we cannot directly study the 
integral terms as in Section \ref{L10} as we did for the two-dimensional case. 
Instead, we are forced to also regularize the initial domain, by a standard 
convolution, with a parameter $\epsilon>0$ fixed independently of $\kappa$, on 
the charts defining it locally, so that the initial regularized domain 
$\Omega_\epsilon=\tilde\Omega$ obtained in this fashion is of class $C^\infty$.
The regularized initial velocity, by a standard convolution, will be denoted 
$u_0(\epsilon)$. We then start at Section \ref{L4} in the same way except that 
the regularity of the functional framework is increased by one degree for each 
quantity. This leaves us with the existence of a solution to (\ref{smoothl}) 
on $(0,T_{\kappa,\epsilon})$, with initial domain $\Omega_\epsilon$ and initial
velocity $u_0(\epsilon)$. We then perform the same asymptotic analysis as 
$\kappa\rightarrow 0$ as we did in Sections \ref{L5} to \ref{L10}, in this new 
framework. We then see that the problematic term is now updated to one which 
can be treated directly by the Sobolev embedding of $H^2$ into $L^\infty$ in 
3d. This leads us to the existence of a solution to a system similar to 
(\ref{euler}) (with $\sigma=0$) with initial domain $\Omega_\epsilon$ on 
$(0,T_\epsilon)$, with 
$\eta_\epsilon\in L^\infty(0,T_\epsilon;H^{\frac{9}{2}}(\Omega_\epsilon))$, 
with initial domain $\Omega_\epsilon$ and initial velocity $u_0(\epsilon)$. 

We then study hereafter the asymptotic behavior of this solution and of 
$T_\epsilon$ as $\epsilon\rightarrow 0$. This will be less problematic than in 
Section \ref{L10} since the convolutions by layers with the parameter $\kappa$ 
do not appear in the problem (\ref{euler}) with smoothed initial data and 
domain. We will denote the dependence on $\epsilon$ this time by a tilde, 
$\v$ standing here for $v_\epsilon$ for instance, and prove that as 
$\epsilon\rightarrow 0$, the time of existence and norms of $\v$ are 
$\epsilon-$independent, which leads to the existence of a solution with 
optimal regularity on the initial data, as stated in Theorem \ref{ltheorem2}. 

Our functional framework will be different than in Sections \ref{L5} to 
\ref{L10}. Our continuous in time energy will be:
\begin{definition}
\begin{align}
\tilde H(t)=&\sup_{[0,t]}\bigl[|\tilde n|_{2,\tilde\Gamma}+\|\v\|_3+\|\v_t\|_\ci+\|\q\|_3+\|\v_{tt}\|_2 \bigr]+1,
\end{align}
where $\tilde n$ denotes the unit exterior normal to $\e(\Omega)$.
\end{definition}
Our condition on $T_\epsilon$ will be that on $(0,T_\epsilon)$,
\begin{subequations}
\label{assume2}
\begin{align}
\ud\le \text{det} \nabla\e &\le \td\ \text{in}\ \tilde \Omega,
\label{assume2.a}\\
\|\e\|_3&\le |\Omega|+1,\ \
\|\q\|_3\le \|q_0\|_3+1,\ \
\|\v\|_\cd\le \|u_0\|_\cd+1,\\
\|\v_t\|_2&\le \|w_1\|_2+1,\\
\forall l\in\{1,...,K\},&\ \bigl|\e\circ\theta_l,_1\times\e\circ\theta_l,_2\bigr|\ge \ud \bigl|\theta_l,_1\times\theta_l,_2\bigr|\ \text{on}\ (0,1)^2\times\{0\},
\label{assume2.d}
\end{align}
\end{subequations} 
where $w_1=-\nabla q_0\in H^\cd(\Omega)$.
We will use a more straightforward approach than in Section \ref{L10}, which is enabled by the fact that we have $\tilde a$ 
instead of the convolution 
by layers $\tak$ in our equation, by defining the following energy:
\begin{definition}
$$\d E^\epsilon(t)= \sum_{l=1}^K \int_{(0,1)^3}\xi\circ\theta_l |D^2(\tilde v_{tt}\circ\theta_l)|^2,$$ where $D^2 f$ stands for any
second space derivative in a horizontal direction, i.e.,
 $f,_{\alpha_1\alpha_2}$, where $\alpha_i\in\{1,2\}$. Summation over all
horizontal derivatives is taken in the expression for $E^\kappa$. 
\end{definition}
\begin{remark}
We also note that this energy is associated with the second time-differentiated 
problem; we thus avoid the use of the curl relation  (\ref{curleta}) for $\e$, 
which necessitates the supplementary condition 
$\curl u_0\in H^\cd(\Omega)$ (which we do not have here). 
\end{remark}
With $\tilde b_l=[\nabla(\tilde\eta\circ\theta_l)]^{-1}$, we have:
$E^\epsilon_t=\sum_{i=1}^9 E_i,$ with
\begin{align*}
E_1(t)&= -\sum_{l=1}^K \int_{(0,1)^3} \xi(\theta_l) D^2([(\tilde b_l)_j^k]_{tt})(\q\circ\theta_l),_kD^2({\tilde v}_{tt}\circ\theta_l)^j,\\
E_2(t)&= -2 \sum_{l=1}^K \int_{(0,1)^3} \xi(\theta_l) D([(\tilde b_l)_j^k]_{tt})D(\q\circ\theta_l),_k D^2({\tilde v}_{tt}\circ\theta_l)^j,\\
E_3(t)&= -\sum_{l=1}^K \int_{(0,1)^3} \xi(\theta_l) [(\tilde b_l)_j^k]_{tt} D^2(\q\circ\theta_l),_kD^2({\tilde v}_{tt}\circ\theta_l)^j,\\
E_4(t)&= -4\sum_{l=1}^K \int_{(0,1)^3} \xi(\theta_l) D[(\tilde b_l)_j^k]_{t}D(\q_t\circ\theta_l),_k D^2({\tilde v}_{tt}\circ\theta_l)^j,\\
E_5(t)&= -2\sum_{l=1}^K \int_{(0,1)^3} \xi(\theta_l) D^2[(\tilde b_l)_j^k]_{t}(\q_t\circ\theta_l),_k D^2({\tilde v}_{tt}\circ\theta_l)^j,\\
E_6(t)&= -2\sum_{l=1}^K \int_{(0,1)^3} \xi(\theta_l) [(\tilde b_l)_j^k]_{t}D^2(\q_t\circ\theta_l),_k D^2({\tilde v}_{tt}\circ\theta_l)^j,\\
E_7(t)&= -\sum_{l=1}^K\int_{(0,1)^3} \xi(\theta_l)D^2(\tilde b)_j^k (\q_{tt}\circ\theta_l),_kD^2{(\tilde v_{tt}\circ\theta_l)}^j,\\
E_8(t)&= -2\sum_{l=1}^K\int_{(0,1)^3} \xi(\theta_l)D(\tilde b)_j^k D(\q_{tt}\circ\theta_l),_kD^2{(\tilde v_{tt}\circ\theta_l)}^j,\\
E_9(t)&= -\sum_{l=1}^K\int_{(0,1)^3} \xi(\theta_l)(\tilde b)_j^k D^2(\q_{tt}\circ\theta_l),_kD^2{(\tilde v_{tt}\circ\theta_l)}^j.
\end{align*}
\subsection{Estimate for $\q_t$, $\q_{tt}$ and $\q_{ttt}$}
From the elliptic system
\begin{align*}
\tilde a_i^j(\tilde a_i^k \q_t,_k),_j&=-[\tilde a_i^j(\tilde a_i^k]_t \q,_k),_j-[\tilde a_j^k\v^i,_k\tilde a_i^l\v^j,_l]_t\ \text{in}\ \tilde\Omega,\\
\q_t&=0\ \text{on}\ \partial\tilde\Omega,
\end{align*}
we infer
\begin{align}
\|\q_t\|_3&\le C [\ \|\v\|_3+\|\q\|_3+\|\e\|_3+\|\v_t\|_2]\le C \tilde H(t).
\label{qt3}
\end{align}
For similar reasons, we also have
\begin{subequations}
\begin{align}
\|\q_{tt}\|_\cd&\le C \tilde H(t),
\label{qtt3}\\
\|\q_{ttt}\|_2&\le C \tilde H(t).
\label{qttt3}
\end{align}
\end{subequations}

\subsection {Estimate for $E_2$, $E_4$, $E_5$, $E_6$, $E_8$} 
We first immediately have thanks to the embedding of $H^1$ into $L^6$ and $H^\ud$ into $L^3$:
\begin{subequations}
\label{H246}
\begin{align}
|E_2(t)|&\le C \|\tilde a_{tt}\|_1 \|\q\|_\cd \|\v_{tt}\|_2 \le C [\ \|\v_t\|_2\|\e\|_3+
\|v\|_3^2]\ \|\q\|_3 \|\v_{tt}\|_2\le C \tilde H(t)^4,\\
|E_4(t)|&\le C \|\tilde a_t\|_2 \|\q_t\|_\cd \|\v_{tt}\|_2\le C \tilde H(t)^3,\label{H246.b}\\
|E_5(t)|&\le  C \|\v\|_3\|\q_t\|_3 \|\v_{tt}\|_2\le C \tilde H(t)^3,\label{H246.c}\\
|E_6(t)|&\le  C \|\v\|_3\|\q_t\|_3 \|\v_{tt}\|_2\le C \tilde H(t)^3,\label{H246.d}\\
|E_8(t)|&\le C \|\tilde a\|_2\|\q_{tt}\|_\cd\|\v_{tt}\|_2\le C \tilde H(t)^3,\label{H246.e}
\end{align}
\end{subequations}
where we have used (\ref{qt3}) for (\ref{H246.c}), (\ref{H246.d}), and (\ref{qtt3}) for (\ref{H246.e}).

\subsection{Estimate for $E_3$} By integrating by parts, and using $[(\tilde b)_j^k],_k=0$, we obtain $E_3=E_3^1+E_3^2$, with
\begin{align*}
E_3^1&= \sum_{l=1}^K\int_{(0,1)^3} \xi(\theta_l)[(\tilde b)_j^k]_{tt} D^2(\q\circ\theta_l) D^2{(\tilde v_{tt}\circ\theta_l)},_k^j,\\
E_3^2&= \sum_{l=1}^K\int_{(0,1)^3} \xi(\theta_l),_k [(\tilde b)_j^k],_{tt} D^2( \q\circ\theta_l) D^2{(\tilde v_{tt}\circ\theta_l)}^j.
\end{align*}
We first have
\begin{align*}
|E_3^2(t)|&\le C \|\tilde a_{tt}\|_1\|\q\|_3\|\v_{tt}\|_2\le C \tilde H(t)^4.
\end{align*}
Next, $E_3^1=\sum_{i=1}^3 E_3^{1i},$ with
\begin{align*}
E_3^{11}&= -\sum_{l=1}^K\int_{(0,1)^3} \xi(\theta_l)D[(\tilde b)_j^k]_{tt} D^2(\q\circ\theta_l) D{(\tilde v_{tt}\circ\theta_l)},_k^j,\\
E_3^{12}&=-\sum_{l=1}^K\int_{(0,1)^3} \xi(\theta_l)[(\tilde b)_j^k]_{tt} D^3(\q\circ\theta_l) D{(\tilde v_{tt}\circ\theta_l)},_k^j,\\
E_3^{13}&=-\sum_{l=1}^K\int_{(0,1)^3} D(\xi(\theta_l))[(\tilde b)_j^k]_{tt} D^2(\q\circ\theta_l) D{(\tilde v_{tt}\circ\theta_l)},_k^j.
\end{align*}
We obviously have $|E_3^{13}(t)|\le C \tilde H(t)^4$. Next, we have by integrating by parts in time:
\begin{align*}
\int_0^t E_3^{12}&=\sum_{l=1}^K\int_0^t\int_{(0,1)^3} \xi(\theta_l)\bigl([(\tilde b)_j^k]_{tt} D^3(\q_t\circ\theta_l)+[(\tilde b)_j^k]_{ttt} D^3(\q\circ\theta_l)\bigr) D{(\tilde v_{t}\circ\theta_l)},_k^j\\
&\ \ \ +\bigl[\sum_{l=1}^K\int_{(0,1)^3} \xi(\theta_l)[(\tilde b)_j^k]_{tt} D^3(\q\circ\theta_l) D{(\tilde v_{t}\circ\theta_l)},_k^j\bigr]_0^t,
\end{align*}
showing, with the continuous embedding of $H^1$ into $L^6$ and of $H^\ud$ into $L^3$:
\begin{align}
\bigl|\int_0^t E_3^{12}\bigr|&\le C t \sup_{[0,t]} [\|\tilde a_{tt}\|_1\|\q_t\|_3\|\v_t\|_\cd+ \|\tilde a_{ttt}\|_1\|\q\|_3\|\v_t\|_\cd]\n\\
&\ \ \ +\bigl[\sum_{l=1}^K\int_{(0,1)^3} \xi(\theta_l)D{(\tilde v_{t}\circ\theta_l)},_k^j  \bigl[  [(\tilde b)_j^k]_{tt}(0) D^3(\q\circ\theta_l)(0) +\int_0^\cdot
[[(\tilde b)_j^k]_{tt} D^3(\q\circ\theta_l)]_t\bigr]
\bigr]_0^t\n\\
&\le C t \sup_{[0,t]} [\|\tilde a_{tt}\|_1\|\q_t\|_3\|\v_t\|_\cd+ \|\tilde a_{ttt}\|_1\|\q\|_3\|\v_t\|_\cd]\n\\
&\ \ \ +\|\v_{t}|_\cd \|\q(0)\|_3 \|\b_{tt}(0)\|_1 + \|\v_t|_\cd t \sup_{[0,t]} [ \|\q_t\|_3\|\b_{tt}\|_1+\|\q\|_3\|\b_{ttt}\|_1] +N(u_0)\n\\
&\le C\delta \tilde H(t)^2+t \tilde H(t)^4+C_\delta N(u_0),
\label{E312}
\end{align}
for any $\delta>0$.
For the remaining term $E_3^{11}$,
\begin{align*}
|E_3^{11}|&\le C \|\tilde a_{tt}\|_\td\|\q\|_3\|\v_{tt}\|_2\le C \tilde H(t)^4.
\end{align*}
Consequently, we have
\begin{align}
\bigl|\int_0^t E_3\bigr|&\le C\delta \tilde H(t)^2+t \tilde H(t)^4+C_\delta N(u_0).
\label{H3}
\end{align}

\subsection{Estimate for $E_7$} By integrating by parts, $E_7=E_7^1+E_7^2$, with
\begin{align*}
E_7^1&= \sum_{l=1}^K\int_{(0,1)^3} \xi(\theta_l)D^2(\tilde b)_j^k (\q_{tt}\circ\theta_l) D^2{(\tilde v_{tt}\circ\theta_l)},_k^j,\\
E_7^2&= \sum_{l=1}^K\int_{(0,1)^3} \xi(\theta_l),_k D^2(\tilde b)_j^k ( \q_{tt}\circ\theta_l) D^2{(\tilde v_{tt}\circ\theta_l)}^j.
\end{align*}
We first have
\begin{align*}
|E_7^2(t)|&\le C \|\tilde a\|_2\|\q_{tt}\|_2\|\v_{tt}\|_2\le C \tilde H(t)^4.
\end{align*}
Next, we notice by integrating by parts in time and space that
\begin{align*}
\int_0^tE_7^1&= \sum_{l=1}^K\int_0^t\int_{(0,1)^3} [D^2(\tilde b_t)_j^k (\xi\q_{tt}\circ\theta_l),_k +D^2(\tilde b)_j^k (\xi\q_{ttt}\circ\theta_l),_k] D^2{(\tilde v_{t}\circ\theta_l)}^j\\
&\ \ \ +\bigl[\sum_{l=1}^K\int_{(0,1)^3} \xi(\theta_l)D^2(\tilde b)_j^k (\xi\q_{tt}\circ\theta_l),_k D^2{(\tilde v_{t}\circ\theta_l)}^j\bigr]_0^t
\end{align*}
In the same fashion as we obtained (\ref{E312}), we then infer
\begin{align}
\bigl|\int_0^t E_7^{1}\bigr|&\le C t \sup_{[0,t]} [\|\b_{t}\|_2\|\q_{tt}\|_2\|\v_t\|_\cd+ \|\b\|_2\|\q_{ttt}\|_2\|\v_t\|_\cd]\n\\
&\ \ \ +\|\v_t\|_\cd \|\q_{tt}(0)\|_2 \|\b(0)\|_2 + \|\v_{t}\|_\cd t \sup_{[0,t]}[\|\q_{ttt}\|_2\|\b\|_2+\|\q_{tt}\|_2\|\b_t\|_2] +N(u_0)\n\\
&\le C\delta \tilde H(t)^2+t \tilde H(t)^4+C_\delta N(u_0),
\label{H51}
\end{align}
where we have used (\ref{qttt3}) for $\q_{ttt}$.
Consequently, we have
\begin{align}
\bigl|\int_0^t E_7\bigr|&\le C\delta \tilde H(t)^2+t \tilde H(t)^4+C_\delta N(u_0).
\label{H5}
\end{align}

\subsection{Estimate for $E_9$}
We notice by integrating by parts in space that 
\begin{align*}
E_9(t)&= \sum_{l=1}^K \int_{(0,1)^3} \xi(\theta_l) [(\tilde b_l)_j^k] D^2(\q_{tt}\circ\theta_l)D^2({\tilde v}_{tt}\circ\theta_l),_k^j\\
&\ \ \ +\sum_{l=1}^K \int_{(0,1)^3} \xi(\theta_l),_k [(\tilde b_l)_j^k] D^2(\q_{tt}\circ\theta_l)D^2({\tilde v}_{tt}\circ\theta_l)^j
\end{align*}
Next by the divergence condition,
\begin{align*}
E_9(t)&= -\sum_{l=1}^K \int_{(0,1)^3} \xi(\theta_l)  D^2[(\tilde b_l)_j^k] D^2(\q_{tt}\circ\theta_l)({\tilde v}_{tt}\circ\theta_l),_k^j\\
&\ \ \ - 2 \sum_{l=1}^K \int_{(0,1)^3} \xi(\theta_l)  D[(\tilde b_l)_j^k] D^2(\q_{tt}\circ\theta_l)D({\tilde v}_{tt}\circ\theta_l),_k^j\\
&\ \ \ +\sum_{l=1}^K \int_{(0,1)^3} \xi(\theta_l),_k [(\tilde b_l)_j^k] D^2(\q_{tt}\circ\theta_l)D^2({\tilde v}_{tt}\circ\theta_l)^j,
\end{align*}
showing that
\begin{align}
|E_9(t)|\le C \|\b\|_2\|\q_{tt}\|_\cd\|\v_{tt}\|_2\le C \tilde H(t)^4.
\label{H7}
\end{align}

\subsection{Estimate for $E_1$}

If $\epsilon^{ijk}$ denotes the sign of the permutation between $\{i,j,k\}$ and $\{1,2,3\}$, if $i,j, k$ are distinct, and is set to zero otherwise, we obtain
\begin{align*}
E_1&=E_1^1+E_1^2,
\end{align*}
with
\begin{align*}
\d E_1^1&=\sum_{l=1}^K\int_{(0,1)^3} \xi(\theta_l),_k \q\circ\theta_l D^2[(\b_l)^k_j]_{tt} D^2(\v_{tt}\circ\theta_l)^j,\\
E_1^2&=\ud\sum_{l=1}^K\epsilon^{mnj}\epsilon^{pqk}\int_{(0,1)^3} \xi\q(\theta_l) D^2[\e\circ\theta_l,_p^m\e\circ\theta_l,_q^n]_{tt} D^2(\v_{tt}\circ\theta_l),_k^j.
\end{align*}
\subsubsection{\bf Estimate of $E_1^1$}
Now, for $E_1^1$, since
\begin{equation*}
D_t\triangle \u+\u,_j^i \u,_{ij}-\nabla[\u,_i^j\u,_j^i]=0,
\end{equation*}
we obtain in $\tilde\Omega$ that
\begin{align*}
\tilde a_i^j(\tilde a_i^k\v,_k),_j(t)=\triangle \u(0)+\int_0^t [\tilde a_n^m (\tilde a_i^k\v,_k^j\tilde a_j^l\v,_l^i),_m]_{n=1}^3-\int_0^t \tilde a_j^k\v,_k^i\tilde a_j^m(\tilde a_i^l\v,_l),_m,
\end{align*}
and thus,
\begin{align}
\tilde a_i^j(\tilde a_i^k\v_t,_k),_j=-[\tilde a_i^j(\tilde a_i^k]_t\v,_k),_j+[\a_n^m (\tilde a_i^k\v,_k^j\tilde a_j^l\v,_l^i),_m]_{n=1}^3- \tilde a_j^k\v,_k^i\tilde a_j^m(\tilde a_i^l\v,_l),_m.
\label{vtinterior}
\end{align}
By elliptic regularity in the interior of $\tilde\Omega$, we infer that 
for any $\omega$ whose closure is contained in $\tilde\Omega$,
\begin{align*}
\|\v_t\|_{3,\omega}\le C_\omega [ \|\a_i^j(\a_i^k\v_t,_k),_j\|_{1,\tilde\Omega}
+\|\v_t\|_2]\le C_\omega \tilde H(t).
\end{align*}
With this estimate and the condition $\xi\circ\theta_l,_k=0$ in a neighborhood of $(0,1)^2\times\{0\}$, we then obtain 
\begin{align}
|E_1^1|&\le C \|\q\|_2 \bigl[ C \tilde H(t) \|\e\|_3+\|\v\|_3^2\bigr] \|\v_{tt}\|_2\n\\
&\le C (\tilde H(t)^3+1).
\label{E11}
\end{align}
\subsubsection{\bf Estimate for $E_1^2$ and the trace regularity}
We now study $E_1^2$, which will be the term bringing the asymptotic regularity of the moving domain $\e(\tilde\Omega)$. We have that
\begin{align*}
E_1^2= \sum_{l=1}^3 E_1^{2l},
\end{align*}
with
\begin{align*}
E_1^{21}&=\sum_{l=1}^K \epsilon^{mnj}\epsilon^{pqk}\int_{(0,1)^3} \xi\q(\theta_l) D^2[\v\circ\theta_l,_p^m\v\circ\theta_l,_q^n] D^2(\v_{tt}\circ\theta_l),_k^j\\
E_1^{22}&=2\sum_{l=1}^K\epsilon^{mnj}\epsilon^{pqk}\int_{(0,1)^3} \xi\q(\theta_l) [D(\v_t\circ\theta_l),_p^m D(\e\circ\theta_l),_q^n] D^2(\v_{tt}\circ\theta_l),_k^j\\
E_1^{23}&=\sum_{l=1}^K\ud\epsilon^{mni}\epsilon^{pqj}\int_{(0,1)^3}\xi(\theta_l)[D^2(\tilde\eta_{tt}\circ\theta_l),_{p}^m D^2(\tilde\eta_{tt}\circ\theta_l),_j^i]_t\q\circ\theta_l(\tilde\eta\circ\theta_l),_q^n.
\end{align*}
We first notice that
\begin{align}
|E_1^{21}|+|E_1^{22}|&\le C \|\q\|_3\|\v\|_3^2\|\v_{tt}\|_2+C  \|\q\|_3\|\v_t\|_\cd\|\e\|_3\|\v_{tt}\|_2\n\\
&\le C \tilde H(t)^4.
\label{E121122}
\end{align}
Now, for the remaining term $E_1^{23}$, an integration by parts in time provides$$\int_0^t E_1^{23}=\ud E_1^{231}+\ud[E_1^{232}]_0^t,$$ with
\begin{align*}
 E_1^{231}&=- \sum_{l=1}^K\epsilon^{mni}\epsilon^{pqj}\int_0^t\int_{(0,1)^3} \xi(\theta_l)D^2(\tilde v_t\circ\theta_l),_{p}^m D^2(\tilde v_t\circ\theta_l),_j^i(\q\circ\theta_l\ (\tilde \eta\circ\theta_l),_q^n)_t,\\
E_1^{232}&=\epsilon^{mni}\epsilon^{pqj}\int_{(0,1)^3} \xi(\theta_l)D^2(\tilde v_t\circ\theta_l),_{p}^m D^2(\tilde v_t\circ\theta_l),_j^i\q\circ\theta_l\ (\tilde \eta\circ\theta_l),_q^n.\\
\end{align*}
First, for the perturbation term $E_1^{231}$, by integrating by parts in space (and using $\q=0$ on $\Gamma$):
\begin{align*}
E_1^{231}&= \int_0^t\sum_{l=1}^K\epsilon^{mni}\epsilon^{pqj}\int_{(0,1)^3} \xi(\theta_l)D^2(\tilde v_t\circ\theta_l),_{pj}^m D^2(\tilde v_t\circ\theta_l)^i(\q\circ\theta_l\ (\tilde\eta\circ\theta_l),_q^n)_t\\
&\ \ \ +\int_0^t\sum_{l=1}^K \epsilon^{mni}\epsilon^{pqj}\int_{(0,1)^3}D^2(\tilde v_t\circ\theta_l),_{p}^m D^2(\tilde v_t\circ\theta_l)^i(\xi\q\circ\theta_l\ (\tilde\eta\circ\theta_l),_q^n)_t,_j.
\end{align*}
For the first integral, we notice that for any $f$, $g$ smooth,
\begin{align*}
\epsilon^{mni}\epsilon^{pqj}  f,_{pj}^m  f^i g,_q^n=\epsilon^{mni}\epsilon^{jqp}  f,_{jp}^m f^i g,_q^n,
\end{align*}
and, since $\epsilon^{pqj}=-\epsilon^{jqp}$, this quantity equals zero, leading to
\begin{align*}
E_1^{231}&= \int_0^t \sum_{l=1}^K \epsilon^{mni}\epsilon^{pqj}\int_{(0,1)^3} D^2(\tilde v_t\circ\theta_l),_{p}^m D^2(\tilde v_t\circ\theta_l)^i(\xi\q\circ\theta_l(\tilde\eta\circ\theta_l),_q^n)_t,_j.
\end{align*}
Consequently, as in (\ref{int}) since the derivatives in $D^2$ are horizontal,
\begin{align}
|E_1^{231}|&\le C \int_0^t \|\tilde v_t\|_\ci^2\| [\ \|\nabla(\q_t\nabla\eta)\|_1+\|\nabla(\q\nabla\v)\|_1]\n\\
&\le C t H(t)^4.
\label{E1231}
\end{align}
Now for $ E_1^{232}$, we will introduce the notation 
$$V(l)=\v\circ\theta_l \text{ and } E(l)=\e\circ\theta_l.$$ 
We then have after a change of variables made in order to get vector fields whose divergence and curl are controlled:
\begin{align*}
E_1^{232}= \sum_{l=1}^K \epsilon^{mni}\epsilon^{pqj}\int_{(0,1)^3}\bigl[&\xi\q\circ\theta_l\ (D^2 V_t(l)\circ E(l)^{-1}),_{p_1}^m(E(l)) E(l),_p^{p_1} \\
& (D^2 V_t(l)\circ E(l)^{-1}),_{j_1}^i(E(l)) E(l),_j^{j_1} \ (E(l)\circ E(l)^{-1}),_{q_1}^n E(l),_q^{q_1}\bigr].
\end{align*}
Now, we notice that any triplet $(i_1,j_1,q_1)$ such that $\text{Card}\{i_1,j_1,q_1\}<3$ will not contribute to this sum. For instance, if $j_1=q_1$, we notice that by relabeling $j$ and $q$, 
\begin{align*}
\epsilon^{pqj} E(l),_p^{p_1} E(l),_j^{j_1} E(l),_q^{j_1}=\epsilon^{pjq} E(l),_p^{p_1} E(l),_q^{j_1} E(l),_j^{j_1},
\end{align*}
where $j_1=q_1$ is fixed in the sum above. Now, since $\epsilon^{pqj} =-\epsilon^{pjq}$, this shows that 
\begin{align*}
\epsilon^{pqj} E(l),_p^{p_1} E(l),_j^{j_1} E(l),_q^{j_1}=0.
\end{align*}
By a similar argument, 
\begin{align*}
\epsilon^{pqj} E(l),_p^{p_1} E(l),_j^{p_1} E(l),_q^{q_1}=0.
\end{align*}
Consequently, only the triplets where $\text{Card}\{i_1,j_1,q_1\}=3$ contribute to $E_1^{232}$, showing that
\begin{align*}
E_1^{232}= \sum_{l=1}^K \epsilon^{mni}\epsilon^{p_1q_1j_1}\int_{(0,1)^3}\bigl[&\xi\q\circ\theta_l\ \epsilon^{pqj}\epsilon^{p_1q_1j_1} E(l),_p^{p_1}E(l),_j^{j_1}  E(l),_q^{q_1}\\
&(D^2 V_t(l)\circ E(l)^{-1}),_{p_1}^m(E(l)) (D^2 V_t(l)\circ E(l)^{-1}),_{j_1}^i(E(l)) \delta_{q_1}^n \bigr].
\end{align*}
Since for each given $(p_1,j_1,q_1)$ we have $\epsilon^{pqj}\epsilon^{p_1q_1j_1} E(l),_p^{p_1}E(l),_j^{j_1}  E(l),_q^{q_1}=\text{det}\nabla E(l)=1$, we then infer
\begin{align*}
E_1^{232}=  \sum_{l=1}^K \epsilon^{mni}\epsilon^{pqj}\int_{\e(\tilde\Omega)}\bigl[\xi\q\circ\e^{-1} 
(D^2 V_t(l)\circ E(l)^{-1}),_{p}^m (D^2 V_t(l)\circ E(l)^{-1}),_{j}^i \delta_{q}^n \bigr].
\end{align*}
By integrating by parts in space, we get by using $\q=0$ on $\tilde\Gamma$:
\begin{align*}
E_1^{232}= - \sum_{l=1}^K \epsilon^{mni}\epsilon^{pqj}\int_{\e(\tilde\Omega)}\bigl[\xi\q\circ\e^{-1} 
(D^2 V_t(l)\circ E(l)^{-1}),_{pj}^m (D^2 V_t(l)\circ E(l)^{-1})^i \delta_{q}^n \bigr]\\
- \sum_{l=1}^K \epsilon^{mni}\epsilon^{pqj}\int_{\e(\tilde\Omega)}\bigl[(\xi\q\circ\e^{-1}),_j 
(D^2 V_t(l)\circ E(l)^{-1}),_{p}^m (D^2 V_t(l)\circ E(l)^{-1})^i \delta_{q}^n \bigr].
\end{align*}
Next, since for any $f$ smooth,
\begin{align*}
\epsilon^{mni}\epsilon^{pqj} f,_{pj}^m \delta_{q}^n=0,
\end{align*}
we infer
\begin{align*}
E_1^{232}&= - \sum_{l=1}^K \epsilon^{mni}\epsilon^{pqj}\int_{\e(\tilde\Omega)}(\xi\q\circ\e^{-1}),_j 
(D^2 V_t(l)\circ E(l)^{-1}),_{p}^m (D^2 V_t(l)\circ E(l)^{-1})^i \delta_{q}^n.
\end{align*}
Now, we notice that for $\delta_n^q\epsilon^{mni}\epsilon^{pqj}\ne 0$, if $p=i$, then necessarily $j=m$. Similarly, if $p\ne i$, then since $p\ne n$, necessarily, $p=m$, and thus $i=j$. Therefore,
\begin{align*}
E_1^{232}&= -\sum_{l=1}^K\epsilon^{mni}\epsilon^{inm}\int_{\e(\tilde\Omega)} (D^2 V_t(l)\circ E(l)^{-1}),_{i}^m (D^2 V_t(l)\circ E(l)^{-1})^i\ (\xi\q\circ\e^{-1}),_m\\
&\ \ \ -\sum_{l=1}^K\epsilon^{mni}\epsilon^{mni}\int_{\e(\tilde\Omega)} (D^2 V_t(l)\circ E(l)^{-1}),_{m}^m (D^2 V_t(l)\circ E(l)^{-1})^i\ (\xi\q\circ\e^{-1}),_i.
\end{align*}
Now, in the same way as we obtained the divergence and curl estimates (\ref{divv}) and (\ref{curlv}), we have the same type of estimates for $V_t(l)$ leading us to
\begin{align*}
\|\sqrt{\xi}(\theta_l)\ [(D^2 V_t(l)\circ E(l)^{-1}),_i^m-(D^2 V_t(l)\circ E(l)^{-1}),_m^i]\|_{H^\ud(\e(\tilde\Omega))'}&\le C,\\
\|\sqrt{\xi}(\theta_l) \text{div}(D^2 V_t(l)\circ E(l)^{-1})\|_{H^\ud(\e(\tilde\Omega))'}&\le C.
\end{align*}
The fact that $D^2$ contains horizontal derivatives once again played a crucial role in these estimates. This implies
\begin{align*}
E_1^{232}&=  \sum_{l=1}^K\sum_{i\ne m}\int_{\e(\tilde\Omega)} (D^2 V_t(l)\circ E(l)^{-1}),_{m}^i (D^2 V_t(l)\circ E(l)^{-1})^i\ (\xi\q\circ\e^{-1}),_m\\
&\ \ \ +\sum_{l=1}^K\int_{\e(\tilde\Omega)} (D^2 V_t(l)\circ E(l)^{-1}),_{i}^i (D^2 V_t(l)\circ E(l)^{-1})^i\ (\xi\q\circ\e^{-1}),_i+R_1\\
&= \sum_{l=1}^K\int_{\e(\tilde\Omega)} (D^2 V_t(l)\circ E(l)^{-1}),_{m}^i (D^2 V_t(l)\circ E(l)^{-1})^i\ (\xi\q\circ\e^{-1}),_m +R_1,
\end{align*}
with $|R_1(t)|\le C t \tilde H(t)$. Consequently,
\begin{align*}
E_1^{232}&=- \frac{1}{2}\sum_{l=1}^K\int_{\e(\tilde\Omega)} 
|D^2 V_t(l)\circ E(l)^{-1}|^2 \triangle(\xi\q\circ\e^{-1})\\
&\ \ \ +\frac{1}{2}\sum_{l=1}^K\int_{\partial\e(\tilde\Omega)} |D^2 V_t(l)\circ E(l)^{-1}|^2  (\xi\q\circ\e^{-1}),_m \tilde n_m+R_1.
\end{align*}
If we write $\d\q=\q(0)+\int_0^t\q_t$, $\d \tilde n=N+\int_0^t \tilde n_t$, and use the fact that $\xi=1$ on $\tilde\Gamma$, we then get
\begin{align}
E_1^{232}=\frac{1}{2} \sum_{l=1}^K\int_{\partial\e(\tilde\Omega)} |D^2 V_t(l)\circ E(l)^{-1}|^2\ \q(0),_m \tilde N_m+R_2,
\label{opt3}
\end{align}
with $|R_2(t)|\le \delta \tilde H(t)^2+C_\delta t \tilde H(t)^4 +C_\delta N(u_0)$.
Together with (\ref{H246}), (\ref{H3}), (\ref{H5}), (\ref{H7}), (\ref{E11}), (\ref{E121122}) and (\ref{E1231}), this provides us on $[0,T_\kappa]$ with
\begin{align}
\sum_{l=1}^K \int_{(0,1)^3}\xi(\theta_l)|D^2(\v_{tt}\circ\theta_l|^2+&\frac{1}{2} \sum_{l=1}^K \int_{\partial\e(\tilde\Omega)} | D^2 V_t(l)\circ E(l)^{-1}|^2\ (-\q(0),_m \tilde N_m)\n\\
&\le  \delta \tilde H(t)^2+C_\delta t \tilde H(t)^4 +C_\delta N(u_0).
\label{opt4}
\end{align}
Similarly as in Section \ref{L10}, from the pressure condition (\ref{lindblad}), this provides us with an estimate of the type:
\begin{align}
\|\v_{tt}\|_2^2+\|\v_t\|^2_\cd\le \delta \tilde H(t)^2+C_\delta t \tilde H(t)^4 +C_\delta N(u_0).
\label{opt5}
\end{align}
The control
\begin{align}
\|\q\|_3^2 \le \delta \tilde H(t)^2+C_\delta t \tilde H(t)^4 +C_\delta N(u_0),
\label{opt6}
\end{align}
is then easy to achieve by elliptic regularity on the pressure system.

Next, we see that (\ref{opt5}) implies that for any $l\in \{1,...,K\}$, $\v_t\circ\theta_l$, and thus $(\q\circ\theta_l),_3 \e\circ\theta_l,_1\times\e\circ\theta_l,_2$, are controlled in $H^2((0,1)^2\times\{0\})$ by the right-hand side of (\ref{opt5}). This implies the same control on $(0,T_\epsilon)$ for
$\d\frac{\q\circ(\theta_l),_3 \e\circ\theta_l,_1\times\e\circ\theta_l,_2}{\bigl|\q\circ(\theta_l),_3 \e\circ\theta_l,_1\times\e\circ\theta_l,_2\bigr|}$ in $H^2((0,1)^2\times\{0\})$, {\it i.e.} that
\begin{equation}
\label{opt7}
|\tilde n|^2_{2,\tilde\Gamma}\le \delta \tilde H(t)^2+C_\delta t \tilde H(t)^4 +C_\delta N(u_0),
\end{equation}
which brings the $H^\se$ regularity of the {\it domain} $\e(\tilde\Omega)$.

Now, for $\v$, we notice that from the identity on $(0,1)^2\times\{0\}$:
\begin{align*}
V_{tt}(l)+\q\circ\theta_l,_3 V(l),_1\times E(l),_2 + \q\circ\theta_l,_3 E(l),_1\times V(l),_2 + \q_t\circ\theta_l,_3 E(l),_1\times E(l),_2=0,
\end{align*}
we infer by taking the scalar product of the above vector by $E(l),_1$ that
\begin{equation*}
|V(l),_1\cdot \tilde n|^2_{\td,(0,1)^2\times \{0\}}\le \delta \tilde H(t)^2+C_\delta t \tilde H(t)^4 +C_\delta N(u_0),
\end{equation*}
which by divergence and curl relations for $\xi(\theta_l) V(l),_1(E(l)^{-1})$ similar in spirit to the ones in Section \ref{L8}, leads to
\begin{equation*}
\|\xi(\theta_l) V(l),_1\|^2_{2,(0,1)^3}\le \delta \tilde H(t)^2+C_\delta t \tilde H(t)^4 +C_\delta N(u_0).
\end{equation*}
In a similar fashion,
\begin{equation*}
\|\xi(\theta_l) V(l),_2\|^2_{2,(0,1)^3}\le \delta \tilde H(t)^2+C_\delta t \tilde H(t)^4 +C_\delta N(u_0).
\end{equation*}
Now, with divergence and curl relations for $\xi(\theta_l) V(l)(E(l)^{-1})$ similar in spirit to the ones in Section \ref{L8}, this leads to
\begin{equation*}
\|\xi \v\|^2_3\le \delta \tilde H(t)^2+C_\delta t \tilde H(t)^4 +C_\delta N(u_0),
\end{equation*}
and consequently, with the control of the divergence and curl of $\v$ inside $\tilde\Omega$ as in Section \ref{L8}, we get
\begin{equation}
\label{opt8}
\|\v\|^2_3\le \delta \tilde H(t)^2+C_\delta t \tilde H(t)^4 +C_\delta N(u_0).
\end{equation}
Now, with the estimates (\ref{opt5}), (\ref{opt6}), (\ref{opt7}) and (\ref{opt8}), we then get similarly as in Section \ref{L11} the existence of a time $T>0$ independent of $\epsilon$ such that on $(0,T)$ the estimates (\ref{assume2}) hold, and such that we have $\tilde H(t)\le  N(u_0)$ on $(0,T)$ for any $\epsilon>0$ small enough. Therefore, we have a solution to the problem with optimal regularity on the initial data and domain as the weak limit as $\epsilon\rightarrow 0$.

\section{Uniqueness}
\label{L13}
Let $(v,q)$ and $(\bar v,\bar q)$ be solutions of \ref{lindblad} on $[0,T]$. We denote $\delta v=v-\bar v$ and $\delta q=q-\bar q$. We then introduce the energy:
$$\d f(t)= \sum_{l=1}^K  \int_{{(0,1)^d}}\xi(\theta_l)|D^2 (v_{tt}\circ\theta_l-\bar v_{tt}\circ\theta_l)|^2,$$ where $D^2 v$ stands for any
second order horizontal space derivative $v,_{i_1i_2}$. By proceeding in the same way as in the previous section, and using the fact that the divergence and curl of $\delta v$ have a transport type structure as well, 
we obtain an energy inequality similar to (\ref{opt5}), without the presence of
$N(u_0)$ (since $\delta v(0)=0$). This establishes uniqueness of solutions.

\section*{Acknowledgments}
DC and SS were supported by the National Science Foundation under 
grant NSF ITR-0313370.  We thank the referee for the major time and effort
expended on the careful reading of our manuscript.

\end{document}